\documentclass[12pt,letterpaper,titlepage]{amsart}
\usepackage{amsmath, amssymb, amsthm, amsfonts,amscd,xr}
\numberwithin{equation}{section}
\newtheorem{thm}{Theorem}[section]
\newtheorem{theorem}[thm]{Theorem}
\newtheorem{cond}[thm]{Conditions}

\newtheorem{conv}[thm]{Convention}
\newtheorem{lemma}[thm]{Lemma}
\newtheorem{lem}[thm]{Lemma}

\newtheorem{cor}[thm]{Corollary}
\newtheorem{ex}[thm]{Example}
\newtheorem{prop}[thm]{Proposition}

\newtheorem{dfn}[thm]{Definition}
\newtheorem{rem}[thm]{Remark}
\newtheorem{con}[thm]{Conjecture}

\makeatletter
\def\section{\@startsection {section}{1}{\z@}{3.5ex plus 1ex minus
    .2ex}{2.3ex plus .2ex}{\large\bf}}
    \def\subsection{\@startsection{subsection}{2}{\z@}{3.25ex plus 1ex minus
 .2ex}{1.5ex plus .2ex}{\bf}}
\makeatother

\def\qed{\hfill $\square$\par\vspace{5pt}}

\def\a{\alpha}
\def\b{\beta}
\def\l{\lambda}
\def\e{\epsilon}
\def\g{\gamma}
\def\i{\iota}
\def\ka{\kappa}
\def\om{\omega}
\def\d{\delta}

\def\cf{\mathfrak{C}}

\def\hf{\mathfrak{H}}
\def\k{{\bf{k}}}

\def\nf{\mathfrak{N}}
\def\nc{\mathcal{N}}

\def\Xf{\mathfrak{X}}
\def\Yf{\mathfrak{Y}}

\def\Hf{\mathfrak{H}}
\def\zf{\mathfrak{Z}}
\def\sf{\mathfrak{S}}
\def\qu{\operatorname{qu}}
\def\N{{\mathbb N}}
\def\R{{\mathbb R}}
\def\C{{\mathbb C}}

\def\Q{{\mathbb Q}}
\def\Z{{\mathbb Z}}
\def\q{{\bf q}}

\def\Bc{\mathcal{B}}
\def\O{\mathcal{O}}

\def\Ze{\mathcal{Z}}
\def\Cc{\mathcal{C}}

\def\M{\mathcal{M}}
\def\I{\mathcal{I}}
\def\A{\mathcal{A}}
\def\B{\mathcal{B}}
\def\H{\mathcal{H}}
\def\F{\mathcal{F}}
\def\J{\mathcal{J}}
\def\K{\mathcal{K}}
\def\L{\mathcal{L}}

\def\P{\mathcal{P}}
\def\V{\mathcal{V}}
\def\W{\mathcal{W}}
\def\Ri{\mathcal{R}}
\def\TW{W_0\backslash T}

\def\cc{C^\infty}
\def\waf{W^{\mathrm{aff}}}
\def\rnr{R_{\mathrm{nr}}}

\renewcommand{\Im}{\mbox{\rm Im}\,}
\makeindex
\begin{document}
\title{On the spectral decomposition of affine Hecke algebras}
%\author{Eric M. Opdam
%\thanks{The author would like to thank Erik van den Ban, Gerrit
%Heckman and Klaas Slooten for many useful remarks and fruitful
%discussions},\\ KdV Institute for Mathematics,\\ Universiteit van
%Amsterdam,\\ Plantage Muidergracht 24\\ 1018 TV Amsterdam, The
%Netherlands\\ email: opdam@wins.uva.nl\\}
\author{Eric M. Opdam}
\address{Korteweg-De Vries Institute for Mathematics,
Universiteit van Amsterdam, Plantage Muidergracht 24, 1018 TV
Amsterdam, The Netherlands.}
\thanks{The author would like to thank Erik van den Ban, Patrick Delorme,
Gerrit Heckman, Mark Reeder, Klaas Slooten and Atsuko Yamamoto for their
helpful comments and remarks}
\email{opdam@science.uva.nl}
\date{\today}

\begin{abstract}
An affine Hecke algebra $\H$ contains a large abelian subalgebra
$\A$ spanned by the Bernstein-Zelevinski-Lusztig basis elements
$\theta_x$, where $x$ runs over (an extension of) the root lattice.
The center $\Ze$ of $\H$ is the subalgebra of Weyl group invariant
elements in $\A$. The natural trace (``evaluation at the
identity'') of the affine Hecke algebra can be written as integral
of a certain rational $n$-form (with values in the linear dual of
$\H$) over a cycle in the algebraic torus $T=spec(\A)$. This cycle
is homologous to a union of ``local cycles''. We show that this
gives rise to a decomposition of the trace as an integral of
positive local traces against an explicit probability measure on
the spectrum $W_0\backslash T$ of $\Ze$. From this result we
derive the Plancherel formula of the affine Hecke algebra.
\end{abstract}
\maketitle \tableofcontents
\section{Introduction}
In this paper we will discuss the spectral decomposition of an affine
Hecke algebra $\H$ defined over $\mathbb{C}$ or, more precisely,
of a natural positive trace $\tau$ defined on $\H$.
In the standard basis of $\H$, $\tau$ is simply defined by
$\tau(T_e)=1$ and $\tau(T_w)=0$ if $w\not=e$.
In addition, $\H$ comes
equipped with the natural $*$-operator $T_w^*=T_{w^{-1}}$.
This defines
a pre-Hilbert structure on $\H$ by $(x,y):=\tau(x^*y)$. The
regular representation $\l\times\rho$ extends to the
Hilbert completion $\mathfrak{H}$ of $\H$,
and by the spectral decomposition of $\tau$ we mean the
decomposition of $\mathfrak{H}$ in irreducible
$*$-representations of $\H\times\H$. By classical results
on the decomposition of traces on $C^*$-algebras of type I
(see for example \cite{dix2}), this is equivalent to
the problem of decomposing the trace $\tau$ as a
superposition of irreducible characters of $*$-representations
of $\H$.
We will call this decomposition the Plancherel
decomposition of $\H$, and the associated positive measure on
the spectrum $\hat\H$ will be called the Plancherel measure.

In the case of a Hecke algebra of finite type we have the well
known decomposition formula
\begin{equation}\label{eq:finite}
\tau=\frac{1}{P}\sum \chi_\pi d_\pi,
\end{equation}
where $P$ denotes the Poincar\'e polynomial of $\H$ (we assume that
$P\not=0$), $\pi$ runs over the finite set of irreducible
representations of $\H$, $\chi_\pi$ denotes the corresponding
character of $\pi$, and $d_\pi$ is the generic degree of $\pi$.
The formula we are going to discuss in the present paper is the
affine analog of equation (\ref{eq:finite}).

This paper is the sequel to \cite{EO},
where we made a basic study of the Eisenstein functionals of an
affine Hecke algebra $\H$. These Eisenstein functionals are
holomorphic functions of a spectral parameter $t\in T$, where $T$
is a complex $n$-dimensional algebraic torus naturally associated
to $\H$. In \cite{EO}, we derived a representation of $\tau$ as
the integral of the
normalized Eisenstein functional times the holomorphic extension
of the Haar measure of the compact form of $T$, against a certain
``global $n$-cycle'' (a coset of the compact form of $T$) in $T$.
The kernel of this integral is a meromorphic $(n,0)$-form on $T$.

The present paper takes off from that starting point, and refines
step-by-step the above basic complex function theoretic
representation formula for $\tau$ until we reach the level of
the spectral decomposition of $\tau$, extended to a tracial state
on the $C^*$-algebra hull $\cf$ of
the regular representation of $\H$ (Main Theorem \ref{thm:mainp}). On
the simpler level of the spherical or the anti-spherical
subalgebra, a similar approach can be found in \cite{Mat} and
\cite{HOH}. In the case of the spherical algebra one should of
course also mention the classical work \cite{Ma}, although the
point of view is different there, and based on analysis on a
reductive $p$-adic group.

\subsubsection{Motivation}
There are various motivations for the study of the spectral
resolution of $\tau$.
A natural application of such a decomposition is the p-adic
analog of the Howlett-Lehrer theory for finite reductive groups,
see for instance \cite{M}, \cite{Lu3}, \cite{HOH}, \cite{Re0}
and \cite{Re}. Here one considers an affine Hecke algebra which
arises as the centralizer algebra of a certain induced
representation
of a p-adic reductive group $G$. The Plancherel measure of $\H$ can
be interpreted as the Plancherel measure of $G$ on a part of
$\hat{G}$ in this situation.
In view of this application it is important that we obtain an
(almost) explicit product formula for the Plancherel measure
(see Main Theorem \ref{thm:mainp}).
In addition we characterize exactly which characters
$W_0r\in W_0\backslash T$
of the spectrum of the center $\Ze$ of $\H$ support a
discrete series representation of $\H$ (see Theorem \ref{thm:support}).
These are the so-called ``residual points''
(see Appendix \ref{sub:defn}). This result was recently applied to the
representation theory of reductive p-adic groups, see \cite{Hei}.

Another motivation for this approach is that it sets the stage
for the definition of a Schwartz-completion $\mathfrak{S}$ of $\H$
(see Subsection \ref{sub:normunif}), and for the subsequent study of
the Fourier transform and its inversion on the level of this
Fr\'echet algebra (joint work with Patrick Delorme, to appear).
This is related to the study of the K-theory and
the cyclic homology of $\H$ and its reduced $C^*$-algebra $\cf$, in
the spirit of \cite{was1}, \cite{was2}.
This point of view is particularly interesting for non-simply laced
cases, since it is natural to expect that the K-theory
does not depend on the parameters $q(s)$ of the Hecke algebra.
On the other hand, in the ``generic case'' these matters seem to be
considerably easier to understand than in the ``natural cases'',
where the logarithms of the parameters have rational relations.
In view of this, it is important that we allow the parameters
$q(s)$ of the affine Hecke algebra to assume any real value $>0$.
\subsubsection{Outline}\label{subsub:out}
It may be helpful to give the reader a rough outline of this paper,
and an indication of the guiding principles in the
various stages. We also refer the reader to Subsection
\ref{sub:out} for a more detailed outline and formulation of the
main results (see in particular \ref{sl2}, \ref{sl3} for
the results on the Plancherel measure).

(0). The starting point of the present paper is the definition of
the Eisenstein functional of the affine Hecke algebra, in \cite{EO}.
These functionals are matrix coefficients of minimal principal
series modules. The study of their intertwining operators
led to a representation of the trace $\tau$ of $\H$, as
an integral of a certain rational kernel over a ``global'' cycle
(see formula (\ref{eq:basic})).

(1). In section \ref{sect:pre} we recall the definition and first
properties of (extended) affine Hecke algebras, we collect some
basic facts from the theory of $C^*$-algebras, and we adapt
certain classical results from the representation theory of reductive
groups to our context. We conclude this section with a discussion
of the properties of the natural map
$p_z:\hat\cf\to\operatorname{Spec}(\Ze)$, the spectrum of the
center $\Ze$ of $\H$, in view of the main results of
this paper.

(2). The study of the residues of the rational kernel for $\tau$
as in formula (\ref{eq:basic}), in Section \ref{sec:loctau}.
This involves a general (but basic) scheme for the calculation of
multivariable residues. After symmetrization over the Weyl group,
the result is a decomposition of $\tau$ as an integral of local
tracial states against an explicit probability measure on the
spectrum $\operatorname{Spec}(Z)=W_0\backslash T$.
The main tools in this process are the {\it positivity} of $\tau$,
and the geometric properties of the collection of residual cosets
(Appendix \ref{sub:defn}).
This step is called the ``localization of the trace $\tau$''.

(3). The local trace (as was mentioned in (2)) defined at an orbit
$W_0t\subset T$, arises as an integral of the Eisenstein kernel
over a ``local cycle'' which is defined in an arbitrarily small
neighborhood of the orbit $W_0t$. This gives a natural extension
of the local trace to localizations of the Hecke algebra itself
(localization as a module over the sheaf of analytic functions on
$W_0\backslash T$).

The analytic localization of the Hecke algebra has a
remarkable structure discovered by Lusztig in \cite{Lu}. This part
of the paper is not self-contained, but draws heavily on the paper
\cite{Lu}. By Lusztig's wonderful structure theorem we can now
investigate the local traces. We find in this way that
everything is organized in accordance with Harish-Chandra
parabolic induction (the philosophy of cusp forms). The local
traces at residual cosets give rise to finite dimensional Hilbert
algebras which we call ``residual algebras''. Their generic structure
reduces (via parabolic induction) to the case of the residual
algebras at ``residual points'' of certain semisimple subquotients
of the Hecke algebra (see Subsection \ref{sub:gene}).

These matters concerning the localization
of $\H$ are studied in Section \ref{sect:loc}, leading to
the main result
Theorem \ref{thm:mainp}. The support of the Plancherel
measure and the Plancherel density are expressed in terms of the
discrete series of Levi subquotient algebras of $\H$, and their
Plancherel masses (formal dimensions).

(4). At this point, two essential problems remain: The
classification of the discrete series representations, and the
determination of their formal dimensions. Regarding the first problem,
we have determined the orbits $W_0r\in W_0\backslash T$ which arise as
the central character of a discrete series representation of $\H$
in Theorem \ref{thm:support}. We have no further information to
offer on this problem in this paper.

Section \ref{sec:inv} is devoted to the second of these problems.
In order to explain our approach, let $\nu(\{W_0r\})$ denote the
Plancherel mass of the central character $W_0r$
with respect to the restriction of the
tracial state $\tau$ to the center $\Ze$ of $\H$.
In Subsection \ref{sub:chiA} we find that the formal dimension of
an irreducible discrete series representation $\d$ whose
central character
is a certain residual point $W_0r$, is equal to the product of
$\nu(\{W_0r\})$ and a certain positive real number
$d_\d>0$ (called ``the residual degree'' of $\d$) depending
on $\d$ (see Corollary \ref{cor:fdim} and
Theorem \ref{thm:nu}).
(These residual degrees are normalized such
that $\sum\operatorname{dim}(\d)d_\d=1$, where the sum runs over
all square integrable $\d$ whose central character is $W_0r$).

{\it The factor $\nu(\{W_0r\})$ is a certain explicit product
(explicit up to a nonzero rational multiple) of rational
functions evaluated at the central character $W_0r$.}
The problem that arises here is that we have not much
information about
the behaviour of the individual ``residual degrees''
$d_\d$ as functions of the parameters $q(s)$.

In Section \ref{sec:inv} we resolve this matter. If we write the
labels $q(s)>0$ in the form $q(s)=\q^{f_s}$ for certain real numbers
$f_s$ and $\q>1$, we prove that the residual Hilbert algebras
are {\it independent} of the base $\q>1$. In other words, the constants
$d_\d>0$ are {independent} of $\q>1$. This proves that
{\it all
irreducible discrete series representations of $\H$ associated
with a central character $W_0r$, have a formal dimension which is
proportional to the mass $\nu(\{W_0r\})$, with a positive real
ratio of proportionality which is independent of $\q>1$}.

In addition we conjecture that the positive reals $d_\d$ are
actually rational numbers (cf. Conjecture \ref{rem:ell}).
This conjecture is subject of joint work in progress with Mark Reeder
and Antony Wasserman.

(5). In Appendix \ref{sub:defn} we study the geometry of the set
of singularities with maximal pole order of the rational $n$-form
\begin{equation}\label{eq:simple}
\frac{dt}{c(t)c(t^{-1})}
\end{equation}
on $T$. This leads to the notion ``residual coset'', which is
crucially important for the understanding of the residues
of the kernel for $\tau$ in (\ref{eq:basic}). It is analogous to
the notion of residual subspace which was introduced in
\cite{HOH0}. The collection of these cosets can be classified,
and from this classification we verify certain important geometric
properties of this collection. These geometric facts are used in
Section \ref{sec:loctau} (especially in Subsection \ref{sub:chiA})
to establish regularity properties of the residues to be
considered in this paper.
\subsubsection{Residue calculus}
Let us make some remarks about the ``residue calculus''
on which much of this paper is ultimately based. At the heart of it
lies the elementary Lemma \ref{thm:resbasic}, which is an adapted version
of Lemma 3.1 of \cite{HOH0}.
This lemma roughly states that on a complex torus $T$,
any linear functional $\tau$ on the ring of Laurent polynomials
$\mathbb{C}[T]$ of the form
\begin{equation}
\tau(f)=\int_{t_0T_u}f\omega
\end{equation}
where $\omega$ is a rational $(n,0)$-form whose pole set is a
union of cosets of codimension $1$ subtori of $T$, can be represented
by a unique collection of ``local distributions'' living on certain
cosets of the compact form $T_u$ of $T$, and satisfying certain support
conditions.

In the context of the representation theory of $\H$, this
lemma becomes remarkably efficient.
We apply the lemma
to linear functionals of the form $a\to\tau(ah)$, where $h\in \H$
and $a\in\A$, a maximal abelian subalgebra of $\H$, using
formula (\ref{eq:basic}).
At this stage we symmetrize the ``local distributions'' for the
action of $W_0$. Using the elementary notion of ``approximating sequence''
(see Lemma \ref{lem:approx}) it is not hard to show that the symmetrized
local distributions inherit the positivity of $\tau$. This implies easily
that these symmetrized distributions are in fact compactly supported measures
on the spectrum $W_0\backslash T$ of the center $\Ze=\A^{W_0}$ of $\H$,
with values in the positive traces on $\H$ (see Corollary \ref{cor:exten}).
This means that all higher order terms in the local distributions
cancel out by the symmetrization by $W_0$.

In addition it follows by positivity that all measures are absolutely
continuous with respect to a scalar measure $\nu$,
the Plancherel measure of the center
$\Ze$ of $\H$.
In fact $\nu$ is obtained by evaluation
of the symmetrized local distributions at $1\in\H$.
Fortunately the poles of
formula (\ref{eq:basic}) simplify to the poles of expression
(\ref{eq:simple}) by this evaluation.
In this way we see that the contributions at non-residual, quasi-residual
cosets must cancel. We can bring into play the geometric
properties of the residual subspaces now, established in Appendix
\ref{sub:defn}, to prove that the Plancherel measure is smooth on
its support, that the local traces are tempered, and that the local traces
at discrete mass points of $\nu$ are finite linear combinations of
discrete series characters.

We get in this way a decomposition of $\tau$ as a superposition of
positive ``local traces'', which is an important step towards
the Plancherel decomposition of $\H$.

At the time of the writing of the paper \cite{HOH0}, working on the
quantum theory of a certain exactly solvable $n$-particle systems,
we were not aware of the already existing results in the spirit of
the above lemma on existence and uniqueness of residue distributions.
But we should certainly mention here the basic work of Langlands
\cite{L}, where residues of Eisenstein series are studied in the
theory of automorphic forms for reductive groupes. We also mention
the work of Arthur \cite{A}, \cite{A2} in this direction.
Langlands' work \cite{L} was elucidated by Moeglin and Waldspurger
in \cite{MW}. Langlands' result on existence of ``residue data''
can be found in Theorem V.2.2, and on uniqueness of ``residue data''
in the formulation of Theorem V.3.13(i) of \cite{MW}.
It should be pointed out however that these results are of a
different nature than our basic Lemma \ref{thm:resbasic}. Lemma
\ref{thm:resbasic} is a (very elementary) general result in distribution
theory, which has nothing to do with group theory.
On the other hand, the above results in \cite{MW} are formulated with
already symmetrized ``residue data'', using intertwining operators.
In order to even formulate a uniqueness property in this setting,
one first needs to show rather deep statements on the holomorphic
continuation of certain residue sums of Eisenstein series (see
V.3.2 of \cite{MW}).

More recently, inspired by the approach in \cite{HOH0},
Van den Ban and Schlichtkrull \cite{BS1} extended the method by allowing
for so-called residue weights. In this generality they applied
the residue calculus in their proof of the Plancherel formula for
semisimple symmetric spaces.
\section{Preliminaries and description of results}\label{sect:pre}
The algebraic background for our analysis was discussed in the
paper \cite{EO}. The main result of that paper is an inversion
formula (see equation (\ref{eq:basic})) which will be the starting
point in this paper. The purpose of this section is to define the
affine Hecke algebra $\H$ and to review the relevant notations and
concepts involved in the above result. Moreover we introduce a
$C^*$-algebra hull $\cf$ of $\H$, which will be the main object of
study in this paper. Finally we will give a more precise outline
of the results in the paper. We refer the reader to \cite{Lu} and
\cite{EO} for a more systematic introduction of the basic
algebraic notions.
\subsection{The affine Weyl group and its root datum}
A reduced root datum is a 5-tuple ${\mathcal
R}=(X,Y,R_0,R_0^\vee,F_0)$
\index{R@$\mathcal R$, root datum},
where $X$\index{X@$X,Y$, lattices}
and $Y$ are free abelian groups with
perfect pairing over $\Z$, $R_0\subset
X$\index{R2@$R_0\subset X$, reduced integral root system} is a
reduced integral root system, $R_0^\vee\subset
Y$\index{R2@$R_0^\vee\subset Y$, coroot system}
is the dual root
system of coroots of $R_0$, and
$F_0\subset R_0$\index{F@$F_0\subset
R_0$, simple roots of $R_0$} is a basis of simple roots. Each
element $\a\in R_0$ determines a reflection
$s_\a\in\operatorname{GL}(X)$ by
\begin{equation}\index{s@$s_\a$, reflection in $\a$}
s_\a(x)=x-x(\a^\vee)\a.
\end{equation}
The group $W_0$\index{W2@$W_0$, Weyl group of $R_0$} in
$\operatorname{GL}(X)$ generated by the $s_\a$ is called the Weyl
group. As is well known, this group is in fact generated by the
set $S_0$\index{S@$S_0$, simple reflections of $W_0$}
consisting of
the reflections $s_\a$ with $\a\in F_0$. The set $S_0$ is called
the set of simple reflections in $W_0$.

By definition the affine Weyl group $W$
\index{W@$W$, affine Weyl
group} associated with a reduced root datum $\mathcal R$ is the
group $W=W_0\ltimes X$. This group $W$ naturally acts on the set
$X$.

We choose once and for all a rational, symmetric, positive
definite, $W_0$-invariant pairing
$\langle\cdot,\cdot\rangle$
\index{<@$\langle\cdot,\cdot\rangle$!a@rational inner product on $X$,$Y$}
on $\mathbb{Q}\otimes Y$. This
defines a $W_0$ pairing on the Euclidean spaces
$\mathfrak{t}:=\mathbb{R}\otimes Y$\index{T4@$\mathfrak{t}$, Lie
algebra $\operatorname{Lie}(T_{rs})=\mathbb{R}\otimes Y$} and its dual
$\mathfrak{t}^*=\mathbb{R}\otimes Y$.
The action of $W$ on $X$ extends to an action of $W$ on
$\mathfrak{t^*}$ by means of isometries.

We can identify the set of
integral affine linear functions on $X$ with $Y\times\Z$ via
$(y,k)(x):=(x,y)+k$. It is clear that $w\cdot
f(x):=f(w^{-1}x)$ defines an action of $W$ on $Y\times \Z$. The
affine root system is by definition the subset
$R^{\mathrm{aff}}=R_0^\vee\times\Z\subset Y\times \Z$
\index{R5@$R^{\mathrm{aff}}$, affine root system}.
Notice that $R^{\mathrm{aff}}$ is a
$W$-invariant set in $Y\times \Z$ containing the set of coroots
$R_0^\vee$. Every element
$a=(\alpha^\vee,k)\in R^{\mathrm{aff}}$
\index{a5@$a=(\alpha^\vee,k)$, affine root}
defines an
affine reflection $s_a\in W$
\index{s@$s_a$, affine reflection in $a$}
, acting on $X$ by
\begin{equation}
s_a(x)=x-a(x)\alpha.
\end{equation}
The reflections $s_a$ with $a\in R^{\mathrm{aff}}$ generate
a normal subgroup
$\waf=W_0\ltimes Q$
\index{W2@$\waf=W_0\ltimes Q\subset W$}
of $W$,
where $Q\subset X$
\index{Q@$Q$, root lattice}
denotes the root
lattice $Q=\mathbb{Z}R_0$. We can choose a basis of simple
affine roots $F^{\mathrm{aff}}$
\index{F@$F^{\mathrm{aff}}$, affine simple roots}
by
\begin{equation}
F^{\mathrm{aff}}:=\{(\alpha^\vee,1)\mid\alpha\in
S^m\}\cup \{(\alpha^\vee,0)\mid \alpha\in F_0\},
\end{equation}
where $S^m$ consists of the set of minimal coroots with respect to
the dominance ordering on $Y$. It is easy to see that every affine
root is an integral linear combination of elements from $F^{\mathrm{aff}}$
with
either all nonnegative or all nonpositive coefficients. The set
$R^{\mathrm{aff}}$ of affine roots is thus a disjoint union of the set of
positive affine roots $R^{\mathrm{aff}}_+$ and the set of negative affine roots
$R^{\mathrm{aff}}_-$. The set $S^{\mathrm{aff}}$\index{S@$S^{\mathrm{aff}}$,
simple reflections of $\waf$}
of simple reflections in $W$ is by definition
the set of reflections in $W$ associated with the fundamental
affine roots. They constitute a set of Coxeter generators for the
normal subgroup $\waf\subset W$.

There exists an Abelian complement to $\waf$ in $W$.
This is best understood by introducing the important length
function $l$ on $W$.
The splitting $R^{\mathrm{aff}}=R^{\mathrm{aff}}_+\cup R^{\mathrm{aff}}_-$
\index{R5@$R^{\mathrm{aff}}_\pm$, positive (negative) affine roots}
described above implies that
$R^{\mathrm{aff}}_+ \cap s_a(R^{\mathrm{aff}}_-)=\{a\}$
when $a\in F^{\mathrm{aff}}$.
Define, as usual, the length of an element
$w\in W$ by
\[
l(w):=|R^{\mathrm{aff}}_+\cap w^{-1}(R^{\mathrm{aff}}_-)|.
\]
\index{l1@$l$, length function on $W$}
It follows that, when $a\in F^{\mathrm{aff}}$,
\begin{equation}
\label{eqn:simlength}
{l}(s_aw)=
\begin{cases}
{l}(w)+1\text{ if }w^{-1}(a)\in R^{\mathrm{aff}}_+.\\
{l}(w)-1\text{ if }w^{-1}(a)\in R^{\mathrm{aff}}_-.\\
\end{cases}
\end{equation}
For any $w\in W$ we may therefore write $w=\omega \tilde w$ with
$\tilde w\in\waf$ and with $l(\omega)=0$ (or equivalently,
$\omega(F^{\mathrm{aff}})=F^{\mathrm{aff}}$). This shows that the set $\Omega$
\index{0Y@$\Omega$, length $0$ elements in $W$}
of elements of
length $0$ is a subgroup of $W$ which is complementary to the
normal subgroup $\waf$, so that we have the decomposition
\[
W=\Omega\ltimes \waf.
\]
Hence $\Omega\simeq W/\waf\simeq X/Q$ is a
finitely generated Abelian group.

Let $m:X\to P$
(where $P$\index{P@$P$, weight lattice} denotes the weight lattice)
denote the homomorphism that is adjoint to the inclusion
$Q^\vee\to Y$. If we
write $Z_X\subset X$
\index{Z@$Z_X$, length $0$ translations in $W$}
for its kernel, then $Z_X\subset\Omega$.
We have
$\Omega/{Z_X}=\Omega_f$
\index{0Y@$\Omega_f=\Omega/Z_X$}
where $Z_X$ is free and
$\Omega_f=m(X)/Q\subset P/Q$ is finite. It is easy to see that
$Z_X$ is the subgroup of elements in $X$ that are central in $W$.
The finite group $\Omega_f$ acts faithfully on $S^{\mathrm{aff}}$
by diagram automorphisms.

The dual cone $X^+$
\index{X2@$X^+\subset X$, dominant cone}
of the cone $Q_+$ spanned by the positive
roots is called the cone of dominant elements of $X$. Thus $x\in
X$ belongs to $X^+$ if and only if $\langle x,\a^\vee\rangle\geq 0$
for all positive
roots $\a\in R_{0,+}$.
Notice that $X^+\cap X^-$ equals the
sublattice $Z_X\subset X$ of translations of length $0$.

Write $v=v_0+v^0$ for the splitting of $v\in\mathfrak{t^*}$
according to the orthogonal decomposition $\mathfrak{t^*}=
\mathfrak{t^*}_0+\mathfrak{t^*}^0$, where
$\mathfrak{t^*}_0=\mathbb{R}\otimes Q$. We define a norm
\begin{equation}\label{eq:bigL}
\mathcal{N}(w)=l(w)+\|w(0)^0\|
\end{equation}
\index{N@$\mathcal{N}$, norm function on $W$}
for $w\in W$. Notice that for all $w,w^\prime\in W$,
$ww^\prime(0)^0=w(0)^0+w^\prime(0)^0$.
Thus $\mathcal{N}(\omega w)=\nc(w\omega)=
l(w)+\mathcal{N}(\omega)$
if $w\in\waf$ and $\omega\in\Omega$.
We also see that for all $\omega\in\Omega$,
$\nc(\omega^k)=k\nc(\omega)$ for
$k\in\mathbb{N}$. It follows easily
that $\omega\in \Omega$ has finite order if and only if
$\nc(\omega)=0$. Finally notice that it also follows that
\begin{equation}
\nc(ww^\prime)\leq\nc(w)+\nc(w^\prime)
\end{equation}
\subsection{Parabolic subsystems}\label{sub:par}
\index{X1a@$X_L\supset R_L$, lattice of $\Ri_L$,
character lattice of $T_L$|(}
\index{Y@$Y_L\supset R_L^\vee$, lattice of $\Ri_L$,
cocharacter lattice of $T_L$|(}
\index{W4@$W_P$, Weyl group of $R_P$,
parabolic subgroup $W_0$|(}
\index{R1@$\Ri_L$, semisimple root datum associated to $L$|(}
\index{R1@$\Ri^L$, root datum associated to $L$|(}
\index{R3@$R_P\subset R_0$, parabolic subsystem,
root system of $\Ri_P$|(}
\index{W4a@$W^P=W_0/W_P$, set of left cosets $wW_P$. If $P\subset F_0$,
identified with shortest length representatives|(}
An important role will be played by parabolic subgroups
of a Weyl group. A root subsystem $R^\prime\subset R_0$ is
called parabolic if $R^\prime=R_0\cap \mathbb{Q}R^\prime$.
Let $P\subset R_{0,+}\cap R^\prime$ be the basis of simple roots.
We then often write $R_P$ instead of $R^\prime$. The subgroup
$W_P:=W(R_P)$ is called the associated parabolic subgroup. If
$P\subset F_0$, we call $R_P$ and $W_P$ standard parabolic. Every
parabolic subgroup is conjugate to a standard parabolic subgroup.
We denote by $W^P$ the set of left cosets $W_0/W_P$. If $W_P$
is standard, we identify this quotient with the set of
distinguished coset representatives of minimal length.

In many instances we obtain a parabolic subsystem $R^\prime$
as the set of roots orthogonal to some subspace
$V^L\subset\mathfrak{t}:=\mathbb{R}\otimes Y$ which has the property
that $V^L=\cap\operatorname{ker}(\a)$ where we take the intersection
over all the roots $\a$ such that $\a(V^L)=0$. By abuse of
notation we usually denote this parabolic subsystem by $R_L$.
Similarly we write $W_L$ and $W^L$. We now denote the basis
of $R_{L,+}$ by $F_L$.

To a parabolic subsystem $R_P\subset R_0$ we associate a root
datum $\Ri^P:=(X,Y,R_P,R_P^\vee,P)$ and a root datum
$\Ri_P:=(X_P,Y_P,R_P,R_P^\vee,P)$ where
$Y_P:=Y\cap\mathbb{Q}R_P^\vee$ and $X_P:=X/(X\cap(R_P^\vee)^\perp)$.
\index{X1a@$X_L\supset R_L$, lattice of $\Ri_L$,
character lattice of $T_L$|)}
\index{Y@$Y_L\supset R_L^\vee$, lattice of $\Ri_L$,
cocharacter lattice of $T_L$|)}
\index{W4@$W_P$, Weyl group of $R_P$,
parabolic subgroup $W_0$|)}
\index{R1@$\Ri_L$, semisimple root datum associated to $L$|)}
\index{R1@$\Ri^L$, root datum associated to $L$|)}
\index{R3@$R_P\subset R_0$, parabolic subsystem,
root system of $\Ri_P$|)}
\index{W4a@$W^P=W_0/W_P$, set of left cosets $wW_P$. If $P\subset F_0$,
identified with shortest length representatives|)}
\subsection{Root labels}\label{sub:rl}
The second ingredient in the definition of $\H$ is a function
$q$\index{q@$q$, $l$-multiplicative function on $W$}
on $S^{\mathrm{aff}}$ with values in the group of invertible elements of a
commutative ring, such that
\begin{equation}\label{cond1}
q(s)=q(s^\prime) {\mathrm{\ if\ }} s\ {\mathrm{and}}\ s^\prime\
{\mathrm{are\ conjugate\ in\ }}W.
\end{equation}
A function $q$ on $S^{\mathrm{aff}}$, satisfying \ref{cond1}, can clearly be
extended uniquely to a length-multiplicative function on
$W$\index{Length multiplicative function}, also
denoted by $q$. By this we mean that the extension satisfies
\begin{equation}
q({ww^\prime})=q(w)q({w^\prime})
\end{equation} whenever
\begin{equation}
{l}(ww^\prime)={l}(w)+ {l}(w^\prime),
\end{equation}
and in addition,
\begin{equation}
\forall \omega\in\Omega: q(\omega)=1.
\end{equation}
Conversely, every length multiplicative function on $W$ restricts
to a function on $S^{\mathrm{aff}}$ that satisfies \ref{cond1}. Another way to
capture the same information is by assigning labels
$q_a$\index{q@$q_a$, affine root label}
to the
affine roots $a\in R^{\mathrm{aff}}$.
These labels are uniquely determined by the
rules
\begin{gather}\label{eqn:afflab}
\begin{split}
(i)&\ q_{wa}=q_a\ \forall w\in W,\mathrm{\ and\ }\\
(ii)&\ q(s_a)=q_{a+1}\ \forall a\in F^{\mathrm{aff}}.\\
\end{split}
\end{gather}
Note that a translation $t_x$ acts on an affine root
$a=(\a^\vee,k)$ by $t_xa=a-\a^\vee(x)$. Hence by $(i)$,
$q_{a}=q_{\a^\vee}$, except when $\a^\vee\in 2Y$, in which
case $q_{a}=q_{(\a^\vee,k(\mathrm{mod}2))}$. This last
case occurs iff $W$ contains direct factors which
are isomorphic to the affine Coxeter group whose
diagram equals $C_n^\mathrm{aff}$.

Yet another manner of labeling will play an important role.
It involves a possibly non-reduced root system
$\rnr$\index{R4@$\rnr$, non reduced root system},
which is defined by:
\begin{equation}
R_{\mathrm{nr}}:=
R_0\cup\{2\alpha\mid
\alpha^\vee\in R_0^\vee\cap 2Y\}.
\end{equation}
Now define labels for the roots $\a^\vee/2$ in
$\rnr^\vee\backslash R_0^\vee$ by:
\[ q_{\alpha^\vee/2}:=
\frac{q_{1+\alpha^\vee}}{q_{\alpha^\vee}}.
\]
\index{q@$q_{\alpha^\vee}$, label for $\a^\vee\in\rnr^\vee$}
This choice is natural, because it implies the formula
\begin{equation}
q(w)=\prod_{\alpha\in R_{\mathrm{nr},+}\cap
w^{-1}R_{\mathrm{nr,-}}} q_{\alpha^\vee},
\end{equation}
for all
$w\in W_0$.

Let $R_L\subset R_0$ be a parabolic root subsystem.
With respect to the root datum $\Ri_L$ we have
$R_{L,\mathrm{nr}}=\mathbb{Q}R_L\cap R_{\mathrm{nr}}
\subset R_{\mathrm{nr}}$.
In this sense we can define a label function denoted
by $q_L$
\index{q@$q_L$, restriction of $q$ to $\Ri_L$}
for the root datum $\Ri_L$, by restriction
from $R_{\mathrm{nr}}^\vee$ to $R_{L,\mathrm{nr}}^\vee$.
Similarly, we define $q^L$ by restriction of $q$ to $\Ri^L$.
\index{q@$q^L$, restriction of $q$ to $\Ri^L$}

We denote by $R_1$\index{R4@$R_1$, system of
long roots in $\rnr$}
the root system of long roots in
$R_\mathrm{nr}$. In other words
\begin{equation}
R_{1}:=
\{\alpha\in R_\mathrm{nr}\mid
2\alpha\not\in R_\mathrm{nr}\}.
\end{equation}
\subsection{The Iwahori-Hecke algebra as a Hilbert algebra}
\label{sub:iwhehil}
Many of the results of this subsection are well known, see
\cite{Mat}.
Let $\Ri$ be a root datum, and let
$\q$\index{q@$\q$, base for the labels $q(s)$}
be a real number with $\q>1$.
We assume that for all $s\in S^{\mathrm{aff}}$ we are
given a {\it real} number
$f_s$\index{f@$f_s=\log_{\q}(q(s))$}.
Throughout this paper we use the convention that
the labels as discussed in the previous subsection are defined by:
\begin{conv}\label{eq:scale} The labels are of the form
\begin{equation}
q(s)=\q^{f_s}\ \forall s\in S^{\mathrm{aff}}.
\end{equation}
\end{conv}

We write $q:=(q(s))_{s\in S^{\mathrm{aff}}}$
for the corresponding label function
on $S^{\mathrm{aff}}$.
The following theorem is well known.
\begin{thm}
There exists a unique complex associative algebra $\H=\H({\mathcal
R},q)$\index{H@$\H$, affine Hecke algebra}
with $\mathbb{C}$-basis
$(T_w)_{w\in W}$\index{T4@$T_w$, basis elements of $\H$}
which satisfy the
following relations:
\begin{enumerate}\label{eq:case}
\item[(a)] If ${l}(ww^\prime)={l}(w)+
{l}(w^\prime)$ then $T_wT_{w^\prime}
=T_{ww^\prime}$.
\item[(b)] If $s\in S^{\mathrm{aff}}$ then
$(T_s+1)(T_s-q(s))=0$.
\end{enumerate}
The algebra $\H=\H(\Ri,q)$ is called the affine Hecke algebra (or
Iwahori-Hecke algebra) associated to $(\Ri,q)$.
\end{thm}
We equip the Hecke algebra $\H$ with an anti-linear
anti-involutive $*$
\index{*@$*$!$h\to h^*$,
conjugate linear anti-involution of $\H$}
operator defined by $$ T_w^*=T_{w^{-1}}. $$
In addition, we define a trace functional
$\tau$
\index{0t@$\tau$, trace functional of $\H$}
on $\H$, by means of
$\tau(T_w)=\delta_{w,e}$. It is a well known basic fact that
\[
\tau(T_w^*T_{w^\prime})=\delta_{w,{w^\prime}}q(w),
\]
implying that $\tau$ is positive and central. Hence the formula
\[
(h_1,h_2):=\tau(h_1^*h_2)
\index{$(\cdot,\cdot)$!inner product on $\H$},
\]
defines an Hermitian inner product satisfying the following rules:
\begin{gather}\label{eqn:inprod}
\begin{split}
%\begin{enumerate}\label{eqn:inprod}
(i)&\ (h_1,h_2)=(h_2^*,h_1^*).\\
(ii)&\ (h_1h_2,h_3)=(h_2,h_1^*h_3).\\
%\end{enumerate}
\end{split}
\end{gather}
The basis $T_w$ is orthogonal for $(\cdot, \cdot)$.
We put
\begin{equation}\label{eq:norm}
N_w:=q(w)^{-1/2}T_w
\end{equation}\index{N@$N_w$, normalized basis elements of $\H$}
for the orthonormal basis of $\H$ that is obtained from the
orthogonal basis $T_w$ by scaling. Let us denote by
$\lambda(h)$
\index{0l@$\lambda(h)$, left multiplication by $h\in\H$}
and
$\rho(h)$\index{0r@$\rho(h)$, right multiplication by $h\in\H$}
the left and right multiplication operators on $\H$
by an element $h\in\H$. Let
$\hf$\index{H8@$\hf$, Hilbert completion of $\H$}
be the Hilbert space obtained
from $\H$ by completion; in other words, $\hf$ is the Hilbert
space with Hilbert basis $N_w$.
The operator $*$ extends to an isometric involution on $\hf$.
Let $B(\hf)$
\index{B@$B(\hf)$, bounded linear operators on $\hf$}
denote the space of bounded
operators on the Hilbert space $\hf$.
\begin{lem}
For all $h\in\H$, both $\l(h)$ and $\rho(h)$
extend to $\hf$ as bounded operators (elements of $B(\hf)$),
with $\Vert\lambda(h)\Vert=\Vert\rho(h)\Vert$.
For a simple reflection $s\in S^{\mathrm{aff}}$,
$\Vert\lambda(N_s)\Vert
=\operatorname{max}\{q(s)^{\pm 1/2}\}$.
\end{lem}
\begin{proof}
We first prove the formula for the norm
of $\Vert\lambda(N_s)\Vert$ ($s\in S^{\mathrm{aff}}$).
For every $w$ such that $l(sw)>l(w)$,
$\lambda(N_s)$ acts on the two-dimensional subspace $V_{w}$ of
$\H$ spanned by $N_w$ and $N_{sw}$ as a self-adjoint operator
with eigenvalues $q(s)^{1/2}$ and $-q(s)^{-1/2}$. Since $\hf$ is
the Hilbert sum of the subspaces $V_{w}$, we see that
$\lambda(N_s)$ extends to $\hf$ as a self-adjoint operator with
operator norm equal to $q(s)^{\pm 1/2}$.
Hence for any $h\in\H$, $\lambda(h)$ extends as a bounded
operator on $\hf$. Finally notice that
$(\l(h)^*(x))^*=\rho(h)(x)$, proving the equality
$\Vert\lambda(h)\Vert=\Vert\rho(h)\Vert$.
\end{proof}
The above lemma shows that $\H$ has the structure of a Hilbert
algebra in the sense of Dixmier \cite{dix1}. Moreover, this
Hilbert algebra is unital, and the Hermitian product is defined
with respect to the trace $\tau$.

We define the operator norm
$\Vert\cdot\Vert_o$
\index{$\Vert\cdot\Vert_o$, operator norm on $\H$}
on $\H$ by
$\Vert h\Vert_o:=\Vert\lambda(h)\Vert=\Vert\rho(h)\Vert$.
The closure of $\H$ with respect to the operator norm
$\Vert\cdot\Vert_o$ is denoted by
$\cf$\index{C@$\cf$, the reduced $C^*$ algebra of $\H$}.
The map $\l$ ($\rho$) extends to an isometry from $\cf$
to the $C^*$-subalgebra $\l(\cf)\subset B(\hf)$
($\rho(\cf)\subset B(\hf)$ resp.), the norm closure
of $\l(\H)$ ($\rho(\H)$ resp.).

We identify $\cf$ with a subset of $\hf$
via the continuous injection $c\to \l(c)(1)$. We
equip $\cf$ with the structure of a unital
$C^*$-algebra by the product
$c_1c_2:=\lambda(c_1)(c_2)=\rho(c_2)(c_1)$
and the $*$-operator coming from $\hf$. Then $\lambda$
($\rho$) is a faithful left (right) representation of $\cf$
in the Hilbert space $\hf$ (note that
we consider $\rho$ as a {\it right}
representation on $\hf$).
\begin{dfn}
We call $\cf$ the reduced $C^*$-algebra of $\H$, and
$\l$ ($\rho$) is called the left (right) regular representation
of $\cf$ on $\hf$.
\end{dfn}
An element $a\in\hf$ is called {\it bounded} if there exists
an element $\lambda(a)\in B(\hf)$ such that for all $h\in\H$,
\begin{equation}
\lambda(a)(h)=\rho(h)(a).
\end{equation}
By continuity we see that $\lambda(a)$ is uniquely determined by
$a=\lambda(a)(1)$.

When $a\in\hf$ is bounded, there also exists a unique
$\rho(a)$ such that for all $h\in\H$,
\begin{equation}
\rho(a)(h)=\lambda(h)(a).
\end{equation}
It is obvious that the elements of $\cf$ are bounded. Let us
denote by $\nf\subset\hf$
\index{N@$\nf$, von Neumann algebra
completion of $\H$} the subspace of bounded elements. We equip
$\nf$ with the involutive algebra structure defined by the product
$n_1n_2:=\lambda(n_1)(n_2)=\rho(n_2)(n_1)$ and the $*$ operator as
before.
\begin{prop}
The subspace $\lambda(\nf):=
\{\lambda(a)\mid a\in\nf\}\subset B(\hf)$ is the
von Neumann algebra completion
of $\lambda(\H)$. In other words, $\lambda(\nf)$ is the closure of
$\lambda(\H)$ in $ B(\hf)$ with respect to the strong topology (defined by
the semi-norms $T\to\Vert T(x)\Vert$ with $x\in\hf$). The analogous
statements hold when we replace $\lambda$ by $\rho$. The centralizing
algebra of $\lambda(\nf)$ is $\rho(\nf)$.
\end{prop}
\begin{proof}
All this can be found in \cite{dix1}, Chapitre I, paragraphe 5.
In general, $\lambda(\nf)$ is a two-sided ideal of the von Neumann
algebra hull of $\lambda(\H)$, but in the presence
of the unit $1\in\H$ the two spaces coincide. In fact, when
$A\in B(\hf)$ and $A$ is in the strong closure of $\lambda(\H)$,
it is simple to see that $A(1)\in\hf$ is bounded.
\end{proof}
The pre-Hilbert structure coming from $\hf$ gives $\nf$
itself the structure of a unital Hilbert algebra.
The algebra $\nf$ can and will be identified with its associated
standard von Neumann algebra $\lambda(\nf)$.
In this situation,
$\nf$ is said to be a saturated Hilbert algebra (with unit element).

Let
$\H^*$\index{H1@$\H^*$, algebraic dual of $\H$}
denote the algebraic dual of $\H$, equipped with its
weak topology. Notice that $\tau$ extends to $\H^*$ by the formula
$\tau(\phi):=\phi(1)$. The $*$-operator can be extended to $\H^*$
by $\phi^*(h):=\overline{\phi(h^*)}$. We have the following chain
of inclusions:
\begin{equation}
\H\subset\cf\subset\nf\subset\hf\subset\cf^\prime\subset\H^*.
\end{equation}
(where
$\cf^\prime$\index{C@$\cf^\prime$,
dual of $\cf$ as topological vector space}
denotes the space of continuous linear
functionals on $\cf$).
\begin{prop}
The restriction of $\tau$ to $\nf$ is central, positive and
finite. It is the natural trace of the Hilbert algebra $\nf$, in
the sense that
\begin{equation}
\tau(a)=(b,b)
\end{equation}
for every positive $a\in\nf$, and $b\in\nf$ such that $a=b^2$.
\end{prop}
\begin{proof}
A square root $b$ is in $\nf$ and is Hermitian (i.e. $b^*=b$).
Then $(b,b)=(1,a)=\tau(a)$.
\end{proof}
\begin{cor}
The Hilbert algebra $\nf$ is finite.
\end{cor}
\begin{cor}\label{cor:regrep}
The tracial state $\tau$ on $\cf$ is finite, and we have
$\l=\l_\tau$ and $\rho=\rho_\tau$, where $\l_\tau$ and $\rho_\tau$
are the representations of $\cf$ naturally associated with the
state $\tau$ (the classical GNS-construction).
\end{cor}
\begin{proof} This is immediate from the definitions,
see \cite{dix2}, Paragraphe 6.7.
\end{proof}
\subsection{Bernstein's description of the center $\Ze$}\label{sub:bern}
By a well known (unpublished) result of J.~Bernstein (see
\cite{Lu}), $\H$ can be viewed as the product $\H_0\A$ or $\A\H_0$
of an abelian subalgebra $\A$ (isomorphic to the group algebra of
the lattice $X$), and the Hecke algebra
$\H_0=\H(W_0,q|_{S_0})$
\index{H2@$\H_0=\H(W_0,q\mid_{S_0})\subset\H$}
of the finite Weyl
group $W_0$. Both product decompositions $\H_0\A$ and $\A\H_0$
give a linear isomorphism of $\H$ with the tensor product
$\H_0\otimes\A$.
The relations between products in $\H_0\A$ and in
$\A\H_0$ are described by the {\it Bernstein-Zelevinski-Lusztig
relations} (see for example \cite{EO}, Theorem 1.10), and with the
above additional description of the structure of $\A$ and $\H_0$
these give a complete presentation of $\H$.

The algebra
$\A$\index{A@$\A$, abelian subalgebra of $\H$}
has a $\C$-basis of invertible elements
$\theta_x$\index{0i@$\theta_x$, basis elements of $\A$}
(with $x\in X$) such that $x\to\theta_x$ is a
monomorphism of $X$ into the group of invertible elements of $\H$.
This monomorphism is uniquely determined by the property that
$\theta_x=N_{t_x}$ (see (\ref{eq:norm})) when $x\in X^+$. As an
important corollary of this presentation of $\H$, Bernstein
identified the center
$\Ze$\index{Z1@$\Ze$, the center of $\H$}
of $\H$ as the space $\Ze=\A^{W_0}$ of
$W_0$-invariant elements in $\A$ (see \cite{EO}, Theorem 1.11).
The following Proposition is well known and easy
(cf. \cite{EO}, Proposition 1.12):
\begin{prop}\label{prop:thetastar}
Let $w_0\in W_0$
\index{wa@$w_0$, longest element of $W_0$} denote
the longest element of $W_0$. Then we have for all $x\in X$:
\begin{equation}
\theta_x^*=T_{w_0}\theta_{-w_0(x)}T_{w_0}^{-1}.
\end{equation}
In particular, $\A\subset\H$ is not a $*$-subalgebra in general.
The center $\Ze\subset\H$ is a Hilbert subalgebra.
\end{prop}
Let $T=\operatorname{Spec}(\A)=\operatorname{Hom}_\Z(X,\C^\times)$
\index{T@$T=\operatorname{Hom}_\Z(X,\C^\times)$, complex algebraic torus}.
This algebraic torus of complex characters of $X$ has a natural
$W_0$-action, and we have
$\operatorname{Spec}(\Ze)\simeq W_0\backslash T$.
\begin{prop}\label{prop:rk}
The Hecke algebra $\H$ is finitely generated over its center
$\Ze$. At a maximal ideal $m=m_t$ (with $t\in T$) of $\Ze$, the
local rank equals $|W_0|^2$ if and only if the stabilizer group
$W_t\subset W_0$
\index{W3d@$W_t$, stabilizer in $W_0$ of $t\in T$}
is generated by reflections.
\end{prop}
\begin{proof}
It is clear that $\H\simeq \H_0\otimes \A$ is finitely generated
over $\Ze=\A^{W_0}$. When $W_t$ is generated by reflections, it is
easy to see that the rank of $m$-adic completion $\hat\A_m$ over
$\hat\Ze_m$ is exactly $|W_0|$ (see Proposition 2.23(4) of
\cite{EO}).
\end{proof}
This fact plays a predominant role in the representation
theory of $\cf$. Let us look at some basic consequences.
\begin{cor}
\begin{enumerate}
\item[(i)] Let $\pi$ be a finite dimensional irreducible representation
of $\H$ with representation space $V$.
The dimension of $V$ is less than or equal to $|W_0|$.
\item[(ii)] In addition, the center $\Ze$ of $\H$ acts by scalars on $V$.
Thus $\pi$ determines a ``central character''
$t_\pi\in\operatorname{Spec}(\Ze)$ such that for all $z\in\Ze$,
$\pi(z)=t_\pi(z)\operatorname{Id}_V$.
\item[(iii)] The characters of any finite set of inequivalent
finite dimensional irreducible representations of $\H$ are linearly
independent.
\item[(iv)] A topologically irreducible $*$-representation $\pi$
of the involutive algebra $\H$ is finite dimensional.
\end{enumerate}
\end{cor}
\begin{proof}
Elementary and well known. Use the Frobenius-Schur theorem for
(iii), Dixmier's version of Schur's lemma for (iv), and
Proposition \ref{prop:rk}.
\end{proof}
%\begin{proof}
%(i)
%There exists a joint eigenvector $v\in V$ for the finitely generated
%commutative subalgebra $\A$. Since $\pi$ is irreducible, this vector is cyclic.
%Consequently, $V=\H v$ has dimension at most $|W_0|$ by Bernstein's
%presentation $\H=\H_0\A$.

%(ii) As above, we can find a joint eigenvector $v\in V$ for the the finitely
%generated commutative subalgebra $\Ze$. Since $v$ is cyclic, this implies
%that the action of $\Ze$ is scalar in $V$.

%(iii)
%By the above, an irreducible finite dimensional character $\chi_\pi$
%factors through the finite dimensional complex algebra $\H/{\Ze_\pi\H}$, where
%$\Ze_\pi$ denotes the maximal ideal of $\Ze$ associated to $t_\pi$. By the
%Frobenius-Schur Theorem, the characters of the irreducible representations
%of $\H/{\Ze_t\H}$ are linearly independent. On the other hand,
%a finite set of distinct characters of $\Ze$ is linearly independent.
%Together these facts imply the
%desired result.

%(iv) In a
%topologically irreducible representation of the involutive algebra $\H$,
%the action of $\Ze$ is scalar and every nonzero vector is topologically
%cyclic. It follows by Proposition \ref{prop:rk} that $V$ has finite
%dimension.
%\end{proof}
\begin{cor}\label{cor:typeI} (See also \cite{Mat})
Restriction to $\H$ induces an injection of the set
$\hat\cf$\index{C@$\hat\cf$, dual (spectrum) of $\cf$}
into the space
$\hat{\H}$\index{H8@$\hat{\H}$, space of irreducible
$*$-representations of $\H$}
of finite dimensional irreducible $*$-representations
of $\H$.
Consequently, the $C^*$-algebra $\cf$ is of finite type I.
\end{cor}
\begin{proof}
Because $\H\subset\cf$ is dense, it is clear that a representation
of $\cf$ is determined by its restriction to $\H$ and that
(topological) irreducibility is preserved. Hence by the previous
Corollary, all irreducible representations of $\cf$ have finite
dimension.
\end{proof}
We equip
$T$ and $\TW$ with the {\it analytic} topology. Given
$\pi\in\hat\cf$ we denote by $W_0t_\pi\in\TW$ the character of
$\Ze$ such that $\chi_\pi(z)=\operatorname{dim}(\pi)z(t_\pi)$
(note that $\pi(\Ze)$ can not vanish identically since $1\in\Ze$).

By Proposition \ref{prop:thetastar}, the $*$-operator on
$\Ze$ is such that
$z^*(t)=\overline{z(\overline{t^{-1}})}$. When $\pi\in\hat\cf$,
we have $\chi_\pi(x^*)=\overline{\chi_\pi(x)}$.  It follows that
$\overline{t_\pi^{-1}}\in W_0t_\pi$ for all $\pi\in\hat\cf$.
Let us denote by $\TW^{herm}$ the closed subset $\{W_0t\in\TW\mid
\overline{t_\pi^{-1}}\in W_0t_\pi\}$ of $\TW$.
\begin{prop}\label{prop:im}
The map
$p_z:\hat\cf\to\TW$
\index{p9@$p_z:\hat\cf\to\operatorname{Spec(\Ze)}$, projection}
defined by $p_z(\pi)=W_0t_\pi$ is
continuous and finite. Its image
$S=p_z(\hat\cf)\subset\TW^{herm}$
\index{S@$S\subset\TW$, image of $p_z$}
is the spectrum $\hat{\overline{\Ze}}$ of
the closure ${\overline{\Ze}}$
\index{Z1@$\overline{\Ze}$, closure of $\Ze$ in $\cf$}
of $\Ze$ in $\cf$. The map $p_z:\hat\cf\to S$
is closed.
\end{prop}
\begin{proof} It is clear that the image is in $\TW^{herm}$ and
that the map is finite (by Proposition \ref{prop:rk}).
Since $W_0\backslash T$ is
Hausdorff and $\hat\cf$ is compact, the map $p_z$ is closed
if it is continuous.

So it remains to show that $p_z$ is continuous.
The closure ${\overline{\Ze}}\subset \cf$ is
a unital commutative $C^*$-subalgebra of $\cf$. By the Gelfand
transform it is isomorphic to the algebra of continuous functions
$C(\hat{\overline{\Ze}})$ on the
compact Hausdorff space $\hat{\overline{\Ze}}$.
Denote
by $\a$ the map $\a:\hat\cf\to\hat{\overline{\Ze}}$ defined by the condition
$\chi_\pi|_{{\overline{\Ze}}}=\operatorname{dim}(\pi)\a(\pi)$. By Proposition
2.10.2 of \cite{dix2}, $\a$ is surjective. In other words, every
primitive ideal $M$ of $\cf$ intersects ${\overline{\Ze}}$ in a maximal ideal
$m$ of ${\overline{\Ze}}$, and all maximal ideals of ${\overline{\Ze}}$ are
of this
form. The corresponding surjective map from the set
$\operatorname{Prim}(\cf)$ of primitive ideals of $\cf$ to the set
of maximal ideals $\operatorname{Max}({\overline{\Ze}})$ is also denoted by
$\a$. Next we claim that $\a$ is continuous. The topologies of
$\hat\cf$ and $\hat{\overline{\Ze}}$ are defined by the Jacobson topologies
on $\operatorname{Prim}(\cf)$ and $\operatorname{Max}({\overline{\Ze}})$.
This means that $U\subset\operatorname{Prim}(\cf)$ is closed if
and only if every $M\in\operatorname{Prim}(\cf)$ which contains
$I(U)=\cap_{u\in U}u$ is in $U$. Let $V\subset
\operatorname{Max}({\overline{\Ze}})$ be closed, and put $U=\a^{-1}(V)$. By
the surjectivity of $\a$ we have $I(U)\cap{\overline{\Ze}}=I(V)$. Hence if
$M\in\operatorname{Prim}(\cf)$ contains $I(U)$, then
$\a(M)=M\cap{\overline{\Ze}}$ contains $I(V)$, implying that $\a(M)\in V$.
Therefore $M\in U$, proving that $U$ is closed as desired.

Next, we consider the injective map $\b:\hat{\overline{\Ze}}\to\TW$ defined
by restriction to $\Ze\subset {\overline{\Ze}}$. Its image
$S\subset\TW^{herm}$ is bounded because for every $z\in\Ze$,
$\Vert z\Vert_o=\max_{\chi\in\hat{\overline{\Ze}}}|\chi(z)|=\max_{s\in
S}|z(s)|$, showing that each $|z|$ with $z\in\Ze$ has a maximum on
$S$. Because $\overline{S}\subset \TW^{herm}$, we see that
$z^*(s)=\overline{z(s)}$ for each $z\in\Ze$ and $s\in
\overline{S}$. By the Stone-Weierstrass theorem, the restriction
to $S$ of a continuous function $f\in C(\TW)$ can be uniformly
approximated by elements in $\Ze$ considered as functions on $S$.
In other words, there exists a $z\in{\overline{\Ze}}$ such that
$f(\b(\chi))=\chi(z)$ for all $\chi\in\hat{\overline{\Ze}}$. Hence $f\circ
\b$ is continuous on $\hat{\overline{\Ze}}$ for all $f\in C(\TW)$, showing that
$\b$ is continuous and $S=\overline{S}$. Since $\hat{\overline{\Ze}}$ is
compact and $S$ is Hausdorff it follows that $\b:\hat{\overline{\Ze}}\to S$
is a homeomorphism. The proposition now follows from the remark
that $p_z=\b\circ\a$.
\end{proof}
\subsection{Positive elements and positive functionals}
\begin{dfn}
We denote by
$\H_+$
\index{H3@$\H_+$, positive elements in $\H$}
the set of Hermitian elements $h\in\H$ such that
$\forall x\in \H: (hx,x)\geq 0$. We call this the set of positive
elements of $\H$.
\end{dfn}
By spectral theory in the Hilbert completion $\Hf\supset\H$, this
is equivalent to saying that $\lambda(h)\in B(\Hf)$ is Hermitian
and has its spectrum in $\R_{\geq 0}$. Thus $\H_+$ is the
intersection of $\H$ with the usual positive cone $\cf_+$ of the
completion $\cf$. It is clear that for all $x\in\H$, $x^*x\in\H_+$
but not every positive element is of this form.
We write
$\H^{re}$
\index{H3@$\H^{re}$, Hermitian (or real) elements in $\H$}
for the real subspace of Hermitian (or real) elements,
i.e. $h\in\H$ such that $h^*=h$.
\begin{lem}\label{lem:easy}
\begin{enumerate}
\item[(i)] If $z\in\Ze_+,\ h\in\H_+$ then $zh\in\H_+$.
\item[(ii)] If $h\in\H^{re}$ and $A\in\R_+$ such that $A\geq\Vert h\Vert_o$,
then $A.1+h\in\H_+$.
\end{enumerate}
\end{lem}
\begin{proof}
A square root $\sqrt{z}\in{\overline{\Ze}}$, the closure of $\Ze$ in $\cf$,
of $z$ has obviously the property that
$\lambda(\sqrt{z})=\rho(\sqrt{z})$. Hence for every $x\in\H$,
$h\in\H_+$:
\[
(zhx,x)=(h\lambda(\sqrt{z})x,\lambda(\sqrt{z})x)\geq 0.
\]
The second assertion follows since $\operatorname{Spec}(\lambda(h))\subset
[-\Vert h\Vert_o,\Vert h\Vert_o]$.
\end{proof}
\begin{dfn}
We call a linear functional $\chi\in\H^*$ {\it positive} if
$\chi(x)\geq 0$ for all $x\in\H_+$.
\end{dfn}
\begin{cor}\label{cor:contrace}
\begin{enumerate}
\item[(i)]
A positive linear functional $\chi\in\H^*$ extends uniquely to a
continuous functional $\chi\in\cf^\prime$ with norm
$\Vert\chi\Vert=\chi(1)$.
\item[(ii)] The character $\chi_\pi$ of an irreducible representation
$\pi$ of $\cf$ is a positive functional $\chi_\pi\in\cf^\prime$.
\item[(iii)]
An irreducible $*$-representation $\pi$ of $\H$ extends to $\cf$
if and only if its character is a positive functional.
\end{enumerate}
\end{cor}
\begin{proof} (i). By the above Lemma \ref{lem:easy}, $|\chi(x)|\leq
\chi(1)\Vert x\Vert_o$ for all
Hermitian $x\in \H$. In addition, the bitrace $(x,y):=\chi(x^*y)$
is a positive semi-definite Hermitian form, and thus satisfies the
Schwarz inequality. Hence for arbitrary $x\in\H$ we have
$|\chi(x)|^2\leq \chi(x^*x)\chi(1)\leq
\chi(1)^2\Vert x^*x\Vert_o=
\chi(1)^2\Vert x\Vert_o^2$, proving the
continuity of $\chi$.

(ii). Because every irreducible representation $\pi$ of $\cf$ is finite
dimensional (Corollary \ref{cor:typeI}),
it is clear that the character $\chi_\pi(x)$ is a well
defined positive functional on $\cf$. It is continuous by (i).

(iii). If the character $\chi_\pi$ is positive, we have by (i)
that $\chi_\pi$ extends to a finite continuous character of $\cf$.
Because $\cf$ is of finite type I (and thus liminal),
the standard construction in
\cite{dix2}, paragraphe 6.7 shows that there is up to equivalence
a unique irreducible representation $\tilde\pi$ of $\cf$ whose
character is $\chi_\pi$. The converse statement follows by (ii).
\end{proof}
\begin{rem}
In general
not all irreducible $*$-representations of $\H$
extend to $\cf$. See for instance Corollary \ref{cor:onedimlext}.
\end{rem}
\subsection{Casselman's criteria}\label{subsect:cass}
For later use, we discuss in this subsection
a suitable version of
Casselman's criteria (see \cite{cas}, Lemma 4.4.1)
to decide whether a representation of $\H$ is tempered (see the
definition below) or is a subrepresentation of $\hf$. See also
\cite{Mat}.

Recall the norm function $\mathcal{N}$ on $W$ as was introduced in Section
\ref{sect:pre}.
\begin{dfn}
A functional $f\in\H^*$ is called tempered if there exists an
$N\in\mathbb{N}$ and constant $C>0$ such that for all $w\in W$,
\begin{equation}
|f(N_w)|\leq C(1+\nc(w))^N.
\end{equation}
Here $N_w=q(w)^{-1/2}T_w$ are the orthonormal basis elements of $\H$
introduced in (\ref{eq:norm}).
\end{dfn}
Let $(V,\pi)$
be a finite dimensional module  over $\H$, and
let $t\in T$. We define $V_t:=\{v\in V\mid \forall a\in\A\ \exists
n\in\N:(a-a(t))^n(v)=0\}$. The nonzero subspaces of the form
$V_t\subset V$ are called the generalized $\A$-weight spaces of $V$.
We call the corresponding elements $t\in T$ the $\A$-weights of $V$.
\begin{lem}\label{lem:cas}(Casselman's criterion for
temperedness).
The following
statements are equivalent:
\begin{enumerate}
\item[(i)] All matrix coefficients of $\pi$ are tempered.
\item[(ii)] The character $\chi$ of $\pi$ is tempered.
\item[(iii)] The weights $t$ of
the generalized $\A$-weight spaces of $V$ satisfy $|x(t)|\leq 1$,
for all $x\in X^+$.
\end{enumerate}
\end{lem}
\begin{proof}
(i)$\Rightarrow$(ii). This is trivial.

(ii)$\Rightarrow$(iii).
If there exists a weight $t$ of $V$ violating the condition, then
there exists a $x\in X^+$ such that $|x(t)|>1$. We may assume that
$|x(t)|\geq|x(t^\prime)|$ for all weights $t^\prime$ of $V$. By
Lemma 4.4.1 of \cite{cas}, the function
\begin{equation}
f_x(n)=|x(t)|^{-n}
\sum_{t^\prime}\operatorname{dim}(V_{t^\prime})x(t^\prime)^n=
|x(t)|^{-n}\chi(\theta_{nx})
\end{equation}
is not summable on $\mathbb{N}$. Hence
for all $\epsilon>0$, $\chi(\theta_{nx})$ can not be bounded by
a constant times $|x(t)|^{n(1-\epsilon)}$.
On the other hand, suppose
that $\chi$ is tempered. Since $\theta_{nx}=N_{nx}$ and
$\nc(nx)=n\nc(x)$, we obtain that $\chi(\theta_{nx})$ is bounded by
a polynomial in $n$, a contradiction.

(iii)$\Rightarrow$(i)
Recall that the elements $T_u \theta_x T_v$ with $x\in
X^+$, $u,v\in W_0$ span the subspace of $\H$ with basis $N_w$,
where $w$ runs over the double coset $W_0xW_0\subset W$ (see the
proof of Lemma 3.1 of \cite{EO}). It is not difficult to see that
in fact we can write, for $w=uxv\in W_0xW_0$ with $x\in X^+$,
\begin{equation}
N_w=\sum_{v^\prime,u^\prime\in W_0}c_{w,(u^\prime,
v^\prime)}T_{u^\prime}\theta_xT_{v^\prime},
\end{equation}
such that the coefficients $c_{uxv,(u^\prime, v^\prime)}$ and
$c_{uyv,(u^\prime, v^\prime)}$ are equal if $x$ and $y$ belong to
the same facet of the cone $X^+$. Moreover, by the length formula
\cite{EO}, equation 1.1, we have
\begin{equation}
l(x)-|R_{0,+}|\leq l(uxv)\leq l(x)+|R_{0,+}|
\end{equation}
with $l(x)=x(2\rho^\vee)$.
Thus we also have
\begin{equation}
\nc(x)-|R_{0,+}|\leq \nc(uxv)\leq \nc(x)+|R_{0,+}|.
\end{equation}
It therefore suffices to show that the
matrix entries of $\pi(\theta_x)$ with $x\in X^+$ are
polynomially bounded in $\nc(x)=x(2\rho^\vee)+\|x^0\|$.
Since $V$ is a direct
sum of generalized $\A$-weight spaces $V_t$, it is enough to
consider the matrix coefficients of
a generalized $\A$-weight space $V_t$ with weight $t$, satisfying
the condition (iii). Observe finally that it is sufficient to
consider the case that $x=x_0+x^0\in Q+Z_X$, a sublattice in $X$ of
finite index.

By Lie's Theorem we can
put the $\pi(\theta_x)$ simultaneously in upper triangular form.
Choose $\Z_+$-generators $x_1, \dots, x_m$ for the cone $Q^+$,
and a basis $x_{m+1}, \dots, x_n$ for the lattice $Z_X$.
The Jordan decomposition $\pi(\theta_{x_i})=D_iU_i$ gives mutually
commuting matrices $D_i$ and $U_i$, with $D_i$ semisimple and $U_i$
unipotent upper triangular.
By conjugation in the group of invertible upper triangular matrices
we may assume that the commuting semisimple matrices $D_i$ are
diagonal.
The strictly upper
triangular matrices $M_i=\log(U_i)$ are commuting and satisfy
\begin{equation}
\pi(\theta_{x_i})|_{V_t}=x_i(t)\exp(M_i)|_{V_t}.
\end{equation}
Hence for $x=\sum_{i=1}^n l_ix_i$, with $l_i\in\mathbb{Z}_+$ when
$1\leq i\leq m$, we have
\begin{equation}
\pi(\theta_{x})|_{V_t}=x(t)\exp(\sum_{i=1}^n l_iM_i|_{V_t}).
\end{equation}
Since $|x(t)|\leq 1$ by assumption, and the exponential map is
polynomial of degree $N:=\max_t\{\dim(V_t)\}-1$ on the space of strictly
upper triangular matrices commuting with $\pi(\A)$, we
see that the matrix entries are bounded by a polynomial of degree
$N$ in the coefficients $l_i$.
Observe that $x_i(2\rho^\vee)\in\Z_{>0}$ when $1\leq i\leq m$.
Since the
coefficients $l_i$ are nonnegative this implies
that for all $1\leq i\leq m$,
$l_i\leq x(2\rho^\vee)$. On the other hand, there exists a
constant $d>0$ independent of $x$ such that
$|l_i|\leq d\|x^0\|$ for all $i>m$.
Thus there exists a constant $d^\prime$ independent of $x$ such
that $|l_i|\leq d^\prime\nc(x)$ for all $i$.
This gives us the desired estimate of the matrix
entries by a polynomial in $\nc(x)$, of degree
$N$.
\end{proof}
\begin{dfn}
A representation $\pi$ of $\H$ satisfying the above equivalent
conditions is
called a tempered representation of $\H$.
\end{dfn}
Along the same lines one proves:
\begin{lem}\label{lem:casds} (Casselman's criterion for discrete
series representations.)
Let $(V,\pi)$ be a finite dimensional representation of $\H$.
The following are equivalent:
\begin{enumerate}
\item[(i)] $(V,\pi)$ is a subrepresentation of $(\hf,\l)$.
\item[(ii)] All matrix coefficients of $\pi$ belong to $\hf$.
\item[(iii)] The character $\chi$ of $\pi$ belongs to $\hf$.
\item[(iv)] The weights $t\in T$ of the generalized
$\A$-weight spaces of $V$ satisfy: $|x(t)|<1$, for all
$0\not=x\in X^+$.
\item[(v)] $Z_X=\{0\}$, and there exist $C>0,\e>0$ such that
the inequality $|m(N_w)|<C\q^{-\e l(x)}$holds for all
matrix coefficients $m$ of $\pi$.
\end{enumerate}
\end{lem}
\begin{proof}
(i)$\Leftrightarrow$(ii) Let $E$ denote the projector of $\hf$ onto
$V$. Then $E=\rho(e)$ for some idempotent of $\nf$, and since
$E$ is open we have $V=\H e\subset\nf$. Choose an orthonormal basis
$v_i$ of $V$.
The corresponding matrix coefficients $(v_i,xv_j)=(v_iv_j^*,x)$ can
be identified with the elements $v_jv_i^*\in\hf$. Conversely, suppose
that, given a basis $v_i$ of $V$ with dual basis $\phi_j$ of $V^*$,
there exist elements $h_{i,j}\in\hf$ such that for all $x\in\H$,
$\phi_i(\pi(x)v_j)=(h_{i,j}^*,x)$. It follows that for each $i$, the map
$v_j\to h_{i,j}$ defines an embedding of $(V,\pi)$ as a
subrepresentation of $(\hf,\l)$.

(ii)$\Rightarrow$(iii)$\Rightarrow$(iv)$\Rightarrow$(v) This is similar
to the proof of Lemma \ref{lem:cas}. For the last implication we first remark
that (iv) implies that $X^+$ can not contain $-x$ for any $x\in X^+$.
Thus $Z_X=\{0\}$ in this case, hence $\nc(x)=l(x)$.

(v)$\Rightarrow$(ii) The number of elements in $W$ with length $l$ grows
polynomially in $l$. Thus, by the exponential decay, it is clear that
$m=\sum_w\overline{m(N_w^*)}N_w$ is convergent in $\hf$,
and moreover $m(x)=(m^*,x)$.
\end{proof}
\begin{cor}
If the $\A$-weights of a finite dimensional representation $(V,\pi)$
of $\H$ satisfy the condition of Lemma
\ref{lem:casds}(iv), then $V$ carries
a Hilbert structure
such that $\pi$ is a $*$-representation of $\H$, and moreover $\pi$
extends to a representation of $\cf$.
\end{cor}
\begin{dfn}
Irreducible representations of $\H$ satisfying the equivalent conditions of
Lemma \ref{lem:casds}
are called discrete series representations.
\end{dfn}
\subsection{The Plancherel measure}
Since $\cf$ is separable, liminal and unital, the spectrum
$\hat\cf$ is a compact $T_1$ space with countable base. Moreover
it contains an open dense Hausdorff subset.

The algebra $\cf$ comes equipped with the tracial state $\tau$,
defining the representations $\l,\rho:\cf\to B(\hf)$ of
$\cf$ (see Corollary \ref{cor:regrep}).
The general theory of the decomposition of a trace on a
separable, liminal $C^*$-algebras (see \cite{dix2}, paragraphe
8.8), asserts that there exists a unique positive Borel measure
$\mu_{Pl}$
\index{0m@$\mu_{Pl}$, Plancherel measure on $\hat\cf$}
on $\hat\cf$ such that
\begin{equation}\label{eq:dec}
\hf\simeq\int^\oplus_{\hat\cf}\operatorname{End}(V_\pi)d\mu_{Pl}(\pi)
\end{equation} and such that
\begin{equation}\label{eq:chadec}
\tau(h)=\int_{\hat\cf}\chi_\pi(h)d\mu_{Pl}(\pi).
\end{equation}
The measure $\mu_{Pl}$ is called the {\it Plancherel measure}.
\begin{thm}\label{thm:supds}  The
support of $\mu_{Pl}$ is equal to $\hat\cf$.
In addition, an irreducible representation $\pi$ of
$\cf$ is a subrepresentation of $(\hf,\l)$ if and only if
$\mu_{Pl}(\pi)>0$
\end{thm}
\begin{proof}
The representation $\l_\tau=\l$ associated with the state $\tau$
is faithful (see Subsection \ref{sub:iwhehil}). The results thus
follow from Proposition 8.6.8 of \cite{dix2}.
\end{proof}

The center $\zf$
\index{Z3@$\zf$, the center of $\nf$}
of $\nf$ will be mapped
onto the algebra of diagonalizable operators
$L^\infty(\hat\cf, \mu_{Pl})$. This is an isomorphism of
algebras, continuous when we give $\zf$ the weak operator
topology and $L^\infty(\hat\cf, \mu_{Pl})$ the weak topology
of the dual of $L^1(\hat\cf, \mu_{Pl})$. It is an isometry.
\begin{rem}
In many cases there exist non $\mu_{Pl}$-negligible subsets $V$ of
$\Hat\cf$ such that $1<|p_z^{-1}(y)|(<\infty)$ for all $y\in
p_z(V)$. For instance, the affine Hecke algebra of type $G_2$
has two discrete series representations $\pi_3$, $\pi_{2,1}$ associated
with the subregular unipotent class $G_2(a1)$ (notation as in
\cite{C}, Section 13.3)
of the complex algebraic group of
type $G_2$ (also see Appendix \ref{KL}).
Then $p_z(\pi_3)=p_z(\pi_{2,1})$ is equal to the $W_0$-orbit of
the weighted Dynkin diagram associated with $G_2(a1)$.
According to Lemma \ref{lem:casds} and
Theorem \ref{thm:supds}, both $\{\pi_3\}$ and $\{\pi_{2,1}\}$
are non-negligible.

In such case, the above remarks imply in particular
that the weak closure of $\Ze$ in
$\nf$ is strictly smaller than $\zf$, the center of $\nf$.
\end{rem}
\subsection{Outline of the main results}\label{sub:out}
This subsection is a continuation of the outline
given in \ref{subsub:out}.

It is our goal to describe the spectral measure of the tracial
state $\tau$ of $\cf$ explicitly. We will not completely succeed,
as was explained in \ref{subsub:out}(4), but we will obtain a
product formula for the density of $\mu_{Pl}$ on each component of
its support, explicit up to a (intractable) positive real constant
factor. An important intermediate step is the description of the
more easily accessible spectral measure $\nu$ of the restriction
of $\tau$ to the closure $\overline{\Ze}\subset\cf$.
\subsubsection{Plancherel measure $\nu$ of $\Ze$}\label{sl}
\index{S@$S\subset\TW$, image of $p_z$|(}
\index{Tempered coset|(}
\index{Residual coset|(}
\index{0n@$\nu$, Plancherel measure of $\overline{\Ze}$
on $W_0\backslash T$|(}
\index{0n@$\nu_L$, smooth measure on $L^{temp}$ such that $\nu=\sum_L\nu_L$|(}
\index{m@$m_L$, density function of $\nu_L/{\overline \ka}_{W_LL}$|(}
\index{m@$m^L$, quotient $m_L/k_L\nu_{\Ri_L,\{r_L\}}(\{r_L\})$|(}
\index{L@$L^{temp}$, tempered residual coset|(}
\index{r9@$r_L$, element of $L\cap T_L$|(}
\index{0kl@${\overline \ka}_{W_0r}$(=${\overline \ka}_{\Ri,W_0r}$),
rational factor in $\nu(\{r\})$|(}
\index{0k@${\overline \ka}_{W_LL}$, rational factor in
$\nu_L$; average of $\ka_L$|(}
The subalgebra $\Ze\subset \H$ is a $*$-subalgebra.
The spectral measure $\nu$ of the restriction of $\tau$ to
$\overline{\Ze}$ is determined in Subsection \ref{sub:chiA}
by the use of the residue calculus.

A residual coset
$L\subset T$ is a coset of a subtorus of $T$ such that the pole
order of the rational function
\begin{equation}
\frac{1}{c(t,q)c(t^{-1},q)}
\end{equation}
along $L$ is equal to $\operatorname{codim}(L)$. Here $c(t,q)$
denotes Macdonald's $c$-function, see equation (\ref{eq:defc}).
In Appendix \ref{sub:defn} this collection of residual cosets is
carefully introduced, classified and studied.
It turns out to be a finite, $W_0$-invariant collection of
cosets, with various good properties which play an important role
in the calculus of residues
(see Subsection \ref{sub:resiprop} of Appendix).
The residual cosets are of the form
(cf. Proposition \ref{prop:red} and \ref{prop:conv})
$L=r_LT^L$,
where $T^L$ is the connected component of the unit element $e$
in the fix point set in $T$ of the Weyl group $W_L$ of
a parabolic subsystem $R_L$ of $R_0$ (a subtorus of $T$),
and where
$r_L\in T_L=\operatorname{Hom}(X_L,\mathbb{C}^\times)\subset T$
is a residual {\it point} with respect to the
root datum $\Ri_L$ and the restriction $q_L$ of $q$
(see Subsection \ref{sub:rl}). This reduces the classification
of these cosets to the case of the residual points. The {\it tempered
form} $L^{temp}$ of $L=r_LT^L$ is the compact form of $L$ defined
by $L^{temp}:=r_LT^L_u$.

Using the identification of the space of $W_0$-invariant
continuous functions on $T$ and the space of continuous functions
on $W_0\backslash T$, $\nu$ can be viewed as a $W_0$-invariant
measure on $T$ supported on $\cup_{L}L^{temp}$ (union over the
residual cosets).
We show (cf. Theorem \ref{thm:nu}, Proposition \ref{prop:par}
and Theorem \ref{thm:support})
that $\nu=\sum\nu_L$, where the sum is over all
residual cosets $L$, and where $\nu_L$ is a {\it smooth} (with
respect to the Haar measure $d^L$ on $L^{temp}$) measure
with support equal to $L$, such that for all $w\in W_0$,
$\nu_{wL}=w_*\nu_L$ (the push forward of $\nu_L$ along
$w:L^{temp}\to (wL)^{temp}$). There is an explicit
(up to a certain rational constant factor ${\overline \ka}_{W_LL}$)
product formula for
$\nu_L$,
compatible with parabolic induction
(Proposition \ref{prop:par}(iv)):
\begin{equation}\label{eq:firstmL}
\nu_L(r_Lt^L)=k_L\nu_{\Ri_L,\{r_L\}}(\{r_L\})m^L(r_Lt^L)d^Lt^L,
\end{equation}
where $k_L=|K_L|$ with $K_L=T_L\cap T^L$, $m^L$ is the
rational function (\ref{eq:m^L}), and
$\nu_{\Ri_L,\{r_L\}}(\{r_L\})$ is the the mass of the residual
{\it point} $\{r_L\}\subset T_L$ with respect to the $W_L$-invariant
spectral measure $\nu_{\Ri_L}$ on $T_L$ determined by
$(\Ri_L,q_L)$.
In the case where $L=\{r\}$ is a residual point
we have (cf. Theorem \ref{thm:nu})
\begin{equation}\label{eq:a}
\nu_{\{r\}}(\{r\})={\overline \ka}_{W_0r}m_{\{r\}}(r),
\end{equation}
where $m_{\{r\}}$ is the given by the product (\ref{eq:m_L}).

The support of $\nu$ is
by definition equal to the image $S$ of the map $p_z$ (cf.
Proposition \ref{prop:im}). Thus we conclude from the above
description of $\nu$ that
$S=W_0\backslash\cup_{L}L^{temp}$ (union over all
residual cosets $L$ with respect to $\Ri$
and root labels $q$).
\index{S@$S\subset\TW$, image of $p_z$|)}
\index{Tempered coset|)}
\index{Residual coset|)}
\index{0n@$\nu$, Plancherel measure of $\overline{\Ze}$
on $W_0\backslash T$|)}
\index{0n@$\nu_L$, smooth measure on $L^{temp}$ such
that $\nu=\sum_L\nu_L$|)}
\index{m@$m_L$, density function of $\nu_L/{\overline \ka}_{W_LL}$|)}
\index{m@$m^L$, quotient $m_L/k_L\nu_{\Ri_L,\{r_L\}}(\{r_L\})$|)}
\index{L@$L^{temp}$, tempered residual coset|)}
\index{r9@$r_L$, element of $L\cap T_L$|)}
\index{0kl@${\overline \ka}_{W_0r}$(=${\overline \ka}_{\Ri,W_0r}$),
rational factor in $\nu(\{r\})$|)}
\index{0k@${\overline \ka}_{W_LL}$, rational factor in
$\nu_L$; average of $\ka_L$|)}
\subsubsection{Separation by central character}
\index{0w@$\chi_t$, local trace of $\H$,
sum (over $c$) of densities $d(\Yf_c^h)/d\nu$ at $t$|(}
Next we deduce the formula (cf. Corollary \ref{cor:exten})
\begin{equation}\label{eq:int:int}
\tau(hz)=\int_T z(t)\chi_t(h)d\nu(t)
\end{equation}
for all $z\in \Ze$ and $h\in\H$. The
$W_0$-invariant
function $t\to\chi_t\in\H^*$
(defined on the support $\cup_L L^{temp}$ of $\nu$, and extended
to $T$ by $0$) has values in
the positive tracial states of $\cf$.
It follows that
$\chi_t$ is a finite positive linear combination of irreducible
characters of $\cf$, which have central character
$W_0t$ (cf. Definition \ref{dfn:resalg} and Theorem
\ref{thm:mainind}). In other words, the state $\chi_t$ is a
positive linear combination of the irreducible characters of
$\cf$ in the (finite) fiber $p_z^{-1}(W_0t)$.
This decomposition of $\chi_t$ is the subject of Section 4, and
will be described below.
\subsubsection{Generic spectrum and residual algebra}\label{sl1}
\index{0D3@$\Delta_{\Ri,W_0r}$, irreducible discrete
series representations of $\H(\Ri,q)$ with central character $W_0r$|(}
\index{d@$d_\d(=d_{\Ri,\d})$, residual degree; degree of $\d$ in the residual
Hilbert algebra $\overline{\H^r}$|(}
\index{H4@$\H^P=\H(\Ri^P,q^P)$, parabolic subalgebra of $\H$|(}
\index{H4@$\H_L:=\H(\Ri_L,q_L)$, semisimple quotient of $\H^L$|(}
\index{0p@$\pi(\Ri_P,W_Pr,\d,t)=\operatorname{Ind}_{\H^P}^\H(\d_t)$,
parabolically induced representation|(}
\index{0w@$\chi_{\Ri_L,W_Lr_L,\d,t^L}$, character of
the induced representation $\pi(\Ri_L,W_Lr_L,\d,t^L)$|(}
The projection $p_z:\hat\cf\to S=\hat{\overline{\Ze}}$ is
complicated near non-Hausdorff points of $\hat\cf$.
Using a variation of techniques introduced in \cite{Lu} we
define an open dense subset $S^{gen}\subset S$ such that the
restriction of $p_z$ to $p_z^{-1}(S^{gen})$ is a covering map
(cf. Theorem \ref{thm:homeom}).

The absolute continuity of $\nu_L$ with respect to
the Haar measure $d^L$ on $L^{temp}$ enables us to disregard
a set of positive codimension, so that we can restrict the domain
of integration in the above integral to (the pull back to $T$ of)
$S^{gen}$.

Given $t\in T$ such that $W_0t\in S^{gen}$, there exists a unique
residual coset $L$ such that $t\in L^{temp}$. Choose $r_L\in T_L\cap
L$, and write $t=r_Lt^L$ with $t^L\in T^L_u$. The results of
\cite{Lu}, suitably adapted, show that in this
situation there exists a bijective correspondence between the
equivalence classes $[\Delta_{\Ri_L,W_Lr_L}]$ of irreducible
discrete series representations of $\H_L:=\H(\Ri_L,q_L)$
with central
character $W_Lr_L$ and the equivalence classes of irreducible
tempered representations of $\H$ with central character $W_0t$.
The correspondence is established by an inflation to
$\H^L:=\H(\Ri^L,q^L)$ using
the induction parameter $t^L\in T^L_u$, and induction from $\H^L$
to $\H$ (Subsection \ref{sub:lus}, in particular
Corollary \ref{cor:short}).

For $t=r_Lt^L\in L^{temp,gen}$ such that $R_L\subset R_0$
is standard parabolic (cf. Theorem \ref{thm:mainind})
we obtain
\begin{equation}\label{eq:pointt}
\chi_t=|W^L|^{-1}\sum_{\d\in\Delta_{\Ri_L,W_Lr_L}}
d_{\Ri_L,\d}\chi_{\Ri_L,W_Lr_L,\d,t^L},
\end{equation}
where $\chi_{\Ri_L,W_Lr_L,\d,t^L}$ is the character of the
representation $\pi_{\Ri_L,W_Lr_L,\d,t^L}$ which is induced
from the irreducible discrete series module $\d$
of $\H_L$
(with central character
$W_Lr_L$)
with induction parameter $t^L$, and where
$d_{\Ri_L,\d}$ denotes the coefficient of the character $\chi_\d$
in the decomposition of the tracial state $\chi_{\Ri_L,r_L}$ of $\H_L$.

The point is that the coefficients in
(\ref{eq:pointt}) are {\it independent} of the induction parameter
$t^L$. They are certain positive real constants, which we
conjecture to be rational, see Conjecture \ref{rem:ell} below.

The positive real $d_{\Ri,\d}$ is the degree
of $\d$ with respect to the finite dimensional
``residual Hilbert algebra'' $\overline{\H^r}$, the quotient
of $\H$ with respect to the radical of the positive semi-definite
form defined by the tracial state $\chi_r$. Hence we have
(cf. Corollary \ref{cor:fdim})
\begin{equation}\label{eq:b}
\mu_{Pl}(\{\d\})=d_{\Ri,\d}\nu(\{W_0r\}).
\end{equation}
Although Conjecture \ref{rem:ell} is out of reach for the methods
used in this paper, there is a weaker statement which is
relatively easy to prove within the context of this paper, and which is
already useful for certain applications (see Subsection \ref{sub:uni}).
We prove in Section 5 that
the real constants $d_{\Ri,\d}$ are {\it independent} of $\q$, where
we assume that $q(s)$ is written in the form Convention
\ref{eq:scale}.
\index{0D3@$\Delta_{\Ri,W_0r}$, irreducible discrete
series representations of $\H(\Ri,q)$ with central character $W_0r$|)}
\index{d@$d_\d(=d_{\Ri,\d})$, residual degree; degree of $\d$ in the residual
Hilbert algebra $\overline{\H^r}$|)}
\index{0w@$\chi_t$, local trace of $\H$,
sum (over $c$) of densities $d(\Yf_c^h)/d\nu$ at $t$|)}
\index{H4@$\H^P=\H(\Ri^P,q^P)$, parabolic subalgebra of $\H$|)}
\index{H4@$\H_L:=\H(\Ri_L,q_L)$, semisimple quotient of $\H^L$|)}
\index{0p@$\pi(\Ri_P,W_Pr,\d,t)=\operatorname{Ind}_{\H^P}^\H(\d_t)$,
parabolically induced representation|)}
\index{0w@$\chi_{\Ri_L,W_Lr_L,\d,t^L}$, character of
the induced representation $\pi(\Ri_L,W_Lr_L,\d,t^L)$|)}
\subsubsection{Plancherel measure $\mu_{Pl}$ and Fourier
transform}\label{sl2}
\index{W8@$\mathcal{W}$, groupoid whose set of objects is
$\mathcal{P}$, with morphisms
$\operatorname{Hom}_\mathcal{W}(P,Q)=
\mathcal{W}(P,Q):=K_Q\times W(P,Q)$|(}
\index{W6a@$W(\O_1,\O_2)=\{n\in W(P_1,P_2)\mid \exists
k\in K_{P_2}:(k\times n)
\in\W(\d_1,\d_2)\}$|(}
\index{W6b@$W(\O)=W(\O,\O)$|(}
\index{0O@${\Xi}=\Lambda\times_{\Gamma}\Delta$|(}
\index{O@$\tilde\O$, connected component of ${\Xi}$|(}
\index{F@$\F_\H$, Fourier transform on $\H$|(}
Equations (\ref{eq:int:int}) and (\ref{eq:pointt})
yield a decomposition of the trace $\tau$ in terms of an
integral over $t\in S^{reg}$, where the integration kernel is a
sum over the finite fiber $p_z(W_0t)$. This integral can be
rewritten more sensibly as an integral over a space of ``standard
induction data'' $\Xi$, invariant for the action of a groupoid $\W$
acting on the standard induction data (cf. Subsection
\ref{sub:pla}).
Next we decompose $\W\backslash\Xi$ in its connected components. This
leads to the final formulation of the spectral decomposition of
$\mathfrak{H}$ in terms of the Fourier isomorphism $\F$
(cf. Theorem \ref{thm:mainp}). This formulation is parallel to
the formulation of the Harish-Chandra Plancherel formula for
$p$-adic groups, cf. \cite{W2}, \cite{D}.

Let $R_P\subset R_0$ be standard parabolic root system,
and $\d$
an irreducible discrete series representation of
$\H_P=\H(\Ri_P,q_P)$ with central character $W_Pr$.
The space $\O$ of all equivalence classes of ``twists of $\d$'',
representations of $\H^P$ of the form $\d_{t^P}$ where
$t^P$ varies over $T^P_u$, is a compact torus of the form
$K_\d\backslash T^P_u$, where $K_\d$ is the isotropy subgroup
of $[\d]$ in
$K_P=T^P\cap T_P$. There exists a principal fiber bundle $\V_\O$
over $\O$ whose fiber at $\om=(\Ri_P,W_Pr_P,\d,K_\d t^P)\in\O$ is
equal to the common representation space $i(V_\d)$ of the
induced representations $\pi(\Ri_P,W_Pr_P,\d,t^P)$. Thanks to the
regularity of certain intertwining operators (see Subsection
\ref{sub:gene}) there exists a natural action of the group
$W(\O)=\{w\in W_0\mid w(R_{P,+})=
R_{P,+},\ \mathrm{and}\ \exists k\in
K_P:\Psi_w(\d)\simeq\Psi_k(\d)\}$
(where $\Psi_w(\d), \Psi_k(\d)$ denote twists of $\d$ by
automorphisms $\psi_w, \psi_k$ of $\H_P$ induced by $w$ and
$k$ respectively)
on the smooth sections of
$\operatorname{End}(\V_\O)$.

The Fourier transform $\F_\H$ is the algebra homomorphism
from $\H$ into the direct sum of the algebras of smooth,
$W(\O)$-equivariant sections of
$\operatorname{End}(\V_\O)$, defined by
$(\F_\H(h))(\om)=\pi(\Ri_P,W_Pr_P,\d,t^P)(h)$ if
$\om=(\Ri_P,W_Pr_P,\d,K_\d t^P)$.

In this terminology the Plancherel measure can be expressed
as follows (cf. Theorem \ref{thm:mainp}):
\begin{equation}
d\mu_{Pl}(\pi(\om))=
\mu_{\Ri_P,Pl}(\{\d\})|K_P\d|m^P(\om)d^\O\om,
\end{equation}
\index{0m@$\mu_{Pl}$, Plancherel measure on $\hat\cf$}
where $K_P\d$ denotes the set of equivalence classes of
discrete series representations of $\H_P$ in the $K_P$-orbit
of $\d$, $d^\O$ is the normalized Haar measure on $\O$,
$m^P(\om)=m^P(r_Pt^P)$ is as in equation (\ref{eq:firstmL}),
and $\mu_{\Ri_P,Pl}(\{\d\})$ is given by (\ref{eq:b})
(applied to $\Ri_P$).
When we equip the space of smooth sections of
$\operatorname{End}(\V_\O)^{W(\O)}$
with the inner product
\begin{equation}
(f_1,f_2)=\sum_{\O/\sim}|W(\O)|^{-1}\int_{\O}
\operatorname{tr}(f_1(\om)^*f_2(\om))d\mu_{Pl}(\pi(\om)),
\end{equation}
then the Fourier transform $\F_\H$ is an isometry,
which extends uniquely to a unitary isomorphism
$\F:\mathfrak{H}\to
L^2(\operatorname{End}(\V_\Xi))^{\W}$.
\index{W8@$\mathcal{W}$, groupoid whose set of objects is
$\mathcal{P}$, with morphisms
$\operatorname{Hom}_\mathcal{W}(P,Q)=
\mathcal{W}(P,Q):=K_Q\times W(P,Q)$|)}
\index{W6a@$W(\O_1,\O_2)=\{n\in W(P_1,P_2)\mid \exists
k\in K_{P_2}:(k\times n)
\in\W(\d_1,\d_2)\}$|)}
\index{W6b@$W(\O)=W(\O,\O)$|)}
\index{0O@${\Xi}=\Lambda\times_{\Gamma}\Delta$|)}
\index{O@$\tilde\O$, connected component of ${\Xi}$|)}
\index{F@$\F_\H$, Fourier transform on $\H$|)}
\subsubsection{Further remarks}\label{sl3}
\index{0l@$\l_\d(=\l_{\Ri,\d}):=
\overline{\ka}_{W_0r}\mid W_0r\mid d_\d$,
constant factor in $\mu_{Pl}(\{\d\})$,
$\d\in\Delta_{W_0r}$|(}
In \cite{DO} we prove that the Fourier isomorphism
restricts to an isomorphism of the Schwartz completion
$\mathfrak{S}$ of $\H$ (cf. \ref{subsub:schwartz}) onto
$C^\infty(\operatorname{End}(\V_\Xi))^{\W}$. Consequently,
$\cf\simeq\F(\cf)=C(\operatorname{End}(\V_\Xi))^{\W}$. In particular,
the connected components of $\hat\cf$ are the closures
$\hat\cf_\O$ of $\pi(\O^{gen})\subset\hat\cf$.
We expect that these results will provide an effective approach towards
the problem of classification of irreducible tempered modules,
using an analog of
the analytic R-group (see for example \cite{A1}) in our context,
granted the classification of the discrete series.

The methods used in this paper are not suitable to obtain a
parametrization of the finite set of discrete series
representations $[\Delta_{W_0r}]$ with central character $W_0r$. If
all the labels of $\H$ are equal this information is contained in
the work of Kazhdan and Lusztig \cite{KL}. They solved this
problem using equivariant K-theory in the case when the labels
$q_i$ are equal, and $X=P$. This result can be extended to
the general equal label case, see \cite{RR}, \cite{Re2}. In
the appendix Section \ref{KL} one can find an account of the
results of Kazhdan and Lusztig, and the relation with residual
cosets.

Let $F$ be a p-adic field and let $\mathcal{G}$ be the group
of $F$-rational points in an adjoint semisimple group over $F$
which splits over an unramified extension of $F$.
Lusztig \cite{Lu2}, \cite{Lu4} solved the above classification
problem in principle for any Hecke algebra which arises as the
centralizer algebra of a representation of $\mathcal{G}$
which is induced from a cuspidal
unipotent representation of (the Levi quotient of) a parahoric
subgroup of $\mathcal{G}$.
The Hecke algebras that are not dealt with by
Lusztig are ``generic'', and these generic algebras are simpler
with respect to this problem of parametrization.

Starting from
the generic case, Slooten \cite{Slooten} formulated an
interesting combinatorics which (among others) conjecturally
parametrizes the irreducible tempered modules with real central
character for all classical root systems (a generalized Springer
correspondence).

For $\d\in\Delta_{W_0r}$
we define $\l_\d:={\overline \ka}_{W_0r}|W_0r|d_{\d}$, so
that
\begin{equation}
\mu_{Pl}(\{\d\})=\l_\d m_{\{r\}}(r)
\end{equation}
(see (\ref{eq:a}) and
(\ref{eq:b})).
This constant $\l_\d$ has
been computed explicitly by Mark Reeder in the cases where the
Hecke algebra arises as the endomorphism algebra of a
representation of a simple p-adic group of exceptional, split
adjoint type which is induced from a cuspidal unipotent
representation of a parahoric subgroup \cite{Re}. He conjectured
an interpretation (in this situation) of $\l_\d$
(see also \cite{Re0}) in terms of the
Kazhdan-Lusztig parameters of $\d$.
In the exceptional cases he
verified this conjecture, using a formula of Schneider and Stuhler
\cite{schstu} for the formal degree of a discrete series
representation of an almost simple $p$-adic group which contains
fixed vectors for the pro-unipotent radical $U$ of a maximal
compact subgroup $K$. This formula of Schneider and Stuhler
however is an alternating sum of terms which does not
explain the product structure of the formal dimension.
In order to rewrite this sum as a product one needs to
resort to a case-by-case analysis (computer aided)
in \cite{Re}.

The method of \cite{Re} is likely to extend to the general case
(joint work with Mark Reeder and Antony Wasserman, in progress).
This would imply the following conjecture:
\begin{con}\label{rem:ell}
The $d_{\d}$ (equivalently, the $\l_\d$) are rational numbers.
\end{con}
\index{0l@$\l_\d(=\l_{\Ri,\d}):=
\overline{\ka}_{W_0r}\mid W_0r\mid d_\d$,
constant factor in $\mu_{Pl}(\{\d\})$,
$\d\in\Delta_{W_0r}$|)}
\section{Localization of $\tau$ on
$\operatorname{Spec}(\Ze)$}\label{sec:loctau}
Recall the
decomposition of $\tau$ we derived in \cite{EO}, Theorem 3.7:
\begin{gather}\label{eq:basic}
\begin{split}
\tau=\int_{t\in t_0T_u}
\left(\frac{E_t}{q(w_0)\Delta(t)}\right)
\omega\\
\end{split}
\end{gather}
where $\omega$ denotes the rational $(n,0)$-form
\begin{equation}\label{eq:omo}
\omega:=\frac{dt}{c(t,q)c(t^{-1},q)}
\end{equation}
\index{0y@$\omega=\frac{dt}{c(t,q)c(t^{-1},q)}$, $(n,0)$-form on $T$}
on $T$.
Let us briefly review the various ingredients of this formula.
First of all, $T_u=\operatorname{Hom}(X,S^1)$
\index{T3@$T_u=\operatorname{Hom}(X,S^1)$, compact form of $T$},
the compact form of the
algebraic torus $T=\operatorname{Hom}(X,\C^\times)$, and $t_0\in
T_{rs}$
\index{T1@$T_{rs}=\operatorname{Hom}(X,\R^\times_+)$,
real split form of $T$},
the real split part of $T$, and should be deep in the negative
chamber $T_{rs,-}$
\index{T2@$T_{rs,-}$, negative chamber in $T_{rs}$}.
The precise conditions will be formulated
below, see equation (\ref{eq:mininf}).

The form $dt$
\index{d@$dt$, holomorphic extension of Haar measure on $T_u$}
denotes the holomorphic $(n,0)$-form on $T$ which
restricts to the normalized Haar measure on $T_u$. It is given
by the formula
\[
dt:=(2\pi i)^{-n}(x_1x_2\dots x_n)^{-1}dx_1\wedge
dx_2\wedge\dots\wedge dx_n
\]
if $(x_1,\dots,x_n)$ is a basis of $X$.

The function
$\Delta(t):=\prod_{\a\in R_{1,+}}\Delta_\a(t)$
\index{0D@$\Delta$, Weyl denominator} with
\begin{equation}\label{eq:defD}
\Delta_\alpha:=1-\theta_{-\alpha}\in\A
\end{equation}
is the Weyl denominator. Here we use the convention to
consider the subalgebra $\A\subset\H$ as the algebra of
regular functions on $T$ via $\theta_x(t):=x(t)$.

The function
$c(t)=c(t,q)$
\index{c@$c=c(t,q)$, Macdonald's $c$-function}
is Macdonald's
$c$-function.
This $c$-function is introduced as an element of
${}_\F\A$
\index{F@${}_\F\A$, field of fractions of $\A$}, the
field of fractions of $\A$, and will be interpreted as a rational
function on $T$ (cf. \cite{EO}, Definition 1.13). Explicitly, we
put
\begin{equation}\label{eq:defc}
c:=\prod_{\a\in R_{0,+}}c_\a=\prod_{\alpha\in R_{1,+}}c_\a.
\end{equation}
Here we define
$c_\a$
\index{c@$c_\a$, rank one $c$-function}
for $\a\in R_1$ by
\begin{equation}
c_\a:= \frac{(1+q_{\alpha^\vee}^{-1/2}\theta_{-\alpha/2})
(1-q_{\alpha^\vee}^{-1/2}q_{2\alpha^\vee}^{-1}\theta_{-\alpha/2})}
{1-\theta_{-\alpha}}\in{}_\F\A.
\end{equation}
If $\a\in R_0\backslash R_1$ then we define
$c_\a:=c_{2\a}$.

\begin{rem}\label{rem:conv}
We have thus associated a $c$-function $c_\a$ to each root
$\a\in R_{nr}$, but $c_\a$ only depends on the direction of
$\a$. This convention is different from the one used in
\cite{EO}. It is handy to write the formulas for the $c$ functions
in the above form, but strictly
speaking incorrect if $\a\in R_1$ and $\a/2\not\in X$.
However, we formally put $q_{2\alpha^\vee}=1$ if
$\alpha/2\not\in R_0$,
and then rewrite the numerator as
$(1-q_{\a^\vee}^{-1}\theta_{-\alpha})$.
Here and below we use this convention.
\end{rem}

The exact inequalities which have to be met by $t_0\in
T_{rs}$ in order to represent the trace functional $\tau$ are
as follows. If $q(s)>1$ for all $s\in S^{\mathrm{aff}}$,
then according to \cite{EO}, Definition 1.4 and Corollary 3.2,
the representation (\ref{eq:basic}) holds if
\begin{equation}
\forall \a\in F_0:\a(t_0)<q_{\a^\vee}^{-1}q_{\a^\vee/2}^{-1/2}.
\end{equation}
It is clear that this representation of $\tau$ remains true
if we vary the parameters $q$ in a connected open set $U$
such that
$\{q\mid \forall s:q(s)>1\}\subset U\subset \{q\mid \forall s:q(s)>0\}$,
as long as the poles of the kernel of the integral for any $q\in U$
do not intersect the integration cycle $t_0T_u$.
It follows that the representation (\ref{eq:basic})
of $\tau$ holds for any $q$ such that
$\forall s\in S^\mathrm{aff}:q(s)>0$, provided that
\begin{equation}\label{eq:mininf}
\forall \a\in F_0:\a(t_0)<\operatorname{min}
\{(q_{\a^\vee}q_{\a^\vee/2}^{1/2})^{\pm 1},q_{\a^\vee/2}^{\pm 1/2}\}.
\end{equation}
Observe that
\begin{equation}
q_{\a^\vee}q_{\a^\vee/2}^{1/2}=q_{\a^\vee}^{1/2}q_{\a^\vee+1}^{1/2};\
q_{\a^\vee/2}^{1/2}=q_{\a^\vee}^{-1/2}q_{\a^\vee+1}^{1/2}
\end{equation}

The expression
$E_t\in\H^*$
\index{E@$E_t$, Eisenstein functional}
is the holomorphic {\it Eisenstein
series} for $\H$, with the following defining properties (cf.
\cite{EO}, Propositions 2.23 and 2.24):
\begin{gather}\label{fundeis}
\begin{split}
(i)&\ \forall h\in\H,\mbox{ the map }T\ni t\to E_t(h)\mbox{ is regular}.\\
(ii)&\ \forall x, y\in X,h\in\H,\ E_t(\theta_x h\theta_y)=t(x+y)E_t(h).\\
(iii)&\ E_t(1)=q(w_0)\Delta(t).\\
\end{split}
\end{gather}
We want to rewrite the integral (\ref{eq:basic}) representing the
trace functional as an integral over the collection of tempered
residual cosets, by a contour shift. It turns out that such a
representation exists and is unique. To find it, we need an
intermediate step. We will first rewrite the integral as a sum of
integrals over a larger set of tempered ``quasi-residual cosets'',
and then we will show that if we symmetrize the result over $W_0$,
all the contributions of non-residual cosets cancel.
\subsection{$\om$-residual cosets}
The basic scheme to compute residues has nothing to do with the
properties of root systems. It is therefore convenient to
formulate everything in a more general setting first. Later we
will consider the consequences that are specific to our context.

Let $T$ be a complex algebraic torus with character lattice $X$.
\begin{dfn}\label{dfn:om}
Let $\om=pdt/q$ be a rational $(n,0)$-form on $T$. Assume that
$p,q$ are of the form
\begin{equation}
q(t)=\prod_{m\in\M}(d_m^{-1}x_m(t)-1);\ p(t)=\prod_{m^\prime\in
\M^\prime}(d_{m^\prime}^{-1}x_{m^\prime}(t)-1),
\end{equation}
where the products are taken over finite multisets $\M,\M^\prime$. The
multisets $\M$ and $\M^\prime$ come equipped with maps
$m\to(x_m,d_m)\in X\times\C^\times$.
For $m\in\M\cup\M^\prime$ we denote by $L_m\subset T$
the codimension $1$ subvariety $L_m=\{t\mid x_m(t)=d_m\}$, and
we denote by $D_\om$
\index{D@$D_\om$, minus the divisor of $\om$ on $T$}
the divisor
$\sum_{m\in\M}L_m-\sum_{m^\prime\in\M^\prime}L_{m^\prime}$
on $T$ of $q/p$.

An $\om$-residual coset $L$
is a connected component of $\cap_{m\in J}L_m$
for some $J\subset \M$, such that the pole order $i_L$ of $\om$
along $L$ satisfies
\[
i_L:= |\{m\in \M\mid L\subset L_m\}|-|\{m^\prime\in \M^\prime\mid L\subset
L_{m^\prime}\}|\geq\operatorname{codim}(L).
\]
\index{i@$i_L$, pole order along $L$}
The collection of $\om$-residual cosets is denoted by $\L^{\om}$
\index{L1@$\L^\om$, collection of $\om$-residual cosets}.
This is a finite, nonempty collection of cosets of subtori of $T$,
which includes by definition $T$ itself (the empty intersection
of the cosets $L_m$).
\end{dfn}
Note that $\om$ as in the above definition is completely
determined by the divisor $D_\om$ on $T$.

Let $\langle\cdot,\cdot\rangle$
\index{<@$\langle\cdot,\cdot\rangle$!a@rational inner product on $X$,$Y$}
be a rational inner product on the vector space $\Q\otimes Y$,
where $Y$ is the cocharacter lattice of $T$.
This defines an isomorphism between
$\Q\otimes X$ and $\Q\otimes Y$, and we also denote by
$\langle\cdot,\cdot\rangle$ the corresponding inner product on
$\Q\otimes X$.
Through the exponential map
$\exp:\mathfrak{t}_\C:=\C\otimes Y\to T$ we obtain a distance
function on $T$. It is defined by taking the distance between
$2\pi i Y$-orbits in $\mathfrak{t}_\C$.
We denote by
$|t|$
\index{t@$|t|$, distance $d(t,e)$ on $T$}
the distance of $t\in T$ to $e\in T$.

Suppose that $L$ is a connected component of the intersection
$\cap_{m\in J}L_m$ for some subset $J\subset \M$. Then $L$ is
a coset for the connected component of $e$ of the subgroup
$\cap_{m\in J}T^m\subset T$, where $T^m:=\{t\in T\mid x_m(t)=1\}$.
We denote this connected component by
$T^L\subset T$
\index{T@$T^L\subset T$,
algebraic subtorus of which $L$ is a coset}.
Its character
lattice
$X^L:=\operatorname{Hom}(T^L,\C^\times)$
\index{X1b@$X^L$, character lattice of $T^L$}
is equal to the
quotient $X^L=X/((\sum_{m\in J}\Q x_m)\cap X)$.
Let
$X_L$
\index{X1a@$X_L\supset R_L$, lattice of $\Ri_L$,
character lattice of $T_L$}
be the quotient $X_L:=X/(\cap_{m\in J}x_m^{\perp}\cap X)$.
Then $T_L:=\operatorname{Hom}(X_L,\C^\times)$
\index{T@$T_L\subset T$, algebraic subtorus orthogonal to $L$}
is an algebraic
subtorus of $T$, the subtorus ``orthogonal to $T^L$''.
The intersection $K_L:=T_L\cap T^L$
\index{K@$K_L$, finite abelian group $T_L\cap T^L$}
is a finite abelian group,
and is canonically isomorphic to character group of the
quotient $X/(X_L+X^L)$. It follows that $L\cap T_L$ is a coset
for the finite subgroup $K_L\subset T_u$.

We denote by $\M_L\subset \M$ the subset $\{m\in \M\mid
x_m(L)=d_m\}$. We choose an element
$r_L=s_Lc_L\in T_L\cap L$
\index{r9@$r_L$, element of $L\cap T_L$}
for
each $L$ so that we can write $L=r_LT^L$. We call
$c_L\in T_{rs}$
\index{c1@$c_L$, center of $L$}
the center of $L$, and note that this center is determined
uniquely by $L$. We write
$c_L=\exp{\g_L}$
\index{0c@$\g_L\in \mathfrak{t}$, logarithm of $c_L$}
with
$\g_L\in\mathfrak{t}_L$. The set of centers of the $\om$-residual
cosets is denoted by
$\Cc^{\om}$
\index{C@$\Cc^{\om}$, set of centers of $\om$-residual cosets}.
The tempered form of a
$\om$-residual $L=r_LT^L$ is by definition
$L^{temp}:=r_LT^L_u$
\index{L@$L^{temp}$, tempered residual coset}
(which is independent of the choice of $r_L$), and such a coset
will be called an $\om$-tempered coset.

Basically, the only properties of the collection $\L^{\om}$ we
will need are
\begin{prop}\label{prop:triv}
\begin{enumerate}
\item[(i)] If $c\in \Cc^{\om}$ then the union
$$S_c:=\cup_{\{L\in \L^{\om}\mid c_L=c\}}L^{temp}\subset cT_u$$
\index{S@$S_c$, support of $\mathfrak{X}_c$}
is a regular support in the sense of \cite{Schw} in $cT_u$. This
means that a distribution on $cT_u$ with support in $S_c$ can be
written as a sum of derivatives of push forwards of measures on
$S_c\subset cT_u$.
\item[(ii)]
If $c=\exp\g\in T_{rs}$, and $L$ is $\om$-residual with $|\g_L|\geq |\g|$
but $\g_L\not=\g$, then there exists a $m\in \M_L$ such that
$f(t)=x_m(t)-d_m$ is non-vanishing on $cT_u$.
\end{enumerate}
\end{prop}
\begin{proof}
The set $S_c$ is a finite union of smooth varieties, obviously
satisfying the condition of \cite{Schw}, Chapitre III, \S 9 for
regularity, proving (i). As for (ii), first note that the
assumption implies that $\g_L\not=0$, hence that $L\not=T$. Thus
the codimension of $L$ is positive, and $\M_L\not=\emptyset$.
Clearly $\g\not\in \g_L+\mathfrak{t}^L=\log(T_{rs}\cap LT_u)$
since $\g_L$ is the unique smallest vector in this affine linear
space. Because $\{x_m\mid m\in \M_L\}$ spans
$\mathfrak{t}_L=(\mathfrak{t}^L)^\perp$, we can find a $m\in \M_L$
such that $x_m(\g)\not=x_m(\g_L)$. This implies the result.
\end{proof}
\subsection{The contour shift and the local contributions}
The following lemma is essentially the same as Lemma 3.1 of
\cite{HOH0}, but because of its basic importance we have included the
proof here, adapted to the present context.
See also \cite{BS1} for a more general method in the same spirit.
\begin{lem}\label{thm:resbasic}
Let $\om$ be as in Definition \ref{dfn:om} and
let $t_0\in T_{rs}\backslash \cup (T_{rs}\cap T_uL_m)$.
Fix an inner product $\langle\cdot,\cdot\rangle$ on $\Q\otimes Y$.
Then there
exists a unique collection of distributions $\{\mathfrak{X}_c\in
C^{-\infty}(cT_u)\}_{c\in \Cc^{\om}}$
\index{X3@$\mathfrak{X}_c$, local contribution to
$\int_{t_0T_u}a\om$}
such that the following
conditions hold:
\begin{enumerate}
\item[(i)] The support of $\mathfrak{X}_c$ satisfies
$\operatorname{supp}(\mathfrak{X}_c)\subset S_c$ (cf. Proposition
\ref{prop:triv}).
\item[(ii)] For every $a\in \A^{an}(T)$
(the ring of analytic functions on $T$) we have
\begin{equation}
\int_{t\in t_0T_u}a(t)\om(t) =\sum_{c\in\Cc^{\om}}
\mathfrak{X}_c(a|_{cT_u}).
\end{equation}
\end{enumerate}
\end{lem}
\begin{proof}
The existence is proved by induction on the dimension $n$ of $T$,
the case of $n=0$ being trivial. Suppose that the result is true
for tori of dimension $n-1$. Choose a smooth path in $T_{rs}$ from
$t_0$ to the identity $e$ which intersects the real projection
$L_{m,r}=T_{rs}\cap T_uL_m$ of the codimension $1$ cosets $L_m$
transversally and in at most one point $t(L_{m,r})$. We may
assume that the intersection points are mutually distinct
with possible exception for the cases $t(L_{m,r})=e$, i.e.
when $e\in L_{m,r}$.
When we move $t_0$ along the curve towards $e$,
then we pick up residues when we pass at a
point $t=t(L_{m,r})\not=e$ on the curve.
We may assume that the cosets $L_m$ are connected (by
factoring the defining equations, and adapting $\M$
accordingly).
Let $L=L_m$ be such that $t\in L_{m,r}$.
For simplicity of notation we write $(x,d)$
instead of $(x_m,d_m)$.
Recall the decomposition $L=r_LT^L=s_Lc_LT^L$ with $r_L\in T_L$.
Let $d^Lt$
\index{d1@$d^Lt$ ($d_Lt$), holomorphic extension of
normalized Haar measure on $T^L_u$ ($T_{L,u}$)}
denote the holomorphic extension to
$T^L$ of the normalized Haar measure on $T^L_u$,
and similarly for $d_Lt$ on $T_L$.
Let $K_L=T^L\cap T_L$, and let
$k_L$
\index{k_L@$k_L$, order of $K_L$}
be its order.
The product homomorphism
$\pi:T^L\times T_L\to T$ has kernel
$\{(k,k^{-1})\mid k\in K_L\}$.
The residue that is picked up on $L$ when we cross at $t$
can be written as follows:
\[
\int_{t^L\in T^L_u}\sum_{k\in K_L}\int_{t_L\in kC}
(ap/q)(ts_Lt^L t_L)d_L(t_L)d^L(t^L),
\]
where $C$ denotes a small circle in $T_L\simeq\C^\times$
around $1$. Using the action of $\operatorname{Ker}(\pi)$
and the invariance of $d^L$ and $d_L$, and in addition
using $r_L$ as a base point of $L$, this equals
\[
k_L\int_{t^L\in tc_L^{-1}T^L_u}\int_{t_L\in C}
(ap/q)(r_Lt^Lt_L)d_L(t_L)d^L(t^L).
\]

Let $x_L\in X_L$ be a generator of $X_L$. Let $D$ be the holomorphic
constant vector field on $T_L$ which is dual to $x_L$. We extend
$D$ to a constant holomorphic vector field on $T$.
We define a $k_L$-th root of $d$ by $x_L(r_L)=d^{1/k_L}$,
so that the pull back
of $d^{-1}x-1$ to $T_L\times T^L$ factors as follows:
\[
(d^{-1}x-1)=\prod_{k\in K_L}(x_L(k^{-1})d^{-1/k_L}x_L-1):=
\prod_{k\in K_L}(d_k^{-1}x_L-1).
\]
With these notations, the above residue contribution is of the form
\[
\int_{t^\prime\in tc_L^{-1}T^L_u}
\left(B_{i_L-1}(D)
(((d_k^{-1}x_L-1)^{i_L}p/q)a)(r_Lt^Lt_L)\right)|_{t_L=1}d^L(t^L),
\]
where $B_{j}(T)\in\Q[T]$ is a certain polynomial of degree $j$.
Note that there may exist other $L_{m^\prime}$ with $t\in L_{m^\prime,r}$.
We pick up similar residues with respect to these cosets as well
when we cross at $t$.

The above integral can be rearranged as follows
\[
\sum_{j=0}^{i_L-1}\int_{ts_LT^L_u}(D^{j}(a)|_L)\om_{j},
\]
where $\om_{j}$ is itself a rational $(n-1,0)$-form on $L$ which
is a linear combination
\[
\om_j=\sum_i f_{j,i}\om_{j,i}
\]
on $L$ with regular holomorphic
coefficients $f_{j,i}$, and $(n-1,0)$-forms
$\om_{j,i}$ which factor as in Definition \ref{dfn:om}.
The forms $\om_{j,i}$ have poles along
the intersections $L^\prime_n=L\cap L_{n}$ (with $n\in \M$) which
are of codimension $1$ in $L$. A simple computation shows that we
can choose this decomposition of $\om_j$ such that for
every ${j,i}$ and every connected component $H$ of an intersection
of cosets of the form $L_n^\prime\subset L$, the index
$i_{\om_{j,i},H}$ of $H\subset L$ satisfies $i_{\om_{j,i},H}\leq
(i_{\om,H}-1)-j$. It follows that the union over all $j,i$ of the
$\om_{j,i}$-residual cosets in $L$ is contained in the collection
of $\om$-residual cosets of $T$. Moreover, when we take $r_L$ as
a base point of $L$, so that we identify $L$ with $T^L$ through
the map $t^L\to r_Lt^L$, then the tempered form of a
$\om_{j,i}$-residual coset in $L$ is equal to its tempered
form as a $\om$-residual coset in $T$.
By the induction hypotheses we can thus rewrite the residue on $L$
in the desired form, where the role of the identity element in the
coset $L$ is now played by $r_L$.

At the identity $e\in T$ itself we have to
take a boundary value of $\om$ towards $T_u$, which defines a
distribution on $T_u$. This proves the existence.

The uniqueness is proved as follows. Suppose that we have a
collection $\{\mathfrak{Y}_c\in C^{-\infty}(cT_u)\}_{c\in
\Cc^{\om}}$ of distributions such that
\begin{enumerate}
\item[(i)] $\operatorname{supp}(\mathfrak{Y}_c)\subset S_c$.
\item[(ii)] $\forall a\in \C[T]:
\sum_{c\in \Cc^{\om}}\mathfrak{Y}_c(a|_{cT_u})=0$.
\end{enumerate}
We show that $\mathfrak{Y}_c=0$ by induction on $|\g=\log(c)|$. Choose $c\in
\Cc^{\om}$ such that $\Yf_{c^\prime}=0$ for all $c^\prime$ with
$|\g^\prime|<|\g|$. For each $L\in\L^{\om}$ with $|\g_L|\geq|\g|$ and
$\g_L\not=\g$ we choose a $l\in \M_L$ such that $x_l(t)-d_l$ does not
vanish on $cT_u$ (Proposition \ref{prop:triv}) and we set
\[
\nu(t):=\prod_{\{L:|\g_L|\geq|\g|\text{\ and\
}\g_L\not=\g\}}(x_l(t)-d_l).
\]
It is clear that for sufficiently large $N\in \N$,
$\Yf_c(\nu^Na)=0$ for all $a\in \C[T]$. On the other hand, by the
theory of Fourier series of distributions on $T_u$,
$\C[T]|_{cT_u}$ is a dense set of test functions on $cT_u$. Since
$\nu^N$ is nonvanishing on $cT_u$, this function is a unit in the
space of test functions in $cT_u$. Thus also $\nu^N\C[T]|_{cT_u}$
is dense in the space of test functions. It follows that
$\Yf_c=0$.
\end{proof}
\subsubsection{Approximating sequences}\label{subsub:approx}
There is an ``analytically dual'' formulation of the result on
residue distributions that will be useful later on. The idea to
deal with the residue distributions in this way was inspired by
the approach in \cite{H} to prove the positivity of certain
residual spherical functions.
\begin{lemma}\label{lem:approx}
For all $N\in \N$ there exists a collection of sequences
$\{a_n^{N,c}\}_{n\in\N}$ ($c\in \Cc^\om$)  in $\A$ with the
following properties:
\begin{enumerate}
\item[(i)] For all $n\in \N$, $\sum_{c\in \Cc^\om}a^{N,c}_n=1$.
\item[(ii)] For every holomorphic constant coefficient differential operator
$D$ of order at most $N$ on $T$, $D(a^{N,c}_n)\to D(1)$ uniformly
on $S_c$ and $D(a^{N,c}_n)\to 0$ on $S_{c^\prime}$ if
$c^\prime\not=c$.
\end{enumerate}
\end{lemma}
\begin{proof}
We construct the sequences with induction on the norm
$|\g=\log(c)|$. We fix $N$ and suppress it from the notation. Let
$c\in\Cc^\om$ and assume that we have already constructed such
sequences $a^{c^\prime}_n$ satisfying (ii) for all $c^\prime$ with
$|\g^\prime|<|\g|$. Consider the function $\nu$ constructed in
the second part of the proof of Lemma \ref{thm:resbasic}. By
Fourier analysis on $cT_u$ it is clear that there exists a
sequence $\{\phi_n\}_{n\in \N}$ in $\C[T]$ such that for each
holomorphic constant coefficient differential operator $D$ of
order at most $N$ there exists a constant $c_D$ such that
\[
\Vert (D(\phi_n)-D(\nu^{-(N+1)}))|_{cT_u}\Vert_\infty<c_D/n
\]
Applying Leibniz' rule to
$\nu^{(N+1)}\phi_n-1=\nu^{(N+1)}(\phi_n-\nu^{-(N+1)})$ repeatedly
we see that this implies that there exists a constant $c^\prime_D$
for each holomorphic constant coefficient differential operator
$D$, such that
\[
\Vert
(D(\nu^{(N+1)}\phi_n)-D(1))|_{cT_u}\Vert_\infty<c_D^\prime/n.
\]
Notice that $D(\nu^{(N+1)}\phi_n)=0$ on all $S_{c^\prime}$ with
$|\g^\prime|\geq|\g|$ but $\g^\prime\not=\g$. On the other
hand, for each holomorphic constant coefficient differential
operator $E$ the function
$E(1-\sum_{\{c^\prime\mid|\g^\prime|<|\g|\}} a^{c^\prime}_{k})$
converges uniformly to zero on each $S_{c^\prime}$ with
$|\g^\prime|<|\g|$. Again applying Leibniz' rule repeatedly we
see that there exist a $k\in \N$ (depending on $n$) such that the
function
\[
a^{c}_n:=\nu^{(N+1)}\phi_n(1-\sum_{\{c^\prime\mid|\g^\prime|<|\g|\}}
a^{c^\prime}_{k})
\]
has the property that
\[
\Vert D(a^{c}_n)|_{\cup S_{c^\prime}}\Vert_\infty<c^\prime_D/n,
\]
where the union is taken over all $c^\prime$ with
$|\g^\prime|<|\g|$. It is clear that the sequence $a^c_n$ thus
constructed satisfies (ii). We continue this process until we have
only one center $c$ left. For this last center we can simply put
\[
a_n^c:=1-\sum_{c^\prime\not=c}a_n^{c^\prime}.
\]
This satisfies the property (ii), and forces (i) to be valid.
\end{proof}
The use of such collections of sequences is the following:
\begin{prop}
In the situation of Lemma \ref{thm:resbasic} and given any
collection of sequences $\{a_n^c\}$ as constructed in Lemma
\ref{lem:approx} we can express the residue distributions as (with
$a\in\A$):
\[
\Xf_c(a)=\lim_{n\to\infty}\tau(a_n^ca),
\]
provided $N$ (in Lemma \ref{lem:approx}) is chosen sufficiently
large.
\end{prop}
\begin{proof}
Because we are working with distributions on compact spaces, the
orders of the distributions are finite. Take $N$ larger than the
maximum of all orders of the $\Xf_c$. By Proposition
\ref{prop:triv} we can thus express $\Xf_{c^\prime}$ as a sum of
derivatives of order at most $N$ of measures supported on
$S_{c^\prime}$. The result now follows directly from the defining
properties of the sequence $a_n^c$.
\end{proof}
\subsubsection{Cycles of integration}\label{subsub:cycle}
Yet another useful way to express the residue distribution is by
means of integration of $a\om$ over a suitable compact n-cycle.
The results of this subsection will be needed later on to compute
certain residue distributions at ``generic'' points of their
support.

In the proposition below we will use the distance function on $T$
which measures the distance between $2\pi i Y$-orbits in
$\mathfrak{t}_\C$. For $\d>0$ and each $L$ which is a connected
component of an intersection of codimension $1$ cosets $L_m\subset
T$ with $m\in \M$, we denote by $\Bc_L(r_L,\d)$
\index{B@$\Bc_L(r_L,\d)$, ball in $T_L$,
center $r_L$ and radius $\d$}
a ball in $T_L$
with radius $\d$ and center $r_L$, and by
$\Bc^L_{rs}(\d)$
\index{B@$\Bc^L_{rs}(\d)$, ball in $T^L_{rs}$,
radius $\d$ and center $e$}
a ball
with radius $\d$ and center $e$ in $T^L_{rs}$. We put
$\M_L\subset\M$ for the $m\in\M$ such that $L\subset L_m$, and
$\M^L\subset \M$ for the $m\in\M$ such that $L_m\cap L$ has
codimension $1$ in $L$. We write $T^m=\{t\mid x_m(t)=1\}$.

Let $U^L(\d)\subset T^L$ be the open set $\{t\in T^L\mid\forall
m\in\M^L: t\overline{\Bc_L(r_L,\d)}\cap L_m=\emptyset\}$. Note
that $U^L(\d_1)\subset U^L(\d_2)$ if $\d_1>\d_2$, and that the
union of these open sets is equal to the complement of union of
the codimension $1$ subsets $r_L^{-1}(L\cap L_m)\subset T^L$ with
$m\in \M^L$.
\begin{prop}\label{prop:cycle} Let $\e>0$ be such that for all
$m\in\M$ and $L\in\L^\om$,
$L_m\cap\Bc_L(r_L,\e)\Bc^L_{rs}(\e)T_u^L\not=\emptyset$ implies
that $L^{temp}\cap L_m\not=\emptyset$. Denote by $\M^{L,temp}$ the
set of $m\in \M^L$ such that $L^{temp}\cap L_m\not=\emptyset$.
There exist
\begin{enumerate}
\item[(i)] $\forall L\in\L^\om$, a point
$\e^L\in\Bc^L_{rs}(\e)\backslash\cup_{m\in\M^{L,temp}}T^m$,
\item[(ii)] a $0<\d<\e$ such that $\forall L\in\L^\om,\ \e^LT^L_u
\subset U^L(\d)$, and
\item[(iii)] $\forall L\in\L^\om$, a compact cycle
$\xi_L\subset\Bc_L(r_L,\d)\backslash\cup_{m\in\M_L}L_m$
\index{0o@$\xi_L$, compact cycle in
$\Bc_L(r_L,\d)\backslash\cup_{m\in\M_L}L_m$}
of dimension
$\operatorname{dim}_{\mathbb{C}}(T_L)$,
\end{enumerate}
such that
\begin{equation}
\forall c\in\Cc^\om, \forall\phi\in
C^\infty(cT_u):\Xf_c(\phi)=\sum_{\{L\mid c_L=c\}}\Xf_L(\phi),
\end{equation}
where
$\Xf_L$
\index{X4@$\Xf_L$, contribution to
$\int_{t_0T_u}a\om$ supported on $L^{temp}$}
is the distribution on $cT_u$ with support $L^{temp}$
defined by
\begin{equation}
\forall a\in\A:\Xf_L(a)=\int_{\e^LT_u^L\times\xi_L}a\om.
\end{equation}
If $\M^{L,temp}=\emptyset$ we may take $\e^L=e$.
\end{prop}
\begin{proof}
We begin the proof by remarking that (i), (ii) and (iii) imply
that the functional $\Xf_L$ on $\A$ indeed defines a distribution
on $c_LT_u$, supported on $L^{temp}$. Consider for $t\in U^L(\d)$
the inner integral
\begin{equation}\label{eq:inner}
\int_{t\xi}a\om:=i(a,t)d^Lt.
\end{equation}
Then $i(a,t)$ is a linear combination of (possibly higher order)
partial derivatives $D_\ka a$ of $a$ at $r_Lt$ in the direction of
$T_L$,
with coefficients in the ring of
meromorphic functions on $T^L$ which are regular outside the
codimension $1$ intersections $r_L^{-1}(L\cap L_m)$:
\begin{equation}
i(a,t)=\sum_{\ka}f_\ka D_{\ka} a.
\end{equation}
Hence $\Xf_L(a)$ is equal to the sum of the boundary value
distributions $\operatorname{BV}_{\e^L,f_\ka}$ of the meromorphic
coefficient functions, applied to the corresponding partial
derivative $D_\ka a$ of $a$, restricted to $L^{temp}$:
\begin{equation}
\Xf_L(a)=\sum_\ka \operatorname{BV}_{\e^L,f_\ka}(D_\ka
a|_{L^{temp}}).
\end{equation}
We see that $\Xf_L$ is a distribution supported in
$L^{temp}\subset c_LT_u$, which only depends on $\xi_L$ and on the
component of $\Bc^L_{rs}(\e)\backslash\cup_{m\in\M^{L,temp}}T^m$
in which $\e^L$ lies.

Hence, by the uniqueness assertion of Lemma \ref{thm:resbasic},
we conclude that it is sufficient to prove that we can choose
$\e^L,\ \d,\ \xi_L$ in such a way that
\begin{equation}
\forall a\in\A: \int_{t_0T_u}a\om=\sum_{L\in\L^\om}\Xf_L(a).
\end{equation}
In order to prove this it is enough to show that we can choose
$\e^L,\ \d,\ \xi_L$ as in (i), (ii) and (iii) for the larger
collection $\tilde\L^\om$ of all the connected components of
intersections of the $L_m$ (with $m\in\M$), such that
\begin{equation}
t_0T_u\sim\cup_{L\in\tilde\L^\om}\e^LT^L_u\times\xi_L.
\end{equation}
Here $\sim$ means that the left hand side and the right hand side
are homologous cycles in $T\backslash \cup_{m\in\M}L_m$. The
desired result follows from this, since the functional $\Xf_L$ is
equal to $0$ unless $L$ is $\om$-residual (because the inner
integral (\ref{eq:inner}) is identically equal to $0$ for
non-residual intersections, by an elementary argument which is
given in detail in the proof of Theorem \ref{thm:nu}).

Let $k\in\{0,1,\dots,n-1\}$. Denote by $\tilde\L^\om(k)$ the
collection of connected components of intersections of the $L_m$
($m\in\M$) such that $\operatorname{codim}(L)<k$. Assume that we
already have constructed points $\e^L,\ \d,\ \xi_L$ satisfying
(i), (ii) and (iii) for all $L\in\tilde\L^\om(k)$ and in addition,
for each $L\in\tilde\L^{\om}$ with $\operatorname{codim}(L)=k$, a
finite set of points $\Omega_L\subset T^L_{rs}$ such that
$\Omega_L T^L_u\subset U^L(\d)$ and a cycle
$\xi_{L,w}\subset\Bc_L(r_L,\d)\backslash\cup_{m\in\M_L}L_m$ for
each $w\in \Omega_L$, such that $t_0T_u$ is homologous to
\begin{equation}\label{eq:induct}
\cup_{L\in\tilde\L^\om(k)}(\e^LT^L_u\times\xi_L)
\cup\cup_{L\in\tilde\L^\om(k+1)\backslash\tilde\L^\om(k)}
\cup_{w\in\Omega_L}(wT^L_u\times\xi_{L,w}).
\end{equation}
This equation holds for $k=0$, with $\Omega_T=\{t_0\}$, which is
the starting point of the inductive construction to be discussed
below. We will construct $\e^L$, $\d_1$ and $\xi_L$ for
$L\in\tilde\L^\om(k+1)\backslash\tilde\L^\om(k)$, and finite sets
$\Omega_L$ for $L\in\tilde\L^\om(k+2)\backslash\tilde\L^\om(k+1)$,
with a cycle $\xi_w$ for each $w\in\Omega_L$ such that equation
(\ref{eq:induct}) holds with $k$ replaced by $k+1$, and $\d$ by
$\d_1$.

First of all, notice that we may replace $\d$ by any
$0<\d^\prime<\d$ in equation (\ref{eq:induct}), because we can
shrink the $\xi_L$ and $\xi_{L,w}$ within their homology class to
fit in the smaller sets
$\Bc_L(r_L,\d^\prime)\backslash\cup_{m\in\M_L}L_m$. Choose
$\d^\prime$ small enough such that for each $L\in
\tilde\L^\om(k+1)\backslash\tilde\L^\om(k)$ there exists a point
$\e^L\in \Bc^L_{rs}(\e)$ with the property that $\e^LT^L_u\subset
U^L(\d^\prime)$.

The singularities of the inner integral are located at codimension
$1$ cosets in $T^L$ of the form $r_L^{-1}N$, where $N$ is a
connected component of $L\cap L_m$ for some $m\in\M^L$. We have
$r_L^{-1}N=r_L^{-1}r_NT^N\subset T^L$, and thus
$c_L^{-1}c_NT^N_{rs}\subset T^L_{rs}$. Choose paths inside
$T^L_{rs}$ from $w\in\Omega_L$ to the point $\e^L$. We choose each
path such that it intersects the real cosets $c_L^{-1}c_NT^N_{rs}$
transversally and in at most one point, and such that these
intersection points are distinct. If $p=\g(x_0)$ is the
intersection point with the path $\g$ from $w\in\Omega_L$ to
$\e^L$ then $p$ is of the form $p=c_L^{-1}c_Nw_{L,N,w}\in
c_L^{-1}c_NT^N_{rs}$ with $w_{L,N,w}\in T^N_{rs}$. Given $N\in
\tilde\L^\om(k+2)\backslash\tilde\L^\om(k+1)$ we denote by
$\Omega_N$ the set of all $w_{L,N,w}$ arising in this way, for all
the $L\in\tilde\L^\om(k+1)\backslash\tilde\L^\om(k)$ such that
$L\supset N$, and $w\in\Omega_L$.

Notice that if $v=w_{L,N,w}\in\Omega_N$ and $vs\in r_N^{-1}(N\cap
L_m)$ for some $m\in\M^N$ and $s\in T_u$, we have that
$c_L^{-1}c_Nv\in c_L^{-1}(c_NT^N_{rs}\cap
c_{N^\prime}T^{N^\prime}_{rs})$ where $N^\prime=L\cap L_m$. Since
$T^{N^\prime}\not=T^N$, this contradicts the assertion that the
intersection points of the paths in $T^L_{rs}$ and the cosets
$c_L^{-1}c_NT^N_{rs}$ are distinct. We conclude in particular that
the compact set $\Omega_NT^N_u$ is contained in the union of the
open sets $U^N(\d^\prime)$. We can thus choose $\d^\prime$ small
enough such that in fact $\Omega_NT^N_u\subset U^N(\d^\prime)$, as
required in equation (\ref{eq:induct}).

Write $T_{N\subset L}$ for the identity component of the
$1$-dimensional intersection $T_N\cap T^L$, and decompose the
torus $T^L$ as the product $T^N\times T_{N\subset L}$. Let
$v=w_{L,N,w}\in\Omega_N$ and put $p=c_L^{-1}c_Nv$ for the
corresponding intersection point in $T^L_{rs}$. Notice that for a
codimension $1$ coset $r_L^{-1}N^\prime\subset T^L$ with
$N^\prime\in\tilde\L^\om$ we have that
\begin{equation}
pT_{N\subset L,u}\cap r_L^{-1}N^\prime=\left\{
\begin{aligned}{}
&\emptyset\text{\ if\ }
c_L^{-1}c_{N^\prime}T^{N^\prime}\not=c_L^{-1}c_NT^{N},
\\ &G_{L,N^\prime,w}\text{\ otherwise}\\
\end{aligned}
\right.
\end{equation}
where $G_{L,N^\prime,w}$ is a coset of the subgroup $T_{N\subset
L}\cap T^N$ of the finite group $K_{N^\prime}=K_N=T_N\cap
T^N\subset T^N_u$, of the form
\begin{equation}
G_{L,N^\prime,w}=(T_{N\subset L}\cap T^N)r_L^{-1}r_{N^\prime}v.
\end{equation}
The cosets $G_{L,N^\prime,w}$ are disjoint. Let $\d(L,w)$ be the
minimum distance of two points in the union of these cosets, and
let $\d(k+1)$ denote the minimum of the positive real numbers
$\d(L,w)$ when we vary over all the $L$ and $w\in\Omega_L$. Choose
$\d_1>0$ smaller than the minimum of $\d^\prime$ and $\d(k+1)$.
Let $\eta$ be a circle of radius $\d_1/2$ with center $e$ in
$T_{N\subset L}$. Next we make $\d^\prime$ sufficiently small so
that $\cup_{N^\prime}G_{L,N^\prime,w}\eta\subset U^L(\d^\prime)$.
For $x_-,x_+$ suitably close to $x_0$ with $x_-<x_0<x_+$ we have
in $U^L(\d^\prime)$:
\begin{equation}\label{eq:point}
\g(x_-)T_{N\subset L,u}\sim\g(x_+)T_{N\subset L,u}\cup
\cup_{N^\prime}G_{L,N^\prime,w}\eta,
\end{equation}
where the union is over all $N^\prime\subset L$ such that
$c_L^{-1}c_{N^\prime}T^{N^\prime}=c_L^{-1}c_NT^{N}$. Define
\begin{equation}
\xi_{L,N^\prime,v}:=r_L^{-1}r_{N^\prime}\eta\times\xi_{L,w}.
\end{equation}
Observe that $T_{N\subset L,u}\times T^N_u$ is a
$|T_{N\subset L}\cap T^N|$-fold covering of $T^L_u$, and that
$g\eta\times vT^N_u\sim g^\prime\eta\times vT^N_u$ if $g,g^\prime\in
G_{L,N^\prime,w}$. We thus have
\begin{equation}\label{eq:subscycle}
\g(x_-)T_u^L\times \xi_{L,w}\sim\g(x_+)T^L_u\times\xi_{L,w}\cup
\cup_{N^\prime}vT_u^N\times \xi_{L,N^\prime,v}.
\end{equation}
By possibly making $\d^\prime$ smaller we get that
$\xi_{L,N,v}\subset \B_N(r_{N},\d_1)$ for all possible choices $N,
L$ and $w$. If $L_m\supset N$ but $L_m\not\supset L$, then, since
$r_L^{-1}r_{N}\eta\subset U^L(\d^\prime)$ and $\xi_{L,w}\subset
\B_L(r_L,\d^\prime)$, we have $\xi_{L,N,v}\cap L_m=\emptyset$. If
on the other hand $L_m\supset L$ then $\xi_{L,N,v}\cap
L_m=r_L^{-1}r_{N}\eta\times(\xi_{L,w}\cap L_m)=\emptyset$. Finally
we put
\begin{equation}
\xi_{N,v}:=\cup_{(L,w)}\xi_{L,N,v},
\end{equation}
where we take the union over all pairs $(L,w)$ with
$L\in\tilde\L^\om(k+1)\backslash\tilde\L^\om(k)$ such that
$L\supset N$ and $w\in\Omega_L$ such that there is an intersection
point $w_{L,N,w}$ with $w_{L,N,w}=v$. We have shown that
\begin{equation}
\xi_{N,v}\subset\B_N(r_{N},\d_1)\backslash\cup_{m\in\M^{N}}L_m,
\end{equation}
as required in equation (\ref{eq:induct}).

Applying equation (\ref{eq:subscycle}) for all the intersections of
all the paths we chose, we obtain equation (\ref{eq:induct}) with
$k$ replaced by $k+1$ and $\d$ by $\d_1$. We thus take
$\xi_L=\cup_{w\in\Omega_L}\xi_{L,w}$ for
$L\in\tilde\L^\om(k+1)\backslash\tilde\L^\om(k)$, and for
$N\in\tilde\L^\om(k+2)\backslash\tilde\L^\om(k+1)$ we take
$\Omega_N$ and $\xi_{N,v}$ as constructed above.

This process continues until we have $k=n-1$ in equation
(\ref{eq:induct}). In the next step we proceed in the same way.
Notice that for $N\in\tilde\L^\om(n+1)\backslash\tilde\L^\om(n)$,
either $\Omega_N=\{e\}$ (if we cross $c_L^{-1}c_N$ with some curve
from $\Omega_L$ to $\e_L$ in $T^L_{rs}$, for one of the one dimensional
residual cosets $L$ containing $N$), or else $\Omega_N=\emptyset$.
The process now stops, since also $\e^N=e$.
This proves the desired result, with
$\d$ equal to the $\d_1$ obtained in the last step of the
inductive construction.
\end{proof}
\begin{rem} The homology classes of the cycles
$\xi_L$ are not uniquely determined by the above algorithm. The
splitting $\Xf_c=\sum_{\{L\mid c_L=c\}}\Xf_L$ is not unique
without further assumptions. However, in our application to
spectral theory of $\cf$, we shall see that the decomposition
$\Xf_c=\sum_{\{L\mid c_L=c\}}\Xf_L$ is such that each $\Xf_L$
is a regular measure
supported on $L^{temp}$, and such a decomposition is of course
unique.
\end{rem}
We list some useful properties of the cycles $\xi_L$. We fix
$\om$, and suppress it from the notation.
\begin{dfn}
Let $L\in\L$. Denote by
$\L^L$
\index{L@$\L^L$, real projections of residual cosets $\supset L$}
the configuration of real cosets
$M^L:=c_LT^M_{rs}$
where $M\in\L$ such that $M\supset L$, $M\not=
T$. The ``dual'' configuration, consisting of the cosets
$M_L:=c_LT_{M,rs}\subset T_L$
with $M\in\L$ such that
$M\supsetneqq L$, is denoted by
$\L_L$
\index{L@$\L_L$, dual configuration of $\L^L$}.
Given an (open) chamber
$C$ in the complement of $\L^L$, we call $C^d=\{c_L\exp(v)\mid
(v,w)<0\forall w\in\overline{\log(c_L^{-1}C)}\backslash\{0\}\}$
the anti-dual cone. This anti-dual cone is the interior of the
closure of a union of chambers of the dual configuration $\L_L$ in
$T_L$.
We denote by
$\L(L)$
\index{L@$\L(L)$, intersection of $T_L$ with residual
cosets $\supset L$}
the residual cosets in $T_L$
with respect to $K_L$-invariant divisor
$\sum_{m\in\M_L}(L_m\cap T_L)-\sum_{m^\prime\in\M_L^\prime}
(L_{m^\prime}\cap T_L)$ on $T_L$.
\end{dfn}
\begin{prop}\label{prop:t0}
\begin{enumerate}
\item[(i)]
If $t_0$ is moved inside a chamber of $\L^L$ we can leave $\xi_L$
unchanged.
\item[(ii)] Let $t_0(L)=T^L_{rs}t_0\cap T_L$.
For each $k\in K_L:=T^L\cap T_L$,
we can choose the cycle $\xi_{kr_L}(L)$ (defined with respect to
the configuration $\L(L)$ in $T_L$ and initial point $t_0(L)\in T_L$)
equal to $k\xi_L$.
\end{enumerate}
\end{prop}
\begin{proof}
(i) If $t_0$ is moved within a chamber of $\L^L$, the path from $t_0$
to $e$ can be chosen equal to the original path up to a path which
only crosses codimension one cosets of the form $L_mT_u\cap
T_{rs}$ which do not contain $c=c_L$. Therefore this does not
change $\xi_L$.

(ii) We may replace $\M$ by $\M_L$ and $\M^\prime$ by
$\M_L^\prime$ and leave $\xi_L$ is unchanged, because the
$L_m\not\supset L$ do not contribute to $\xi_L$ in the procedure
of the proof of Proposition \ref{prop:cycle}.
By (i) we may also replace $t_0$ by $t_0(L)$
without changing $\xi_L$. We apply Proposition \ref{prop:cycle}
in this situation in $T$. Then we intersect with $T_L$ and use the
formula $T_L\cap(T^L_u\times\xi_L)=\sum_{k\in K_L}k\xi_L$.
\end{proof}
\begin{prop}\label{prop:regequiv}
Write $L=r_LT^L=c_Ls_LT^L$ as usual, and let $M\in\L$. Then
$L^{temp}\subset M^{temp}$ if and only if $L\subset M$ and
$e\in M_L$.
In particular, $L^{temp}$ is maximal in the collection of
$\om$-tempered cosets if and only if $e$ is regular with respect
to the configuration $\L_L$.
\end{prop}
\begin{proof}
If $M\in\L$, then $L^{temp}\subset M^{temp}\Leftrightarrow
L\subset M\mathrm{\ and\ }c_L=c_M$
(since then $s_L\in(c_M^{-1}M)\cap T_u=s_MT_u^M$,
implying that $r_L\in M^{temp}$). Now $c_L=c_M\Leftrightarrow
c_L\in T_M\Leftrightarrow e\in M_L$.
\end{proof}
\begin{prop}\label{prop:antidual}(cf. \cite{HOH0}, Lemma 3.3.)
If $e$ is not in the closure of the anti-dual cone of the chamber
of $\L^L$ in which $t_0$ lies, we can take $\xi_L=\emptyset$.
\end{prop}
\begin{proof}
By Proposition \ref{prop:t0} it is sufficient to show this
in the case where $L=r_L$ is a residual point.

We identify $T_{rs}$ with the real vector space $\mathfrak{t}$ via
the map $t\to \log(c_L^{-1}t)$, and we denote by
$\langle\cdot,\cdot\rangle$
\index{<@$\langle\cdot,\cdot\rangle$!Euclidean inner product on
$T_{rs}$}
the Euclidean inner product thus
obtained on $T_{rs}$. Notice that the role of the origin is played
by $c_L$. The sets $M^L=c_MT^M_{rs}$ with $M\in\L(L)$ satisfy
$c_MT^M_{rs}\ni c_L$, and are equipped with the induced Euclidean
inner product.

By the assumption and Proposition \ref{prop:t0} we can choose
$t_0$ within its chamber such that $\langle t_0, e\rangle>0$.
Assume by induction that in the $k$-th step of the inductive
process of Proposition \ref{prop:cycle} we have, $\forall
N\in\L(L)$ with $\operatorname{codim}(N)=k$ and $\forall
w\in\Omega_N$, that
\begin{equation}\label{eq:posi}
\langle c_Nw,c_N \rangle>0
\end{equation}
(see equation (\ref{eq:induct}) for the meaning of $\Omega_N$).
Notice in particular that this implies that $\Omega_N=\emptyset$
if $c_N=c_L$. By choosing $\e$ sufficiently small, we therefore have
$\langle c_N\e^N,c_N \rangle>0$ when $\Omega_N\not=\emptyset$.
Take the path $\gamma$ in $T^N_{rs}$
from $w\in\Omega_N$ to $\e^N$ equal to
$c_N^{-1}[c_Nw,c_N\e^N]$, where $[c_Nw,c_N\e^N]$ denotes the
(geodesic) segment from $c_Nw$ to $c_N\e^N$ in the Euclidean space
$c_NT_{rs}^N$. Consequently, we have $\langle x,c_N \rangle>0$ for all
$x\in[c_Nw,c_N\e^N]$. Let $M\subset\L(L)$ with
$\operatorname{codim}(M)=k+1$ and $M\subset N$. If $\gamma$
intersects $c_N^{-1}c_MT^M_{rs}$ in $c_N^{-1}c_Mw_{N,M,w}$, then
we have $0<\langle c_Mw_{N,M,w},c_N\rangle=\langle
c_Mw_{N,M,w},c_M\rangle$. By induction on $k$ this proves that we
can perform the contour shifts in such a way that (\ref{eq:posi})
holds for each $k\in\{0,\dots,\operatorname{codim}(L)\}$. This
implies that $\Omega_L=\emptyset$, and thus that
$\xi_L=\emptyset$.
\end{proof}
In the next proposition we view the constants $d_m$ as variables.
We choose a continuous path $[0,1]\ni \sigma\to(d_m(\sigma))_{m\in\M}$
from $(d_m)_{m\in\M}$ to $(d_m^\prime)_{m\in\M}$, and consider the
resulting deformation of $\om$ and $\L$. The end point of the path
corresponds to the form $\om^\prime$ and its collection of
$\om^\prime$-residual cosets, denoted by $\L^\prime$. Recall that
$\M_L$ denotes the multiset of $m\in\M$ such that $L_m\supset L$.
Assume that $\cap_{m\in\M_L}L_m(\sigma)\not=\emptyset$ for all $\sigma$.
In this situation there exists a continuous path $\sigma\to r_L(\sigma)$
such that $L(\sigma):=r_L(\sigma)T^L$ is a connected component of
$\cap_{m\in\M_L}L_m(\sigma)$. We may take $r_L(\sigma)\in T_L\cap L(\sigma)$.
We put $L^\prime=L(1)$. Assume that $\{m\in\M\mid
x_m(L^\prime)=d_m^\prime\}=\{m\in\M\mid x_m(L)=d_m\}$.
\begin{prop}\label{prop:cyinv}
Assume that $e(\sigma):=c_Lc_{L(\sigma)}^{-1}$ stays within a facet of
$\L_L$ for all $\sigma$, and $t_0(\sigma):=t_0c_Lc_{L(\sigma)}^{-1}$ stays
within a chamber of $\L^L$. With these assumptions we can take
$\xi_{L^\prime}=r_L^{-1}r_{L^\prime}\xi_L$.
\end{prop}
\begin{proof}
As above, we may assume that in fact $L=r_L$ is a point.
The only contributions to $\xi_L$ come from contour shifts inside
residual cosets of the configuration $\L(L)$ as in the proof of
Proposition \ref{prop:antidual}. Likewise, for the construction of
$\xi_{L^\prime}$ we only need to consider the translated
configuration $r_L^{-1}r_{L^\prime}\L(L)$. By the assumption on
$t_0$ and Proposition \ref{prop:t0} we can construct
$r_Lr_{L^\prime}^{-1}\xi_{L^\prime}$ by working with $\L(L)$ and
$t_0$, but with the center $e$ of $T$ replaced by
$e^\prime:=e(1)$.

We now follow the deformations of the centers $c_M(\sigma)$ with
$\sigma\in [0,1]$ and $M\in\L(L)$. The assumption on $e(\sigma)$ implies
that $c_M\not=c_L\Leftrightarrow\forall \sigma:\ c_M(\sigma)\not=c_L$.
This implies we can use $\xi_L$ also as the cycle associated with
$L$ relative to the center $e^\prime$.
\end{proof}
\begin{rem}
Note that for some $\sigma\in (0,1)$ there may be additional
$M\in\L(\sigma)$ such that $L(\sigma)\subset M$. It may also happen that
for some values of $\sigma\in [0,1]$, $L(\sigma)^{temp}$ contains smaller
tempered cosets. We may need to adjust $\e^{L(\sigma)}$ accordingly.
\end{rem}
\subsection{Application to the trace functional} We will now
apply the above results to the integral (\ref{eq:basic}). We thus
use the rational $(n,0)$-form
\begin{equation}\label{eq:om}
\eta(t):=\frac{dt}{q(w_0)^2\Delta(t)c(t,q)c(t^{-1},q)}
\end{equation}
and define the notion of {\it quasi}-residual coset
\index{Quasi residual coset}
as the
$\eta$-residual cosets introduced above. We write
$\L^{\operatorname{qu}}$
\index{L1@$\L^{\operatorname{qu}}$, collection of quasi-residual cosets}
for the collection of these
$\eta$-residual spaces,
$\Cc^{\operatorname{qu}}$
\index{C@$\Cc^{\operatorname{qu}}$, centers of quasi-residual cosets}
for their
centers etc. Note: the collection of residual cosets of Appendix
\ref{sub:defn} is {\it strictly} included in this collection.

Apply Lemma \ref{thm:resbasic} to $\eta$ of equation
(\ref{eq:om}), with $t_0$ such that (\ref{eq:mininf}) is satisfied.
Denote the resulting local distributions by $\Xf_{\eta,c}$.
\begin{prop}\label{prop:dfn}
The collection
$\{\Xf_c^h\}_{c\in\Cc^{\operatorname{qu}},h\in\H}$
\index{X5@$\Xf_c^h$, local contribution to $a\to\tau(ah)$ at $c$}
of distributions $\Xf_c^h\in C^{-\infty}(cT_u)$ defined by
$\Xf_c^h(a):=\Xf_{\eta,c}(\{t\to a(t)E_t(h)\})$ satisfies
\begin{enumerate}
\item[(i)] $\operatorname{supp}(\Xf_c^h)\subset S_c^{\qu}$
\index{S@$S_c^{\qu}$, support of $\Xf_c^h$}.
\item[(ii)] $\forall a\in \A:\
\tau(ah)=\sum_{c\in\Cc^{\qu}}\Xf_c^h(a)$ (where $\Xf_c(a)$ means
$\Xf_c(a|_{cT_u})$).
\item[(iii)] The application $h\to\Xf_c^h$ is $\C$-linear.
\item[(iv)] $\forall a,b\in\A,h\in\H:\Xf_c^{ah}(b)=\Xf_c^h(ab)$.
\end{enumerate}
\end{prop}
\begin{proof}
These properties are simple consequences of \ref{fundeis}.
\end{proof}
\subsubsection{Symmetrization and positivity}\label{subsub:pos}
The main objects of this section
are the $W_0$-symmetric versions of the local distributions
$\Xf_c^h$.
\begin{dfn}\label{dfn:Y}
Let $\Cc_-^{\qu}$
\index{C@$\Cc_-^{\qu}$, quasi-residual centers in
$\overline{T_{rs,-}}$}
denote the set of elements in $\Cc^{\qu}$ which
lie in the closure of the negative chamber
$T_{rs,-}=\{t\in T_{rs}\mid
\forall\a\in R_{0,+}:\a(t)< 1\}$.
For $h\in \H$, $a\in \A$, and $c\in\Cc_-^{\qu}$ put:
\begin{equation}
\Yf_c^h(a):=\sum_{c^\prime\in W_0c}\Xf^h_{c^\prime}(\bar{a}),
\end{equation}
\index{Y@$\Yf_c^h$, symmetrized local contribution to
$a\to\tau(ah)$ at $c$}
where $\bar{a}:=|W_0|^{-1}\sum_{w\in W_0}a^w$. Then $\Yf_c^h$ is a
$W_0$-invariant distribution on $\cup_{c^\prime\in W_0c}c^\prime T_u$, with
support in $W_0S_c^{\qu}$, such that for all $z\in \Ze$:
\begin{equation}\label{eq:restz}
\tau(zh)=\sum_{c\in\Cc_-^{\qu}}\Yf_c^h(z).
\end{equation}
\end{dfn}
It is elementary to compute the distribution $\Yf_c^h$ when $c=e$.
Recall that (\cite{EO}, Corollary 2.26) we have the following
identity for the character of the minimal principal series $I_t$:
\begin{equation}
\chi_{I_t}=q(w_0)^{-1}\sum_{w\in W_0}\Delta(wt)^{-1}E_{wt}.
\end{equation}
Hence we can write for all $z\in \Ze$:
\begin{equation}
\begin{split}
\Yf_e^h(z)&=\int_{T_u}z(t)E_t(h)\eta(t)\\
&=\int_{W_0\backslash T_u}z(t)\chi_{I_t}(h)d\mu_T(t),
\end{split}
\end{equation}
where $\mu_T$ is the positive measure on $T_u$ given by
\begin{equation}\label{eq:prin}
d\mu_T(t):=\frac{dt}{q(w_0)c(t)c(t^{-1})}.
\end{equation}
Here we used the $W_0$-invariance of $c(t)c(t^{-1})$, and the fact
that for $t\in T_u$ we have
\begin{equation}
c(t)c(t^{-1})=c(t)c(\bar{t})=|c(t)|^2.
\end{equation}
We see that $h\to\Yf_e^h(1)$ is the integral of the function
$T_u=S_e\ni t\to\chi_{I_t}(h)$ against a positive measure on
$T_u$. Moreover, for every $t\in T_u$, the function
$h\to\chi_{I_t}(h)$ is positive and central, and is a
$\Ze$-eigenfunction with character $t$. Our first task will be to
prove these properties for arbitrary centers $c\in\Cc^{\qu}$. The
main tools we will employ are the approximating sequences.
\subsubsection{Positivity and centrality of the kernel}
Let us choose, for a suitably large $N$, approximating sequences
$a_n^c$ for the distributions $\Xf_{\eta,c}$. We remark that the
group $\pm W_0$ acts on the collection of quasi-residual
subspaces. In addition, complex conjugation also leaves this
collection stable. We define an action $\cdot$ on $\A$ of the
group $G$ of automorphisms of $T$ generated by $W_0$,
$\operatorname{inv}:t\to t^{-1}$ and $\operatorname{conj}:t\to
\bar{t}$. For elements $g\in\pm W_0$ this action is given by
$g\cdot a:=a^g$, and $(\operatorname{conj}\cdot
a)(t):=\overline{a(\overline{t})}$.
\begin{lemma}\label{lem:equiv}
We can choose the $a_n^c$ in a $G$-equivariant way, i.e. such that
$\forall g\in G:\ a_n^{gc}=g\cdot(a_n^c)$.
\end{lemma}
\begin{proof}
Just notice that for any given collection of approximating
sequences $A:=\{a^c_n\}$ and any $g\in G$, $g\cdot A=\{g\cdot
a_n^{g^{-1}c}\}$ is also a collection of approximating sequences
for the distributions $\Xf_{\eta,c}$, and this defines an action
of $G$ on the set of collections of approximating sequences for
the $\Xf_{\eta,c}$. Hence we can take the average over $G$.
\end{proof}
For $c\in\Cc^{\qu}_-$ we now define
\begin{equation}
z^c_n:=\sum_{c^\prime\in W_0c}a_n^{c^\prime}.
\end{equation}
Then these sequences in the center $\Ze$ of $\H$ have the property
that for all $c\in \Cc_-^{\qu}$, $z\in \Ze$ and $h\in\H$:
\begin{equation}
\Yf_c^h(z)=\lim_{n\to\infty}\tau(z_n^czh).
\end{equation}

It is easy to see that the map $h\to\Yf_c^h$ is central:
\begin{prop}\label{prop:cent}
For all $c\in\Cc_-^{\qu}$, we have $\Yf_c^h=0$ if $h$ is a
commutator.
\end{prop}
\begin{proof}
We compute $\Yf_c^h(z)=\lim_{n\to\infty}\tau(z^c_nzh)=0$, because
$z^c_nzh$ is also a commutator and $\tau$ is central.
\end{proof}
We define an anti-holomorphic involutive map
$t\to t^*$
\index{*@$*$!$t\to t^*$, anti-holomorphic involution on $T$}
on $T$ by
$t^*:=\overline{{t}^{-1}}$. In view of the action of conjugation
on $\A$, we see that for all $z\in\Ze$,
$z^*(t)=\overline{z(t^*)}$. By Lemma \ref{lem:equiv} we have, for
each $c\in \Cc^{\qu}_-$,
\begin{equation}
z_n^{c^*}(t)=z_n^c(t^{-1})=(z_n^c)^*(t).
\end{equation}

Now we embark on the proof that the distributions $\Yf_c^h$ are in
fact (complex) measures.
\begin{lem}\label{lem:supp}
\begin{enumerate}
\item[(i)] If $c^*\not\in W_0c$ then $\Yf^h_c=0$.
\item[(ii)] Let $c^*\in W_0c$, and $cs\in S_c^{\qu}$ such that
$(cs)^*=c^{-1}s\not\in W_0(cs)$. Then
$cs\not\in\operatorname{Supp}(\Yf_c^h)$.
\end{enumerate}
\end{lem}
\begin{proof}
(i). Any $h\in\H$ can be decomposed as $h=h_r+ih_i$ with
$h_r^*=h_r$ and $h_i^*=h_i$, so it suffices to prove the assertion
for $h\in \H^{re}$. Thus by Lemma \ref{lem:easy} it is sufficient
to prove the assertion for a positive element $h\in \H_+$.
Similarly $z\in\Ze$ is a linear combination of positive central
elements, so that it is sufficient to show that
$\Yf_c^h(z)=\Yf_c^{zh}(1)=0$ for each positive central element
$z$. By Lemma \ref{lem:easy} this reduces our task to proving that
$\Yf_c^h(1)=0$ for an arbitrary element $h\in\H_+$. Then
\begin{equation}
%\begin{array}{ll}
0\leq\lim_{n\to\infty}\tau(h(z^{c^*}_n+uz^c_n)^*(z^{c^*}_n+uz^c_n))=
u\Yf_c^h(1)+\overline{u}\Yf_{c^*}^h(1).
%\end{array}
\end{equation}
It follows easily that $\Yf_{c^*}^h(1)=\Yf_c^h(1)=0$.

(ii). This is essentially the same argument that we used to prove
(i). Since $(cs)^*\not\in W_0(cs)$, we can find an open
neighborhood $U\ni cs$ in $cT_u$ such that $W_0U\cap
U^*=\emptyset$. Let $\phi\in C_c^{\infty}(W_0U)^{W_0}$. Then
$\phi^*\phi=0$, where $\phi^*(x):=\overline{\phi(x^*)}$. We want
to prove that $\Yf_c^h(\phi)=0$ for $h\in\H_+$. By Fourier
analysis on $cT_u$ we can find a sequence $f_n\in\A^{W_c}$ such
that $D(f_n)$ converges uniformly to $D(\phi)$ on $cT_u$ for every
holomorphic constant coefficient differential operator $D$ on $T$
of order at most $N$ on $T$. We can then find a sequence $g_n$ of
the form $g_n=f_na_{k(n)}^c$ such that $D(g_n)$ converges
uniformly to $D(\phi)$ on $S_c^{\qu}$, and to $0$ on
$S_{c^\prime}^{\qu}$ for every $c^\prime\not=c$. Hence if we put
\[
\phi_n=\sum_{w\in W^c}g_n^{w}\in\Ze,
\]
then for each holomorphic constant coefficient differential
operator $D$ on $T$ of order at most $N$, $D(\phi_n)\to D(\phi)$
uniformly on $W_0S_c^{\qu}$, and $D(\phi_n)\to 0$ uniformly on
$S_{c^\prime}^{\qu}$ for $c^\prime\not\in W_0c$. Hence $\forall
h\in\H_+, u\in\C$,
\begin{equation}
\begin{split}
0&\leq\lim_{n\to\infty}\tau(h(uz_n^c+\phi_n)^*(uz_n^c+\phi_n))\\&=
|u|^2\Yf^h_c(1)+u\Yf^h_c(\phi^*)+\overline{u}\Yf_c^h(\phi)
\end{split}
\end{equation}
If we divide this inequality by $|u|$ and send $|u|$ to $0$, we get that
$\forall\e\in\C$ with $|\e|=1$,
\begin{equation}
0\leq \e\Yf^h_c(\phi^*)+\overline{\e}\Yf_c^h(\phi)
\end{equation}
It follows that $\forall h\in\H_+$, $\Yf_c^h(\phi)=\Yf_c^h(\phi^*)=0$. Hence
the same is true for arbitrary $h\in \H$.
\end{proof}
\begin{cor}\label{cor:pos}
If $h\in\H_+$, the distribution $\Yf_c^h$ is a $W_0$-invariant
positive Radon measure on $W_0cT_u$, supported on $W_0S_c^{\qu}$.
\end{cor}
\begin{proof}
It suffices to show that $\Yf_c^h$ is a positive distribution.
Assume that $\phi\in C^\infty(W_0cT_u)^{W_0}$ and that $\phi>0$.
Then the positive square root $\sqrt{\phi}$ is also in
$C^\infty(W_0cT_u)^{W_0}$. Using the approximating sequences as we
did before, we can find a sequence $f_n\in\Ze$ such that
$D(f_n)\to D(\sqrt{\phi})$, uniformly on $W_0S_c^{\qu}$, and to
$0$ on $S_{c^\prime}^{\qu}$ for $c^{\prime}\not=c$. By Lemma
\ref{lem:supp}, the support of $\Yf_c^h$ is contained in
$W_0S^{\operatorname{herm}}_c:=W_0S_c^{\qu}\cap
T^{\operatorname{herm}}$, where $T^{\operatorname{herm}}:=\{t\in
T\mid t^*\in W_0t\}$. This is itself a regular support for
distributions. On $W_0S^{\operatorname{herm}}_c$, the sequence
$\phi_n:=f_n^*f_n\in\Ze_+$ converges uniformly to $\phi$ up to
derivatives of order $N$. Hence
\begin{equation}
0\leq\lim_{n\to\infty}\tau(hf_n^*f_n)=\Yf_c^h(\phi).
\end{equation}
This proves the desired inequality.
\end{proof}
\begin{cor}\label{cor:cont}
Put $\nu_c:=\Yf_c^1$
\index{0n@$\nu_c=\Yf_c^1$, positive $W_0$-invariant measure on $T$}.
This is a positive Radon measure, with support
in $W_0S_c^{\qu}$, and for all
$h\in\H$, $\Yf^h_c$ is absolutely continuous with respect to $\nu_c$.
\end{cor}
\begin{proof}
It is enough to prove this for $h\in \H$ which are Hermitian, i.e. such that
$h^*=h$. By Lemma \ref{lem:easy} and Corollary \ref{cor:pos} we see that for
positive functions $\phi\in C^\infty(W_0cT_u)^{W_0}$,
\begin{equation}\label{eq:cont}
-\Vert h\Vert_o\nu_c(\phi)\leq\Yf_c^h(\phi)\leq\Vert h\Vert_o\nu_c(\phi).
\end{equation}
\end{proof}
\begin{dfn}\label{dfn:chi}
Let $\nu:=\sum_{c\in\Cc^{\qu}_-}\nu_c$
\index{0n@$\nu$, Plancherel measure of $\overline{\Ze}$
on $W_0\backslash T$}. By equation (\ref{eq:restz}), this
is the spectral measure on $\hat{\overline\Ze}$ of the restriction
of $\tau$ to $\Ze$ (the ``Plancherel measure'' of $\Ze$).
For $h\in\H$ we define a measurable, essentially bounded,
$W_0$-invariant function $t\to\chi_t(h)$
\index{0w@$\chi_t$, local trace of $\H$,
sum (over $c$) of densities $d(\Yf_c^h)/d\nu$ at $t$}
on $T$ by
\begin{equation}
\sum_{c\in\Cc^{qu}_-}\Yf_c^h(\phi|_{W_0S_c^{\qu}})=
\int_T\phi(t)\chi_t(h)d\nu(t)
\end{equation}
for each $\phi\in C_c(T)^{W_0}$. For $t$ outside the support of
$\nu$ we set $\chi_t(h)=0$.
\end{dfn}
\begin{cor}\label{cor:exten}
The function $t\to \chi_t\in\H^*$
\index{0w@$\chi_t$, local trace of $\H$,
sum (over $c$) of densities $d(\Yf_c^h)/d\nu$ at $t$}
satisfies
\begin{enumerate}
\item[(i)] The support of $t\to\chi_t$ is the support of $\nu$.
\item[(ii)] $\chi_t\in\H^*$ is a positive,
central functional such that $\chi_t(1)=1$, $\nu$ almost
everywhere on $T$.
\item[(iii)] For $h\in\H,\ z\in\Ze:\chi_t(zh)=z(t)\chi_t(h)$, $\nu$ almost
everywhere on $T$.
\item[(iv)] $\chi_t$ extends, for $\nu$-almost all $t$,
to a continuous tracial state of the $C^*$-algebra $\cf$.
\item[(v)] We have the following decomposition of $\tau$ for all
$h\in\H$,
\begin{equation}
\tau(h)=\int_T\chi_t(h)d\nu(t).
\end{equation}
\end{enumerate}
\end{cor}
\begin{proof}
Everything is clear. Assertion (iv) follows from Corollary
\ref{cor:contrace} by (ii).
\end{proof}
For $t\in\operatorname{Supp}(\nu)$, we define the positive semi-definite
Hermitian form $(x,y)_t:=\chi_t(x^*y)$ associated to the tracial state
$\chi_t$ of $\H$. It is clear that
the maximal ideal $\mathcal{I}_t\subset \Ze$ of elements vanishing
at $t$ is contained in the radical $\operatorname{Rad}_t$ of
$(\cdot,\cdot)_t$. Hence the radical is a cofinite two-sided ideal
of $\H$. Consequently the GNS-construction produces a finite dimensional
Hilbert algebra associated with $\chi_t$:
\begin{dfn}\label{dfn:resalg}
The algebra
$\overline{\H^t}:=\H/\operatorname{Rad}_t$
\index{H7@$\overline{\H^t}$, residual Hilbert algebra at $t$}
is a
finite dimensional Hilbert algebra with trace $\chi_t$. We will
refer to this Hilbert algebra as the residual algebra at $t$.

Let $\{e_i\}_{i=1}^{l_t}$
\index{e1@$e_i$, minimal central idempotent of $\overline{\H^t}$}
denote the set of minimal central
idempotents of $\overline{\H^t}$, and
$\chi_{t,i}$
\index{0w@$\chi_{t,i}$, irreducible character of $\overline{\H^t}$}
the associated
irreducible characters given by
\begin{equation}
\chi_{t,i}(x)=\operatorname{dim}(e_i\overline{\H^t})^{1/2}
\chi_t(e_i)^{-1}\chi_t(e_ix)
\end{equation}
We define $d_{t,i}:=\operatorname{dim}(e_i\overline{\H^t})^{-1/2}
\chi_t(e_i)\in\R_+$
\index{d@$d_{t,i}(=d_{W_0t,i})$, residual degree; degree of $\chi_{t,i}$ in
$\overline{\H^t}$}
\index{d@$d_{t,i}(=d_{W_0t,i})$, residual degree; degree of $\chi_{t,i}$ in
$\overline{\H^t}$|see{$d_\d$}},
so that
\begin{equation}
\chi_t=\sum_{i=1}^{l_t}\chi_{t,i}d_{t,i}.
\end{equation}
Note that everything in sight depends on the orbit $W_0t$ rather
than on $t$ itself. We will sometimes use the notation
$d_{W_0t,i}$ etc. in order to stress this.
(This notation and parametrization for the irreducible
characters of $\overline{\H^t}$ is provisional. We return to these
matters in a systematic way in Section \ref{sect:loc}
(see e.g. Theorem \ref{thm:mainind}).)
\end{dfn}
\subsection{The Plancherel measure $\nu$ of $\Ze$, and the $\A$-weights of
$\chi_t$}\label{sub:chiA}
The results in this subsection are based on the fact that the
Eisenstein kernel of (\ref{eq:basic}) simplifies considerably when
restricted to the subalgebra $\A\supset\Ze$ of $\H$.
This means that the $(n,0)$-form $\eta$ (see (\ref{eq:om})) can be
replaced by the better behaved $(n,0)$-form $\omega$ (cf.
(\ref{eq:omo}) and subsection \ref{sub:quick}) in the residue calculus.
This has as an
important consequence (see below) that the support of
the measure $\nu$
can be identified as the union of the tempered residual cosets,
which is only a small subcollection
of the tempered quasi residual cosets, and very well behaved
(see Subsection \ref{sub:resiprop} of Appendix \ref{sub:defn}).
Since we have derived that $\Yf^h$ is absolutely continuous with
respect to $\nu$ for general $h\in\H$(see Corollary \ref{cor:cont}),
we conclude that the support of the density function $t\to\chi_t$
is the union of the tempered residual cosets.

The  probability measure $\nu$ can be computed almost explicitly,
due to the
good properties of residual cosets.
We will exploit these
facts here to study the behaviour of the states $\chi_t$ on $\A$.
\begin{thm}\label{thm:nu}
The $W_0$-invariant probability measure $\nu$
has a decomposition $\nu=\sum_{L}\nu_L$,
where $L$ runs over the collection of residual cosets
as defined in Appendix \ref{sub:defn}, and
where $\nu_L$ is the push forward
to $T$ of a smooth measure on $L^{temp}$.
Let $d^L$ denote the
normalized Haar measure on $T^L_u$, transported to the coset
$L^{temp}$ by translation. The measure
$\nu_L$
\index{0n@$\nu_L$, smooth measure on $L^{temp}$ such that $\nu=\sum_L\nu_L$}
is given by a density function
${\overline \ka}_{W_LL}m_L(t):=\frac{d\nu_L(t)}{d^Lt}$
\index{m@$m_L$, density function of $\nu_L/{\overline \ka}_{W_LL}$},
where ${\overline \ka}_{W_LL}\in\Q$
is a constant, and where $m_L$ is of the form
\begin{equation}\label{eq:m_L}
m_L(t)=q(w_0)\frac{\prod^\prime_{\a\in
R_1}(\alpha(t)-1)}{\prod^\prime_{\a\in
R_1}(q_{\alpha^\vee}^{1/2}{\alpha(t)^{1/2}}+1)
\prod^\prime_{\a\in R_1}(q_{\alpha^\vee}^{1/2}
q_{2\alpha^\vee}\alpha(t)^{1/2}-1)}.
\end{equation}
Here we used the convention of Remark \ref{rem:conv}. The constant
${\overline \ka}_{W_LL}$ is independent of $\q$ if we assume $q$ to be as in
Convention \ref{eq:scale}. The notation $\prod^\prime$ means that we omit the
factors which are identically equal to $0$
on $L$. The density
$m_L$ is a smooth function on $L^{temp}$.
\end{thm}
\begin{proof}
We know already that $\nu$ is a $W_0$-invariant measure supported on
the union of the tempered quasi residual cosets.
We apply Proposition \ref{prop:cycle} to the integral
\[
\tau(a)=\int_{t_0T_u}a\omega
=\sum_{c\in\Cc^{\qu}}\Xf_c^1(a)
\]
(cf. equations (\ref{eq:basic}), (\ref{eq:omo}) and \ref{fundeis}).
Choose $\e>0$. For a suitably small $\d>0$ we can find, for each
quasi residual subspace $L$, an $\e^L\in T_{rs}^L$ in an $\e$
neighborhood of $e$, and a cycle
$\xi_L\subset\Bc_L(r_L,\d)\backslash\cup_{L_m^\prime\supset
L}L_m^\prime$, where $\Bc_L(r_L,\d)\subset T_L$ denotes a ball of
radius $\d>0$ centered around $r_L$, such that
\begin{equation}\label{eq:roughly}
\Xf_c^1(a):=\sum_{L:c_L=c}
k_L\int_{t\in\e^LT^L_u}\left\{\int_{\xi_L}a(tt^\prime)\frac{d_L(t^\prime)}
{q(w_0)c(tt^\prime)c({(tt^\prime)}^{-1})}\right\}d^L(t).
\end{equation}
Here $d^L(t)$ is the holomorphic extension to $L$ of $d^L$, and
$d_L{t^\prime}$ denotes the Haar measure on $T_{L,u}$, also
extended as a holomorphic form on $T_L$. We assume that $\d$ is
small enough to assure that $\log$ is well defined on
$\Bc_L(r_L,\d)$. For the inner integral we use a basis $(x_i)$
of $X\cap\Q R_L$ as coordinates on $\log{(\Bc_L(r_L,\d))}$,
shifted so that the coordinates are
centered at $\log(r_L)$. We can then
write the integration kernel as:
\begin{equation}
t^\prime\to
a(tt^\prime)m_L(tt^\prime)(1+f_t(t^\prime))\om_L(t^\prime)
\end{equation}
where $\om_L$ is a rational homogeneous $(l:=\dim(T_L),0)$-form
(independent of $t$) in the $x_i$,
and $f_t$ is a power series in $x_i$ such
that $f_t(0)=0$. In fact, the form $\om_L$ is easily seen to be
(including the factor $k_L$ of (\ref{eq:roughly}))
\begin{equation}\label{eq:tochhandig}
\om_L(x)=\frac{\prod_{\a\in
R_{L}^z}\a(x)}{(2\pi i)^l
\prod_{\b\in R_L^p}\b(x)}dx_1\wedge dx_2\dots\wedge dx_l.
\end{equation}
By Corollary \ref{cor:simpledefres} it follows that the form
$\om_L$ has homogeneous degree $\geq 0$ if $L$ is residual in the
sense of Definition \ref{dfn:ressub}. A homogeneous closed
rational form of positive homogeneous degree is exact. Hence the
inner integral will be nonzero only if $L$ is in fact a residual
coset. In that case the inner integral will have value
\begin{equation}\label{eq:inint}
\ka_La(r_Lt)m_L(r_Lt)
\end{equation}
with
\begin{equation}
\ka_L=\int_{\xi_L}\omega_L.
\end{equation}
\index{0k@$\ka_L$, rational number $\int_{\xi_L}\omega_L$}
We note that $\ka_L\in\Q$, since $\om_L$ defines a rational
cohomology class. Let us therefore assume that $L$ is residual
from now on. Write $r_L=sc$. By Theorem \ref{thm:ster} we know
that $r_L^*=sc^{-1}=w_s(r_L)$ with $w_s\in W(R_{L,s,1})$. When
$t\in L^{temp}$, the expression $m_L(t)$ can be rewritten as
\begin{equation}\label{eq:ml}
q(w^L)m_{\Ri_L,\{r_L\}}(r_L)\prod_{\a\in R_{1,+}\backslash
R_{L,1,+}}\frac{|1-\a(t)|^2}{|1+q_{\a^\vee}^{1/2}\a(t)^{1/2}|^2
|1-q_{\a^\vee}^{1/2}q_{2\a^\vee}\a(t)^{1/2}|^2}.
\end{equation}
Here we used that if $t=cu\in L^{temp}$ with $u\in sT^L_u$, we
have $w_sc=c^{-1}$, $w_su=u$, and $w_s(R_{1,+}\backslash
R_{L,1,+})=R_{1,+}\backslash R_{L,1,+}$. By the same argument as
was used in Theorem 3.13 of \cite{HOH0} we see that this
expression is real analytic on $L^{temp}$. This implies that we
can in fact take $\e^L=e$ for all residual $L$ in equation
(\ref{eq:roughly}) after we evaluate the inner integrals. This leads
to
\begin{equation}
\Xf_c^1(a)=\sum_{L:c_L=c}\ka_L\int_{L^{temp}}a(t)m_L(t)d^L(t).
\end{equation}
where the sum is taken over residual cosets only. When we combine
terms over $W_0$ orbits of residual cosets we find the desired
result. Let $W_0L$ denote the set of residual cosets in the orbit
of $L$. We have to take
\begin{equation}\label{eq:denk}
{\overline \ka}_{W_LL}=\frac{1}{|W_0L|}\sum_{L^\prime\in W_0L}\ka_{L^\prime}.
\end{equation}
\index{0k@${\overline \ka}_{W_LL}$, rational factor in
$\nu_L$; average of $\ka_L$}
When we now define a measure $\nu_L$ on $L^{temp}$ by
\begin{equation}
\int_{t\in L^{temp}}f(t)d\nu_L(t):={\overline \ka}_{W_LL}\int_{t^L\in
T^L_u}f(r_Lt^L)m_L(r_Lt^L)d^L(t^L)
\end{equation}
then we have the equality $\nu=\sum_L\nu_L$ (sum over the residual subspaces).

We note in addition that
$\ka_L=k_L\ka_{\Ri_{L},\{r_L\}}$
\index{0k@$\ka_{\Ri_L,\{r_L\}}$|see{$\ka_L$}},
because the cycle $\xi_L$ is constructed inside $T_L$, entirely in
terms of the root system $R_{L}$
(see Proposition \ref{prop:t0}) (the factor $k_L$ comes from
the facorization $dt=k_Ld^Lt^Ld_Lt_L$, see (\ref{eq:roughly})).
Also, it is clear that
$m_{\Ri_L,\{r_L\}}(r_L)$ is independent of the choice of $r_L$,
because the finite group $K_L=T_L\cap T^L$ is contained in the
simultaneous kernel of the roots of $R_{L}$. Finally, the
independence of $\q$ is clear from Proposition \ref{prop:cyinv}.
When we apply a scaling transformation $\q\to\q^\e$, the point
$c_L$ moves such that the facet of the dual configuration
containing $e$ does not change. Hence $r_L^{-1}\xi_L$ and $\om_L$
will be independent of $\e$.
\end{proof}
\begin{rem}\label{rem:smoothnest}
We note that the smoothness of $m_L$ implies Theorem
\ref{thm:nonnest}, similar to \cite{HOH0}, Remark 3.14.
\end{rem}
\begin{prop}\label{prop:par}
For $L$ residual consider the root datum
$\Ri_L=(X_L,Y_L,R_L,R^\vee_L,F_L)$ (see Subsection \ref{sub:par})
associated with the parabolic root subsystem $R_L\subset R_0$.
Let $q_L$ be the restriction of the label function $q$ to $\Ri_L$.
Then $\{r_L\}\subset T_L$ is a $(\Ri_L,q_L)$ residual point.
Assume that $R_L$ is a standard parabolic subsystem of roots, and
thus that $F_L\subset F_0$. Denote by $W_L$ the standard parabolic
subgroup $W_L=W(R_L)$ of $W_0$, and let
$W^L$
denote the set of
minimal length representatives of the left $W_L$ cosets in $W_0$.
\begin{enumerate}
\item[(i)] When $w\in W^L$, we may take $\xi_{wL}=w(\xi_L)$.
Consequently, $\ka_L=\ka_{wL}$ if $w\in W^L$.
\item[(ii)] Put
\begin{align}\label{eq:m^L}
m^L(t)&=q(w^L)^{-1}\prod_{\a\in R_1\backslash R_{L,1}}c_\a(t)^{-1}\\
&=q(w^L)\prod_{\a\in R_{1,+}\backslash
R_{L,1,+}}\frac{|1-\a(t)|^2}{|1+q_{\a^\vee}^{1/2}\a(t)^{1/2}|^2
|1-q_{\a^\vee}^{1/2}q_{2\a^\vee}\a(t)^{1/2}|^2}.
\end{align}
\index{m@$m^L$, quotient $m_L/k_L\nu_{\Ri_L,\{r_L\}}(\{r_L\})$}
Then $m^L$ and $m_L$ are $\operatorname{Aut}(W_0)$-equivariant,
i.e. $m^L(t)=m^{gL}(gt)$ and $m_L(t)=m_{gL}(gt)$ for every
$g\in\operatorname{Aut}(W_0)$. In particular, $m^L$ and $m_L$ are
invariant for the stabilizer $N_L$ of $L$ in $W_0$.
\item[(iii)] We have $\ka_L=k_L\ka_{\Ri_L,\{r_L\}}$,
${\overline \ka}_{W_LL}=k_L{\overline \ka}_{\Ri_L,W_Lr_L}$.
\item[(iv)] For $z\in\Ze$, we have
\begin{equation}
\frac{1}{|W_0L|}\int_T zd\nu_{W_0L}=k_L\nu_{\Ri_L,\{r_L\}}(\{r_L\})
\int_{L^{temp}}z(t)m^L(t) d^L(t).
\end{equation}
\item[(v)] Assuming that $q$ is expressed as in Convention
\ref{eq:scale} with $f_s\in2\mathbb{Z}$. Then
$\nu_{\Ri_L,\{r_L\}}(\{r_L\})=
{\overline \ka}_{\Ri_L,W_Lr_L}m_{\Ri_L,\{r_L\}}(r_L)$
is of the form $d\q^nf(\q)$, where $d\in\mathbb{Q}$,
$n\in\Z$, and where $f$ is a quotient of products of cyclotomic
polynomials in $\q$.
\end{enumerate}
\end{prop}
\begin{proof}
(i). We note that for $w\in W^L$, $t_0$ and $w^{-1}t_0$ are in the
same chamber of $\L^L$. Hence, by application of
Proposition \ref{prop:t0}, we
may replace $\xi_{wL}$ by $w(\xi_L)$.

(ii). This is trivial.

(iii). The formula $\ka_L=k_L\ka_{\Ri_L,\{r_L\}}$ was explained
in the proof of Theorem \ref{thm:nu}.
Let $W_L=W(R_L)$. Let $N_{T^L}$ be the stabilizer of $T^L$
in $W_0$. Observe that $N_L\subset N_{T^L}$ and
$W_L\vartriangleleft N_{T^L}$. If we define $\Gamma_L=N_{T^L}\cap
W^L$ then $\Gamma_L$ is a complementary subgroup of $W_L$ in
$N_{T^L}$. Using (i), (ii) and the remark
$\ka_L=k_L\ka_{\Ri_L,\{r_L\}}$ we see that
\begin{equation}
\begin{split}
{\overline \ka}_{W_LL}&=
\frac{1}{|W_0L|}\sum_{L^\prime\in W_0L}\ka_{L^\prime}\\
&=\frac{|W_0T_L|}{|W_0L|}\sum_{L^\prime\in
N_{T^L}L}\ka_{L^\prime}\\
&=\frac{|W_0T_L||N_{T^L}L|}{|W_0L|}
k_L{\overline \ka}_{\Ri_L,W_Lr_L}=
k_L{\overline \ka}_{\Ri_L,W_Lr_L}
\end{split}
\end{equation}
\index{0kl@${\overline \ka}_{W_0r}$(=${\overline \ka}_{\Ri,W_0r}$),
rational factor in $\nu(\{r\})$}
Using Theorem \ref{thm:nu} and equation (\ref{eq:ml}) the result
follows.

(iv). Follows easily from (iii).

(v). Since equation (\ref{eq:m_L}) involves only roots in $R_0$,
it is sufficient to prove the statement for $R_0$ indecomposable
and $X=Q$.
Notice that for all $\a\in R_0$, $\a(s)$ is a root of unity and,
by Theorem \ref{thm:ster}(iii), $\a(c)$ is an integral power of $\q$.
Looking at the explicit formula (\ref{eq:m_L}), we see that
it remains to show that this expression has rational coefficients
if $L=r=sc$ is a residual point.
Let $k$ be the extension of $\mathbb{Q}$ by the values
of $\a(s)$, where $\a$ runs over $R_0$. In the case where $\mathcal{R}$ is of
type $C_n^{\text{aff}}$ it follows by Lemma \ref{lem:order2} that $k=\mathbb{Q}$,
and we are done. For the other classical cases it follows from the result of
Borel and de Siebenthal \cite{BS} that the order of $s$ is at most two, and
hence that $k=\mathbb{Q}$.
Next let $\mathcal{R}$ be of exceptional type, and
$\sigma\in\operatorname{Gal}(k/\mathbb{Q})$. Define a character $\sigma(s)$
of $X=Q$ by $Q\ni x\to \sigma(x(s))=:x(\sigma(s))$. By Lemma \ref{lem:conj}
we see that there exists a $w_1\in W_0$ such that $\sigma(s)=w_1s$. Moreover,
$w_1:R_{s,0}\to R_{s,0}$ acts as an automorphism and $c$ is an $R_{s,0}$-residual
point. If $F_{s,0}$ contains isomorphic components then these are of type $A$,
which has only one real residual point up to the action of $W(R_{s,0})$.
Hence by Theorem \ref{thm:ster}(i), there exists a $w_2\in W(R_{s,0})$
such that $w_1(c)=w_2(c)$. Put $w=w_2^{-1}w_1$, so that $wr=c\sigma(s)$.
By the $W_0$-equivariance of $m_{\{r\}}(r)$ we see that (with the
action of $\sigma$ being extended to $k[\q,\q^{-1}]$ by its
action on the coefficients)
$\sigma(m_{\{r\}}(r))=m_{\{wr\}}(wr)=m_{\{r\}}(r)$, whence the desired
rationality.
\end{proof}
The next proposition is a direct consequence of (the proof of)
Theorem \ref{thm:nu} and the definition of $\chi_t$.
\begin{prop}\label{chiA}
Let $r=sc\in T$ be a residual point, and let $a\in \A$. Then
\begin{equation}
\nu(W_0r)\chi_r(a)=m_{\{r\}}(r)\sum_{r^\prime\in
W_0r}\ka_{\{r^\prime\}}a(r^\prime).
\end{equation}
\end{prop}
\begin{theorem}\label{thm:support}
The support of $\nu$ is exactly equal to the union of the tempered
residual cosets. In other words, $S=W_0\backslash\cup_{L\
\mathrm{residual}}L^{temp}$.
\end{theorem}
\begin{proof}
The equality $S=W_0\backslash\operatorname{Supp}(\nu)$ was
explained in \ref{sl}, so it suffices to show that the
support of $\nu$ is equal to the union of the tempered residual cosets.
By Theorem \ref{thm:nu} we know that $\nu$ is supported on this
set, so we need only to show that $W_0L^{temp}$
is contained in the support for each tempered residual coset $L$.

By Proposition \ref{prop:par} this reduces to the case of a
residual point $r=sc$. By Proposition \ref{chiA} it is enough to
show that there exists at least one $r^\prime=wr\in W_0r$ such
that $\ka_{\{r^\prime\}}\not=0$. In other words, using
Proposition \ref{chiA} we single out the point residue at
$r^\prime$. In particular, we ignore all residues at residual
cosets which do not contain $r^\prime$ and thus do not contribute
to $\ka_{\{r^\prime\}}$ in the argument below.

By the $W_0$-invariance of $\om$, we can formulate the problem as
follows. Recall from the proof of Theorem \ref{thm:nu} that
\begin{equation}
\ka_{\{r\}}m_{\{r\}}(r)=\int_{\xi}\om,
\end{equation}
where $\xi$ is the residue cycle at $r$, which is obtained from
Proposition \ref{prop:cycle} if we use the $n$-form
\begin{equation}
\om(t)=\frac{dt}{c(t)c(t^{-1})}
\end{equation}
and a base point $t_0\in T_{rs}$ such that $\forall \a_i\in F_0:\
\a_i(t_0)<q(s_i)$. By definition, $m_{\{r\}}(r)\not=0$. For
$r^\prime=wr$ we have
\begin{equation}\label{eq:niet}
\ka_{\{r^\prime\}}m_{\{r\}}(r)=\int_{\xi(w)}\om,
\end{equation}
where $\xi(w)$ is the cycle near $r$ which we obtain in
Proposition $\ref{prop:cycle}$ when we replace $t_0$ by
$w^{-1}t_0$. Hence we have to show that there exists a proper
choice for $t_0$ such that when we start the contour shift
algorithm from this point, the corresponding point residue at $r$
will be nonzero. The problem we have to surmount is possible
cancellation of nonzero contributions to $\ka_{\{r^\prime\}}$.
We will do this by showing that there exists at least one chamber
such that the residue at $r$ consists only of one nonzero
contribution.

We consider the real arrangement $\L^{\{r\}}$ in $T_{rs}$, and
transport the Euclidean structure of $\mathfrak{t}$ to $T_{rs}$ by
means of $t\to\log(c^{-1}t)$. Then $\L^{\{r\}}$ is the lattice of
intersections of a central arrangement of hyperplanes with center
$c$. We assign indices $i_L$ to the elements of $\L^{\{r\}}$ by
considering the index of the corresponding complex coset
containing $r$, and we note that by Corollary
\ref{cor:simpledefres},
$i_{\{r\}}=n:=\operatorname{codim}(\{r\})$. From Corollary
\ref{cor:simpledefres} we further obtain the result that there
exist full flags of subspaces $c_LT^L_{rs}\in\L^{\{r\}}$ such that
$i_L=\operatorname{codim}(L)$. In particular, there exists at
least one line $l$ through $r$ with $i_l=n-1$.

By Theorem \ref{thm:nonnest} we see that the centers $c_L$,
$c_{L^\prime}$ of two ``residual subspaces'' $c_LT^L\subset
c_{L^\prime}T^{L^\prime}$ (i.e. $\operatorname{codim}(T^L)=i_L$
and $\operatorname{codim}(T^{L^\prime})=i_{L^\prime}$) in
$\L^{\{r\}}$ satisfy $c_{L^\prime}\not=c_{L}$ unless
$c_LT^L= c_{L^\prime}T^{L^\prime}$.
Hence $d(e,c_{L^\prime})\leq d(e,c_{L})$ (where $d$ denotes the
distance function), with equality only if
$c_LT^L= c_{L^\prime}T^{L^\prime}$.
In the case of a residual line
$l\in \L^{\{r\}}$, $cT^l$ is divided in two half
lines by $c$, and $c_l$ lies in one of the two halves
(i.e. does not coincide with $c$).

We want to find a chamber for $t_0$ such that the corresponding
point residue $\ka_{\{r^\prime\}}m_{\{r\}}(r)$ at $r$ is nonzero.
We argue by induction on the rank. If the rank of $R_0$ is $1$,
obviously we get $\ka_{\{r^\prime\}}\not=0$ if we choose $t_0$
in the half line not containing $e=c_T$,
because we then have to pass a simple pole
of $\om$ at $r$ when moving the contour $t_0T_u$ to $T_u$ (since
$t_0$ and $e=c_T$ are separated by $c$). Assume by induction that
for any residual point $p$ of a rank $n-1$ root system, we can
choose a chamber for $t_0$ such that $\ka_{\{p\}}\not= 0$. Let
$S\subset T_{rs}$ be a sphere centered at $r$ through $e$, and
consider the configuration of hyperspheres in $\L^{\{r\}}\cap S$.
Let us call $e\in S$ the north pole of $S$.
If $L_S=c_LT^L_{rs}\cap S$ with $\operatorname{dim}(T^L)>1$, we
denote by $c_{L\cap S}$ the intersection of the half line through
$c_L$ beginning in $c$ (recall that $c\not=c_L$) and $L_S$.
By the above remarks,
$c_{L\cap S}$ is in the northern hemisphere for all
residual $L\supset r$ of dimension $>1$.
We call this point $c_{L\cap S}$ the center of $L_S$.

In the case when
$\operatorname{dim}(T^L)=1$, $L_S$ is disconnected and
consists of two antipodal points $c_{L\cap S}$ (north) and
$\overline{c_{L\cap S}}$ (south), its opposite.
In this case of residual lines through $r$, both of these
antipodal points are considered as centers of $\L^{\{r\}}\cap S$.
We call $c_{L\cap S}$ the northern center, and its opposite
is called the southern center.
All centers of $\L^{\{r\}}\cap S$ lie in the northern hemisphere,
with
the exception of the southern centers of the residual lines
through $r$.

Consider a closed (spherical) ball $D\subset S$ centered at $e$
such that $D$ contains a southern center $p$ in its boundary,
but no southern centers in its interior. Since $e$ is regular
with respect to $\L^{\{r\}}$
(a trivial case of Theorem \ref{thm:nonnest}, as $e$ is the
center of $T$), we have $D\not=S$.

We take $t_0$
in $S$, and we apply the algorithm as described in the proof of
Proposition \ref{prop:cycle}, {\it but now on the sphere $S$, and
with respect to the sets $L_S$ and their centers}.

By the induction hypothesis, we can take $t_0\in S$ close to $p$
in a chamber of the configuration $\L^{\{r\}}\cap S$ which
contains $p$ in its closure, such that a nonzero residue at $l$ is
picked up in $p$. Denote by $\L^p\cap S$ the central
subarrangement of elements of $\L^{\{r\}}\cap S$ containing $p$.
Consider any alternative ``identity element'' $\tilde{e}$ which
belongs to the same chamber of the {\it dual} configuration of
$\L^p\cap S$ as the real (original) identity element $e$.

As was explained in (the proof of) Proposition \ref{prop:cyinv},
when we apply the contour shifts as in (the proof of) Proposition
\ref{prop:cycle} to $\L^{\{r\}}\cap S$, the residue at $p$ only
depends on the dual chamber which contains the identity element.
In other words, we may use the new identity $\tilde{e}$ instead of
$e$ without changing the residue at $p$. We can and will choose
$\tilde{e}$ close to $p$, and in the interior of $D$. By
Proposition \ref{prop:cycle} we can replace the integral over
$t_0T_u$ by a sum of integrals over cosets of the form
$\tilde{c}_{L\cap S}\tilde{s}_LT_u^L$ (for some $\tilde{s}_L\in
T_u$) of the residue kernel $\tilde{\ka}_Lm_L$ (cf. equation
(\ref{eq:inint})) on $L$. As was mentioned above, we are only
interested in such contributions when $r\in L$, which means that
we may take $\tilde{s}_L=s$. The new ``centers'' $\tilde{c}_{L\cap
S}$ with respect to the new identity element $\tilde{e}$ are in
the interior of $D$.

Next we apply the algorithm of contour shifts as in
Proposition
\ref{prop:cycle} to move the cycles $\tilde{c}_{L\cap S}sT_u^L$ to
${c}_{L\cap S}sT_u^L$. Since both the new centers
$\tilde{c}_{L\cap S}$ and the original centers ${c}_{L\cap S}$
belong to the interior of $D$, and since the intersection of $D$
with $L_S$ is connected if $\operatorname{dim}(L_S)>0$, we can
choose every path in the contour shifting algorithm inside the
interior of $D$. Thus, the centers $c_{L\cap S}$ of the residual
cosets $L$ that arise in addition the one southern center $c_l=p$ in
the above process are in the interior of $D$. In particular, with
the exception of $psT_u^l$, the one dimensional cosets of
integration which show up in this way, all have a {\it northern}
center.

Finally, in order to compute the residue
$\ka_{r^\prime}m_{\{r\}}(r)$ at $r^\prime$, we now have to move
the center $c_{L\cap S}\in S$ of $L_S$ to the corresponding center
$c_L\in T_{rs}$ of $L$, for each residual coset $L$ which contains
$r$ and which contributes to $\int_{t_0T_u}\om$. The only such
center of $\L^{\{r\}}\cap S$ which will cross $c$ is the southern
center $p$. Since $m_l$ has a simple pole at $r=sc$, we conclude
that this gives a nonzero residue at $r$. Hence with the above
choice of $t_0$ we get $\ka_{\{r^\prime\}}\not=0$, which is
what we wanted to show.
\end{proof}
\subsection{Discrete series}
In this subsection we show that the irreducible characters
$\chi_{r,i}$ (see Definition \ref{dfn:resalg}) associated to a residual point
are in fact discrete series characters.
\begin{cor}(of Theorem \ref{thm:support})\label{cor:cas}
For every residual point $r=sc$, the sum
${\overline \ka}_{W_0r}|W_0r|=\sum_{r^\prime\in
W_0r}\ka_{r^\prime}\not=0$, and for all $a\in\A$:
\begin{equation}
\chi_r(a)=\frac{1}{{\overline \ka}_{W_0r}|W_0r|}\sum_{r^\prime\in
W_0r}\ka_{\{r^\prime\}}a(r^\prime).
\end{equation}
Moreover, $\ka_{\{r^\prime\}}=0$ unless $\forall x\in
X^+\backslash\{0\}:|x(r^\prime)|<1$ (where $X^+$ denotes the set
of dominant elements in $X$).
\end{cor}
\begin{proof}
This is immediate from Proposition \ref{chiA} and Theorem
\ref{thm:support},
except for the last assertion. This fact
follows from Proposition \ref{prop:antidual}. We know that $e$ is
regular in $\L_{\{r^\prime\}}$ by Theorem \ref{thm:nonnest}. On
the other hand, $t_0$ lies in $c^\prime T_{rs,-}$, which is
clearly a subset of a chamber of $\L^{\{r^\prime\}}$. The
anti-dual of the chamber of $\L^{\{r^\prime\}}$ containing $t_0$
is thus a subset of $c^\prime T_{rs}^+$, with
$T_{rs}^+:=\{t\in T_{rs}\mid \forall x\in X^+\backslash\{0\}:x(t)> 1\}$.
Thus when
$e$ is contained in the anti-dual chamber we have $c^\prime\in
T_{rs}^-$
\index{T2@$T_{rs}^-$, anti-dual of the positive chamber $T_{rs,+}$}
as desired.
\end{proof}
We introduce the
notation $\Delta_{\Ri}(=\Delta_{\Ri,q})$
\index{0D3@$\Delta_{\Ri}$(=$\Delta_{\Ri,q}$), irreducible discrete
series of $\H(\Ri,q)$}
for a complete set of representatives of the finite set of
equivalence classes of the
irreducible discrete series representations of $\H(\Ri,q)$, and
$\Delta_{\Ri,W_0r}(=\Delta_{\Ri,W_0r,q})$
\index{0D3@$\Delta_{\Ri,W_0r}$, irreducible discrete
series representations of $\H(\Ri,q)$ with
central character $W_0r$} for the representatives of the
classes of irreducible discrete
series of $\H(\Ri,q)$ with central character $W_0r$.
(We sometimes drop $\Ri$ from the notation if no
confusion is possible, and write $\Delta_{W_0r}$.)
\begin{lem}
$\Delta_{W_0r}$ is nonempty if and only if $r$ is residual.
If $r$ is residual, $\Delta_{W_0r}$ is in bijective
correspondence with the collection $\{\d_{r,i}\}$
of irreducible characters of $\overline{\H^r}$.
In particular, $\H(\Ri,q)$ has
at most finitely many discrete series representation.
\end{lem}
\begin{proof}
We have
\begin{equation}
\chi_{r,i}(a)=\sum_{r^\prime\in
W_0r}\dim(V_{r,i}^{r^\prime})a(r^\prime).
\end{equation}
Hence from $d_{r,i}>0$,
\begin{equation}\label{eq:packet}
\chi_{r}(a)=\sum_{i}\chi_{r,i}(a)d_{r,i},
\end{equation}
and Corollary \ref{cor:cas} we conclude that the generalized
weight spaces of $V_{r,i}$ indeed satisfy the Casselman criterion
Lemma \ref{lem:casds} for
discrete series.

Conversely, if $\d$ is a discrete series representation,
Theorem \ref{thm:supds} implies that $\mu_{Pl}(\d)>0$. By Corollary
\ref{cor:exten},
the central character $W_0r$ of $\d$ is such that $\nu(\{r\})>0$.
Theorem \ref{thm:nu} implies that such points $r$ are
necessarily residual.
\end{proof}
In view of the above, we adapt the notations of
Definition \ref{dfn:resalg} accordingly,
i.e. we write
$d_{\Ri,\d}$ (or simply $d_{\d}$)
\index{d@$d_\d(=d_{\Ri,\d})$, residual degree; degree
of $\d$ in the residual Hilbert algebra $\overline{\H^r}$}
instead of
$d_{r,i}$
\index{d@$d_{t,i}(=d_{W_0t,i})$, residual degree; degree of $\chi_{t,i}$ in
$\overline{\H^t}$}
if $\d\in\Delta_{\Ri,W_0r}$, and its character
$\chi_\d$
\index{0w@$\chi_\d$, character of $\d$}
descends to
$\chi_{r,i}$
\index{0w@$\chi_{t,i}$, irreducible character of $\overline{\H^t}$}
on
$\overline{\H^r}$
etc.
\begin{cor}\label{cor:fdim}
Let $\d\in\Delta_{W_0r}$.
The formal dimension
$\mu_{Pl}(\d)$ of $\d$ equals
\begin{equation}\label{eq:fdim}
\mu_{Pl}(\d)=\operatorname{fdim}(\d)={d_\d\nu(\{W_0r\})}
=|W_0r|{\overline{\ka}}_{W_0r}
d_\d m_{\{r\}}(r)
\end{equation}
\end{cor}
\begin{proof}
Combine equation (\ref{eq:chadec}), Corollary \ref{cor:exten},
and Theorem \ref{thm:nu}.
\end{proof}
\begin{cor}\label{rem:cas}
For a residual point $r$ there exist constants $C,
\e>0$ such that
\begin{equation}
|\chi_r(N_w)|\leq C\exp(-\e l(w)).
\end{equation}
\end{cor}
\begin{cor}
The residual degrees $d_{\d}>0$ of the irreducible characters $\chi_{\d}$ of
the residual algebra $\overline{\H^r}$ (with $r$ a residual point)
satisfy the following system of linear equations.
\begin{equation}\label{eq:rats}
\sum_{\d\in[\Delta_{\Ri,W_0r}]} \dim(V^{r^\prime}_\d)d_\d=
\frac{\ka_{\{r^\prime\}}}{{\overline \ka}_{W_0r}|W_0r|}.
\end{equation}
(with $V_\d^{r^\prime}$ the generalized $r^\prime$-weight space
in the space $V_\d$ of $\d$).
In particular we conclude that the nonzero $\ka_{\{r^\prime\}}$
all have the same sign (equal to the sign of $m_{\{r\}}(r)$).
\end{cor}
\begin{rem}
We note in addition that if the restrictions $\chi_\d|_{\A}$ to
$\A$ of the characters $\chi_\d$ are linearly independent, it
follows from the equations (\ref{eq:rats}) that $d_\d\in\Q$ for
all $\d\in [\Delta_{\Ri}]$.
I did not find any argument in favor of this
linear independence. However, we do conjecture that the constants
$d_{\d}$ are rational, see Conjecture \ref{rem:ell}.
\end{rem}
\subsection{Temperedness of the traces $\chi_t$}
In this subsection we discuss the tempered growth behaviour of the
$\chi_t$ on the orthonormal basis $N_w$ of $\H$, as a corollary of
the analysis of the $\A$-weights of $\chi_t$.
\begin{prop}\label{cor:indA}
Let $L$ be residual such that $W_L$ is a standard parabolic
subgroup of $W_0$. For $t\in L^{temp}$ we write $t=r_Lt^L$, with
$t^L\in T^L_u$. We consider $\chi_t|_{\A}$ as a formal linear
combination of elements of $T$. Likewise, let $\A_L=\C[X_L]$ be
the ring of regular functions on $T_L\subset T$. We consider
$\chi_{\Ri_L,\{r_L\}}|_{\A_L}$ as a formal linear combination of
elements of $T_L$. In this sense we have, $\nu_L$-almost
everywhere on $L^{temp}$,
\begin{equation}\label{eq:div}
\chi_t|_{\A}=\frac{1}{|W^L|}\sum_{w\in
W^L}w(t^L\chi_{\Ri_L,\{r_L\}}|_{\A_L}).
\end{equation}
Hence $\nu$-almost everywhere, $\chi_t$ is a nonzero tempered
functional on $\H$.
\end{prop}
\begin{proof}
Equation (\ref{eq:div}) follows by a straightforward computation
similar to Proposition \ref{chiA}, using Proposition
\ref{prop:par} and the definition of $\chi_t$. Since $\chi_t$ is a
positive combination of the irreducible characters of the residual algebra
$\overline{\H^t}$, it follows that the weights
$t^\prime\in W_0t$ of the generalized $\A$-eigenspaces of the irreducible
characters of $\overline{\H^t}$ all satisfy the condition
$\forall x\in X^+:\ |x(t^\prime)|\leq 1$.
This shows, by Casselman's criterion Lemma \ref{lem:cas},
that $\chi_t$ is a
tempered functional on $\H$.
\end{proof}
\section{Localization of the Hecke algebra}\label{sect:loc}
We have obtained thus far a decomposition of the trace $\tau$ as
an integral of positive, finite traces $\chi_t$ against an
explicit probability measure $\nu$ on $T$, such that each $\chi_t$
is a finite positive linear combination of finite dimensional,
irreducible characters of $\cf$. This is an important step towards
our goal of finding the Plancherel decomposition, but it is not
yet satisfactory because we know virtually nothing about the
behavior of the decomposition of $\chi_t$ in irreducible
characters at this stage, neither as a function of $t$,
nor as a function of $\q$.
In particular, the residual degrees $d_{t,i}\in \mathbb{R}_+$
of the residual algebras are obscure at this point, and these
degrees are involved in the Plancherel measure $\mu_{Pl}$.

In the remaining part of the paper we will formulate the Plancherel
theorem, and also remedy to some extend the above problems. The
support $S$ of $\nu$ (viewed as a $W_0$-invariant measure on $T$)
decomposes as a union of the closed
sets $L^{temp}$ (see \ref{sl}). For each $L$ we show that, up to
isomorphism of Hilbert algebras, the residual algebras
$\overline{\H^t}$ are independent of $t$, $\nu$-almost everywhere
on $L^{temp}$.

The above is based on ideas of Lusztig \cite{Lu} about
completions of the affine Hecke algebra. Lusztig describes the
$\I_t$-adic completion of $\H$,
where $\I_t$ is a maximal ideal of $\Ze$. It is not
hard to see that Lusztig's arguments can be adapted to
(analytic) localization with respect to suitably small open
neighborhoods $U\supset W_0t$ of orbits of points in $T$, and this
will be discussed in present section.

When $s=s_\alpha\in S_0$ (with $\alpha\in F_1$), we define an
intertwining element $\i_s$ as follows:
\begin{equation}
\begin{split}
\i_s&=(1-\theta_{-\alpha})T_s+((1-q_{\alpha^\vee}q_{2\alpha^\vee})
+q_{\alpha^\vee}^{1/2}(1- q_{2\alpha^\vee})\theta_{-\alpha/2})\\
&=T_s(1-\theta_{\alpha})+((q_{\alpha^\vee}q_{2\alpha^\vee}-1)
\theta_{\alpha}
+q_{\alpha^\vee}^{1/2}(q_{2\alpha^\vee}-1)\theta_{\alpha/2})\\
\end{split}
\end{equation}
\index{0j@$\iota_s$, intertwining element of $\H$}
We remind the reader of the convention of Remark \ref{rem:conv}.
These elements are important tools to study the Hecke algebra. We
recall from \cite{EO}, Theorem 2.8 that these elements satisfy the
braid relations, and they satisfy (for all $x\in X$)
\[
\i_s\theta_x=\theta_{s(x)}\i_s,
\]
and finally they satisfy
\[
\i_s^2=(q_{\alpha^\vee}^{1/2}+\theta_{-\alpha/2})
(q_{\alpha^\vee}^{1/2}+\theta_{\alpha/2})
(q_{\alpha^\vee}^{1/2}q_{2\alpha^\vee}-\theta_{-\alpha/2})
(q_{\alpha^\vee}^{1/2}q_{2\alpha^\vee}-\theta_{\alpha/2}).
\]
(where we have again used the convention of Remark
\ref{rem:conv}!). Suitably normalized versions of the $\i_s$
generate a group isomorphic to the Weyl group $W_0$. In order to
normalize the intertwiners, we need to tensor $\H$ by the field of
fractions $\F$ of the center $\Ze$. So let us introduce the
algebra
\begin{equation}
{}_\F\H:=\F\otimes_\Ze\H
\end{equation}
\index{H6@${}_\F\H:=\F\otimes_\Ze\H$, $\F$ field of fractions of $\Ze$}
with the multiplication defined by $(f\otimes h) (f^\prime\otimes
h^\prime):=ff^\prime\otimes hh^\prime$. Notice that this an
algebra over $\F$ of dimension $|W_0|^2$. The subalgebra
${}_\F\A=\F\otimes_\Ze\A$
\index{A@${}_\F\A:=\F\otimes_\Ze\A$, $\F$ field of fractions of $\Ze$}
is isomorphic to the field of fractions
of $\A$. The field extension $\F\subset{}_\F\A$ has Galois group
$W_0$, and we denote by $f\to f^w$ the natural action of $W_0$ on
the field of rational functions on $T$. The elements $T_w$ with
$w\in W_0$ form a basis for ${}_\F\H$ for multiplication on the
left or multiplication on the right by ${}_\F\A$, in the sense
that
\begin{equation}
{}_\F\H=\oplus_{w\in W_0}{}_\F\A T_w=\oplus_{w\in W_0}T_w{}_\F\A.
\end{equation}
The algebra structure of ${}_\F\H$ is determined by the
Bernstein-Zelevinski-Lusztig relations as before: when
$f\in{}_\F\A$ and $s=s_{\a}$ with $\a\in F_1$, we have
\begin{equation}
fT_s-T_sf^s=((q_{2\alpha^\vee}q_{\alpha^\vee}-1)
+q_{\alpha^\vee}^{1/2} (q_{2\alpha^\vee}-1)\theta_{-\alpha/2})
\frac{f-f^s} {1-\theta_{-\alpha}}
\end{equation}
We have identified $\A$ with the algebra of regular functions on
$T$ in the above formula.

Let us introduce
\begin{equation}
\begin{split}
n_\alpha &:= q(s_\alpha)\Delta_\alpha c_\a \\
&=(q_{\alpha^\vee}^{1/2}+\theta_{-\alpha/2})
(q_{\alpha^\vee}^{1/2}q_{2\alpha^\vee}-\theta_{-\alpha/2})\in
\A,\\
\end{split}
\end{equation}
\index{n@$n_\a$, numerator of $c_\a$}
where we used the Macdonald $c$-function introduced in equation
(\ref{eq:defD}) and (\ref{eq:defc}).

The normalized intertwiners are now defined by (with $s=s_\a$,
$\a\in R_1$):
\begin{equation}\label{eq:defint}
\i^0_s :=n_\a^{-1}\i_s\in {}_\F\H.
\end{equation}
\index{0j@$\iota_w^0$, normalized intertwining element of ${}_\F\H$}
By the properties of the intertwiners listed above it is clear
that $(\i_s^0)^2=1$. In particular, $\i_s^0\in{}_\F\H^\times$, the
group of invertible elements of ${}_\F\H$. From the above we have
the following result:
\begin{lemma}
The map $S_0\ni s\to \i^0_s\in {}_\F\H^\times$ extends (uniquely)
to a homomorphism $W_0\ni w\to \i^0_w\in {}_\F\H^\times$. Moreover,
for all $f\in {}_\F\A$ we have that $\i_w^0f\i_{w^{-1}}^0=f^w$.
\end{lemma}
Lusztig (\cite{Lu}, Proposition 5.5) proved that in fact
\begin{theorem}\label{thm:ind}
\begin{equation}
{}_\F\H=\oplus_{w\in W_0}\i^0_w{}_\F\A=\oplus_{w\in W_0}{}_\F\A
\i^0_w
\end{equation}
\end{theorem}
Let $U\subset T$ be a nonempty, open, $W_0$-invariant subset.
We denote by
$\Ze^{an}(U)$
\index{Z2@$\Ze^{an}(U)$, ring of $W_0$-invariant holomorphic
functions on $U\subset T$}
the ring of $W_0$-invariant holomorphic
functions of $U$.  Consider the algebras
$\A^{an}(U):=\Ze^{an}(U)\otimes_\Ze\A$
\index{A1@$\A^{an}(U):=\Ze^{an}(U)\otimes_\Ze\A$, ring of
holomorphic functions on $U\subset T$}
and
$\H^{an}(T):=\Ze^{an}(T)\otimes_\Ze\H$.
The algebra structure on
$\H^{an}(T)$ is defined by $(f\otimes h) (f^\prime\otimes
h^\prime):=ff^\prime\otimes hh^\prime$ (similar to the definition
of ${}_\F\H$).

Let us first remark that the finite dimensional representation
theory of the ``analytic'' affine Hecke algebra $\H^{an}(T)$ is
the same as the finite dimensional representation theory of $\H$.
Every finite dimensional representation $\pi$ of $\H$ determines a
co-finite ideal $J_\pi\subset \Ze$, the ideal of central elements
of $\H$ which are annihilated by $\pi$. Denote by $J_\pi^{an}$ the
ideal of $\Ze^{an}(T)$ generated by $J_\pi$. Because of the
co-finiteness we have an isomorphism
\begin{equation}
\Ze/J_\pi \tilde{\rightarrow} \Ze^{an}(T)/J^{an}_\pi(T).
\end{equation}
This shows that $\pi$ can be uniquely lifted to a representation
$\pi^{an}$
\index{0p1@$\pi^{an}$, $\pi$ extended to $\H^{an}$}
of $\H^{an}(T)$ whose restriction to $\H$ is $\pi$. The
functor
$\pi\to\pi^{an}$
defines an equivalence between the
categories of finite dimensional representations of $\H$ and
$\H^{an}(T)$ (with the inverse given by restriction).

For any $W_0$-invariant nonempty open set $U\subset T$ we define
the localized affine Hecke algebra
\begin{equation}
\H^{an}(U):=\Ze^{an}(U)\otimes_{\Ze}\H.
\end{equation}
\index{H61@$\H^{an}(U):=\Ze^{an}(U)\otimes_\Ze\H$, the Hecke algebra with
coefficients in $\Ze^{an}(U)$}
This defines a presheaf of $\Ze^{an}$-algebras on $W_0\backslash
T$, which is finitely generated over the analytic structure sheaf
$\Ze^{an}$ of the geometric quotient $W_0\backslash T$.

A similar argument as above shows that
\begin{prop}\label{prop:anequiv}
The category $\operatorname{Rep}(\H^{an}(U))$ of finite
dimensional modules $\pi^{an}_U$
\index{0p2@$\pi^{an}_U$, $\pi$ extended to $\H^{an}(U)$}
over $\H^{an}(U)$ is equivalent to
the category
$\operatorname{Rep}_U(\H)$
\index{Rep@$\operatorname{Rep}_U(\H)$, category of finite dimensional
representations of $\H$ whose $\Ze$-spectrum is contained in $U$}
of finite dimensional
modules $\pi$ over $\H$ whose $\Ze$-spectrum is contained in $U$.
\end{prop}
\begin{lem}\label{lem:chinese}
For every $W_0$-invariant nonempty open set $U$ in $T$,
we have the isomorphism $\A^{an}(U)\simeq\Ze^{an}(U)\otimes_\Ze
\A$, where $\A^{an}(U)$ denotes the ring of analytic functions on
$U$.
\end{lem}
\begin{proof}
Both the left and the right hand side are finitely generated
modules over $\Ze^{an}(U)$, and we have a natural morphism from
the right hand side to the left hand side (product map).
In order to prove that this map is an isomorphism it suffices to show
this in the stalks of the corresponding
sheaves at each point of $W_0\backslash U$. Let $\I_t$
denote the maximal ideal in $\Ze$ corresponding to $W_0t$, and let
$\hat\Ze_{t}$ denote the $\I_t$-adic completion.
Because
$\hat\Ze_{t}$ is faithfully flat over $\Ze_t^{an}$
(the stalk at $W_0t$ of the sheaf $\Ze^{an}$), it suffices to
check that for each $t\in U$, we have
\begin{equation}\label{imp}
\hat\Ze_{t}\otimes_{\Ze^{an}_{t}}\A^{an}_{W_0t}\simeq
\hat\Ze_{t}\otimes_{\Ze} \A,
\end{equation}
where $\A^{an}_{W_0t}=\oplus_{t^\prime}\A^{an}_{t^\prime}$ denotes
the space of analytic germs at the set $W_0t$.
Let $m_t$ denote the maximal ideal of $\A$ at $t\in T$, and
let $\I_t\A=\prod_{t^\prime\in W_0t}j_{t^\prime}$ with
$j_{t^\prime}=\I_t\A\cap m_{t^\prime}$.
For all $t^\prime\in W_0t$ we have
$\widehat{\A^{an}_{t^\prime}}_{j_{t^\prime}\A^{an}_{t^\prime}}
=\hat{\A}_{j_{t^\prime}}$.
Since $A^{an}_t\cap\I_tA^{an}_{W_0t}=j_t\A^{an}_t$,
the left hand side of \ref{imp} is equal to
$\oplus_{t^\prime\in W_0t}\hat{\A}_{j_{t^\prime}}$, the sum of
the completions of $\A$ with respect to $j_{t^\prime}$.
The right hand side of \ref{imp}
is equal to the completion $\hat\A_{\I_t\A}$.
By the
Chinese remainder theorem,
$\hat\A_{I_t\A}\simeq
\oplus_{t^\prime\in W_0t}\hat\A_{j_{t^\prime}}$,
finishing the proof.
\end{proof}
\begin{prop}
The algebra $\H^{an}(U)$ is a free $\A^{an}(U)$ module of rank
$|W_0|$, with basis $T_w\otimes 1$ ($w\in W_0$). When
$f\in\A^{an}(U)$ and $s=s_{\a}$ with $\a\in F_1$ we have again the
Bernstein-Zelevinski-Lusztig relation
\begin{equation}
fT_s-T_sf^s=((q_{2\alpha^\vee}q_{\alpha^\vee}-1)
+q_{\alpha^\vee}^{1/2} (q_{2\alpha^\vee}-1)\theta_{-\alpha/2})
\frac{f-f^s} {1-\theta_{-\alpha}}.
\end{equation}
This describes the multiplication in the algebra $\H^{an}(U)$. The
center of $\H^{an}(U)$ is equal to $\Ze^{an}(U)$.
\end{prop}
Similarly we have the localized meromorphic affine Hecke algebra
$\H^{me}(U)$, which is defined by
\begin{equation}
\H^{me}(U):=\F^{me}(U)\otimes_{\Ze}\H,
\end{equation}
\index{H6@$\H^{me}(U):=\F^{me}(U)\otimes_{\Ze}\H$,
localized Hecke algebra with meromorphic coefficients}
where
$\F^{me}(U)$
\index{F@$\F^{me}(U)$, quotient field of $\Ze^{an}(U)$}
it the quotient field of $\Ze^{an}(U)$. We
write
$\A^{me}(U):=\F^{me}(U)\otimes_{\Ze}\A$
\index{A3@$\A^{me}(U):=\F^{me}(U)\otimes_{\Ze}\A$,
ring of meromorphic functions on $U\subset T$}.
It
is the ring of meromorphic functions on $U$.
\begin{theorem}\label{thm:gralg}
\begin{equation} \H^{me}(U)=\oplus_{w\in
W_0}\A^{me}(U)\i^0_w=\oplus_{w\in W_0} \i_w^0\A^{me}(U)
\end{equation}
\end{theorem}
\begin{proof}
This is clear from Theorem \ref{thm:ind} by the remark that
$\H^{me}$ arises from the $\F$-algebra ${}_\F\H$ by extension of
scalars according to
\begin{equation}
\begin{split}
\H^{me}(U)&=\F^{me}(U)\otimes_\Ze\H\\
&=\F^{me}(U)\otimes_\F{}_\F\H.\\
\end{split}
\end{equation}
\end{proof}
\subsection{Lusztig's structure theorem and parabolic
induction}\label{sub:lus} We shall investigate the structure of
the tracial states $\chi_t$, using Lusztig's technique of
localization of $\H$ as discussed above. The results in the
present subsection are substitutes for the usual techniques of
parabolic induction for reductive groups. The results in this
subsection are closely related to the results on parabolic
induction in the paper \cite{BM}.

We use in fact a slight variation of the results of Lusztig \cite{Lu}.
There are two main differences. First of
all we work with analytic localization at suitably small
neighborhoods, instead of Lusztig's use of adic completion.
In addition we have
replaced the root system of the localized algebra which Lusztig
has defined by something slightly different. Lusztig's construction
only works with the additional assumption in Convention \ref{eq:scale}
that $f_s\in\mathbb{N}$, and this assumption is not natural in our
context. We have therefore adapted the construction.

We define a function
\begin{equation}
T\ni t\to R_{P(t)}\subset R_0, {\mathrm{\ a\ parabolic\ subsystem}}
\end{equation}
\index{R6@$R_{P(t)}=R_{P(\varpi)}\subset R_0$, parabolic subsystem
associated with $t\in\varpi\subset T$}
by putting $R_{P(t)}:=R_0\cap\mathfrak{t^*}_{<t>}$, with
$\mathfrak{t^*}_{<t>}\subset\mathfrak{t^*}=
\mathbb{R}\otimes_{\mathbb{Z}}X$ the subspace spanned by the roots
$\a\in R_0$ for which one of the following properties holds
\begin{enumerate}
\item $c_\a\not\in\mathcal{O}_t^\times$
(the invertible holomorphic germs at $t$).
\item $\a(t)=1$,
\item $\a(t)=-1$ and $\a\not\in 2X$.
\end{enumerate}
We let
$P(t)\subset R_{P(t),+}:=R_{P(t)}\cap R_{0,+}$
\index{P@$P(\varpi)=P(t)$, basis of simple roots in
$R_{P(t),+}$ where $\varpi=W_{P(t)}t$}
be the basis of simple root for $R_{P(t),+}$.
%Notice that in the situation $\a\in R_0$, $2\a\in R_1$,
%$q_{\a^\vee/2}=1$ (thus R=B_n, X=Q, \a=e_n\in Q, not in 2Q)
%and $\a(t)=-1$, $\a$ satisfies (iii) but not (i).
%If $a(t)=1$ and $q_{a^\vee/2}^{1/2}q_{\a^\vee}=1$
%then we have similarly that (ii) holds, but not (i)
We have the following easy consequences of the definition:
\begin{prop}\label{prop:eltprop}
\begin{enumerate}\label{def:two}
\item $t\to R_{P(t)}$ is lower semi-continuous with respect to
the Zariski-topology of $T$ and the ordering of subsets of
$R_0$ by inclusion.
\item $t\to R_{P(t)}$ is equivariant: for all $w\in W$ we have
$R_{P(wt)}=w(R_{P(t)}$.
\end{enumerate}
\end{prop}
We denote by
$W_{P(t)}$
\index{W4@$W_{P(t)}=W(R_{P(t)})$, parabolic subgroup
associated with $t\in T$}
the parabolic subgroup of $W_0$ generated
by the reflections $s_\a$ with $\a\in R_{P(t)}$.
We say that $t_1,t_2\in W_0t$ are equivalent if there exists
a $w\in W_{P(t_1)}$ such that $t_2=w(t_1)$. To see that this is
actually an equivalence relation, observe that
$R_{P(t_2)}=R_{P(t_1)}$ for all $t_2\in W_{P(t_1)}t_1$. The equivalence
classes are the orbits $\varpi=W_{P(t)}t$.
This gives a partition of $W_0t$ in a collection equivalence
classes which are denoted by $\varpi\subset W_0t$.
If $t\in \varpi$ we sometimes write $P(\varpi)$,
$W_{P(\varpi)}$ etc. instead of $P(t)$, $W_{P(t)}$ etc.
Note that $W_0$ acts transitively on the set of equivalence
classes and that for each equivalence class $\varpi$,
$W_{P(\varpi)}$ acts transitively on $\varpi$.

Let $\varpi\subset W_0t$ be the equivalence class of $t$.
We define:
\begin{equation}
W_\varpi:=\{w\in W_0\mid w(\varpi)=\varpi\}.
\end{equation}
\index{W4a@$W_\varpi$, stabilizer in $W_0$ of $\varpi=W_{P(t)}t$}
By Proposition \ref{prop:eltprop} it is clear that
$W_{P(\varpi)}\vartriangleleft W_\varpi$, and that this normal
subgroup is complemented by the subgroup
\begin{equation}
W(\varpi):=\{w\in W_\varpi\mid w(P(\varpi))=P(\varpi)\}
\end{equation}
\index{W4b@$W(\varpi)=
\{w\in W_\varpi\mid w(P(\varpi))=P(\varpi)\}$,
complement of $W_{P(\varpi)}$ in
$W_\varpi$}
\begin{lem}\label{lus:lem8.2b} For $\a\in R_0$ we have:
$\a\in R_{P(\varpi)}\Longleftrightarrow
s_\a\in W_\varpi$.
\end{lem}
\begin{proof}
We only need to show that $s_\a\in W_\varpi$ implies that
$\a\in R_{P(\varpi)}$ (the other direction being obvious). Notice that
if $t\in\varpi$ we have
\begin{equation}
t^{-1}\varpi\subset\mathbb{Z}R_{P(\varpi)}^\vee\otimes\mathbb{C}^\times.
\end{equation}
If $s_\a\in W_\varpi$ then $s_\a(t)\in\varpi$, and thus
\begin{equation}\label{eq:lusz}
\a^\vee\otimes\a(t)\in\mathbb{Z}
R_{P(\varpi)}^\vee\otimes\mathbb{C}^\times.
\end{equation}
By Proposition \ref{prop:eltprop}, we have
$s_\a(R_{P(\varpi)})=R_{P(\varpi)}$.
Since $R_{P(\varpi)}$ is parabolic this implies that
either $\a\in R_{P(\varpi)}$ or that $\a(R_{P(\varpi)}^\vee)=0$.
In the first case we are done, so let us assume the second case.
By (\ref{eq:lusz}) it follows that
$1=\a(\a^\vee\otimes\a(t))=\a(t)^2$.
If $\a(t)=1$ we have $\a\in R_{P(\varpi)}$ by definition,
contradicting the assumption.
If $\a(t)=-1$ and $\a\not\in2X$ then, by definition,
$\a\in R_{P(\varpi)}$, contrary to the assumption.
If $\a(t)=-1$ and $\a=2x$ for some $x\in X$ then
(\ref{eq:lusz}) implies $1=x(\a^\vee\otimes\a(t))=\a(t)=-1$,
again a contradiction.
We conclude that the second case does not arise altogether,
and we are done.
\end{proof}
Consider the algebra
$\H^{P(t)}:=\H(X,Y,R_{P(t)},R_{P(t)}^\vee,P(t))$.
Note that $W(\varpi)$ acts by means of
automorphisms on $\Ri^{P(t)}=(X,Y,R_{P(t)},R_{(t)}^\vee,P(t))$,
compatible with
the root labels $q$. Thus we may define an action of
$\g\in W(\varpi)$ on $\H^{P(t)}$ by $\g(T_w\theta_x)=T_{(\g
w\g^{-1})}\theta_{\g x}$. In this way we form the algebra
$\H^\varpi:=\H^{P(t)}[W(\varpi)]$
\index{H7@$\H^\varpi$, cross product of $\H^{P(\varpi)}$ by $W(\varpi)$},
with its product being defined
by $(h_1\g_1)(h_2\g_2)=h_1\g_1(h_2)\g_1\g_2$.

By Proposition \ref{prop:eltprop}(i) it is obvious that for any $t\in T$
there exists an open ball
$B\subset\mathfrak{t}_\C$ centered around the origin such that
the following conditions are satisfied:
\begin{cond}\label{cond}
\begin{enumerate}
\item[(i)] $\forall \a\in R_0,b\in B:|\Im(\a(b))|<\pi$.
In particular, the map $\exp:\mathfrak{t}_\C\to T$ restricted to
$B$ is an analytic diffeomorphism onto its image $\exp(B)$ in $T$.
\item[(ii)] If $w\in W_0$ and $t\exp(B)\cap
w(t\exp(B))\not=\emptyset$ then $wt=t$.
\item[(iii)] For all $t^\prime\in t\exp(B)$, we have
$R_{P(t^\prime)}\subset R_{P(t)}$.
\end{enumerate}
\end{cond}

Let $t\in T$. We take $B\subset\mathfrak{t}_\C$ as above and we put
$U=W_0t\exp(B)$. Concerning the analytic localization $\H^{an}(U)$ we
have the following analog of Lusztig's first reduction theorem
(see \cite{Lu}):
\begin{thm}\label{thm:lusind}
For $\varpi\subset W_0t$ an equivalence class, we put
$U_\varpi:=\varpi\exp(B)$. We define $1_\varpi\in\A^{an}(U)$ by
$1_\varpi(u)=1$ if $u\in U_\varpi$ and $1_\varpi(u)=0$ if
$u\not\in U_\varpi$. The elements
$1_\varpi$
\index{>@$1_\varpi$, idempotent in $\A^{an}(U)$ with support
$U_\varpi$}
are mutually
orthogonal idempotents. Let $t\in\varpi$.
\begin{enumerate}
\item[(i)] We have
$\H^{\varpi,an}(U_\varpi):=\H^{P(\varpi),an}(U_\varpi)[W(\varpi)]\simeq
1_\varpi\H^{an}(U)1_\varpi$.
\item[(ii)] We can define linear isomorphisms
\begin{equation}
\Delta_{\varpi_1,\varpi_2}:\H^{\varpi,{an}}(U_\varpi)\to
1_{\varpi_1}\H^{an}(U)1_{\varpi_2}.
\end{equation}
\index{0D4@$\Delta_{\varpi_1,\varpi_2}:\H^{\varpi,{an}}(U_\varpi)\to
1_{\varpi_1}\H^{an}(U)1_{\varpi_2}$, linear isomorphisms}
such that
$\Delta_{\varpi_1,\varpi_2}(h)\Delta_{\varpi_3,\varpi_4}(h^\prime)=
\Delta_{\varpi_1,\varpi_4}(hh^\prime)$ if $\varpi_2=\varpi_3$, and
$\Delta_{\varpi_1,\varpi_2}(h)\Delta_{\varpi_3,\varpi_4}(h^\prime)=0$
else.
\item[(iii)] The center of $\H^{\varpi,an}(U_\varpi)$ is
$\Ze^{\varpi,an}(U_\varpi):=(\A^{an}(U_\varpi))^{W_\varpi}$. This
algebra is isomorphic to $\Ze^{an}(U)$ via the map $z\to 1_\varpi
z$, and this gives $\H^{\varpi,an}(U_\varpi)$ the structure of a
$\Ze^{an}(U)$-algebra.
\item[(iv)] Let $N$ denote the number of equivalence classes
in $W_0t$. There exists an isomorphism
$\H^{an}(U)\simeq(\H^{\varpi,an}(U_\varpi))_N$, the algebra of
$N\times N$ matrices with entries in
$1_\varpi\H^{an}(U)1_\varpi\simeq\H^{\varpi,an}(U_\varpi)$. It is
an isomorphism of $\Ze^{an}(U)$-algebras.
\end{enumerate}
\end{thm}
\begin{proof}
The difference with Lusztig's approach is that he works with the
$\I_t$-adic completions of the algebras instead of the
localizations to $U$. In addition, we have a different definition
of the root system $R_{P(t)}$.

Using Lemma \ref{lem:chinese} we can copy the arguments of \cite{Lu},
because of the Conditions
\ref{cond} for $B$ and because of Lemma \ref{lus:lem8.2b} (which
replaces in our situation Lemma 8.2b of \cite{Lu}).
By this we see that the function $c_\a$ is analytic and
invertible on $U_\varpi\cup U_{s_\a\varpi}$ for all $\a\in R_0$
such that $s_\a\not\in W_\varpi$ (compare \cite{Lu}, Lemma 8.9),
and this is the crucial point of the construction.
\end{proof}
\begin{cor}
The functor $V\to V_\varpi:=1_\varpi V$ defines an equivalence
between the category of finite dimensional representations
of $\H^{an}(U)$ and the category of finite dimensional
representations of
$\H^{\varpi,an}(U_\varpi)=\H^{P(t),an}(U_\varpi)[W(\varpi)]$. We
have $\dim(V)=N\dim(V_\varpi)$ where $N$ denotes the
number of equivalence classes in $W_0t$. \qed
\end{cor}
\begin{dfn}\label{defn:generic}
Let $R_P\subset R_0$ be a parabolic root subsystem,
with $P\subset R_{P,+}:=R_P\cap R_{0,+}$ its basis
of simple roots. We denote the
corresponding parabolic subgroup of $W_0$ by $W_P:=W(R_P)$. We
call $t\in T$ an {\it $R_P$-generic point} if
$W_\varpi\subset W_P$ for $\varpi=W_{P(t)}t$.
\end{dfn}
\begin{cor}\label{cor:gen}
If $t$ is $R_P$-generic we have $R_{P(t)}\subset R_P$.
\end{cor}
\begin{proof}
This is immediate from the definitions.
\end{proof}

We define for any parabolic subsystem $R_P\subset R_0$ with
basis $P$ of $R_{P,+}$ the parabolic subalgebra
$\H^P=\H(X,Y,R_P,R_P^\vee,P)\subset \H$
with root labels $q^P$\index{H4@$\H^P=\H(\Ri^P,q^P)$, parabolic
subalgebra of $\H$}.

Assume that $B$ satisfies the Conditions \ref{cond}. Notice that if
$t$ is $R_P$-generic, then every $t^\prime\in t\exp(B)$ is
$R_P$-generic. Indeed, let $\varpi^\prime=W_{P(t^\prime)}t^\prime$
and $\varpi=W_{P(t)}t$.
If $w\in W_{\varpi^\prime}$, then there exists a
$w^\prime\in W_{P(t^\prime)}\subset W_{P(t)}$ (by condition \ref{cond}(iii))
such that $w^\prime t^\prime= w t^\prime$ (since the equivalence class
of $t^\prime$ is a $W_{P(t^\prime)}$-orbit). By condition \ref{cond}(ii),
also $w^\prime t=wt$. Hence $w\in W_\varpi\subset W_P$, as required.

We now put $U=W_0t\exp(B)$,
$U_P=W_Pt\exp(B)$
\index{U@$U_P$, certain $W_P$-invariant open set in $T$}
and consider the
localization $\H^{P,an}(U_P)$.
\begin{cor}\label{cor:struct}
Assume that $t\in T$ is $R_P$-generic. We have
$\H^{an}(U)\simeq(\H^{P,an}(U_P))_{|W^P|}$, where
$W^P=W_0/W_P$. Moreover, when we define
$1_P:=\sum_{\varpi\subset W_Pt}1_\varpi$
\index{>@$1_P$, idempotent in $\A^{an}(U)$, with support $U_P$}
then
$\H^{P,an}(U_P)\simeq 1_P\H^{an}(U)1_P$. These are isomorphisms of
$\Ze^{an}(U)$-algebras.
\end{cor}
\begin{proof}
The fact that $t$ is $R_P$-generic implies that the
$W_0$-equivalence classes of the elements of $W_Pt$ are equal to
the $W_P$-equivalence classes of these elements. Therefore we
have, by the above theorem,
$\H^{P,an}(U_P)\simeq(\H^{\varpi,an}(U_\varpi))_{n_P}$, where
$n_P$ is the number of equivalence classes $\varpi^\prime$ in the
orbit $W_Pt$. And for each $w\in W_0$, $wW_Pt\subset W_0t$ is a
union of $n_P$ distinct $W_0$-equivalence classes. The orbit
$W_0t$ is the disjoint union of $|W^P|$ subsets of the form
$wW_Pt\subset W_0t$, since the stabilizer of $t$ is contained in
$W_P$ (because $t$ is $R_P$-generic).  Each
subset $wW_{P(t)}t$ in $W_0t$ is partitioned
into $n_P$ equivalence classes, and the result follows.
\end{proof}
Recall that, by Proposition \ref{prop:anequiv}, a finite
dimensional representation $(V,\pi)$ of $\H$ with its
$\Ze$-spectrum contained in $U$ extends uniquely to a
representation $(V^{an},\pi^{an})$ of $\H^{an}(U)$.
\begin{cor}\label{cor:restr}
In the situation of Corollary \ref{cor:struct}, there exists an
equivalence $(V,\pi)\to (V_P,\pi_P)$ between
$\operatorname{Rep}_U(\H)$ and
$\operatorname{Rep}_{U_P}(\H^P)$,
characterized by $V_P^{an}=1_PV^{an}$
\index{V@$(V_P,\pi_P)$, representation of $\H^P$ with $V_P=1_PV$}.
We have
$\dim(V)=|W^P|\dim(V_P)$, and the inverse functor is given by
$V_P\to\operatorname{Ind}_{\H^P}^\H(V_P)=\H\otimes_{\H^P}V_P$. The
character $\chi^P$ of the module $(V_P,\pi_P)$ of $\H^P$ is given
in terms of the character $\chi_{\pi}$ of $(V,\pi)$ by the formula
$\chi^P(h)=\chi_{\pi}(1_Ph)$.
\end{cor}
\begin{proof}
We localize both the algebras $\H$ and $\H^P$ and use Proposition
\ref{prop:anequiv} and Corollary \ref{cor:struct}. Using Corollary
\ref{cor:struct} we see that the functor $V\to 1_PV^{an}|_{\H^P}$
is the required equivalence. The relation between the dimensions
of $V$ and $V_P$ is obvious from this definition. Conversely,
again using Corollary \ref{cor:struct}, we have
\begin{equation}
\begin{split}
1_P(\operatorname{Ind}_{\H^P}^\H V_P)^{an}&=
1_P(\H\otimes_{\H^P}V_P)^{an}\\
&=1_P(\sum_{P^\prime,P^{\prime\prime}}1_{P^\prime}
\H^{P,an}(U_P)1_{P^{\prime\prime}})\otimes_{\H^{P,an}(U_P)}1_PV^{an}\\
&=1_PV^{an}=V_P^{an},
\end{split}
\end{equation}
finishing the proof.
\end{proof}
\begin{prop}\label{luspardef}
Let $P\subset F_0$ be a subset, and let $R_P\subset R_0$ be the
corresponding {\it standard} parabolic subsystem. We define the
subtori $T_P$, $T^P$ and the lattices $X_P$, $Y_P$ as in
Proposition \ref{prop:red}. Put
$\Ri_P=(X_P,Y_P,R_P,R_P^\vee,F_P)$, and let $t\in T^P$. There
exists a surjective homomorphism
$\phi_{t}:\H^P\to\H_P$
\index{0v@$\phi_t:\H^P\to\H_P$, surjective homomorphism}
which is
characterized by (1) $\phi_t$ is the identity on the finite
dimensional Hecke algebra $\H(W_P)$, and (2)
$\phi_{t}(\theta_x)=x(t)\theta_{\overline{x}}$, where
$\overline{x}\in X_P$ is the natural image of $x$ in
$X_P=X/{}^PX=X/(X\cap{Y_P}^\perp)$.
\end{prop}
\begin{proof}
We have to check that $\phi_{t}$ is compatible with the
Bernstein-Zelevinski-Lusztig relations. Let $s=s_\a$ with $\a\in
P\subset F_0$. Then
\begin{gather}
\begin{split}
\theta_x&T_s-T_s\theta_{s(x)}=\\ &=\left\{
\begin{array}{ccc}
&(q_{\alpha^\vee}-1)\frac{\theta_x-\theta_{s(x)}}
{1-\theta_{-\alpha}}\ &{\rm if}\ 2\alpha\not\in\rnr.\\
&((q_{\alpha^\vee/2}q_{\alpha^\vee}-1) +q_{\alpha^\vee/2}^{1/2}
(q_{\alpha^\vee}-1)\theta_{-\alpha}) \frac{\theta_x-\theta_{s(x)}}
{1-\theta_{-2\alpha}}\ &{\rm if}\ 2\alpha\in\rnr.\\
\end{array}
\right.\\
\end{split}
\end{gather}
Since $s$ acts trivially on $T^P$, we have $x(t)=sx(t)$. This
implies the result.
\end{proof}
\begin{dfn}\label{dfn:ind}
Let $P\subset F_0$ be a subset. In this case we identify the
algebra $\H^P=\H(X,Y,R_P,R_P^\vee,P)$ with the subalgebra in
$\H$ generated by $\H(W_P)$ and $\C[X]$. Let $(V,\d)$ be a
representation of $\H_P$ with central character $W_Pr\in
W_P\backslash T_P$, and let $t\in T^P$. Denote by $\d_t$ the
representation $\d_t=\d\circ\phi_t$ of $\H^P$. We define a
representation $\pi(\Ri_P,W_Pr,\d,t)$ of $\H$ by
$\pi(\Ri_P,W_Pr,\d,t)=\operatorname{Ind}_{\H^P}^\H(\d_t)$
\index{0p@$\pi(\Ri_P,W_Pr,\d,t)=\operatorname{Ind}_{\H^P}^\H(\d_t)$,
parabolically induced representation}.
We refer
to such representations as parabolically induced representations.
\end{dfn}
\begin{cor}\label{cor:short}
Let $W_0t\in W_0\backslash T$, and let $R_P$ be a standard
parabolic subsystem of $R_0$. Suppose that there exists an $r\in
T_P$ and $t^P\in T^P$ such that $rt^P\in W_0t$ is an $R_P$-generic
point. The map $\delta\to\pi(\Ri_P,W_Pr,\d,t^P)$ gives an
equivalence between the representations of $\H$ with central
character $W_0t$ and the representations of $\H_P$ with central
character $W_Pr$.
\end{cor}
\begin{proof}
By Corollary \ref{cor:restr}, the induction functor from
representations of $\H^P$ to $\H$ gives rise to an equivalence
between the representations of $\H^P$ with $R_P$-regular central
character $W_Pt$ and the representations of $\H$ with central
character $W_0t$. If $\pi$ is a representation of $\H^P$ with
central character $W_Pt$, then it is easy to see that the
annihilator of $\pi$ contains the kernel of the homomorphism
$\phi_{t^P}$. Thus $\pi$ is the lift via $\phi_{t^P}$ of a
representation $\delta$ of $\H_P$. This
gives an equivalence between the category of representations of
$\H_P$ with central character $W_Pr$ and the representations of
$\H^P$ with central character $W_Pt$.
\end{proof}
The following proposition describes the induced modules analogous
to the ``compact realization'' of parabolically induced representations of
real reductive groups.
\begin{prop}\label{prop:unit}
Let $(V,\d)$ be an irreducible representation of $\H_P$
with central character
$W_Pr\in W_P\backslash T_P$ as before.
Suppose that $(V,\d)$ is unitary
with respect to an Hermitian inner product
$(\cdot,\cdot)$, and that $t\in T^P_u$.
We identify the underlying representation space $V_\pi$ of
$\pi:=\pi(\Ri_P,W_Pr,\d,t)$ with
\begin{equation}
V_\pi:=\H(W^P)\otimes V,
\end{equation}
where $\H(W^P)\subset H(W_0)$ denotes the subspace spanned by the
elements $T_w$ with $w\in W^P$.
Then $\pi$ is unitary with respect
to the Hermitian
inner product $\langle\cdot,\cdot\rangle$
\index{<@$\langle\cdot,\cdot\rangle$!inner product on $V_\pi$}
defined on $V_\pi$ by
(with $x,y\in \H(W^P)$, and $u,v\in V$):
\begin{equation}
\langle x\otimes u,y\otimes v\rangle:=
\tau(x^*y)(u,v).
\end{equation}
\end{prop}
\begin{proof}
The above form is clearly Hermitian and positive definite.
It remains to show that the inner product satisfies
\begin{equation}\label{eq:desire}
\langle \pi(h)m_1,m_2\rangle=\langle m_1,\pi(h^*)m_2\rangle
\end{equation}
for each $m_1,m_2\in V_\pi, h\in\H$. To this end we
recall Theorem 2.20 of \cite{EO}. Let
$i_s:\H\to\operatorname{End}(\H_0)$ denote the
minimal principal series representation induced from
$s\in T$. Then the nondegenerate sesquilinear pairing
defined on $\H_0\times \H_0$ by
\begin{equation}\label{eq:natinprod}
(x,y):=\tau(x^*y)
\end{equation}
satisfies the property
\begin{equation}\label{eq:starcomp}
(i_{s}(h)x,y)=(x,i_{s^*}(h^*)y)
\end{equation}
(see Theorem \ref{thm:ster} for the definition of $s^*$).
We have $\H_0=\H(W^P)\otimes \H(W_P)$, and the pairing
(\ref{eq:natinprod}) on $\H_0$ factors as the tensor product of the
pairings on $\H(W^P)$ and on $\H(W_P)$ which are also defined by
equation (\ref{eq:natinprod}) but with $x,y$ both in $\H(W^P)$
or both in $\H(W_P)$.

We choose $r\in T_P$ such that $V$ contains a simultaneous
eigenvector $v$ for $X_P$ with eigenvalue $r$.
Via $\d_t$, the vector $v\in V$ has eigenvalue $rt\in T$ with
respect to $X$. Thus there is a surjective
morphism of $\H^P$-modules $\alpha: \H(W_P)\twoheadrightarrow V$, where
$\H(W_P)$ is the minimal principal series module $i_{rt}^P$
for $\H^P$, and $V$ is the representation space of $\d_t$.
By the above, applied to $\H^P$, we have the adjoint injective
morphism $\a^*:V\hookrightarrow \H(W_P)$, where the action on
$\H(W_P)$ is via $i^P_{r^*t}$ (since $(rt)^*=r^*t$, because
$t\in T^P_u$).
By the exactness and the transitivity of induction we get
morphisms of $\H$-modules
$\operatorname{Ind}(\a):i_{rt}\twoheadrightarrow\pi$
and $\operatorname{Ind}(\a^*):\pi\hookrightarrow
i_{r^*t}$. Notice that $\operatorname{Ind}(\a)=
\operatorname{Id}_{\H(W^P)}\otimes \a$ and similarly,
$\operatorname{Ind}(\a^*)=
\operatorname{Id}_{\H(W^P)}\otimes \a^*$.
By the factorization of the pairing
(\ref{eq:natinprod}) we see that $\operatorname{Ind}(\a)$ and
$\operatorname{Ind}(\a^*)$ are adjoint with respect to
the pairings $\langle\cdot,\cdot\rangle$ on $V_\pi$ and
(\ref{eq:natinprod}) on $\H_0$. This, the surjectivity
of $\operatorname{Ind}(\a)$ and (\ref{eq:starcomp}) gives the
desired result
(\ref{eq:desire}).
\end{proof}
\begin{prop}\label{prop:indtemp}
With the notations as above, assume that $(V,\d)$ is a
tempered representation with central character
$W_Pr\in W_P\backslash T_P$ and that $t\in T^P_u$.
Then $\pi:=\pi(\Ri_P,W_Pr,\d,t)$ is a tempered
representation of $\H$.
\end{prop}
\begin{proof}
Recall that we have the identification
\begin{equation}
V_\pi:=\H(W^P)\otimes V,
\end{equation}
where $\H(W^P)\subset H(W_0)$ denotes the subspace spanned by the
elements $T_w$ with $w\in W^P$.
Recall from the proof of Lemma \ref{lem:cas} that
we can find a basis $(v_j)$ of $V$ such that $X_P$ acts by upper
triangular matrices with respect to this basis. By Casselman's
criterion, the diagonal entries are characters $r_{j,j}\in W_Pr$
of $X_P$ which satisfy $|x(r_{i,i})|\leq 1$ for $x\in X_{P}^+$.
Let $(w_i)$ denote an ordering of
the set $W^P$ such that the length $l(w_i)$ increases with $i$. We
take the tensors $T_{w_i}\otimes v_j$, ordered lexicographically,
as a basis for the representation space of $\pi$.
From a direct application of the Bernstein-Zelevinski-Lusztig
relations we see that the $\theta_x$ are simultaneously upper
triangular in this basis, and that the diagonal entries are
the elements $w_i(tr_{j,j})$. Since $w_i\in W^P$, $t\in T^P_u$
and since the vector part of $r_{j,j}$ is an element of the cone
generated by the negative roots of $R_P^\vee$, it follows from
a well known characterization of $W^P$ that the vector parts of
these diagonal entries are in the antidual of the positive
chamber. Again using Lemma \ref{lem:cas} we conclude that $\pi$ is
tempered.
\end{proof}
\subsection{The tracial states $\chi_t$ and parabolic
induction}\label{sub:states}
In this subsection we will compute the states $\chi_t$ on
$W_0$-orbits of tempered residual cosets of positive dimension in
terms of the characters of unitary representations which are
induced from discrete series characters of parabolic subalgebras,
as was discussed in Subsection \ref{sub:lus}.

Let $L$ be a residual coset such that $R_L\subset R_0$ is a
standard parabolic subsystem. In other words, $F_L\subset F_0$.
Let us denote by $\H_L$ the affine Hecke algebra with root datum
$\Ri_L:=(X_L,Y_L,R_L,R_L^\vee,F_L)$ and root labels $q_L$.
\index{H4@$\H_L:=\H(\Ri_L,q_L)$, semisimple quotient of $\H^L$}
(see Proposition
\ref{prop:red}). Let $r_L=c_Ls_L\in T_L$ be the corresponding
residual point of $\Ri_L$.
\begin{lem}\label{lem:eisext}
Let $U\subset T$ be a nonempty $W_0$-invariant open subset. Let
$t\in U$. There exists a unique extension of the Eisenstein
functional (cf. equations \ref{fundeis}) $E_t$ (which we will also denote by
$E_t$) to the localization $\H^{an}(U)$, such that
$E_t(fh)=E_t(hf)=f(t)E_t(h)$ for all $f\in\A^{an}(U)$.
\end{lem}
\begin{proof}
The functional $E_t$ factors to a functional of the finite
dimensional $\C$-algebra
$\H^t:=\H/\I_t\H$
\index{H7@$\H^t=\H/\I_t\H$, where $\I_t$ is the maximal ideal
of $W_0t$ in $\Ze$},
where $\I_t$ is the maximal
ideal in $\Ze$ corresponding to $W_0t$. We have
$\H/\I_t\H=\H^{an}(U)/\I^{an}_t(U)\H^{an}(U)$ for $t\in U$, and
this defines the extension with the required property uniquely.
\end{proof}
\begin{lem}\label{lem:pareis}
Let $L$ be such that $R_L\subset R_0$ is a standard parabolic
subset of roots, and let $t_0\in T$ be $R_L$-generic.
Set $U=W_0t_0\exp(B)$ with $B$
satisfying the conditions \ref{cond} (i), (ii), and (iii).
As before, we put
$U_L=W_Lt_0\exp(B)$.

We denote by
$E^L_t$
\index{E1@$E^L_t$, Eisenstein functional of $\H^L$}
the Eisenstein functional of the subalgebra
$\H^L\subset \H$. For $t_L\in T_L$ we write
$E_{L,t_L}$
\index{E1@$E_{L,t_L}$, Eisenstein functional of $\H_L$}
to denote
the Eisenstein functional at $t_L\in
T_L=\operatorname{Hom}(X_L,\C^\times)$ of the algebra $\H_L$. Let
$t=t^Lt_L\in U$ with $t^L\in T^L$ and $t_L\in T_L$. Recall $1_L$
is the characteristic function of $U_L$. We have,
for all $h\in \H^L$:
\begin{enumerate}
\item[(i)] $E_t(1_Lh1_L)={q(w^L)1_L(t)\Delta^L(t)}E^L_t(h)$.
\item[(ii)] $E_t^L(h)=E_{L,t_L}(\phi_{t^L}(h))$.
\end{enumerate}
\end{lem}
\begin{proof}
Because these are both equalities of holomorphic functions of $t\in
U$ it suffices to check them for $t$ regular, and outside the
union of the residual cosets (in other words,
$c(t)c(t^{-1})\not=0$).

(i). By the defining properties \ref{fundeis} and \cite{EO},
2.23(4) we need only to show that the left hand side satisfies the
properties $E_t(1_Lxh1_L)= E_t(1_Lhx1_L)=t(x)1_L(t)E_t(h)$ and
$E_t(1_L)=q(w_0)1_L(t)\Delta(t)$. These facts follow from
Lemma \ref{lem:eisext}.

(ii). We see that
\begin{align}\nonumber
E_{L,t_L}(\phi_{t^L}(\theta_xh))&=E_{L,t_L}(x(t^L)
\theta_{\overline{x}}\phi_{t^L}(h))\\\nonumber
&=x(t^L)\overline{x}(t_L)E_{L,t_L}(\phi_{t^L}(h))\\\nonumber
&=x(t)E_{L,t_L}(\phi_{t^L}(h)),
\end{align}
and similarly for $E_{L,t_L}(\phi_{t^L}(h\theta_x))$. The value at
$h=1$ is equal to $q(w_L)\Delta_L(t)$ on both the left and the
right hand side. As in the proof of (i), this is enough to prove
the desired equality.
\end{proof}
\begin{theorem}\label{thm:mainind}
Let $L$ be a residual coset such that $R_L\subset R_0$ is a
standard parabolic subset of roots. Let $t^L\in T^L_u$ be such
that $t:=r_Lt^L\in L^{temp}$ is $R_L$-generic. Notice that this
condition is satisfied outside a finite union of real codimension
one subsets in $T^L_u$. Let $\Delta_{\Ri_L,W_Lr_L}$ be complete set
of inequivalent irreducible representations of the residual algebra
$\overline{\H_L^{r_L}}$, and let
$\chi_{\H_L,r_L}=\sum_{\d\in\Delta_{\Ri_L,W_Lr_L}}
\chi_{\d}d_{\Ri_L,\d}$
be the corresponding decomposition in irreducible discrete series
characters of the tracial state $\chi_{\H_L,r_L}$ of ${\H_L}$.
\begin{enumerate}
\item[(i)] For all $\d$, $\pi(\Ri_L,W_Lr_L,\d,t^L)$ is
irreducible, unitary and tempered with central character $W_0t$.
These representations are mutually inequivalent.
\item[(ii)] We have
\begin{equation}|W^L|\chi_t=
\sum_{\d\in\Delta_{\Ri_L,W_Lr_L}}\chi_{\Ri_L,W_Lr_L,\d,t^L}
d_{\Ri_L,\d},
\end{equation}
where $\chi_{\Ri_L,W_Lr_L,\d,t^L}$
\index{0w@$\chi_{\Ri_L,W_Lr_L,\d,t^L}$, character of
the induced representation $\pi(\Ri_L,W_Lr_L,\d,t^L)$}
denotes the character of
$\pi(\Ri_L,W_Lr_L,\d,t^L)$. In particular, the
constants $d_{t,i}$ as in Definition \ref{dfn:resalg} are
independent of $t$.
\item[(iii)] For all (not necessarily $R_L$-generic)
$t=r_Lt^L\in L^{temp}$,
the character $\chi_{\Ri_L,W_Lr_L,\d,t^L}$ is a positive trace on
$\H$. Consequently, the irreducible subrepresentations of
$\pi(\Ri_L,W_Lr_L,\d,t^L)$ extend to $\cf$.
\end{enumerate}
\end{theorem}
\begin{proof}
(i). This is a direct consequence of Corollary \ref{cor:short},
Proposition \ref{prop:unit} and Proposition \ref{prop:indtemp}.

(ii). Recall the definition of the states $\chi_t$. Recall that
the support of the measure $\nu$ is the union of the tempered
residual cosets. We combine Definition \ref{dfn:Y}, Proposition
\ref{prop:dfn}, and Proposition \ref{prop:cycle} to see that (with
$N_L$ the stabilizer of $L$ in $W_0$)
\begin{align}\label{inteq}
\frac{|W_0|}{|N_L|}\int_{t\in
L^{temp}}z(t)&\chi_t(h)d\nu_L(t)\\\nonumber&=\sum_{M\in
W_0L}\int_{t^M\in T^M_u}\int_{t^\prime\in t^M\e^M\xi_M}z(t^\prime)
E_{t^\prime}(h)\eta(t^\prime)\\\nonumber
\end{align}
for all $h\in\H$ and $z\in\Ze$. Rewrite the right hand side as
\begin{equation}\label{eq:av}
\frac{k_L}{|N_L|}\int_{t^L\in T^L_u}\sum_{w\in
W_0}J_{w,L}(\e^{wL}r_{wL}w(t^L))d^L(t^L)
\end{equation}
where the inner integrals equal, with $s\in wT^L$,
\begin{equation}\label{eq:J}
J_{w,L}(r_{wL}s)=\int_{t^\prime\in s\xi_{wL}}z(t^\prime)
\frac{E_{t^\prime}(h)}{q(w_0)\Delta^{wL}(t^\prime)}
m^{wL}(t^\prime)\frac{1}{\Delta_{wL}(t^\prime)}\om_{wL}(t^\prime),
\end{equation}
where
\begin{equation}
\om_{wL}(t):=
\frac{d_{wL}(t_{wL})}{q(w_L)c_{R_L}(w^{-1}t)c_{R_L}(w^{-1}t^{-1})}.
\end{equation}
Hence $J_{w,L}(r_{wL}s)$ is a linear combination of
(possibly higher order) partial derivatives
(in the direction of $wT_L$) of the kernel
\begin{equation}\label{eq:kern}
z(t^\prime) \frac{E_{t^\prime}(h)}{q(w_0)\Delta^{wL}(t^\prime)}
m^{wL}(t^\prime),
\end{equation}
evaluated at $r_{wL}s$. The $N_L$-invariant
measure on $L^{temp}$ on the left hand side of \ref{inteq}
is thus obtained by taking
the boundary values $\e^{wL}\to 1$ of the
$J_{w,L}(\e^{wL}r_{wL}w(t))$, and then sum over the Weyl group as
in equation (\ref{eq:av}). Notice that the collection of
$R_L$-generic points in $L$ is the complement of a union of
algebraic subsets of $L$ of codimension $\geq 1$. The kernel
(\ref{eq:kern}) is regular at such points of $L$. The boundary
values at $R_L$-generic points are therefore computed simply by
specialization at $\e^{wL}=1$. We thus have
\begin{equation}\label{eq:sumJ}
z(t)\chi_t(h){\overline \ka}_{W_LL}m_L(t)=\frac{k_L}{|W_0|}\sum_{w\in W_0}
J_{w,L}(r_{wL}w(t^L)).
\end{equation}
The expression on the left hand side can be extended uniquely to
$z\in \Ze^{an}(U)$ and $h\in \H^{an}(U)$. By equation (\ref{eq:J}),
each summand in the expression on the right hand can also be
extended uniquely to such locally defined analytic $z$ and $h$.

Take $U=W_0t\exp(B)$. We restrict both sides to
$1_L\H^L1_L\subset\H^{L,an}(U_L)$. Substitute $h$ by $1_Lh1_L$
with $h\in \H^L$. On the left hand side we get, by Corollary
\ref{cor:restr},
\begin{equation}\label{eq:left}
\frac{1}{|W^L|}z(t)\chi_t^L(h){\overline \ka}_{W_LL}m_L(t)
\end{equation}
where $\chi_t^L$ is a central functional on $\H^L$, normalized by
$\chi_t^L(e)=1$.

On the other hand, by Lemma \ref{lem:pareis}, if $h=1_Lh1_L$ with
$h \in \H^L$ then
\begin{equation}
J_{w,L}(w(r_{L}t^L))=\int_{t^\prime\in w(t^L)\xi_{wL}}z(t^\prime)
\frac{E_{L,t^\prime_L}(\phi_{t^{\prime,L}}(h))}{q(w_L)\Delta_L(t^\prime_L)}
m^{wL}(t^\prime)\om_{wL}(t^\prime)
\end{equation}
if $w(r_Lt^L)\in U_L$, and $J_{w,L}(w(r_{L}t^L))=0$ otherwise.

Observe that, because of condition (ii) for $B$, $wt\in U_L$
implies that there exists a $w^\prime\in W_L$ such that
$wt=w^\prime t$. Since $t=r_Lt^L$ is $R_L$-generic, we see in
particular that the stabilizer of $t$ is contained in $W_L$. Thus
$wt\in U_L$ implies that $w\in W_L$, and hence that
$wt=w(r_L)t^L$.

Therefore the sum at the right hand side of equation (\ref{eq:sumJ})
reduces, if $h$ is of the form $1_Lh1_L$ with $h\in \H^L$, to
\begin{equation}
\frac{k_L}{|W_0|}\sum_{w\in W_L}\int_{t^\prime\in
\xi_{wL}}z(t^Lt^\prime)m^{wL}(t^Lt^\prime)
\frac{E_{L,t^\prime}(\phi_{t^{L}}(h))}{q(w_L)\Delta_L(t^\prime)}
\om_{wL}(t^\prime)
\end{equation}
The function $wr_L\exp(B\cap\mathfrak{t}_L)\ni t^\prime\to
m^{wL}(t^Lt^\prime)$ is $W_L$-invariant on
$W_Lr_L\exp(B\cap\mathfrak{t}_L)=(t^L)^{-1}U_L\cap T_L$, because
$m^L(t)$ is $W_L$-equivariant (i.e. $m^{wL}(wt)=m^L(t)$
when $w\in W_L$). In other words, this function is in the center
of $\H_L^{an}((t^L)^{-1}U_L\cap T_L)$. By Definition \ref{dfn:Y},
Corollary \ref{cor:pos}, Definition \ref{dfn:chi}, and Theorem
\ref{thm:nu} applied to $\H_L$ therefore, this sum reduces to
\begin{equation}
\frac{|W_L|}{|W_0|}z(t)\chi_{\H_L,r_L}(\phi_{t^L}(h))
m^L(t){\overline \ka}_{W_LL}m_{\Ri_L,r_L}(r_L),
\end{equation}
which we can rewrite as
\begin{equation}
\frac{|W_L|}{|W_0|}z(t)\chi_{\H_L,r_L}(\phi_{t^L}(h))
{\overline \ka}_{W_LL}m_L(t).
\end{equation}
Comparing this with (\ref{eq:left}) we see that,
in view of equation (\ref{eq:packet}), this implies that for $h\in
\H^L$,
\begin{equation}
|W^L|\chi_t(1_Lh)=\sum_{\d\in\Delta_{\Ri_L,W_Lr_L}}
\chi_{\d}(\phi_{t^L}(h))d_{\Ri_L,\d}=\chi^L_t(h).
\end{equation}
Applying Corollary \ref{cor:restr}, Definition \ref{dfn:ind}, and
Corollary \ref{cor:short} we obtain (ii).

(iii). By Corollary \ref{cor:exten}, $\chi_t$ is $\nu$-almost
everywhere a positive trace. On the set of $R_L$-generic points
$t\in L^{temp}$, we have expressed $\chi_t$ as a positive linear
combination of the irreducible induced characters
$\chi_{\Ri_L,W_Lr_L,\d,t^L}$.
This gives the decomposition
of $\chi_t$ in terms of irreducible characters of the finite
dimensional algebra $\H^t:=\H/J_t$, where $J_t$ denotes the two-sided
ideal of $\H$ generated by the maximal ideal $\I_t$ of the
elements of the center $\Ze$ which vanish at $W_0t$.

On the other
hand, we have the decomposition of $\chi_t$ in irreducible
characters of the finite dimensional Hilbert algebra
$\overline{\H^t}$, as in Definition \ref{dfn:resalg}. This algebra
is a quotient algebra of $\H^t$. Because $\H^t$ is finite
dimensional, there is no distinction between topological and
algebraic irreducibility. We therefore have two decompositions of
$\chi_t$ in terms of irreducible characters of $\H^t$. The
irreducible characters of $\H^t$ are linearly independent, and
thus the two decompositions are necessarily the same. This implies
that the induced characters $\chi_{\Ri_L,W_Lr_L,\d,t^L}$ are
characters of the Hilbert algebra $\overline{\H^t}$. In
particular, the characters are positive traces for all
$R_L$-generic $t\in L^{temp}$.
These induced characters are regular functions of $t\in L^{temp}$.
Hence by continuity, they are positive traces for all $t^L$.

By Corollary \ref{cor:contrace}(i), $\chi(\Ri_L,W_Lr_L,\d,t^L)$
extends to a continuous trace of $\cf$ for all $t\in L^{temp}$.
According to the construction in \cite{dix2}, Paragraphe 6.6.
this character is the trace
of a (obviously finite dimensional) representation of $\cf$,
quasi-equivalent with $\pi(\Ri_L,W_Lr_L,\d,t^L)$ when restricted
to $\H$. Hence all subrepresentations of $\pi(\Ri_L,W_Lr_L,\d,t^L)$
extend to $\cf$.
\end{proof}
\begin{cor}\label{cor:xcont}
For all $x\in\H$, the $\nu_L$-measurable function $L^{temp}\ni
t\to\chi_t(x)$ can be defined by the restriction to $L^{temp}$ of
a regular function on $L$. For $x\in\cf$, the function
$t\to\chi_t(x)$ is continuous on $L^{temp}$.
\end{cor}
\begin{proof}
The first assertion was shown in the proof of Theorem
\ref{thm:mainind}(iii). If $x\in\cf$ there exists a sequence
$(x_i)$ with $x_i\to x$ and $x_i\in \H$. By Corollary
\ref{cor:contrace}(i), the function $t\to\chi_t(x)$ is a uniform
limit of the functions $t\to\chi_t(x_i)$, proving the continuity.
\end{proof}
The next Theorem basically is the Plancherel decomposition of
$\hf$. (In the next subsection we will refine the formula by adding
more details about the spectrum of $\cf$.)
\begin{thm}\label{thm:thisisathm}
We have the following isomorphism of Hilbert spaces:
\begin{equation}\label{eq:intdec}
\Hf=\int_{W_0\backslash T}^\oplus\overline{\H^t}|W_0t|d\nu(t).
\end{equation}
The support of the probability measure $\nu$ is the union of the
tempered residual cosets. If $t=r_Lt^L\in L^{temp}$ is
$R_L$-generic, then the residue  algebra $\overline{\H^t}$ has the
structure
\begin{equation}\label{eq:struct}
\overline{\H^t}\simeq(\overline{\H^{r_L}_L})_{|W^L|}.
\end{equation}
Finally, the residue  algebra $\overline{\H^{r}}$ at a residual
point $r\in T$ is of the form
\begin{equation}
\overline{\H^{r}}=\bigoplus_{\d\in\Delta_{\Ri,W_0r}}\operatorname{End}(V_{\d})
\end{equation}
with the Hermitian form on the summand
$\operatorname{End}(V_{\d})$ given by
\begin{equation}\label{eq:innerp}
(A,B)=d_{\d}\operatorname{trace}(A^*B),
\end{equation}
where the positive real numbers $d_{\d}$ are defined as in
Definition \ref{dfn:resalg} (with the notational convention that
$d_{\d}=d_{r,i}$ if $\d=\d_i$ is the irreducible representation
of $\overline{\H^{r}}$ corresponding to the central idempotent
$e_i$).
\end{thm}
\begin{proof}
The Hilbert space $\hf$ is the completion of $\H$ with respect to
the positive trace $\tau$. In  Corollary \ref{cor:exten}(v) we
have written $\tau$
as a positive superposition of positive traces $\chi_t$, with
$t\in W_0\backslash T$.
In Theorem \ref{thm:mainind}
we established that, outside a set of measure zero, $\chi_t$ is a
finite linear combination of traces of irreducible representations
$\pi_{\Ri_L,W_Lr_L,\d,t^L}$ which extend to $\cf$. Thus Corollary
\ref{cor:exten}(v) is a positive decomposition of $\tau$ in
terms of traces of elements of $\hat{\cf}$. This allows us to apply
8.8.5 and 8.8.6 of \cite{dix2} in order to obtain (\ref{eq:intdec}).
The statements about the residual algebra of a residual point follow
directly from the Definition \ref{dfn:resalg}. Finally, equation
(\ref{eq:struct}) follows from Theorem \ref{thm:mainind} in
combination with the factorization Proposition \ref{prop:unit} of
the inner product $\langle\cdot,\cdot\rangle$ on the induced
representations.
\end{proof}
\subsection{The generic residue of the Hecke algebra}\label{sub:gene}
In this subsection we will use Theorem \ref{thm:mainind} in order
to compute explicitly the local traces $\chi_t$ when $t=r_Lt^L$ is
an $R_L$-generic element of $L^{temp}$. Here we assume that $L$ is
a residual coset such that $R_L$ is a standard parabolic
subset of $R_0$ with basis $F_L$ of simple roots.
Since we assume that $t$ is $R_L$-generic, we
have $P(t)=F_L$. As before, we put $W_Lt=\varpi$,
the equivalence class of $t$. By the genericity of $t$,
$W(\varpi)=1$.

Observe that the residual algebra $\overline{\H^t}$
(see Definition \ref{dfn:resalg}) is of the form
$\overline{\H^t}=\H^t/\operatorname{Rad}_t$,
where $\H^t$ is the quotient algebra
$\H^t:=\H/\I_t\H$ (with $\I_t$ the maximal ideal of $\Ze$
corresponding to $t$), and where
$\operatorname{Rad}_t$ is the radical of the positive
semi-definite form
$(x,y)_t:=\chi_t(x^*y)$
\index{$(x,y)_t=\chi_t(x^*y)$, semi-definite Hermitian form on $\H^t$}
(viewed as
a form on $\H^t$).

By Lusztig's Structure Theorem \ref{thm:lusind}, $\H^t$ has the
following decomposition in the case where $t$ is $R_L$-generic:
\begin{equation}
\H^t=\bigoplus_{u,v\in W^L}\i^0_ue_\varpi\H^L \i^0_{v^{-1}}
\end{equation}
where $e_\varpi=1_\varpi\ \operatorname{mod}(\I_t\H^{an}(U))$,
the image of $1_\varpi$ in $\H^t$
\index{e1@$e_\varpi$, image of $1_\varpi$ in $\H^t$}
(in the notation of
Section \ref{sect:loc}). We remark that this is {\it not} an
orthogonal decomposition with respect to $(x,y)_t$.

The subspace $\i^0_ue_\varpi\H^L \i^0_{v^{-1}}$ is equal to
$e_{u\varpi}\H^t e_{v\varpi}$. If $u=v$ then this is
a subalgebra of $\H^t$. If $u=v=e$ then this subalgebra is
isomorphic to $\H^{L,t}$ via the isomorphism
$\H^{L,t}\ni x\to e_\varpi x\in e_\varpi\H^L$.

Let $P,Q\subset F_0$ be subsets. We denote by
$W(P,Q)$ the following set of Weyl group elements:
$W(P,Q):=\{w\in W_0\mid w(P)=Q\}$
\index{W5a@$W(P,Q)=\{w\in W_0\mid w(P)=Q\}$,
with $P,Q\subset F_0$}.
If $P=Q\subset F_0$ then we simply write
$W(P)=W(P,P)$
\index{W5b@$W(P)$ for the stabilizer in $W_0$ of $P\subset F_0$}

Let $n\in W_0$ be such that $n(F_L)=F_M\subset F_0$, in other
words, let $n\in W(F_L,F_M)$. Then the map
\begin{align}\label{eq:delta}
\Delta_{n\varpi,n\varpi}: e_\varpi\H^L&\to e_{n\varpi}\H^M\\
\nonumber        x&\to \i^0_n x \i^0_{n^{-1}}
\end{align}
\index{0D4@$\Delta_{\varpi_1,\varpi_2}:\H^{\varpi,{an}}(U_\varpi)\to
1_{\varpi_1}\H^{an}(U)1_{\varpi_2}$, linear isomorphisms}
is an isomorphism. By the results of Lusztig \cite{Lu}, section 8,
it satisfies (with $h\in \H^L$):
\begin{equation}\label{eq:psi}
\Delta_{n\varpi,n\varpi}(e_\varpi h)=e_{n\varpi} \psi_n(h),
\end{equation}
where $\psi_n:\H^L\to \H^M$
\index{0x@$\psi_g:\H^L\to\H^M$, isomorphism for
$g\in K_M\times W(F_L,F_M)$}
is the isomorphism of algebras defined
by (with $w\in W_L$ and $x\in X$):
\begin{equation}
\psi_n(T_w)=T_{nwn^{-1}},
\psi_n(\theta_x)=\theta_{nx}.
\end{equation}

Recall that
Theorem \ref{thm:mainind}(ii) implies that for all $h\in\H^L$,
\begin{equation}\label{eq:bel}
|W^L|\chi_t(e_\varpi h)=\chi_{\H_L,r_L}\bigl(\phi_{t^L}(h)\bigr).
\end{equation}
\begin{cor}\label{cor:con}
Write $n(t)=s=r_M^\prime s^M$.
\begin{enumerate}
\item
\begin{equation}
\chi_{\H_M,r_M^\prime}\bigl(\phi_{s^M}(\psi_n(h))\bigr)=
\chi_{\H_L,r_L}\bigl(\phi_{t^L}(h)\bigr).
\end{equation}
\item Let
$\Psi_n:\Delta_{\Ri_L,W_Lr_L}\to\Delta_{\Ri_M,W_Mr_M^\prime}$
\index{0X@$\Psi_g:\Delta_{\Ri_L,W_Lr_L}\to\Delta_{\Ri_M,W_Mr_M^\prime}$,
bijection induced by $\psi_g$}
be the bijection such that $\Psi_n(\d)\simeq \d\circ \psi_n^{-1}$.
Then $\Psi_n$ respects the residual degree:
$d_{\Ri_L,\d}=d_{\Ri_M,\Psi_n(\d)}$.
\end{enumerate}
\end{cor}
\begin{proof}
(i) This follows from the above text, and the fact that $\chi_t$
is a central functional.

(ii) It is clear from Casselman's criteria that $\Psi_n$ indeed
defines a bijection between the sets of discrete series representations
$\Delta_{\Ri_L,W_Lr_L}$ and $\Delta_{\Ri_M,W_Mr_M^\prime}$.
The result follows from (i) and Theorem \ref{thm:mainind}(ii).
\end{proof}
\begin{prop}\label{prop:constterm}
For all $h,g\in \H^L$ we have
\begin{equation}
|W^L|m^L(t)\chi_t\bigl((e_\varpi h)^*(e_\varpi g)\bigr)
=\chi_{\H_L,r_L}\bigl(\phi_{t^L}(h^\sharp g)\bigr),
\end{equation}
where $\sharp$
\index{$\sharp$, the $*$-operator of $\H^L$}
denotes the $*$-operator of $\H^L$ (thus
$T^\sharp_w=T_{w^{-1}}$ if $w\in W_L$, and $\theta_x^\sharp=
T_{w_L}\theta_{-w_Lx}T_{w_L}^{-1}$ where
$w_L$
\index{wb@$w_P$, longest element of $W_P$}
denotes the longest
element of $W_L$).
\end{prop}
\begin{proof}
This will be proved by the computation in the proof of Theorem
\ref{thm:unitint}.
\end{proof}
\begin{cor}\label{cor:corr}
We equip $e_\varpi\H^t e_\varpi$ with the positive semi-definite
sesquilinear pairing obtained
by restriction of the pairing
$|W^L|m^L(t)(x,y)_t=|W^L|m^L(t)\chi_t(x^*y)$ defined on $\H^t$. The
residual algebra $\overline{\H^{L,t}}$ is isomorphic as a Hilbert
algebra to the quotient of $e_\varpi \H^t e_\varpi$ by the radical
of this pairing, via the map $x\to e_\varpi x$.
\end{cor}
%\begin{proof}
%This follows immediately from Proposition \ref{prop:constterm}.
%\end{proof}
\begin{cor}\label{cor:emb}
Recall the notation of Proposition \ref{prop:unit}.
We consider
$V_\pi$ as a module over $\H^t$. Put
$V_{\pi,\varpi}=\pi(e_\varpi)(V_\pi)$. Choose an isometric
embedding $\overline{i}:V\to\overline{\H_{L}^{r_L}}$
(as $\H_L$-modules).
Let $\overline{j}$ denote the unique module map
$\overline{j}:V_\pi\to\overline{\H^t}$
such that
$\overline{j}(1\otimes v)=e_\varpi
(\overline{\phi}_{t^L})^{-1}(\overline{i}(v))$,
where $\overline{\phi}_{t^L}:
\overline{\H^{L,t}}\to\overline{\H_L^{r_L}}$
denotes the isometric isomorphism
determined by the homomorphism $\phi_{t^L}$
(cf. Proposition \ref{luspardef}).
For any $v\in V$ we denote by $i(v)$ any lift of
$\overline{i}(v)$,
and similarly for $j$.
We have:
\begin{enumerate}
\item $V_{\pi,\varpi}=1\otimes V$, and
$\overline{j}:V_{\pi,\varpi}
\stackrel{\sim}{\to}e_\varpi
(\overline{\phi}_{t^L})^{-1}(\overline{i}(V))$.
\item The positive definite Hermitian form
$\langle\cdot,\cdot\rangle$ on $V_\pi$ (see Proposition
\ref{prop:unit}) can be expressed by:
\begin{equation}
\langle v,w\rangle=|W^L|m^L(t)\chi_t(j(v)^*j(w)).
\end{equation}
\end{enumerate}
\end{cor}
\begin{proof} (i) is straightforward by observing
that $\phi_{t^L}$ descends to $\H^{L,t}$ and so defines
an isometric isomorphism  $\overline{\phi}_{t^L}$
by (\ref{eq:bel}), applied to
$\Ri^L$ instead of $\Ri$ (thus $e_\om=1$, and $|W^L|=1$).
Since $V_\pi$ is irreducible, it is enough to compare the
inner products on $V_{\pi,\varpi}$ in order to prove (ii).
Apply Proposition
\ref{prop:constterm} and Corollary \ref{cor:corr}.
\end{proof}
Assume that $R_M$ and $R_L$ are associate standard parabolic
subsystems. Let $\varpi_1=\varpi$ and $\varpi_2$ be equivalence
classes inside $W_0t$ such that $\varpi_1=W_Lt$ and
$\varpi_2=W_Ms$.
\begin{thm} \label{thm:unitint}
Let $n\in W(F_M,F_L)$ be such that $n(\varpi_2)=\varpi_1$.
Let $h\in\H^L$ and $h^\prime\in\H^M$ such that
$e_{\varpi_2}h^\prime=\i^0_{n^{-1}}e_{\varpi_1}h \i^0_n\in
e_{\varpi_2}\H^M$. We have
\begin{align*}
\chi_t\bigl((he_{\varpi_1} \i^0_n)^*(he_{\varpi_1} \i^0_n)\bigr)&=
\chi_t\bigl((he_{\varpi_1})^*(he_{\varpi_1})\bigr)\\&=
\chi_t\bigl((h^\prime e_{\varpi_2})^*(h^\prime e_{\varpi_2})\bigr)\\
\end{align*}
\end{thm}
\begin{proof}
Before we embark on this computation we establish some
useful relations. First recall that (Theorem \ref{thm:ster})
$t^*:=\overline{t^{-1}}\in W_Lt$.
Also recall Proposition \ref{prop:thetastar}.
We see that
\begin{align}\label{eq:ms1}
e_{\varpi_1}^*&=T_{w_0}e_{w_0\varpi_1}T_{w_0}^{-1}\\
\nonumber&=T_{w^L}e_{w^L\varpi_1}T_{w^L}^{-1},
\end{align}
where $w^L$
\index{wc@$w^P$, longest element of $W^P$}
denotes the longest element of set of minimal coset
representatives $W^L$.
Next we observe that for any $w\in W_0$,
\begin{equation}\label{eq:ms2}
(\i^0_w)^*=T_{w_0}\left(
\prod_{{\a>0}\atop{w^\prime(\a)<0}}
\left(\frac{c_\a}{c_{-\a}}\right)\i^0_{{w^\prime}^{-1}}\right)T_{w_0}^{-1},
\end{equation}
where ${w^\prime}:=w_0 w w_0$. This formula follows in a
straightforward way from Definition (\ref{eq:defint}).

If $s$ is a simple reflection and $\varpi\subset W_0t$ an
equivalence class, we check that (recall that $t$ is
$R_L$-generic)
\begin{equation}\label{eq:ms3}
e_{s\varpi}T_se_\varpi=
\left\{
\begin{aligned}
&e_\varpi T_s \hskip70pt\text{\ if\ }s\varpi=\varpi\\
&e_{s\varpi}q(s)c_\a \i^0_s \hskip32pt\text{\ else.\ }\\
\end{aligned}
\right.
\end{equation}
Since we assume that $t$ is $R_L$-generic, we have
$u(\varpi_1)\not=w(\varpi_1)$
for all $w\in W^L$ and all $u\in W_0$ such that $l(u)<l(w)$. From
this, (\ref{eq:ms3}) and induction to the length of $l(w)$
we see that
\begin{align}\label{eq:ms4}
e_{w\varpi_1}T_we_{\varpi_1}&=
q(w)\Bigl(\prod_{{\a>0}\atop{w^{-1}(\a)<0}}c_\a\Bigr)
\i^0_w e_{\varpi_1}\\
\nonumber&=e_{w\varpi_1}
q(w)\Bigl(\prod_{{\a>0}\atop{w^{-1}(\a)<0}}c_\a\Bigr)
\i^0_w
\end{align}
for all $w\in W^L$.
Observe that we also have
\begin{equation}\label{eq:ms5}
e_{w\varpi_1}T_{w^{-1}}^{-1}e_{\varpi_1}=
\Bigl(\prod_{{\a>0}\atop{w^{-1}(\a)<0}}c_\a\Bigr)
\i^0_w e_{\varpi_1}.
\end{equation}
We note that $w_0=w^Lw_L$. Since $w_L$ and $w_0$ are involutions,
this implies that $(w^L)^{-1}=w^{L^\prime}$, where
$R_{L^\prime}=w_0(R_L)$ (also a standard parabolic subsystem).

Let $h=T_w\theta_x\in \H^L$, with $w\in W_L$ and $x\in X$.
Keeping in mind the above preliminary remarks, and the fact that
$\chi_t$ is central, we now compute (with $x^\prime:=-w_0(x)$
and $w_0\varpi_1=\varpi_1^\prime$):
\begin{gather}\label{eq:long}
\begin{split}
\chi_t&\bigl((he_{\varpi_1}\i_n^0)^*(he_{\varpi_1}\i_n^0)\bigr)\\
&=\chi_t\Bigl(T_{w_0}
\prod_{{\a>0}\atop{n^\prime(\a)<0}}
\left(\frac{c_\a}{c_{-\a}}\right)\i^0_{{n^\prime}^{-1}}
e_{\varpi_{1}^{\prime}}\theta_{x^\prime}T_{w_0}^{-1}T_{w^{-1}}T_w
\theta_xe_{\varpi_1} \i^0_n\Bigr)\\
&=\chi_t\Bigl(T_{w^{M^\prime}}
\prod_{{\a>0}\atop{n^\prime(\a)<0}}
\left(\frac{c_\a}{c_{-\a}}\right)
\i^0_{{n^\prime}^{-1}}
e_{\varpi_{1}^{\prime}}\theta_{x^\prime}T_{w^{L^\prime}}^{-1}
T_{w_L}^{-1}T_{w^{-1}}T_w
\theta_x e_{\varpi_1} \i^0_nT_{w_M}\Bigr)\\
&=\chi_t\Bigl(e_{\varpi_2}T_{w^{M^\prime}}
e_{\varpi_{2}^{\prime}}
\Bigl({{\prod_{{\a>0}\atop{w^{M^\prime}}(\a)<0}}{c_{\a}}}\Bigr)
\Bigl({{\prod_{{\a>0}\atop{w^{L^\prime}}(\a)<0}}
{c_{{n^\prime}^{-1}\a}^{-1}}}\Bigr)\\
&\hskip101pt \i^0_{{n^\prime}^{-1}}
\theta_{x^\prime}e_{\varpi_{1}^{\prime}}T_{(w^L)^{-1}}^{-1}e_{\varpi_1}
T_{w_L}^{-1}T_{w^{-1}}T_w
\theta_x e_{\varpi_1} \i^0_nT_{w_M}\Bigr)\\
&=q(w^{M^\prime})\chi_t\Bigl(e_{\varpi_2}
\i^0_{w^{M^\prime}}
\Bigl(\prod_{{\a>0}\atop{w^{M^\prime}\a<0}}c_{-\a}\Bigr)
\Bigl({{\prod_{{\a>0}\atop{w^{M^\prime}}(\a)<0}}{c_{\a}}}\Bigr)
\Bigl({{\prod_{{\a>0}\atop{w^{L^\prime}}(\a)<0}}
{c_{{n^\prime}^{-1}\a}^{-1}}}\Bigr)\\
&\hskip97pt \i^0_{{n^\prime}^{-1}}
\theta_{x^\prime}
\Bigl(\prod_{{\a>0}\atop{w^{L^\prime}}\a<0}c_\a\Bigr)
\i^0_{w^L}
T_{w_L}^{-1}T_{w^{-1}}T_w
\theta_x e_{\varpi_1} \i^0_nT_{w_M}\Bigr)\\
&=q(w^L)\chi_t\Bigl(e_{\varpi_2}
\i^0_{n^{-1}}\i^0_{w^{L^\prime}}
\Bigl(\prod_{{\a>0}\atop{w^{M^\prime}\a<0}}c_{-n^\prime\a}\Bigr)
\Bigl({{\prod_{{\a>0}\atop{w^{M^\prime}}(\a)<0}}{c_{n^\prime\a}}}\Bigr)
\Bigl({{\prod_{{\a>0}\atop{w^{L^\prime}}(\a)<0}}{c_\a^{-1}}}\Bigr)\\
&\hskip124pt
\Bigl(\prod_{{\a>0}\atop{w^{L^\prime}}\a<0}c_\a\Bigr)
\theta_{x^\prime}
\i^0_{w^L}T_{w_L}^{-1}T_{w^{-1}}T_w
\theta_x e_{\varpi_1} \i^0_nT_{w_M}\Bigr)\\
&=q(w^L)\chi_t\Bigl(e_{\varpi_2}
\i^0_{n^{-1}}
\Bigl(\prod_{{\a>0}\atop{w^{M^\prime}\a<0}}c_{-w^{L^\prime}n^\prime\a}\Bigr)
\Bigl({{\prod_{{\a>0}\atop{w^{M^\prime}}(\a)<0}}
{c_{w^{L^\prime}n^\prime\a}}}\Bigr)\\
&\hskip187pt \theta_{-w_Lx}T_{w_L}^{-1}T_{w^{-1}}T_w
\theta_x e_{\varpi_1} \i^0_nT_{w_M}\Bigr)\\
&=q(w^L)\chi_t\Bigl(e_{\varpi_1}
\Bigl(\prod_{\a\not\in R_M}c_{n\a}\Bigr)
e_{\varpi_1}\theta_{-w_Lx}
T_{w_L}^{-1}T_{w^{-1}}T_w
\theta_xT_{w_L}\Bigr)\\
&=q(w^L)
\Bigl(\prod_{\a\not\in R_L}c_{\a}(t)\Bigr)
\chi_t\Bigl(e_{\varpi_1}
T_{w_L}\theta_{-w_Lx}
T_{w_L}^{-1}T_{w^{-1}}T_w
\theta_x\Bigr)\\
&=m^L(t)^{-1}
\chi_t\bigl(e_{\varpi_1}
h^\sharp h\bigr)\\
&=|W^L|^{-1}m^L(t)^{-1}
\chi_{\H_L,r_L}
\bigl(\phi_{t^L}(h^\sharp h)\bigr).\\
\end{split}
\end{gather}
The result is independent of $n$,
implying the first equality of the theorem.
The second equality follows because $\chi_t$ is central.
Indeed, this implies that we have
\begin{equation}\label{eq:short}
\chi_t(e_{\varpi_1}h^\sharp h)
=\chi_t(e_{\varpi_2}(h^\prime)^\sharp h^\prime),
\end{equation}
where the second $\sharp$ of course refers to the
$*$-structure on $\H^M$.

In equation (\ref{eq:long}) we have used the evaluation
\begin{equation}
\Bigl(\prod_{\a\not\in R_L}c_{\a}\Bigr)e_{\varpi_1}=
\Bigl(\prod_{\a\not\in R_L}c_{\a}(t)\Bigr)e_{\varpi_1}.
\end{equation}
This is allowed because the element
\begin{equation}
\Bigl(\prod_{\a\not\in R_L}c_{\a}\Bigr)\in\A
\end{equation}
is $W_L$-invariant, and thus belongs to the center
of $\H^L$.

At several places in equations (\ref{eq:long}) and (\ref{eq:short})
we have freely used formulae of Lusztig \cite{Lu}
(see Theorem \ref{thm:lusind}) for
the structure of $\H^t$. For example,
\begin{equation}
\i^0_{n^{-1}} e_{\varpi_1} h \i^0_{n}=
e_{\varpi_2}\psi_{n^{-1}}(h)=e_\varpi h^\prime
\end{equation}
when $h\in \H^L$.
By this we easily see that for all $h\in \H^L$,
\begin{equation}
\i^0_{n^{-1}} e_{\varpi_1}h^\sharp
\i^0_{n}=e_{\varpi_2}(h^\prime)^\sharp,
\end{equation}
and hence we may conclude by equation (\ref{eq:long}) that
\begin{equation}
\chi_t\bigl((he_{\varpi_1})^*(he_{\varpi_1})\bigr)=
\chi_t\bigl((h^\prime e_{\varpi_2})^*(h^\prime
e_{\varpi_2})\bigr).
\end{equation}
The proof is finished.
\end{proof}
\begin{cor}\label{cor:hehe}
Let $L,M_1,M_2$ be residual cosets such that $F_L,F_{M_1}$ and
$F_{M_2}$
are associate subsets of $F_0$,
and let $n_i\in W(L,M_i)$ ($i=1,2$).
The map
\begin{align}
\Delta_{n_1\varpi,n_2\varpi}:e_\varpi\H^L&\to
\i^0_{n_1} e_\varpi\H^L
\i^0_{n_2^{-1}}=e_{n_1\varpi}\H e_{n_2\varpi}\\
\nonumber x&\to \i^0_{n_1}x \i^0_{n_2^{-1}}
\end{align}
is a partial isometry with respect to the natural
positive semi-definite pairing on $\H^t$ given
by $(x,y)_t:=\chi_t(x^*y)$.
\end{cor}
\begin{proof}
We have, with $\varpi^\prime:=n_1\varpi$,
$\Delta_{n_1\varpi,n_2\varpi}=
\Delta^{\varpi^\prime}_{\varpi^\prime,n_2n_1^{-1}\varpi^\prime}
\circ\Delta_{n_1\varpi,n_1\varpi}$, where
$\Delta^{\varpi^\prime}_{\varpi^\prime,n_2n_1^{-1}\varpi^\prime}$
is defined by
\begin{align}
\Delta^{\varpi^\prime}_{\varpi^\prime,n_2n_1^{-1}\varpi^\prime}
:e_{\varpi^\prime}\H^{M_1}&\to
e_{\varpi^\prime}\H^{M_1}\i^0_{n_1n_2^{-1}}\\
x&\to x\i^0_{n_1n_2^{-1}}
\end{align}
Both these respect the pairing $(\cdot,\cdot)_t$, by Theorem
\ref{thm:unitint}.
\end{proof}
\subsection{Unitarity and regularity of intertwining operators}
Let $L,M$ be associate residual subspaces such that $R_L, R_M$
are standard parabolic subsystems of $R_0$. Let $n\in W_0$ be
such that $n(R_{L,+})=R_{M,+}$.
As before we let $\psi_n:\H^L\to\H^M$ denote
the isomorphism defined by $\psi_n(T_w)=T_{nwn^{-1}}$
and $\psi_n(\theta_x)=\theta_{nx}$.
Let $(V,\d)$ be an irreducible discrete series
representation of $\H_{L}^{r_L}$ and let $t=r_Lt^L$ be
an $R_L$-generic point of $r_LT^L$. Let
$s=n(t)=r_M^\prime s^M$, and
let $(V^\prime,\d^\prime)$ be a realization of the discrete
series representation $\d^\prime=\Psi_n(\d)$.

Choose a {\it unitary} isomorphism $\tilde\delta:V\to V^\prime$
such that
\begin{equation}\label{eq:inttwist}
\tilde\d(\d_t(h)v)=\d^\prime_s(\psi_n(h))(\tilde\d(v)).
\end{equation}
Recall that $V_\pi$ with
$\pi=\pi(\Ri_L,W_Lr_L,\d,t^L)$ is isomorphic to
\begin{equation}
V_\pi\simeq \H^{an}(U)\otimes_{\H^{L,an}(U_\varpi)}V_{t^L}
\end{equation}
(see Subsection \ref{sub:lus}), where $V_{t^L}$ denotes the
representation space $V$ with $\H^L$ action defined by
$h\to\d(\phi_{t^L}(h))$. Put
$\pi^\prime=\pi(\Ri_M,W_Mr_M^\prime,\d^\prime,s^M)$.
\begin{dfn}
For $t^L\in T^L$ such that $r_Lt^L$ is $R_L$-generic,
we define an intertwining
isomorphism $A(n,\Ri_L,W_Lr_L,\d,t^L):V_\pi\to V_{\pi^\prime}$ by
\begin{align}
A(n,\Ri_L,W_Lr_L,\d,t^L):\H\otimes_{\H^L}V_{t^L}
&\to\H\otimes_{\H^M}V^\prime_{s^M}\\
\nonumber  h\otimes v&\to h \i^0_{n^{-1}}\otimes \tilde\d(v)
\end{align}
\index{A@$A(n,\Ri_L,W_Lr_L,\d,t^L)$, unitary intertwining operator
($n\in W(F_L,F_M)$)}
\end{dfn}
It is easy to check that this is well defined and that
this map intertwines the $\H$ actions.
\begin{thm}\label{thm:ms}
Recall the compact realization $V_\pi=\H(W^L)\otimes V$, with its
inner product $\langle\cdot,\cdot\rangle_\pi$ (see
Proposition \ref{prop:unit}).
\begin{enumerate}
\item In the ``compact realization'', the intertwining map
\begin{equation}
A(n,\Ri_L,W_Lr_L,\d,t^L):\H(W^L)\otimes V\to\H(W^M)\otimes V^\prime
\end{equation}
is rational as a function of induction parameter $t^L$, and
regular outside the set of zeroes of the functions
$t^L\to\Delta_\a c_\a(u(r_L)t^L)$, where $\a$ runs over the set
of positive roots in $R_1$ such that $n(\a)<0$,
and $u(r_L)$ (with $u\in W_L$) runs over the set of
$X_L$-weights in $V$.
\item When $t^L\in T^L_u$ and
$A(n,\Ri_L,W_Lr_L,\d,t^L)$ is regular at $t^L$, then in fact
$A(n,\Ri_L,W_Lr_L,\d,t^L)$ is unitary with respect to the inner products
$\langle\cdot,\cdot\rangle_\pi$ and
$\langle\cdot,\cdot\rangle_{\pi^\prime}$.
\item With respect to these inner products we have
\begin{equation}
A^*(n,\Ri_L,W_Lr_L,\d,t^L)=
A(n^{-1},\Ri_M,W_Mr_M^\prime,\Psi_n(\d),n(t^L)).
\end{equation}
\end{enumerate}
\end{thm}
\begin{proof}
(i) The representation $\pi$ is cyclic and generated by
$1\otimes v$, with $v\not=0$ an arbitrary vector in $V$.
By the intertwining property it is therefore enough
to show that
$A(n,\Ri_L,W_Lr_L,\d,t^L)(1\otimes v)\in \H(W^L)\otimes V$ is meromorphic
in $t^L$, and regular outside the indicated set.
Using equation (\ref{eq:defint}), we have
\begin{equation}
A(n,\Ri_L,W_Lr_L,\d,t^L)(1\otimes v)=\pi(\i_{n^{-1}})
\pi(\prod_{\a>0,n(\a)<0}\Delta_\a c_\a)^{-1}(1\otimes \tilde\d(v)).
\end{equation}
Since $\pi(h)$ is a regular function on $T^L$ for
all $h\in \H$, this is a rational expression. The generalized
$X$-weights in $1\otimes V_{t^L}$ are of the form
$u(r_L)t^L$. So the inverse
of $\pi(\prod_{\a>0,n(\a)<0}\Delta_\a c_\a)$ can have poles
only at the indicated set.

(ii) In order to see the unitarity, we first note that by
Corollary \ref{cor:emb}(ii) and Theorem \ref{thm:unitint},
the statement is equivalent to the unitarity with respect
to the inner products on $V_\pi$ and $V_{\pi^\prime}$
defined by the embedding of these spaces in $\overline{\H^t}$
as in Corollary \ref{cor:emb}.
Choose an embedding $\overline{i}:V\to\H_{L}^{r_L}$ as in
Corollary \ref{cor:emb}. By Theorem \ref{thm:unitint} and
Corollary \ref{cor:con}, the map
$\Delta_{n\varpi,n\varpi}|_{\overline{j}(V_{\pi,\varpi})}$
is an isometry. By equation (\ref{eq:psi}) we see that this
isometry satisfies, for $h\in H^L$,
$\psi_n(h)\cdot\Delta_{n\varpi,n\varpi}(\overline{j}(1\otimes v))
=\Delta_{n\varpi,n\varpi}(\overline{j}(1\otimes\d_{t^L}(h)(v)))$. Hence
if we identify $V_{t^L}$ with $\overline{j}(V_{\pi,\varpi})$, we
can define $V^\prime_{n(t^L)}=\Delta_{n\varpi,n\varpi}(V_{t^L})$.
Then the map $\tilde\d=\Delta_{n\varpi,n\varpi}$ defines a unitary
map satisfying (\ref{eq:inttwist}).

Now it is clear, in the notation of Corollary \ref{cor:emb},
that we can identify the space
$\overline{j^\prime}(V_{\pi^\prime,\varpi^\prime})$ with
$\overline{j}(V_{\pi, \varpi})\i^0_{n^{-1}}$, and the map
$A(n,\Ri_L,W_Lr_L,\d,t^L)$ is then identified with the right multiplication
with $\i^0_{n^{-1}}$, thus with $\Delta_{\varpi,n\varpi}$.
This is unitary on $\overline{j}(V_{\pi, \varpi})$,
by Corollary \ref{cor:hehe}.
By the irreducibility of $V_\pi$ and $V_{\pi^\prime}$
this concludes the proof of (ii).

(iii) This last assertion of the Theorem is now obvious, since these
maps are clearly inverse to each other.
\end{proof}
The next Corollary is an important classical application of the
unitarity of the intertwiners, see \cite{BCD}, Th\'eor\`eme 2.
\begin{cor}\label{cor:hol}
The intertwining map $t^L\to A(n,\Ri_L,W_Lr_L,\d,t^L)$
extends holomorphically
to an open neighborhood of $T^L_u$ in $T^L$.
\end{cor}
\begin{proof}
By the unitarity on $T^L_u$, the meromorphic matrix entries of
$A(n,\Ri_L,W_Lr_L,\d,t^L)$
are uniformly bounded for $t^L\in T^L_u$ in the open set
of $T^L_u$ where $A(n,\Ri_L,W_Lr_L,\d,t^L)$ is well defined
and regular. This is the complement of the collection of
real codimension $1$ cosets in $T^L_u$ as described in
Theorem \ref{thm:ms}. This implies that the singularities of the
matrix entries which meet $T^L_u$ are actually removable.
\end{proof}
\subsection{The Plancherel decomposition of
the trace $\tau$}\label{sub:pla}
In this section we rewrite the decomposition
Theorem \ref{thm:thisisathm} as a decomposition of $\tau$
in terms of characters of irreducible tempered representations
induced from cuspidal representations of the subalgebras $\H^P$.

Using the results of the previous section, we show that
the corresponding Fourier homomorphism maps $\H$ into a certain space
of smooth sections defined over orbits of irreducible cuspidal
representations of the subalgebras $\H^P$, equivariant with respect
to the natural actions of intertwining operators.

This final formulation of the results (Theorem \ref{thm:mainp})
is inspired by and parallel to the
notations used in the theory of the Harish-Chandra
Plancherel formula for p-adic groups, as treated in
\cite{W2} and \cite{D}.

We need to develop some notations.
Let $P\subset F_0$
\index{P@$\mathcal{P}$, power set of $F_0$} denote the power
set of $F_0$, and let $\Gamma$
\index{0C@$\Gamma$, set of all pairs $(\Ri_P,W_Pr)$ with $P\in\mathcal{P}$
and $W_Pr$ an orbit of residual points in $T_P$}
denote the set of all pairs $\gamma=(\Ri_P,W_Pr)$ with
$P\in\mathcal{P}$, $\Ri_P$ the associated parabolic root datum,
and $W_Pr$ an orbit of residual points in $T_P$.
We consider the disjoint union of the set of all triples
$\Lambda=\{(\Ri_P,W_Pr,t)\}$
\index{0L@$\Lambda$, disjoint union
of all triples of the form $\{(\Ri_P,W_Pr,t)\}$},
where $(\Ri_P,W_Pr)\in\Gamma$ and $t\in T^P_u$. Let
$\Lambda_{\Ri_P,W_Pr}=\Lambda_\gamma$
\index{0L2@$\Lambda_{\Ri_P,W_Pr}=\Lambda_\gamma$,
triples in $\Lambda$ with
$\gamma=(\Ri_P,W_Pr)\in\Gamma$ fixed}
be the
subspace of such triples with $\gamma=(\Ri_P,W_Pr)\in\Gamma$ fixed.
Hence
for all $\gamma\in\Gamma$, $\Lambda_\gamma$ is a copy of $T^P_u$ and
$\Lambda=\cup\Lambda_\gamma$ (disjoint union).
Therefore $\Lambda$ is a
disjoint union of finitely many compact tori, which gives
$\Lambda$ the structure of a compact Hausdorff space. In addition,
each $\Lambda_\gamma$ comes with its (normalized) Haar measure, thus
defining a measure on $\Lambda$. We denote by
$\Lambda_\gamma^{gen}$
\index{0L3@$\Lambda_\gamma^{gen}$, triples in
$\Lambda_\gamma$ with $rt$ $R_P$-generic}
the open, dense subset of triples
$(\Ri_P,W_Pr,t)$ such that $(\Ri_P,W_Pr)=\gamma$ and
$rt$ is $R_P$-generic. We put
$\Lambda^{gen}=\cup \Lambda_\gamma^{gen}$
\index{0L1@$\Lambda^{gen}=\cup \Lambda_\gamma^{gen}$}
(disjoint union over $\gamma\in\Gamma$).

Define a map $m:\Lambda\to S\subset W_0\backslash T$ by
\begin{equation}
m(\Ri_P,W_Pr,t)=W_0(rt)
\end{equation}
\index{m@$m$!$m:\Lambda\to S\subset W_0\backslash T$, projection}
By Theorem \ref{thm:support}, $m$ is surjective, and obviously
$m$ is continuous and finite.

Let $P,Q\in\mathcal{P}$. Recall the set $W(P,Q)\subset W_0$
defined by
$W(P,Q):=\{n\in W_0\mid n(P)=Q\}$.
We put $W(P)=W(P,P)$, which is a subgroup of $W_0$
and acts on $R_P$ through diagram automorphisms.
Observe that $W(P)\subset N_{W_0}(W_P)$
\index{N@$N_{W_0}(W_P)$, normalizer of $W_P$ in $W_0$}
is a subgroup
which is complementary to the normal subgroup
$W_P\subset N_{W_0}(W_P)$.
Moreover, $W(P,Q)$ is a
left $W(P)$ coset and a right $W(Q)$ coset.

The action of $n\in W(P,Q)$ on $T$ restricts to isomorphisms
$T_P\to T_Q$ and $T^P\to T^Q$.
Recall that $K_P=T_P\cap T^P$, so that $n\in W(P,Q)$ gives
rise to an isomorphism $n:K_P\to K_Q$.

Consider the groupoid $\mathcal{W}$
\index{W8@$\mathcal{W}$, groupoid whose set of objects is
$\mathcal{P}$, with morphisms
$\operatorname{Hom}_\mathcal{W}(P,Q)=
\mathcal{W}(P,Q):=K_Q\times W(P,Q)$}
whose set of objects is $\mathcal{P}$,
with morphisms $\operatorname{Hom}_\mathcal{W}(P,Q)
=\mathcal{W}(P,Q):=K_Q\times W(P,Q)$
\index{W8a@$\mathcal{W}(P,Q)=\operatorname{Hom}_\mathcal{W}(P,Q)$}
and
the composition defined by
$(k_1\times n_1)\circ(k_2\times n_2)=k_1n_1(k_2)\times (n_1\circ
n_2)$. We denote by $\mathcal{W}(P)$ the group
$\mathcal{W}(P)=\mathcal{W}(P,P)$.
\index{W8b@$\mathcal{W}(P)=\mathcal{W}(P,P)$}

If $k\times n\in\mathcal{W}(P,Q)$ we define for
$\gamma=(\Ri_P,W_Pr)\in\Gamma_P$:
\begin{equation}
(k\times n)(\gamma)=(k\times n)(\Ri_P,W_Pr):=(\Ri_Q,W_Q(k^{-1}n(r))).
\end{equation}
This defines a left action of $\mathcal{W}$ on $\Gamma$.
If $t\in T^P_u$ then $(\gamma,t)\in\Lambda_\gamma$, and we define
\begin{equation}
(k\times n)(\gamma,t):=((k\times n)(\gamma),kn(t)).
\end{equation}
This defines a left action of $\mathcal{W}$ on $\Lambda$.
With these definitions we obviously have
\begin{equation}\label{eq:minv}
m(g(\lambda))=m(\lambda)
\end{equation}
for all $g\in \operatorname{Hom}(\mathcal{W})$ and $\lambda\in\Lambda$
such that $g(\lambda)$ is defined. In other words,
$m$ is $\mathcal{W}$-invariant.
\begin{lem}\label{lem:free}
The action of $\mathcal{W}$ on
$\Lambda^{gen}$ is free.
\end{lem}
\begin{proof}
Let $\lambda=(\Ri_P,W_Pr,t)\in\Lambda_{\Ri_P,W_Pr}^{gen}$ and
let $g=k\times n\in \mathcal{W}(P,Q)$ be
such that $g(\lambda)=\lambda$.
Then $Q=P$, $g$ fixes $W_Pr$, and $kn(t)=t$.
We have $g(W_Pr):=W_P(k^{-1}n(r))$, thus
$n(r)=kw(r)$ for some $w\in W_P$. Hence
$n(rt)=n(r)n(t)=w(r)(kn(t))=w(rt)$.
Since $rt$ is $R_P$-generic this implies that $w^{-1}n\in W_P$.
Hence $n=e$, and thus also $k=e$.
\end{proof}
\begin{lem}\label{lem:orb}
Let $L$ be a residual subspace, and let $t=r_Lt^L\in L^{temp}$ be
$R_L$ generic. Then $m^{-1}(W_0t)$ is a $\mathcal{W}$-orbit in
$\Lambda$.
\end{lem}
\begin{proof}
By making a suitable choice of $t$ in the orbit $W_0t$ we may
assume that $R_L=R_P$ for some $P\in\mathcal{P}$. We write
$r_P$ instead of $r_L$ and $t^P$ instead of $t^L$.
Thus it is assumed that $t=r_Pt^P\in L^{temp}=r_PT^P_u$ is $R_P$-generic.
Define $\lambda:=(\Ri_P,W_Pr_P,t^P)\in\Lambda_{\Ri_P,W_Pr_P}^{gen}$.
Clearly, $\mathcal{W}\cdot\lambda\subset m^{-1}(W_0t)$ by the
$\mathcal{W}$-invariance of $m$.

Conversely, suppose that $\mu=(\Ri_Q,W_Qr_Q,t^Q)\in m^{-1}(W_0t)$.
Hence there exists a $w\in W_0$ such that
$r_Qt^Q=wt=w(r_P)w(t^P)$. This is an element of the tempered residual
subspace $L_Q^{temp}:=r_QT^Q_u$, so that
$R_Q\subset R_{P(wt)}$. Since $t\in L^{temp}$ is
$R_P$ generic,
we have $R_{P(t)} = R_P$ by Corollary \ref{cor:gen}.
Because $R_{P(wt)}=w(R_{P(t)})$, we obtain
$R_Q\subset w(R_P)$. This implies that
$w(L^{temp})=w(t)w(T^P_u)\supset r_QT^Q_u=L_Q^{temp}$. By Theorem
\ref{thm:nonnest} we see that $w(L^{temp})=L_Q^{temp}$. Hence we
have $w(R_P)=R_Q$, $w(T_P)=T_Q$ and $w(T^P)=T^Q$. We conclude that
$r_Q^{-1}w(r_P)=k\in K_Q$. There exists a unique $u\in
W_Q$ such that $uw\in W(P,Q)$. One easily checks that
$\mu=(k\times uw)(\lambda)$.

Note that it follows that the intersection $m^{-1}(W_0t)\cap
\Lambda_{\Ri_P,W_Pr_P}$ is contained in
$\Lambda_{\Ri_P,W_Pr_P}^{gen}$ (for any choice of $P\in\mathcal{P}$
and $W_Pr_P$).
\end{proof}
\begin{cor}\label{cor:romeo}
We form the quotients $\Sigma=\mathcal{W}\backslash\Lambda$
\index{0S@$\Sigma$, quotient space of $\Lambda$ for the action
of $\mathcal{W}$} and $\Sigma^{gen}=\mathcal{W}\backslash\Lambda^{gen}$.
The map $m$ factors through $\Sigma$, and defines a homeomorphism
(also denoted by $m$) from $\Sigma^{gen}$ onto the open dense set
$S^{gen}:=m(\Lambda^{gen})\subset S$.
\index{m@$m$!$m:\Sigma^{gen}\to S^{gen}:=m(\Lambda^{gen})\subset S$,
homeomorphism}
\end{cor}
\begin{proof}
By equation (\ref{eq:minv}), $m$
is well defined on $\Sigma$, and thus
$m(\Sigma)=m(\Lambda)=S$. By the previous lemma, the
set $\Lambda^{gen}$ is $m$-saturated. Since $m$ is closed,
this implies that $S^{gen}=m(\Lambda^{gen})\subset S$ is open
(and obviously dense) in $S=m(\Sigma)$. Finally,
again by the previous lemma, $m$ is injective on $\Sigma^{gen}$.
Thus, being a closed map, $m$ is homeomorphic onto its image.
\end{proof}

$\Sigma$ can be realized as a disjoint union of
orbifolds as follows.
Choose a complete set $\Gamma_a$
\index{0C@$\Gamma_a$, complete set of representatives
of the association classes in $\Gamma$}
of representatives for the association classes
(the orbits of $\mathcal{W}$)
of elements in $\Gamma$. Put
\begin{equation}
\Sigma_\gamma:=\mathcal{W}(\gamma)\backslash
\Lambda_\gamma,
\end{equation}
where $\W(\g)$
\index{W8c@$\W(\gamma)=\{g\in\W\mid g\g=\g\}$}
denotes the isotropy group of $\gamma\in\Gamma$ in
$\W$.
Then
\begin{equation}\label{eq:union}
\Sigma\simeq\cup\Sigma_\gamma,
\end{equation}
where the (disjoint) union is taken over the set of
$\gamma\in\Gamma_a$.
\subsubsection{Groupoid $\W_{\Xi}$ of standard induction data}
Recall the complete set of representatives
$\Delta_\gamma$ ($\gamma=(\Ri_P,W_Pr_P)\in\Gamma$) of the irreducible
discrete series representations with central character $W_Pr_P$
of $\H_P$. We denote by $\Delta=\cup\Delta_\gamma$
\index{0D@$\Delta=\cup_{\gamma\in\Gamma}\Delta_\gamma$}
the disjoint
union of these sets over all $\gamma\in\Gamma$.
The composition $\Delta\to\Gamma\to\P$ gives a surjection of
$\Delta\to\P$, whose fibers are denoted by $\Delta_P$.
\index{0D2@$\Delta_P$, fiber at $P\in\P$
of the surjection $\Delta\to\P$}
\index{0D1@$\Delta_\g$, fiber at $\g\in\Gamma$ of
the surjection $\Delta\to\Gamma$}

There is a natural
left action $\Psi$ of $\mathcal{W}$ on $\Delta$ as follows:
When $k\in K_P=T^P\cap T_P$, we have an automorphism
$\psi_k:{\H_P}\to{\H_P}$ defined by
$\psi_k(\theta_xT_w)=k(x)\theta_xT_w$. This induces an isomorphism
$\psi_k:\overline{\H_P^{r_P}}\to \overline{\H_P^{k^{-1}r_P}}$.
We define a bijection
$\Psi_k$ from $\Delta_{\Ri_P,W_Pr_P}$ to
$\Delta_{\Ri_P,k^{-1}W_Pr_P}$ by $\Psi_k(\d)\simeq \d\circ \psi_k^{-1}$.

Let $Q\in F_0$ be associate to $P$, and $n\in W(P,Q)$.
Then $n$ induces an isomorphism of root data
and labels $(\Ri_P,q)\to(\Ri_Q,q)$,
thus inducing an isomorphism $\psi_n$ on $\H_P\to\H_Q$.
Recall that $\psi_n$ induces an isomorphism
$\psi_n:\overline{\H_P^{r_P}}\to\overline{\H_Q^{n(r_P)}}$
(Corollary \ref{cor:con}), and thus a bijection
$\Psi_n:\Delta_{\Ri_P,W_Pr_P} \to \Delta_{\Ri_Q,W_Qn(r_P)}$ by
$\Psi_n(\d)\simeq \d\circ\psi_n^{-1}$.
One easily checks that these definitions combine to define a
left action $\Psi$ of $\mathcal{W}$ on $\Delta$
\index{0X@$\Psi$, left action of $\mathcal{W}$ on $\Delta$},
compatible with the surjection
$\Delta\to\P$ mentioned above.
\index{0x@$\psi_g:\H^L\to\H^M$, isomorphism for
$g\in K_M\times W(F_L,F_M)$}
\index{0X@$\Psi_g:\Delta_{\Ri_L,W_Lr_L}\to\Delta_{\Ri_M,W_Mr_M^\prime}$,
bijection induced by $\psi_g$}

Consider the product ${\Xi}:=\Lambda\times_{\Gamma}\Delta$
\index{0O@${\Xi}=\Lambda\times_{\Gamma}\Delta$}. This set comes
equipped with a natural surjection ${\Xi}\to\P$
and compatible left action
of $\W$ (the diagonal action).
We form the cross product $\W_{\Xi}:=\W\times_\P{\Xi}$
\index{W8d@$\W_{\Xi}:=\W\times_\P{\Xi}$,
groupoid of standard induction data},
which has itself a natural groupoid structure with
$\operatorname{Obj}(\W_{\Xi}):={\Xi}$, and
$\operatorname{Hom}_{\W_{\Xi}}(\xi_1,\xi_2):=\{w\in\W\mid
w(\xi_1)=\xi_2\}$. The composition maps are defined by the composition
in $\W$. We will refer to this structure as the groupoid
of standard induction data of $\H$. Its set of objects ${\Xi}$
has the structure of a disjoint union of compact tori, and with
this structure $\W_{\Xi}$ is obviously a smooth compact groupoid.

Recall that we associated to each
$\xi=\l\times\d=(\Ri_P,W_Pr_P,t)\times(\Ri_P,\d)\in{\Xi}$ (i.e.
$\d$ is an irreducible discrete series representation of
$\H(\Ri_P,q)$
with central character equal to $W_Pr_P$) a tempered, unitary
representation $\pi(\xi)=\pi(\Ri_P,W_Pr_P,\d,t)$ of
$\H$ with central character $m(\l)=W_0(r_Pt)$
and representation space $V_{\pi(\xi)}=\H(W^P)\otimes V_\d$
(the compact realization)
(cf. Definition \ref{dfn:ind}, Proposition \ref{prop:unit} and
Proposition \ref{prop:indtemp}).

For every $(g,\xi)\in\W_{\Xi}$ with source $\xi=\l\times\d$,
we choose a unitary isomorphism
$\tilde\d_g: V_\d\to V_{\Psi_g(\d)}$
\index{0dZ@$\tilde\d_g: V_\d\to V_{\Psi_g(\d)}$, unitary isomorphism
intertwining $\d\circ\psi_g$ and $\Psi_g(\d)$}
so that
we have
\begin{equation}
\Psi_g(\d)(h)\circ\tilde\d_g=\tilde\d_g\circ\d(\psi_g^{-1}h)
\end{equation}
for $h\in \H_P$
(where $P=P(\delta)$).

This enables us to define intertwining operators (depending on the
choices of the isomorphisms $\tilde\d_g$)
\begin{equation}
A(g,\xi)\in\operatorname{Hom}_\H(V_{\pi(\xi)},V_{\pi(g(\xi))})
\end{equation}
as follows:

For $k\in K_P$ and $h\in\H^P$ we have
$\phi_{kt}(h)=\psi_k(\phi_t(h))$, so that
if $h\in \H^P$ we have that $\tilde\d_k\circ\d(\phi_{t}(h))=
\Psi_k(\d)(\phi_{kt}(h))\circ\tilde\d_k$.
With this notation we have for each
$\d\in\Delta_{\Ri_P,W_Pr_P}$, in view of Proposition \ref{prop:unit},
a unitary intertwining isomorphism
\begin{equation}\label{eq:iso}
\operatorname{Id}\otimes\tilde\d_k:\pi(\xi)\to
\pi(k(\xi)).
\end{equation}
We denote this unitary intertwining operator by $A(k,\Ri_P,W_Pr,\d,t)$
or more simply $A(k,\xi)$. Notice that it is constant, i.e.
independent of $t$.

For $n\in W(P,Q)$ (with $P,Q\in\P$ associate subsets) we defined
(cf. Theorem \ref{thm:ms} and Corollary \ref{cor:hol})
an intertwining isomorphism (depending on the choice of
$\tilde\d_n$)
\begin{equation}
A(n,\Ri_P,W_Pr_P,\d,t):\pi(\Ri_P,W_Pr_P,\d,t)\to
\pi(\Ri_Q,W_Qn(r_P),\Psi_n(\d),n(t)),
\end{equation}
which is rational in $t$, well defined and regular in a neighborhood
of $T^P_u$,
and unitary for $t\in T^P_u$. We now denote this isomorphism
by $A(n,\xi)$.

The above intertwining isomorphisms combine (as one easily checks
directly from the definitions) to a functor
\begin{equation}
\W_{\Xi}\to {PRep}_{unit,temp}(\H)
\end{equation}
where ${PRep}_{unit,temp}(\H)$ denotes the category
of finite dimensional, tempered, unitary modules over $\H$, with
morphisms $\operatorname{Hom}_{PRep}(\pi_1,\pi_2)=
PU_\H(V_{\pi_1},V_{\pi_2})$ (the space of unitary intertwiners
modulo the action of scalars).

Summarizing the above we have:
\begin{thm}\label{cor:proact}
There exists an induction functor
\begin{equation}
\pi:\W_{\Xi}\to{PRep}_{unit,temp}(\H)
\end{equation}
\index{0p@$\pi$, induction functor on $\W_{\Xi}$}
such that for
$\xi=(\Ri_P,W_Pr_P,t)\times(\Ri_P,\d)\in {\Xi}$
and $(g,\xi)\in\W_{\Xi}$ (thus $g\in\W$ with source $P(\xi)=P$),
$\pi(\xi):=\pi(\Ri_P,W_Pr_P,\d,t)$ and
$\pi(g,\xi):=A(g,\xi)=A(n,\Ri_P,W_Pr_P,\d,t)$.
\end{thm}
\subsubsection{Generic spectrum}
Consider the natural projection
\begin{equation}
p_\Sigma:\W_{\Xi}\backslash{\Xi}
=\W\backslash(\Lambda\times_\Gamma\Delta)
\to\Sigma=\W\backslash\Lambda.
\end{equation}

Since the action of $\W$ is free on the set of generic points
$\Lambda^{gen}$, we obtain a
finite covering
\begin{equation}
p_\Sigma:\W_{\Xi}\backslash{\Xi}^{gen}\to\Sigma^{gen},
\end{equation}
where ${\Xi}^{gen}:=\Lambda^{gen}\times_\Gamma\Delta$
\index{0O@${\Xi}^{gen}=\Lambda^{gen}\times_\Gamma\Delta$}.

By what was said in the previous subsection
and Corollary \ref{cor:short}, it is clear that the map
(see Proposition \ref{prop:im} for the definition of $p_z$):
\begin{equation}
\begin{split}
[\pi]:
{\Xi}^{gen}&\to p^{-1}_z(S^{gen})\\
\xi&\to[\pi(\xi)]
\end{split}
\end{equation}
factors through the quotient $\W_{\Xi}\backslash{\Xi}^{gen}$
\index{0p1@$[\pi(\xi)]$, class of $\pi(\xi)$ modulo equivalence}.
We thus have the following
commutative diagram:
\begin{equation}
\begin{CD}\label{CD?}
\W_{\Xi}\backslash{\Xi}^{gen}@>[\pi]>>p^{-1}_z(S^{gen})\\
@V{p_\Sigma}VV @VV{p_z}V\\
\Sigma^{gen}@>>m>S^{gen}
\end{CD}
\end{equation}
\begin{thm}\label{thm:homeom}
The map $[\pi]$ in the diagram \ref{CD?} is a homeomorphism.
\end{thm}
\begin{proof}
The topology on $\hat\cf$ is second countable since $\cf$ is
separable. Thus, in order to check the continuity of $[\pi]$, it
suffices to check that $[\pi]$ maps a converging sequence
$\l_i\times\d\to\l\times\d\in\Lambda_{\Ri_P,W_Pr}$ to a
converging sequence in $\hat{\cf}$.
We check this using the Fell-topology
description of the topology of $\hat{\cf}$ (see \cite{Fell}).
By \cite{Fell}, Proposition 1.17, restriction to the dense
subalgebra $\H\subset\cf$ is a homeomorphism with respect to the Fell
topologies.
Let $V_{\l\times\d}=\H(W^P)\otimes V$ be the representation space of
$\pi(\l\times\d)$ (with $\l\in\Lambda_{\Ri_P,W_Pr}$).
We equip $V_{\l\times\d}$ with the inner product
$\langle\cdot,\cdot\rangle$
of Proposition \ref{prop:unit}
(which is independent of $\l\in\Lambda_{\Ri_P,W_Pr}$),
and we choose an orthonormal basis $(e_i)$ of $V_{\l\times\d}$.
In order to check that
$\pi(\l_i\times\d)\to\pi(\l\times\d)$
in the Fell topology with respect to $\H$, we need to check
that for all $h\in\H$,
$\pi(\l_i\times\d)(h)_{k,l}\to\pi(\l\times\d)(h)_{k,l}$ for all
matrix coefficients. This is clear since the matrix coefficients
are regular functions of the induction parameter.

To see that the map $[\pi]$ is closed, assume that we have a
sequence $\rho_i=[\pi](\l_i\times\d_i)$ converging to
$\rho\in p_z^{-1}(S^{gen})$. Since $m$ is a homeomorphism and
$\Sigma^{gen}$ is a finite quotient of $\Lambda^{gen}$, we may
assume that $\l_i$ converges, to $\l_0\in\Lambda^{gen}$ say, by
possibly replacing the sequence by a subsequence. Since
$\Delta_{\Ri_P,W_Pr}$ is finite for each $R_P$ and $W_Pr$, we may
assume that $\forall i: \d_i=\d$, again by taking a subsequence.
Then $d=\operatorname{dim}(\rho_i)$ is independent of $i$, and
lower semi-continuity of $\operatorname{dim}$ on $\hat\cf$ implies
that $\operatorname{dim}(\rho)\leq d$. Choose an orthonormal basis
$B$ for $\rho$. Convergence in the Fell-topology means that there
exists for all $i$ an orthonormal subset $B_i$ of size
$\operatorname{dim}(\rho)$ in the representation space
$V_{\rho_i}=\H(W^P)\otimes V$ of
$\rho_i$, such that the matrix coefficients of $\rho_i$ with
respect to $B_i$ converge to the matrix coefficients of $\rho$
with respect to $B$. By the independence of the inner product
$\langle\cdot,\cdot\rangle$ of the induction parameter
(Proposition \ref{prop:unit}) we may assume, by further
restricting to a subsequence, that the sets $B_i$
converge in $\H(W^P)\otimes V$ to an orthonormal set
$B_0$. It follows that the matrix of
$\rho(x)$ with respect to $B$ equals a principal block of the
matrix of $\pi(\l_0\times\d)(x)$ with respect to a
suitable orthonormal basis $\tilde{B}$ of $\H(W^P)\otimes V$ for
$\pi(\l_0\times\d)$. Since $\pi(\l_0\times\d)$ is
irreducible it is easy to see that this is impossible unless
$\rho\simeq\pi(\l_0\times\d)$.

The map $[\pi]$ is injective by
Corollary \ref{cor:short} and Lemma \ref{lem:orb}.

Finally, by Theorem \ref{thm:nu}, Theorem \ref{thm:support},
Theorem \ref{thm:mainind} and
Theorem \ref{thm:thisisathm} we see that the complement of
$[\pi](\W_{\Xi}\backslash{\Xi}^{reg})$ has measure $0$ in
$\hat\cf$ with respect to
the Plancherel measure of the representation $\hf$ of $\cf$.
The support of the Plancherel measure is equal to $\hat\cf$,
since $\hf$ is a faithful representation of
$\cf$ (by definition of $\cf$).
Thus the closure of $[\pi](\W_{\Xi}\backslash{\Xi}^{reg})$
is $\hat\cf$.
But
$[\pi](\W_{\Xi}\backslash{\Xi}^{reg})\subset p_z^{-1}(S^{reg})$
is closed as we have
seen above, implying that $[\pi]$ is surjective.
\end{proof}
\begin{cor}
The restriction of the map $p_z$ of
Corollary \ref{prop:im} to $p_z^{-1}(S^{reg})$
is a covering map.
\end{cor}
\subsubsection{Fourier transform}
Let $\tilde\O\subset{\Xi}$
\index{O@$\tilde\O$, connected component of ${\Xi}$}
be a connected component of ${\Xi}$. Thus there exists a
$\d\in\Delta$
such that
$\tilde\O=\Lambda_\gamma\times\{\d\}:=\tilde\O_\d$, where
$\g=\g(\d)\in\Gamma$. Explicitly, if $\g(\d)=(\Ri_P,W_Pr_P)$ then
$\tilde\O_\d$ is a copy of the subtorus $T^P_u\subset T_u$.

The representation space $V_{\pi(\xi)}$ of $\pi(\xi)$
is equal to $V_{\pi(\xi)}=\H(W^P)\otimes V_\d$ for
$\xi\in\tilde\O_\d$ with $\d\in\Delta_P$.
In particular, $V_{\pi(\xi)}$ depends
only on the connected component $\tilde\O_\d$ of $\Xi$
containing $\xi$, and not on the choice of $\xi\in\tilde\O_\d$.
We will use the notation
$i(V_\d):=\H(W^P)\otimes V_\d=V_{\pi(\xi)}$ for any choice of
$\xi\in\tilde\O_\d$ (where $P=P(\d)\in\P$).
\index{i@$i(V_{\d})=\H(W^P)\otimes V_\d$ if $\d\in\Delta_P$}

We form the trivial fiber bundle
$\V_{\tilde\O}=\tilde\O_\d\times i(V_\d)$
\index{V1@$\V_{\tilde\O}=\tilde\O_\d\times i(V_\d)$,
trivial fiber bundle over $\tilde\O=\tilde\O_\d$}
over
$\tilde\O=\tilde\O_\d$, and put
\begin{equation}
\V_{\Xi}:=\cup_{\d\in\Delta}\V_{\tilde\O_\d}.
\end{equation}
\index{V3@$\V_{\Xi}$, trivial fiber bundle over $\Xi$}
We identify the connected component
$\tilde\O_\d$ of $\Xi$ with the compact torus
$T^P_u$ ($P=P(\d)$).
This allows us to define the function spaces
$\operatorname{Pol}(\Xi)$
\index{Pol@$\operatorname{Pol}(\Xi)$, space of Laurent polynomials
on $\Xi$} (Laurent polynomials on $\Xi$) and
$\cc(\Xi)$
\index{C@$\cc(\Xi)$, space of $\cc$-functions on
$\Xi$}. We also introduce the space
$\operatorname{Rat}^{reg}(\Xi)=
\oplus_{\d\in\Delta}\operatorname{Rat}^{reg}(\tilde\O_d)$,
\index{RatX@$\operatorname{Rat}^{reg}(\Xi)$,
regular rational functions on $\Xi$}
\index{RatO@$\operatorname{Rat}^{reg}(\tilde\O)$,
rational functions on $\tilde\O$, regular in an
open neighborhood of $\tilde\O\simeq T^P_u\subset T^P$}
where $\operatorname{Rat}^{reg}(\tilde\O_d)$ denotes
the space of restrictions to $T^P_u$ (which we identify
with $\tilde\O_d$)
of rational functions on $T^P$ which are regular
in a open neighborhood of $T^P_u$.
The corresponding spaces of (global) sections are denoted
by
$\operatorname{Pol}(\V_{\Xi})=
\operatorname{Pol}(\Xi)\otimes\V_{\Xi}$
\index{PolV@$\operatorname{Pol}(\V_{\Xi})=
\operatorname{Pol}(\Xi)\otimes\V_{\Xi}$},
$\cc(\V_{\Xi})=
\cc(\Xi)\otimes\V_{\Xi}$
\index{C@$\cc(\V_{\Xi})=
\cc(\Xi)\otimes\V_{\Xi}$}, and
$\operatorname{Rat}^{reg}(\V_{\Xi})=
\operatorname{Rat}^{reg}(\Xi)\otimes\V_{\Xi}$
\index{RatV@$\operatorname{Rat}^{reg}(\V_{\Xi})=
\operatorname{Rat}^{reg}(\Xi)\otimes\V_{\Xi}$} respectively.

Recall that $\pi(g,\xi)\in PU_\H(i(V_\d),i(V_{gd}))$
(with $\xi=\l\times\d=(\Ri_P,W_Pr,t)\allowbreak\times\d$) is
rational and regular for $t\in T^P$
in a neighborhood of $T^P_u$ (Corollary \ref{cor:hol}).
We define
\begin{align*}
\operatorname{Pol}(&\operatorname{End}(\V_\Xi))^\W\\=
&\{f\in\operatorname{Pol}(\operatorname{End}(\V_\Xi))
\mid \forall(g,\xi)\in\W_\Xi:\
\pi(g,\xi)f(\xi)=f(g\xi)\pi(g,\xi)\}\\\simeq
&\bigoplus_{\tilde\O}\operatorname{Pol}(\V_{\tilde\O})^{\W(\tilde\O,\tilde\O)}
\end{align*}
\index{Pol@$\operatorname{Pol}(\operatorname{End}(\V_\Xi))^\W$,
space of $\W_\Xi$-equivariant sections in
$\operatorname{Pol}(\operatorname{End}(\V_\Xi))$}
where the direct sum runs over a complete set of representatives
of connected components $\tilde\O$ for the action of $\W$, and
$\W(\tilde\O_1,\tilde\O_2)$ denotes the set of
$w\in\W$ such that $w(\tilde\O_1)=\tilde\O_2$.
We define the space of $\W_\Xi$-equivariant sections
in other spaces of sections of $\operatorname{End}(\V_\Xi)$ similarly.
\begin{dfn}
The Fourier transform is the algebra homomorphism
\begin{align*}\label{eq:FT}
\F_\H:\H&\to\operatorname{Pol}(\operatorname{End}(\V_\Xi))^\W\\
h&\to\{\xi\to\pi(\xi)(h)\}
\end{align*}
\end{dfn}
\index{F@$\F_\H$, Fourier transform on $\H$}

We would like to replace $\Xi$ by the set of
equivalence classes of cuspidal representations
of the standard parabolic subalgebras $\H^P$.
This can be done as follows.
Consider the subgroupoid $\K\subset \W$
\index{K@$\K$, normal subgroupoid of $\W$}
of $\W$, with set of objects $\P$, and
$\K(P_1,P_2)=\emptyset$ if $P_1\not=P_2$, and
$\K(P,P)=K_P$.
This subgroupoid is normal in the sense that
$gK_Pg^{-1}=K_Q$ if $g\in\W(P,Q)$. The quotient groupoid
$\W/\K=\W/\K$ has $\P$ as set of objects, and $\W/\K(P,Q)=W(P,Q)$.

Suppose that $\d_t\simeq\d^\prime_s$ with
$\d,\d^\prime\in\Delta_{\Ri_P}$ and $s,t\in T^P_u$.
Then $W_Prt=W_Pr^\prime s$, and thus $s=kt$ for
some $k\in K_P$, and $\d^\prime_s=\Psi_k(\d)_{kt}$.
Conversely, in view of the text above (\ref{eq:iso}),
$\xi\simeq k(\xi)$ for every $k\in K_P$ and $\xi\in\Xi_P$,
viewed as representation of $\H^P$.

The connected components $\O$
\index{O@$\O$, orbit of twists of cuspidal representations}
of the quotient
$\K\backslash\Xi$ are called ``orbits of twists of cuspidal
representations'' of the parabolic subalgebras $\H^P$.
Such a component can be viewed as the
collection of mutually inequivalent representations of $\H^P$
of the form $\d_t$. It is isomorphic to a smooth quotient
$\O\simeq\K(\tilde\O,\tilde\O)\backslash\tilde\O$, a finite
quotient of the subtorus $T^P_u\subset T_u$.

We have $\W\backslash\Xi=
(\K\backslash\W)\backslash(\K\backslash\Xi)$. For $\O$ a connected
component of $\K\backslash\Xi$, we choose a connected component
$\Xi\supset\tilde\O\to\O$ covering $\O$. Let
$\tilde\O=\Lambda_\g\times\d$ and write $P(\g)=P$. The isotropy
group $\{k\in K_P\mid k(\tilde\O) =\tilde\O\}$ equals the isotropy
group $K_\d$ \index{K@$K_\d\subset K_P$, isotropy subgroup of
$[\d]\in [\Delta_P]$}. Notice that $K_\d$ is independent of the
choice of $\tilde\O\to\O$. We define the principal fiber bundle
$\V_\O:=\tilde\O\times_{\K_\d}i(V_\d)$
\index{V2@$\V_\O:=\tilde\O_\d\times_{\K_\d}i(V_\d)$} over $\O$.
This fiber bundle is not necessarily trivial. We put
\begin{align}
\operatorname{Pol}(\operatorname{End}(\V_\O))&
=\left(\bigoplus_{\tilde\O:\tilde\O\to\O}\operatorname{Pol}
(\operatorname{End}(\V_{\tilde\O}))\right)^{K_P}\\\nonumber
&\simeq\operatorname{Pol}
(\operatorname{End}(\V_{\tilde\O}))^{K_\d}
\end{align}
\index{Pol@$\operatorname{Pol}(\operatorname{End}(\V_\O))$,
polynomial sections in fiber bundle $\operatorname{End}(\V_\O)$}
and
\begin{equation}
\operatorname{Pol}(\operatorname{End}(\V_{(\K\backslash\Xi)}))=
\bigoplus_{\O\mathrm{\ orbit}}\operatorname{Pol}(\operatorname{End}(\V_\O))
\end{equation}
\index{Pol@$\operatorname{Pol}(\operatorname{End}(\V_{(\K\backslash\Xi)}))$,
polynomial sections in fiber bundle
$\operatorname{End}(\V_{(\K\backslash\Xi)})$}

The quotient $\W/\K$ acts on $\K\backslash\Xi$ and thus
also on the set of orbits.
Given orbits $\O_1,\O_2$ with $P(\O_i):=P_i$ and
$O_i=\Lambda_{\g_i}\times \d_i$, we have
$\W/\K(\O_1,\O_2)=\{n\in \W/\K(P_1,P_2)
\mid n(\O_1)=\O_2\}=
\{n\in W(P_1,P_2) \mid \exists k\in K_{P_2}:k\times n
\in\W(\d_1,\d_2)\}$. We denote this set by
$W(\O_1,\O_2)$
\index{W6a@$W(\O_1,\O_2)=\{n\in W(P_1,P_2)\mid \exists k\in K_{P_2}:(k\times n)
\in\W(\d_1,\d_2)\}$}. We also put $W(\O):=W(\O,\O)$
\index{W6b@$W(\O)=W(\O,\O)$}.

In this way we get
\begin{align}\label{eq:FTbis}
\operatorname{Pol}(\operatorname{End}(\V_\Xi))^{\W}&=
\operatorname{Pol}(\operatorname{End}(\V_{(\K\backslash\Xi)}))
^{\K\backslash\W}\\\nonumber
&\simeq
\bigoplus_\O\operatorname{Pol}(\operatorname{End}(\V_\O))^{W(\O)},
\end{align}
where the direct sum runs over a complete set of representatives
of orbits $\O$ modulo the action of $\W/\K$ (association classes
of orbits).

We use similar notations for spaces of sections with coefficients
in other types of functions spaces
(e.g. continuous, $\cc$, etc.) in
$\operatorname{End}(\V_\Xi)$ and
$\operatorname{End}(\V_{(\K\backslash\Xi)})$.
\subsubsection{Averaging projections}\label{subsub:av}
Consider a function space $F$ on $\Xi$ which is
a module over $\operatorname{Rat}^{reg}(\Xi)$.
Due to the regularity of the intertwining operators
(cf. Corollary \ref{cor:hol}), there exists a
natural averaging projection $f\to\overline{f}$
\index{$f\to\overline{f}$, average of a section of
$\operatorname{End}(\V_\Xi)$}
from $F(\operatorname{End}(\V_\Xi))$ (sections of $\operatorname{End}(\V_\Xi)$
with coefficients in $F$) to $F(\operatorname{End}(\V_\Xi))^\W$.
It is defined by (where $\W_\xi=\{g\in\W\mid (g,\xi)\in\W_\Xi\}$)
\begin{equation}
\overline{f}(\xi)=|W_\xi|^{-1}\sum_{g\in\W_\xi}
\pi(g,\xi)^{-1}f(g(\xi))\pi(g,\xi).
\end{equation}
Notice that the function space  $F=\operatorname{Pol}(\Xi)$ is too
small; in general the average of
$f\in\operatorname{Pol}(\operatorname{End}(\V_\Xi))$ will be in
$\operatorname{Rat}^{reg}(\operatorname{End}(\V_\Xi))^{\W}$.

There is a similar averaging procedure $f\to \overline{f}_{\K}$
which sends the space of sections $F(\operatorname{End}(\V_\Xi))$ to
$F(\operatorname{End}(\V_{(\K\backslash\Xi)}))$ (in this case
$F$ should be a module over $\C$).
\subsubsection{Plancherel formula}
We now define the Plancherel measure on $\W\backslash\Xi$.
The following proposition says that the natural action of
$\W_\Xi$ (via $\psi$) on the residual algebras is through
isomorphisms of Hilbert algebras.
\begin{prop}\label{prop:impinv}
Let $\d\in\Delta_{\Ri_P,W_Pr}$ and let $g=(k\times
n)\in\mathcal{W}_{\Ri_P,W_Pr}$. We have (in the notation of Theorem
\ref{thm:mainind}) $d_{\Ri_P,\Psi_g(\d)}=d_{\Ri_P,\d}$.
\end{prop}
\begin{proof}
This is a simple extension of Corollary \ref{cor:con}, with
a similar proof.
\end{proof}
Let $\tilde\O=\Lambda_\g\times\d$ and let
$\O=K_\d\backslash\tilde\O$. If $P=P(\g)$ then
$\O$ is a copy of the subquotient torus $\K_\d\backslash T^P_u$.
For $\om\in\O$ we put $d^\O\om$ for the normalized Haar measure
on $\O$. Let $\g=(\Ri_P,W_Pr)$ and let
$\om=(\Ri_P,W_Pr,\d,K_\d t^P)$ be $R_P$-generic.
Let $L^{temp}=rT^P_u$ denote a residual subspace
underlying $\tilde\O$.
We define
\begin{equation}\label{eq:plameasure}
\begin{split}
d\mu_{Pl}(\pi(\om))&=|W_0(rt^P)||W^P|^{-1}d_{\Ri_P,\d}
d\nu_L(rt^P)\\
&=\frac{|W_P|}{|W_P\cap W_{r}|}
\nu_{\Ri_P}(\{r\})d_{\Ri_P,\d}k_Pm^P(rt^P)d^P(t^P)\\
&=\mu_{\Ri_P,Pl}(\{\d\})|K_P\d|m^P(\om)d^\O\om,
\end{split}
\end{equation}
\index{0m@$\mu_{Pl}$, Plancherel measure on $\hat\cf$}
where $d_{\Ri_P,\d}>0$ is the residual degree of $\d$ in
the residual algebra $\overline{\H^{r}_L}$,
$\mu_{\Ri_P,Pl}$ is given in Corollary
\ref{cor:fdim},
$m^P(\om)=m^L(\om)$ is the common value
of $m^L$ (as defined in Proposition
\ref{prop:par}) on the $K_\d$ orbit $\om$,
and $k_P:=|K_P|$.
We have used that the isotropy subgroup $W_{rt^P}$
equals $W_P\cap W_r$ if $rt^P$ is $R_P$-generic.

Recall that by definition we have
\begin{equation}
\sum_{\d\in
\Delta_{\Ri_P,W_Pr}}\operatorname{dim}(\d)
d_{\Ri_P,\d}=1
\end{equation}
Recall Conjecture \ref{rem:ell} stating that
$d_{\Ri_P,\d}\in\mathbb{Q_+}$.

We define an Hermitian inner product on
$\operatorname{Pol}(\operatorname{End}(\V_{(\K\backslash\Xi)}))$
as follows:
\begin{equation}\label{eq:hermien}
(f_1,f_2)=\sum_{\O}|W(\O)|^{-1}\int_{\O}
\operatorname{tr}(f_1(\om)^*f_2(\om))d\mu_{Pl}(\pi(\om)),
\end{equation}
\index{$(\cdot,\cdot)$!inner product on
$\operatorname{Pol}(\operatorname{End}(\V_{(\K\backslash\Xi)}))$}
where the sum runs over a complete set of
representatives for the association classes of orbits $\O$
(an association classes is an orbit under the action
of $\W/\K$).
Note that $f_1,f_2$ are in fact $K_\d$-equivariant sections over
the covering $\tilde\O\to\O$. The expression
$\operatorname{tr}(f_1(\xi)^*f_2(\xi))$ is independent of a
choice of $\xi\in\tilde\O$ such that $K_\d\xi=\om$. The common
value is denoted by $\operatorname{tr}(f_1(\om)^*f_2(\om))$.
\begin{thm}\label{thm:mainp} (Main Theorem)
\begin{enumerate}
\item[(i)] Let $\O$ be an orbit (a connected component of
$\K\backslash\Xi$). We put $\hat\cf_\O^{gen}=
[\pi](W(\O)\backslash\O^{gen})\subset\hat\cf$,
and we denote its closure by
$\cf_\O\subset\hat\cf$
\index{C@$\hat\cf_\O$, component of $\hat\cf$, the closure of
$[\pi](W(\O)\backslash\O^{gen})\subset\hat\cf$}.
Then $\cf_{\O_1}^{gen}\cap\cf_{\O_2}^{gen}=\emptyset$ unless
$\O_1$ and $\O_2$ are in the same $\W/\K$-orbit, and
\begin{equation}
{\hat\cf}^{gen}:=\cup\hat\cf_\O^{gen}
\end{equation}
(union over a complete set of representatives
for the association classes of orbits) is a dense set in $\cf$,
whose complement has measure zero.
\item[(ii)]
The Plancherel measure of $\cf$ (i.e. the measure on
$\hat\cf$ determined by the tracial state $\tau$ of $\cf$)
is given on $\hat\cf_\O$ by equation (\ref{eq:plameasure}).
The decomposition of $\tau$ in irreducible,
mutually distinct characters of $\cf$ is given by
\begin{equation}
\tau=\sum_{\O}
\int_{\om\in W(\O)\backslash\O}
\chi_{\pi(\om)}d\mu_{Pl}(\pi(\om))
\end{equation}
(sum over a complete set of representatives for the association
classes of orbits).
\item[(iii)]
Equivalently, the algebra homomorphism $\F_\H$ (see (\ref{eq:FT})
and (\ref{eq:FTbis})) is an isometry with respect to the inner
product (\ref{eq:hermien}), and extends uniquely to an isomorphism
of $\cf\times\cf$ modules
\begin{equation}
\F:\mathfrak{H}\stackrel{\sim}
{\to}L^2(\operatorname{End}(\V_{\Xi}))^\W\simeq
\bigoplus_{\O}L^2(\operatorname{End}(\V_{\O}))^{W(\O)}
\end{equation}
(sum over a complete set of representatives for the association
classes of orbits).
\end{enumerate}
\end{thm}
\begin{proof}
(i) See Theorem \ref{thm:homeom}. The complement of
${\hat\cf}^{gen}$ has measure zero by the argument in the last
part of the proof of that theorem. The density follows since
$\hat\cf$ is the support of the Plancherel measure (cf.
Theorem \ref{thm:supds}).

(ii)
By formula of Proposition \ref{prop:dfn}(v) and
Corollary \ref{cor:romeo}
we have
\begin{equation}
\tau=\int_{\W\backslash\Lambda^{reg}}\chi_{m(\l)}d\nu(m(\l))
\end{equation}
We decompose $\chi_{m(\l)}$ according to  Theorem
\ref{thm:mainind}(ii) to obtain
\begin{equation}
\tau=\int_{\l\in\W\backslash\Lambda^{reg}}
|W^{P(\l)}|^{-1}\sum_{\d\in\Delta_{\g(\l)}}d_{\Ri(\l),\d}
\chi_{\pi(\l\times\d)}d\nu(m(\l))
\end{equation}
By Corollary \ref{cor:short} and Theorem \ref{thm:homeom} we
have $\{[\pi](\l\times \d)\}_{\d\in\Delta_{\g}}=
[\pi](p_\Sigma^{-1}(\W\l))$. Thus (by Theorem \ref{thm:homeom})
we can rewrite the integral as integral over
$\W\backslash\Xi^{reg}$. When we use parameters and notations as
explained in equation
(\ref{eq:plameasure}), and we express $d\nu$ according to
Proposition \ref{prop:par}, we obtain
\begin{equation}
\tau=\int_{\xi\in\W\backslash\Xi^{reg}}
|W^{P(\l)}|^{-1}|W_0(rt^P)|d_{\Ri(\xi),\d(\xi)}
\chi_{\pi(\xi)}
d\nu_L(rt^P)
\end{equation}
According to our definition of $\mu_{Pl}$ this is equal to
\begin{equation}
\tau=\sum_{\O}\int_{\om\in W(\O)\backslash \O^{reg}}
\chi_{\pi(\om)}d\mu_{Pl}(\pi(\om))
\end{equation}
This is a decomposition of $\tau$ in characters of
inequivalent irreducible representations of $\cf$
(see Theorem \ref{thm:mainind}(iii)). Hence this
uniquely determines the Plancherel measure (by \cite{dix2},
Th\'eor\`eme 8.8.6) on $\hat\cf$. We conclude that
$\mu_{Pl}$ is equal to the Plancherel measure of $\cf$.

(iii) The equivalence of (ii) and (iii) is well known, see
the proof of Theorem \ref{thm:thisisathm}. It is allowed to
use the formulation with $W(\O)$-equivariant sections because
of the unitarity and the regularity of intertwining operators
(Theorem \ref{thm:ms}) and by Proposition \ref{prop:impinv}.
\end{proof}
\begin{rem}
In \cite{DO} it is shown that the $\hat\cf_\O$ are the components
of $\hat\cf$. Moreover, $\F(\mathfrak{S})$ (see \ref{subsub:schwartz}
for the definition of $\mathfrak{S}$) and $\F(\cf)$ are determined
in \cite{DO}.
\end{rem}
\begin{cor}\label{cor:av}
Let
$\J:L^2(\operatorname{End}(\V_{(\K\backslash\Xi)}))\to\mathfrak{H}$
denote the adjoint of $\F$. Then $\J\F=\operatorname{Id}$
and $\F\J(f)=\overline{f}$ (see subsection \ref{subsub:av}).
\end{cor}
\begin{proof} By the isometry property of $\F$,
$(\J\F(x),y)=(x,y)$ for all $x,y\in\mathfrak{H}$. Whence
the first assertion. It is clear that $\J(f)=\J(\overline{f})$.
If $g\in L^2(\operatorname{End}(\V_\Xi))^\W$ then $g=\F(x)$ for
some $x\in\mathfrak{H}$. Thus $\F\J(g)=\F\J\F(x)=\F(x)=g$ for
$\W$-equivariant $g$. Hence $\F\J(f)=\F\J(\overline{f})=\overline{f}$.
\end{proof}
\section{Base change invariance of the residual algebra}
\label{sec:inv} Thus far we have found the spectral
decomposition for $\H$ in terms of the ``residual degrees''
$d_{\Ri_L,\d}$
of the residual algebras $\overline{\H_L^r}$. We
prove in this section that the residual algebras are {\it
independent} of $\q$ (using the Convention \ref{eq:scale}),
up to isomorphism of Hilbert algebras.
\subsection{Scaling of the root labels} Let $r=sc\in T$ be fixed,
with $s\in T_u$ and $c=\exp(\gamma)$ with $\gamma\in
\mathfrak{t}$. Assume that $B\subset \mathfrak{t}_\C$ is an open
ball centered around the origin such that the conditions
\ref{cond} (with respect to $r\in T$) are satisfied.

The second condition implies that each connected component of the
union $U:=W_0(r\exp(B))$ contains a unique element of the orbit
$W_0r$. Given $u\in U$ there is a unique $r^\prime =s^\prime
c^\prime\in W_0r$ such that $u\in r^\prime\exp(B)$. By (i) there
is a unique $b\in B$ such that $u=s^\prime
c^\prime\exp(b)=s^\prime\exp(b+\gamma^\prime)$. Now let $\e\in
(0,1]$ be given. We define an analytic map $\sigma_\e$ on $U$ by
\begin{equation}
\sigma_\e(u):=s^\prime\exp(\e\log((s^\prime)^{-1}u))=
s^\prime\exp(\e(b+\gamma^\prime)).
\end{equation}
\begin{lemma}
The map $\sigma_\e$ is an analytic, $W_0$-equivariant
diffeomorphism from $U$ onto $U_\e:=W_0(sc^\e\exp(\e B))$. The
inverse of $\sigma_\e$ will be denoted by $\sigma_{1/\e}$.
\end{lemma}
\begin{proof}
On the connected component $r^\prime\exp(B)$ the map $\sigma_\e$
is equal to $\sigma_\e=\mu_{s^\prime}\circ \exp\circ M_\e \circ
\log\circ\mu_{(s^\prime)^{-1}}$ where $\mu_{s^\prime}$ is the
multiplication in $T$ by $s^\prime$, and $M_\e$ is the
multiplication in $\mathfrak{t}_\C$ by $\e$. These are all
analytic diffeomorphisms, because of condition (i). The $W_0$
equivariance follows from the fact that $\log$ is well defined
(and thus equivariant, since $\exp$ is equivariant) from
$W_0\exp(B+\gamma)$ to $W_0(B+\gamma)$, and that $M_\e$ is
$W_0$-equivariant. This implies that for $w\in W_0$,
$w\exp(\e\log((s^\prime)^{-1}u))=\exp(\e\log((ws^\prime)^{-1}wu))$.
It follows that
\begin{equation}
\begin{split}
\sigma_\e(wu)&=ws^\prime\exp(\e\log((ws^\prime)^{-1}wu))\\
&=ws^\prime w\exp(\e\log((s^\prime)^{-1}u))\\ &=w(\sigma_\e(u)).
\end{split}
\end{equation}
\end{proof}
\begin{lemma}\label{lem:ana}
Denote by $q^\e$ the label function $q^\e(s)=q(s)^\e=\q^{\e f_s}$,
and denote by $\H_{q^\e}$ the affine Hecke algebra with root datum
$\Ri$ (same as the root datum of the affine Hecke algebra
$\H=\H_q$), but with the labels $q$ replaced by $q^\e$. Let
$c_{\a,\e}\in{}_\F\A_{q^\e}\subset{}_\F\H_{q^\e}$ be the
corresponding Macdonald $c$-functions. For every root $\a\in R_1$
we have:
\begin{equation}
U\ni u\to (c_{\a,\e}(\sigma_\e(u))c_\a(u)^{-1})^{\pm 1}\in
\A^{an}(U).
\end{equation}
\end{lemma}
\begin{proof}
For $u$ in the connected component $r^\prime\exp(B)$ of $U$ we
write $u=s^\prime v$ with $v\in c^\prime\exp(B)$. We have
\begin{equation}\label{eq:ana}
\begin{split}
c_{\a,\e}(\sigma_\e(u))&c_\a(u)^{-1}
=\frac{(1+q_{\alpha^\vee}^{-\e/2}\alpha(v)^{-\e/2}\alpha(s^\prime)^{-1/2})}
{(1+q_{\alpha^\vee}^{-1/2}\alpha(v)^{-1/2}\alpha(s^\prime)^{-1/2})}\\
&\times\frac{(1-q_{\alpha^\vee}^{-\e/2}q_{2\alpha^\vee}^{-\e}
\alpha(v)^{-\e/2}\alpha(s^\prime)^{-1/2})
(1-\alpha(v)^{-1}\alpha(s^\prime)^{-1})}
{(1-q_{\alpha^\vee}^{-1/2}q_{2\alpha^\vee}^{-1}\alpha(v)^{-1/2}
\alpha(s^\prime)^{-1/2})(1-\alpha(v)^{-\e}\alpha(s^\prime)^{-1})}
\end{split}
\end{equation}
We remind the reader of the convention Remark \ref{rem:conv}; in
particular, the expression $\a(s^\prime)^{1/2}$ occurs only if
$\a/2\in R_0$, in which case this expression stands for
$(\a/2)(s^\prime)$. If $\a/2\not\in R_0$, we should reduce formula
(\ref{eq:ana}) to
\begin{equation}
c_{\a,\e}(\sigma_\e(u))c_\a(u)^{-1}
=\frac{(1-q_{\a^\vee}^{-\e}\a(v)^{-\e}\alpha(s^\prime)^{-1})
(1-\alpha(v)^{-1}\alpha(s^\prime)^{-1})}
{(1-q_{\a^\vee}^{-1}\a(v)^{-1}\alpha(s^\prime)^{-1})
(1-\alpha(v)^{-\e}\alpha(s^\prime)^{-1})}
\end{equation}
By conditions (i) and (iii) it is clear that poles and zeroes of
these functions will only meet $U$ if $\a(s^\prime)=1$ when $\a\in
R_0\cap R_1$ or $\a(s^\prime)=\pm 1$ if $\a\in 2R_0$. In these
cases the statement we want to prove reduces to the statement that
the function
\begin{equation}
f(x):=\frac{1-\exp(-\e x)}{1-\exp(-x)}
\end{equation}
is holomorphic and invertible on the domain $x\in
p+\a(\gamma^\prime+B)$, where $p$ is a real number and $\a\in
R_0$. By condition (i) both the denominator and the numerator of
$f$ have a zero in this domain only at $x=0$ (if this belong to
the domain), and this zero is of order $1$ both for the numerator
and the denominator. The desired result follows.
\end{proof}
Recall Theorem \ref{thm:gralg}. This result tells us that the
structure of the algebra with coefficients in the locally defined
meromorphic functions on $U$ is independent of the root labels. We
will now show that the subalgebra with analytic coefficients
(defined locally on $U$) is invariant for scaling transformations.
\begin{theorem}
The map
\begin{equation}
\begin{split}
j_\e:\H^{me}(U)&\mapsto\H_{q^\e}^{me}(U_\e)\\ \sum_{w\in
W_0}f_w\i_w^0&\mapsto \sum_{w\in W_0}(f_w\circ
\sigma_{1/\e})\i_{w,\e}^0
\end{split}
\end{equation}
\index{j@$j_\e:\H^{me}(U)\mapsto\H_{q^\e}^{me}(U_\e)$,
``scaling'' isomorphism of localized Hecke algebras}
defines an isomorphism of $\C$-algebras, with the property that
$j_\e(\F^{me}(U))\allowbreak =\F_{q^\e}^{me}(U_\e)$ and
$j_\e(\A^{me}(U))=\A^{me}_{q^\e}(U_\e)$. Moreover (and most
significantly), $j_\e(\H^{an}(U))=\H_{q^\e}^{an}(U_\e)$.
\end{theorem}
\begin{proof}
The map $j_\e$ as defined above is clearly a $\C$-linear
isomorphism by Theorem \ref{thm:gralg}. It is an algebra
homomorphism because we have
\begin{equation}
\begin{split}
j_\e(\sum_{u\in W_0}f_u\i_u^0\sum_{v\in W_0}g_v\i_v^0)&=
j_\e(\sum_{u,v\in W_0}f_ug^u_v\i_{uv}^0)\\ &= \sum_{u,v\in
W_0}(f_u\circ \sigma_{1/\e})(g^u_v\circ
\sigma_{1/\e})\i_{uv,\e}^0\\&= \sum_{u,v\in W_0}(f_u\circ
\sigma_{1/\e})(g_v\circ \sigma_{1/\e})^u\i_{uv,\e}^0\\&=
\sum_{u,v\in W_0}(f_u\circ \sigma_{1/\e})\i_{u,\e}^0(g_v\circ
\sigma_{1/\e})\i_{v,\e}^0\\&= j_\e(\sum_{u\in
W_0}f_u\i_u^0)j_\e(\sum_{v\in W_0}g_v\i_v^0)
\end{split}
\end{equation}
What remains is the proof that
$j_\e(\H^{an}(U))=\H^{an}_{q^\e}(U_\e)$. Notice that $\H^{an}(U)$
is the subalgebra generated by $\A^{an}(U)$ and the elements $T_s$
where $s=s_\a$ with $\a\in R_1$. The $j_\e$-image of $\A^{an}(U)$
equals $\A^{an}_{q^\e}(U_\e)$ since $\sigma_\e$ is an analytic
diffeomorphism. To determine the image of $T_s$ we use formula
Lemma 2.27(2) of \cite{EO}, applied to the situation
$W_0=\{e,s\}$. This tells us that
\begin{equation}
(1+T_s)=q_{\a^\vee} q_{2\a^\vee}c_\a(1+\i_s^0).
\end{equation}
Hence we see that
\begin{equation}
\begin{split}
j_\e(T_s)&=q_{\a^\vee} q_{2\a^\vee}(c_\a\circ
\sigma_{1/\e})(1+\i_{s,\e}^0)-1\\ &=q^{1-\e}_{\a^\vee}
q^{1-\e}_{2\a^\vee}(c_\a\circ \sigma_{1/\e})c_{\a,\e}^{-1}
(1+T_{s,\e})-1.
\end{split}
\end{equation}
By Lemma \ref{lem:ana} it is clear that this is indeed in
$\H^{an}_{q^\e}(U_\e)$, and that these elements together with
$\A^{an}_{q^\e}(m_{\e}(U))$ generate $\H^{an}_{q^\e}(U_\e)$.
\end{proof}
\subsection{Application to the residual algebras} In
order to prove that the residual algebras $\overline{\H^t}$ are
invariant for the scaling transformation $\q\to\q^\e$ it suffices
to consider the case $\overline{\H^r}$ for a residual point $r\in
T$. This follows from Theorem \ref{thm:mainind}, expressing
$\chi_t$ in terms of characters induced from discrete series
characters of proper parabolic subalgebras.

When $r=sc\in T$ is a residual point, the state $\chi_r$ has a
natural extension to the localized algebras $\H^{an}(U)$ where
$U=W_0r\exp(B)$, with $B$ an open ball in $\mathfrak{t}_\C$
satisfying the conditions \ref{cond} with respect to the point
$r\in T$. Because the radical $\operatorname{Rad}_r^{an}(U)$ of
the bitrace $(x,y)_r:=\chi_r(x^*y)$ on $\H^{an}(U)$ is contained
in the maximal ideal $\I_r^{an}(U)$ of functions in the center
$\Ze^{an}(U)$ which vanish in the orbit $W_0r$, we clearly have
\begin{equation}
\overline{\H^r}=\H^{an}(U)/\operatorname{Rad}_r^{an}(U).
\end{equation}
The structure of this algebra as a Hilbert algebra is given by the
bitrace defined by $\chi_r$. Therefore, we need to prove
independence of $\chi_r$ for the scaling transformation. We start
with a simple lemma:
\begin{lem} Let $h\in\H^{an}(U)$
We have
\begin{equation}
\frac{E_{q^\e,\sigma_\e(t)}(j_\e(h))}{q^\e(w_0)\Delta(\sigma_\e(t))}
=\frac{E_{t}(h)}{q(w_0)\Delta(t)}
\end{equation}
\end{lem}
\begin{proof}
For all $x\in X$ we have
\begin{align}
E_{q^\e,\sigma_\e(t)}(j_\e(\theta_xh))
&=E_{q^\e,\sigma_\e(t)}(j_\e(\theta_x)j_\e(h))\\\nonumber
&=(x\circ \sigma_{1/\e})(\sigma_\e(t))
E_{q^\e,\sigma_\e(t)}(j_\e(h))\\\nonumber
&=x(t)E_{q^\e,\sigma_\e(t)}(j_\e(h)),
\end{align}
showing that the left hand side has the correct eigenvalue for
multiplication of $h$ by $\theta_x$ on the left. For the
multiplication of $h$ by $\theta_x$ on the right a similar
computation holds. This shows, in view of Lemma \ref{lem:eisext}
and \cite{EO}, Proposition 2.23(3) that, for regular $t$ and
outside the union of all residual cosets, the left and the right
hand side are equal up to normalization. But both the left and the
right hand side are equal to $1$ if $h=T_e=1$. Hence generically
in $t$, we have the desired equality. Since both expressions are
holomorphic in $t$, the result extends to all $t\in T$.
\end{proof}
\begin{lem}
Let $\e\in(0,1]$ be given. We have, for all $h\in \H^{an}(U)$,
\begin{equation}
\chi_{q^\e,\sigma_\e(r)}(j_\e(h))=\chi_r(h).
\end{equation}
\end{lem}
\begin{proof}
Take a neighborhood $U=W_0r\exp(B)$ with $B$ satisfying conditions
\ref{cond}
relative to $r$. Let ${\cup\xi}\in\H_n(U)$ denote the $n$-cycle
defined by
${\cup\xi}=\cup_{r^\prime\in W_0r}\xi_{r^\prime}$. In view of
Proposition \ref{prop:cycle}, Definition \ref{dfn:Y} and
Definition \ref{dfn:chi} we see that, for all $h\in \H$,
\begin{equation}
\nu(\{W_0r\})\chi_r(h)=\int_{\cup\xi}
\left(\frac{E_{t}(h)}{q(w_0)\Delta(t)}\right)
\frac{dt}{q(w_0)c(t)c(t^{-1})}
\end{equation}
Let $r^\prime\in W_0r$. The scaling operation sends the root
labels $q$ to $q^\e$, and follows the corresponding path of $\e\to
\sigma_\e(r^\prime)$ of the residual point $r^\prime$. Obviously
the position of $t_0$ (in equation (\ref{eq:basic})) relative to
$\L^{\{\sigma_\e(r^\prime)\}}$ is independent of $\e$. And also,
the position of $e$ relative to the facets of the dual
configuration $\L_{\{\sigma_\e(r^\prime)\}}$ is independent of
$\e$, since the effect of the scaling operation on
$\L_{\{r^\prime\}}\subset T_{rs}$ simply amounts to the
application of the map $c\to c^\e$. In view of Proposition
\ref{prop:t0} and Proposition \ref{prop:cyinv}, we can take the
cycle $\sigma_\e({\cup\xi})\in H_n(\sigma_\e(U))$ in order to define the
state $\chi_{\sigma_\e(r)}$ of $\H_{q^\e}^{an}(\sigma_\e(U))$. In
other words, we have, for $h\in \H^{an}(U)$,
\begin{align}\label{eq:start}
\nu_{q^\e}(\{&W_0\sigma_\e(r)\})\chi_{q^\e,\sigma_\e(r)}(j_\e(h))\\
\nonumber
&=\int_{\sigma_\e({\cup\xi})}
\left(\frac{E_{t}(j_\e(h))}{q^\e(w_0)\Delta(t)}\right)
\frac{dt}{q^\e(w_0)c_{\e}(t)c_{\e}(t^{-1})}\\\nonumber
&=\int_{{\cup\xi}}\left(\frac{E_{\sigma_\e(t)}(j_\e(h))}{q^\e(w_0)
\Delta(\sigma_\e(t))}\right)
\frac{d(\sigma_\e(t))}{q^\e(w_0)c_{\e}
(\sigma_\e(t))c_{\e}(\sigma_\e(t^{-1}))}\\\nonumber
&=\int_{{\cup\xi}}\left(\frac{E_{t}(h)}{q(w_0)\Delta(t)}\right)
\phi_\e(t)\frac{dt}{q(w_0)c(t)c(t^{-1})},
\end{align}
where
\begin{equation}
\phi_\e(t):=\frac{\e^nq(w_0)c(t)c(t^{-1})}
{q^\e(w_0)c_{\e}(\sigma_\e(t))c_{\e}(\sigma_\e(t^{-1}))}.
\end{equation}
By Lemma \ref{lem:ana}, the function $t\to\phi_\e$ extends, for
all $\e\in(0,1]$, to a regular holomorphic function on $U$.
Clearly, $\phi_\e$ is $W_0$-invariant. In other words, $\phi_\e$
is an element of $\Ze^{an}(U)$. Its value in $W_0r$ can be
computed easily, if we keep in mind that the index
$i_{\{r^\prime\}}=n$ (by Theorem \ref{thm:equal} applied to the
residual coset $r^\prime$). We obtain, by a straightforward
computation:
\begin{equation}
\phi_\e(W_0r)=\frac{m_{q^\e,\{\sigma_\e(r)\}}(r)}{m_r(r)}
=\frac{\nu_{q^\e}(\{W_0\sigma_\e(r)\})}{\nu(\{W_0r\})}.
\end{equation}
We now continue the computation which we began in equation
(\ref{eq:start}), using the fact that $\phi_\e\in\Ze^{an}(U)$ and
the fact that $\chi_r$ extends uniquely to $\H^{an}(U)$ in such a
way that for all $\phi\in\Ze^{an}(U)$ and $h\in \H^{an}(U)$,
$\chi_r(\phi h)=\phi(r)\chi(h)$. We get
\begin{align}
\nu_{q^\e}(\{W_0\sigma_\e(r)\})\chi_{q^\e,\sigma_\e(r)}(j_\e(h))&=
\nu(\{W_0r\})\chi_r(\phi_\e h)\\\nonumber
&=\nu_{q^\e}(\{W_0\sigma_\e(r)\})\chi_r(h).
\end{align}
This gives the desired result.
\end{proof}
\begin{thm}\label{thm:bch}
The ``base change'' isomorphism $j_\e$ induces an isomorphism
\begin{equation}
\overline{j}_\e:\overline{\H^r}
{\tilde{\longrightarrow}}\overline{\H^{\sigma_\e(r)}_{q^\e}}
\end{equation}
of Hilbert algebras. In particular, the positive constants
$d_{\Ri_P,\d}$ (in the notation of Theorem
\ref{thm:mainind}, see also
equation (\ref{eq:plameasure}))
(in Corollary \ref{cor:fdim} these constant were denoted
by $d_{r,i}$) are independent of $\q$.
\end{thm}
\begin{proof}
This is an immediate consequence of the previous lemma.
\end{proof}
\section{Applications and closing remarks}\label{app}
\subsection{Formation of L-packets of unipotent
representations}\label{sub:uni}
Let $F$ be a nonarchimedean local field,
and let $G$ be a split simple algebraic
group of adjoint type defined over $F$.
We denote by $\mathcal{G}$
the group of $F$-rational points in $G$.
The finite set of
irreducible unipotent discrete series representations of
$\mathcal{G}$ is
by definition the disjoint
union of the discrete series
constituents of
induced representations of the form
$\sigma_{\mathcal{P}}^{\mathcal{G}}$,
where $\mathcal{P}$ is a
parahoric subgroup of $\mathcal{G}$,
$\sigma$ is a cuspidal unipotent representation of the Levi
quotient $L:=\mathcal{P}/{\mathcal{U}_\mathcal{P}}$ of
$\mathcal{P}$, and the union is taken over a complete set
of representatives of conjugacy classes of pairs
$(\mathcal{P},\sigma)$.

The formal dimension is an effective tool to partition
these unipotent discrete series representations into L-packets.
This observation is due to
Reeder \cite{Re0}. He conjectured that the formal dimensions of
the unipotent
discrete series representations of $\mathcal{G}$ within one unipotent
discrete series L-packet are proportional, with a rational ratio
of proportionality {\it independent} of $F$, and used this in
\cite{Re0} to form the unipotent L-packets for groups of small
rank.

It is known \cite{M} that the endomorphism algebra
$\H(\mathcal{G},\mathcal{P},\sigma)$ of
$\sigma_{\mathcal{P}}^{\mathcal{G}}$ has the structure of an
affine Hecke algebra, whose root datum and root labels depend only
on $\mathcal{P}$. The root labels are integral powers $\q^{n_a}$ of the
cardinality $\q$ of the residue field of $F$ (cf. \cite{Lu3}), and
are explicitly known. Moreover, if we define a trace functional
$\operatorname{Tr}$ on $H(\mathcal{G},\mathcal{P},\sigma)$ by
$\operatorname{Tr}(f):=\operatorname{Tr}(f(e),V_\sigma)$, this
corresponds to the trace $\tau$ studied in this paper by the
formula
\begin{equation}\label{eq:tr}
\operatorname{Tr}=\operatorname{Vol}(\mathcal{P})^{-1}
\operatorname{dim}(V_\sigma)\tau.
\end{equation}
Thus there is a bijection between the set of discrete series
representations of $H(\mathcal{G},\mathcal{P},\sigma)$
(in the sense of this paper)
and the unipotent discrete series representations arising from the
pair $(\mathcal{P},\sigma)$. The formal dimension of such a discrete series
representation of $\mathcal{G}$ is then equal to the formal dimension of
the corresponding discrete series representation of the affine
Hecke algebra $H(\mathcal{G},\mathcal{P},\sigma)$, but with its trace
$\operatorname{Tr}$ normalized by equation (\ref{eq:tr}).

In \cite{HOH} we computed the formal dimension of the discrete
series representations of the ``anti-spherical'' subalgebra of the
affine Hecke algebra, i.e. the commutative subalgebra
$e_-\H e_-=e_-\Ze$, where $e_-$ denotes the idempotent of $\H_0$
corresponding to the sign representation
(see Subsection \ref{sub:range}). The formula we obtained
was expressed entirely in terms of the central character of the
representation, the root datum and the root labels. We conjectured
in \cite{HOH} that our formula would also hold for the full affine
Hecke algebra. For the group of type $E_8$, we showed that this,
in combination with Reeder's conjecture, leads to a partitioning
of the unipotent discrete series L-packets which is in agreement
with Lusztig's conjecture \cite{Lu0} for the Langlands parameters
of the members of these packets.
In other words, the formal
dimension seems to be a sufficient criterion to separate the
L-packets of unipotent representations in the case $E_8$.
%This means that in the case $E_8$, the unipotent, discrete series
%L-packets of $\mathcal{G}$ consist of
%the irreducible constituents of various
%$\sigma^{\mathcal{G}}_{\mathcal{P}}$, corresponding to discrete
%series representations of
%$\H(\mathcal{G},\mathcal{P},\sigma)$ with proportional
%formal dimensions (as functions of $\q$).
Theorem
\ref{thm:mainp} of this paper proves the conjecture in \cite{HOH}
mentioned above.

Reeder \cite{Re} proves an exact formula for the formal dimension
of the unipotent discrete series representation of all split
exceptional groups, based on a result of Schneider and Stuhler
\cite{schstu}.
In this approach one first represents the formal
dimension by an alternating sum of rational functions
(depending on the $K$-types (cf. Subsection \ref{sub:types})),
rather than the product formula which we have obtained.
On the other hand, there are no intractable constants
such as the constant $d_\d$ in our formula.
Using his previous work on non-standard intertwining operators for
affine Hecke algebra modules and Theorem \ref{thm:nu} of the
present paper, Reeder gave the precise partitioning of the unipotent
discrete series for exceptional groups into L-packets, in complete
agreement with Lusztig's conjecture mentioned above.

Recently, Lusztig \cite{Lu4} established the partitioning of
unipotent discrete series representations into L-packets
if $G$ is split over an unramified extension of $F$.
This is based on a different approach. It is worth
mentioning that this classification includes a geometric
parametrization of the set $\Delta_{W_0r}$ of discrete series
representations with central character $W_0r$ if the affine Hecke
algebra arises as an endomorphism algebra of an induced
representation of the form $\sigma_{\mathcal{P}}^{\mathcal{G}}$.
\subsection{Operator norm estimate and the Schwartz
completion}\label{sub:normunif}
\subsubsection{Uniform norm estimate}
We know that the generators $N_i:=q(s_i)^{-1/2}T_i$ satisfy
\begin{equation}
\Vert N_i\Vert_o=\operatorname{max}\{q(s_i)^{\pm 1/2}\}
\end{equation}
Therefore we have the trivial estimate $\Vert N_w\Vert_o\leq
\operatorname{max}\{q(w)^{\pm 1/2}\}$.
By the spectral decomposition it is easy to see that
the operator norm $\Vert N_w\Vert_o$ is actually bounded by a
{\it polynomial} in $\nc(w)$:
\begin{theorem}\label{thm:normunif} Let $\q>1$ be fixed.
\begin{enumerate}
\item[(i)] There exist constants $C\in{\mathbb R}_+$ and $d\in\N$ such
that for all $w\in W$,
\begin{equation}
\Vert N_w \Vert_o\leq C(1+\nc(w))^d
\end{equation}
\item[(ii)] For a residual point $r\in T$, let
$\Delta_{W_0r}:=\Delta_{\Ri,W_0r}$ denote the collection of all
discrete series representations with central character $W_0r$.
Denote by $\Vert x \Vert_{ds}$ the operator norm of the left
multiplication by $x\in\H$, restricted to the finite dimensional
subspace
\begin{equation}
\hf_{ds}:=\oplus_{W_0r}
\oplus_{\pi\in\Delta_{W_0r}}\operatorname{End}(\pi)\subset\hf.
\end{equation}
Then there exist constants $C,\e>0$ such that
\begin{equation}
\Vert N_w \Vert_{ds} \leq C\q^{-\e l(w)}.
\end{equation}
\end{enumerate}
\end{theorem}
\begin{proof}
(i) As was explained in the proof of Lemma \ref{lem:cas}, it is
sufficient to prove this statement for $w=x\in (Z_X+Q)\cap X^+$.
Recall that
in this case $N_x=\theta_x$. We have
\begin{equation}
\Vert \theta_x\Vert^2_o=\Vert \theta_x^*\theta_x\Vert_o=
\sup\{\sigma(\pi(T_{w_0}\theta_{-w_0x}T_{w_0}^{-1}\theta_x)) \mid
\pi\in\hat\cf\},
\end{equation}
where $\sigma(A)$ denotes the spectral radius of $A$. According to
Theorem \ref{thm:mainp}, the spectrum $\hat\cf$ equals the union
of the compact sets $\hat\cf_{\Ri_P,W_Pr,\d}$. This set is by
definition the closure in $\hat\cf$ of the set
$\pi(\Gamma_{\Ri_P,W_Pr,\d}^{gen})$.
It is well known that for any
$h\in\H$, the map $\pi\to\Vert \pi(h)\Vert_o$ is lower
semi-continuous as a function of $\pi\in\hat\cf$ (cf. \cite{Fell},
VII, Proposition 1.14). Since there are only finitely many
triples $(R_P,W_Pr,\d)$, it is sufficient to show that
there exist constants $C$, $d$ such that the spectral radius of
\begin{equation}
\pi(\Ri_P,W_Pr,\d,t^P)(T_{w_0}\theta_{-w_0x}T_{w_0}^{-1}\theta_x)
\end{equation}
is bounded by $C(1+\nc(x))^{d}$, uniformly in $t^P$.

The roots of a monic polynomial are bounded by the sum of the
absolute values of the coefficients of the equation (including the
top coefficient $1$). Hence the spectral radius of an $m\times m$
matrix $A$ is bounded by a polynomial of degree $m$ in
$\operatorname{max}(A):=\max\{|a_{i,j}|\mid 1\leq i,j\leq m\}$.
Since $m\leq |W_0|$ it is sufficient to show that there exists a
suitable basis for the parameter family $\pi(\Ri_P,W_Pr,\d,t^P)$ (with
$t^P\in T^P_u$) of representations, in which the matrix
coefficients of the $\theta_x$ (with $x\in (Z_X+ Q)\cap X^+$) are uniformly
bounded by $C(1+\nc(x))^{d}$ for suitable constants $C$ and $d$.

As in the proof of Proposition \ref{prop:indtemp},
there exists a basis
$T_{w_i}\otimes (v_j)$ of $V_\pi=\H(W^P)\otimes V$,
the representation space of $\pi(\Ri_P,W_Pr,\d,t^P)$,
such that the $\theta_x$ ($x\in X$) simultaneously act by means
of upper triangular matrices in this basis.
Moreover, by Proposition
\ref{prop:indtemp} it is clear that the diagonal elements are
bounded in norm by $1$ when $x\in X^+$.
By the compactness of $T^P_u$ we conclude
that there exists, for each $x\in X^+$, an unipotent upper
triangular matrix $U_x$ with positive coefficients such that every
matrix coefficient of $M_x(t^P):=\pi(\Ri_P,W_Pr,\d,t^P)(\theta_x)$ in
the above basis is uniformly (in $t^P$) bounded by the
corresponding matrix coefficient of $U_x$.

Let us denote by $P$ the set of $n\times n$ matrices with
non-negative entries, and introduce the notation
$|A|=(|A_{i,j}|)_{i,j}$ for complex matrices $A$. Introduce a
partial ordering in $P$ by defining $A\leq B$ if and only if
$B-A\in P$. Since $P$ is a semigroup for matrix multiplication, it
is clear that if $A,B$ and $C$ are in $P$ and $A\leq B$, then
$AC\leq BC$. In addition we have the rule $|AB|\leq |A||B|$ for
arbitrary complex matrices $A$ and $B$.

Let $x_1,\dots, x_m$ denote a set of $\mathbb{Z}_+$-generators for
the cone $Q^+$,
and let moreover $x_{m+1},\dots,x_N$ be a basis of $Z_X$. Put
$M_{i,\e}(t^P):=M_{\e x_i}(t^P)$ for $1\leq i\leq N$, $\e=\pm 1$
with $\e=1$ if $i\leq m$.
We can thus find an upper triangular
unipotent matrix $U\in P$ such that  for all $i$, $\e$, $t^P$:
$|M_{i,\e}(t^P)|\leq U$. If we write $x=\sum_il_ix_i$ with $l_i\geq 0$
if $i\leq m$, then
\begin{align*}
\nc(x)&=x(2\rho^\vee)+\Vert \sum_{i>m}l_ix_i\Vert\\
&=\sum_{i\leq m}l_ix_i(2\rho^\vee)+\Vert \sum_{i>m}l_ix_i\Vert\\
\end{align*}
with $x_i(2\rho^\vee)\geq 1$ if $i\leq m$. From this we see that
there exists a constant $K$ independent of $x$ such that
$\a:=\sum |l_i|\leq K\nc(x)$.
Thus, with $\log(U)$ the nilpotent logarithm of $U$,
and $l_i=\e(i)|l_i|$:
\begin{equation}\label{eq:est}
\begin{split}
\operatorname{max}(M_x(t^P))&=
\operatorname{max}(|M_1^{l_1}(t^P)\dots M^{l_N}_N(t^P)|)\\ &\leq
\operatorname{max}(|M_{1,\e(1)}^{|l_1|}(t^P)|
\dots|M^{|l_N|}_{N,\e(N)}(t^P)|)\\ &\leq
\operatorname{max}(U^\a)\\ &\leq
\operatorname{max}(\exp(\a\log(U)))\\ &\leq \sum_{i=0}^d
\operatorname{max}(\log(U)^i)\a^i/i!\\ &\leq c_U(1+\nc(x))^d\\
\end{split}
\end{equation}
where $c_U$ is a constant depending on $U$ only, and $d$ is the
degree of the polynomial function $\a\to\exp({\a\log(U)})$. This
finishes the proof.

(ii) As in the proof of (i), but we restrict ourselves to the
(finitely many) discrete series representations. This implies that
we can find $U\in P$ unipotent and $\e>0$ such that for all $i$,
$|M_i|\leq \q^{-2\e} U$. Inserting this in the inequalities
(\ref{eq:est}) we find that the matrix entries of $M_x$ are bounded
by $C\q^{-\e l(x)}$, with $C$ independent of $x\in X^+$. Hence the
spectral radius of $\q^{2\e l(x)}M_x^*M_x$ is uniformly bounded
for $x\in X^+$, proving the desired estimate.
\end{proof}
\begin{cor}\label{cor:onlytemp}
$\hat\cf$ consists only of tempered representations.
\end{cor}
\begin{proof}
The character $\chi_\pi$ of $\pi\in\hat\cf$ is a positive trace, and
thus satisfies the inequality $|\chi_\pi(x)|\leq \chi_\pi(1)\Vert
x\Vert_o$ by Corollary \ref{cor:contrace}. By Casselman's
criterion
Lemma \ref{lem:cas}
and Theorem \ref{thm:normunif} this implies that $\pi$
is tempered.
\end{proof}
\begin{prop}\label{prop:onediml}
The trivial representation $\pi_{triv}(T_w)=q(w)$ is tempered if
and only if
the point $r_{triv}\in T_{rs}$ defined by
$\forall\a\in F_0:\a(r_{triv})=q_{\a^\vee/2}^{1/2}q_{\a^\vee}$
satisfies $r_{triv}\in\overline{T_{rs}^-}$.
It is discrete series if and only if $r_{triv}\in T_{rs}^-$.
The Steinberg
representation $\pi_{St}(T_w)=(-1)^{l(w)}$ is tempered if and
only if $r_{triv}^{-1}:=r_{St}\in\overline{T_{rs}^-}$, and
discrete series if and only if $r_{St}\in T_{rs}^-$.
\end{prop}
\begin{proof}
This is well known, and follows easily by the remark that
the restriction to $\{\theta_x\mid x\in X\}$ of the trivial
representation is equal the square root of the Haar modulus $\d$
(see e.g. \cite{EO}, Corollary 1.5):
$\pi_{triv}(\theta_x)=q(x)^{1/2}=\delta^{1/2}(x):=x(r_{triv})$.
Now apply the Casselman criteria Lemma \ref{lem:cas}. Similar
remarks apply to the case of the Steinberg representation.
\end{proof}
\begin{cor}\label{cor:onedimlext}
If the trivial representation extends to $\cf$ then $r_{triv}\in
\overline{T_{rs}^-}$. If the Steinberg representation extends to
$\cf$ then $r_{triv}\in\overline{T_{rs}^+}$.
\end{cor}
\begin{proof}
Use Corollary \ref{cor:onlytemp} and Proposition
\ref{prop:onediml}.
\end{proof}
\subsubsection{The Schwartz completion of $\H$}\label{subsub:schwartz}
Using Theorem \ref{thm:normunif} we now define a Fr\'echet completion
of $\H$.
For all $n\in\mathbb{N}$ we define a norm $p_n$ on $\H$ by
\begin{equation}
p_n(h)=\max_{w\in W}|(N_w,h)|(1+\nc(w))^n
\end{equation}
Here $\nc$ denotes the norm function on $W$ which was defined by
equation (\ref{eq:bigL}).
\begin{thm}
The functions $\tau$, $*$ and the multiplication $\cdot$ of $\H$ are
continuous with respect to the family of norms $p_n$.
\end{thm}
\begin{proof}
The continuity of $\tau$ and $*$ is immediate from the
definitions. So let us look at the multiplication.
Let us write
\begin{equation}
N_uN_v=\sum_w c_{u,v}^w N_w
\end{equation}
It is easy to see that $w(0)^0=u(0)^0+v(0)^0$ and
that $l(w)\leq l(u)+l(v)$ if $c_{u,v}^w\not=0$. Therefore
\begin{equation}\label{eq:1}
c_{u,v}^w\not=0\Rightarrow \nc(w)\leq \nc(u)+\nc(v),
\end{equation}
and by Theorem \ref{thm:normunif}, there exist constants $C,d$
such that for all $u,v$ and $w$:
\begin{equation}\label{eq:2}
|c_{u,v}^w|\leq C\min\{(1+\nc(u))^d,(1+\nc(v))^d\}
\end{equation}
We put $D_w=\{(u,v)\in W\times W\mid c_{u,v}^w\not=0\}$.
It is easy to see that there exists a $b\in\mathbb{N}$ such that
\begin{equation}
\sum_{u\in W}\frac{1}
{(1+\nc(u))^b}=\mu<\infty
\end{equation}
converges. By (\ref{eq:1})
we have that
\begin{equation}
(1+\nc(u))(1+\nc(v))\geq (1+\nc(w))
\end{equation}
for all $(u,v)\in D_w$.
Given $n\in \mathbb{N}$, let $k=\max\{b+d,n\}$.
Using these remarks we see that for all $0\not=x=\sum x_u N_u$ and
$0\not=y=\sum y_v N_v$ in $\H$ the following holds:
\begin{align*}
\frac{p_n(xy)}{p_{2k}(x)p_{2k}(y)}
&=\frac{1}{p_{2k}(x)p_{2k}(y)}\max_w |(xy,N_w)|(1+\nc(w))^n\\
&\leq\max_w\sum_{u,v\in D_w}\frac{|x_u||y_v|}{p_{2k}(x)p_{2k}(y)}
|c_{u,v}^w|(1+\nc(w))^n\\
&\leq C\sup_w\sum_{u,v\in D_w}\frac{\min\{(1+\nc(u))^d,(1+\nc(v))^d\}
(1+\nc(w))^n}{(1+\nc(u))^{2k}(1+\nc(v))^{2k}}\\
&\leq C\sup_w\sum_{u,v\in D_w}\frac{(1+\nc(u))^d(1+\nc(v))^d
(1+\nc(w))^{n-k}}{(1+\nc(u))^{k}(1+\nc(v))^{k}}\\
&\leq C\sum_{u,v}
(1+\nc(u))^{d-k}(1+\nc(v))^{d-k}\\
&\leq \mu^2 C
\end{align*}
This finishes the proof.
\end{proof}
Notice that, by Theorem \ref{thm:normunif}, $\|x\|_o\leq Cp_d(x)$
for all $x\in\H$. Therefore the completion of $\H$ with respect to
the family of norms $p_n$ will be a subspace of $\cf$.
\begin{dfn}\label{dfn:schw}
We define the Schwartz completion $\mathfrak{S}$ of $\H$ by
\begin{equation}
\mathfrak{S}:=\{x=\sum_w x_w N_w\in\H^*\mid p_n(x)<\infty\ \forall
n\in\mathbb{N}\}.
\end{equation}
We have $\H\subset\mathfrak{S}\subset{\cf}$, and $\mathfrak{S}$ is
a $*$-subalgebra of $\cf$. $\mathfrak{S}$ is a nuclear
Fr\'echet algebra with respect to the topology defined by the family
of norms $p_n$.
It comes equipped with continuous trace $\tau$ and anti-involution $*$.
\end{dfn}
\index{S@$\sf$, the Schwartz completion of $\H$}
\begin{cor} (of definition)
The topological dual $\mathfrak{S}^\prime$ is the space of
tempered linear functionals on $\H$.
\end{cor}
\subsection{A Hilbert algebra isomorphism; abelian subalgebras}
\label{sub:range}
There exists a trace preserving $*$-algebra isomorphism
\begin{align*}
i:\H(\Ri,q)&\to\H(\Ri,q^{-1})\\
N_w&\to(-1)^{l(w)}N_w.
\end{align*}
(see \cite{HOH}).
Clearly
this induces a $*$-algebra isomorphism, and $\tau$ is respected.
Thus $i$ induces an isomorphism of $C^*$-algebras
$i:\cf(\Ri,q)\to\cf(\Ri,q^{-1})$, also respecting the traces.
The corresponding homeomorphism
$\hat{i}:\hat{\cf}(\Ri,q^{-1})\to\hat{\cf}(\Ri,q)$ is
therefore Plancherel measure
preserving.
Note that $i$ restricts to a (Plancherel measure preserving)
$*$-isomorphism from the subalgebra
$e_+\H(\Ri,q)e_+$ (the spherical subalgebra)
to the subalgebra $e_-\H(\Ri,q^{-1})e_-$
(the anti-spherical subalgebra) (see \cite{HOH}).
\subsubsection{The Plancherel measure of the center}
The commutative subalgebras $e_\pm\H(\Ri,q)e_\pm$
are both isomorphic as algebras to the center $\Ze$
via the Satake isomorphism $\Ze\ni z\to e_\pm z\in
e_\pm\H(\Ri,q)e_\pm$. These subalgebras are
commutative Hilbert subalgebras of $\H(\Ri,q)$.
They are
in general not isomorphic as Hilbert algebras.
The Hilbert algebra isomorphism $i$ restricts to
an isomorphism $\Ze(\Ri,q)\to\Ze(\Ri,q^{-1})$.
Recall the $\nu$ is the Plancherel measure of
$\overline{\Ze}\subset \cf$.
The above is reflected by the symmetry
\begin{cor}
$\nu(t,q)=\nu(t,q^{-1})$,
\end{cor}
\noindent{which can be verified directly}
(see Theorem \ref{thm:nu}
and Proposition \ref{prop:par}).

The spherical algebra $e_+\H(\Ri,q)e_+$ with $X=P$ (weight lattice)
and $q(s)>1$ has a very important
basis, uniquely defined by orthogonality and by a triangularity
requirement with respect to the standard monomial basis
$e_+m_\l$ (with $\l\in P^+$ and
$m_\l=\sum_{\mu\in W_0\l}t^\mu\in\Ze$). In type $A$ these
are the Hall-Littlewood polynomials.
It would be interesting to study such orthogonal,
triangular bases for the center $\Ze$ as well.
\subsection{Central idempotents of $\cf$ and $\mathfrak{S}$}
\label{sub:idem}
\index{C@$\hat\cf_\O$, component of $\hat\cf$, the closure of
$[\pi](W(\O)\backslash\O^{gen})\subset\hat\cf$|(}
Recall that $\hat\cf=\cup\hat\cf_\O$ (union
over a complete set of representatives of the association classes
orbits). In Theorem \ref{thm:mainp} we have shown
that two distinct closed subsets in $\hat\cf$ of the form
$\hat\cf_{\O_i}$ ($i=1,2$) intersect in a subset of measure $0$.
In fact more is true: according to \cite{DO}, these closed
subsets are the components of $\hat\cf$.

There is a bijection $I\to\hat{I}$ between the closed
two-sided ideals of $\cf$ and the open subsets of $\hat\cf$. Hence
the decomposition of $\hat\cf$ into components $\hat\cf_\O$
corresponds to the
decomposition of $1\in\cf$ as a sum of minimal central orthogonal
idempotents $e_\O$ of $\cf$.

If $\O=\K_\d\backslash(\Lambda_\g\times\d)$, $e_\O\in\cf$
is determined by
\begin{equation}
\F(e_\O)(\pi)=
\begin{cases}
\operatorname{Id_{i(V_\d)}}\text{\ if\
}[\pi]\in\hat\cf_\O\\ 0\text{\ else\ },\\
\end{cases}
\end{equation}
\index{e1@$e_\O\in\sf$, central idempotent associated with $\O$}
where $\F$ is the isomorphism of Theorem \ref{thm:mainp}.

In fact, the results of \cite{DO} on smooth
wave packets even imply that $e_\O\in\mathfrak{S}$.
We thus have the following decomposition of the unit element in
central, Hermitian, mutually orthogonal, minimal idempotents
of $\sf$:
\begin{equation}
1=\sum_{\O} e_\O.
\end{equation}
\begin{theorem}\label{thm:idemexp}
\begin{enumerate}
\item[(i)] Let $\O=\K_\d\backslash(\Lambda_\g\times\d)$. Then
\begin{equation}
\begin{split}
(e_\O,e_\O)=|W^P|
\dim(\d)\mu_{Pl}(\hat\cf_\O).
\end{split}
\end{equation}
\item[(ii)]
These idempotents have the following expansion with respect to any
orthonormal basis $B$ of $\hf$:
\begin{equation}
e_\O=\sum_{b\in B}\chi_\O(b^*)b
\end{equation}
with
\begin{equation}
\chi_\O(b):=
\int_{\pi\in\hat\cf_\O}\chi_\pi(b)d\mu_{Pl}(\pi).
\end{equation}
In particular this holds with respect to the orthonormal basis
$(N_w)_{w\in W}$.
\end{enumerate}
\end{theorem}
\begin{proof}
(i). The dimension of $\pi\in\hat\cf_\O$ equals
$|W^P|\dim(\d)$ on an open dense subset, and the measure
$\mu_{Pl}$ is absolutely continuous with respect to the Haar
measure on $\O$. Hence, using the fact that $\F$ is an
isometry, we find
\begin{equation}
\begin{split}
(e_\O,e_\O)&=\int_{\pi\in\hat\cf}
\operatorname{Tr}_{V_\pi}(\F(e_\O)(\pi))d\mu_{Pl}(\pi)\\
&=|W^P|\dim(\d)\mu_{Pl}(\hat\cf_\O).\\
\end{split}
\end{equation}

(ii). As in (i) we get
\begin{equation}
\begin{split}
\chi_\O(b)&=(e_\O,b)\\
&=\int_{\pi\in\hat\cf}\operatorname{Tr}_{V_\pi}
(\F(e_\O)(\pi)\F(b)(\pi))d\mu_{Pl}(\pi)\\
&=\int_{\pi\in\hat\cf_\O}
\operatorname{Tr}_{V_\pi}(\pi(b))d\mu_{Pl}(\pi)\\
&=\int_{\pi\in\hat\cf_\O}\chi_\pi(b)d\mu_{Pl}(\pi).\\
\end{split}
\end{equation}
\end{proof}
The above depends on the results of \cite{DO},
but in the special case of isolated points in $\hat\cf$
these facts are more elementary:
\begin{prop}
When $\pi\in\Delta_{\Ri}$,
put $\O_\pi=\{[\pi]\}$ for the corresponding component
of $\hat\cf$ (an isolated point). Put $e_\pi$ for the
corresponding central idempotent of $\cf$.
The expansion
\begin{equation}
e_\pi=\mu_{Pl}(\{\pi\})\sum_{b\in B}\chi_{\pi}(b^*)b
\end{equation}
is convergent in $\sf$.
\end{prop}
\begin{proof}
The expansion follows as in the above Theorem. It is
convergent in $\sf$ because of Corollary \ref{rem:cas}
and Definition \ref{dfn:schw}.
\end{proof}
\begin{cor}\label{cor:normeq}
Let $B$ be a Hilbert basis of $\hf$. For any residual point $r$
and discrete series representation $\pi\in\Delta_{W_0r}$ we have
\begin{equation}
\sum_{b\in B}|\chi_\pi(b)|^2=\frac{\dim(\pi)}
{|W_0r|{\overline \ka}_{W_0r}d_{\pi}m_{\{r\}}(r)}
\end{equation}
where the constant ${\overline \ka}_{W_0r}\in\mathbb{Q}$ is defined by
(\ref{eq:denk}), the constant $d_{\pi}\in\mathbb{R}_+$ by
Definition \ref{dfn:resalg} and $m_{\{r\}}(r)$ by Theorem
\ref{thm:nu}.
\end{cor}
\begin{proof}
This follows from $\mu_{Pl}(\pi)=|W_0r|{\overline \ka}_{W_0r}
d_{\pi}m_{\{r\}}(r)$. Note that $d_{\pi}$ is indeed constant
(i.e. independent of $\q$) by Theorem \ref{thm:bch}. Also note
Conjecture \ref{rem:ell}.
\end{proof}
\index{C@$\hat\cf_\O$, component of $\hat\cf$, the closure of
$[\pi](W(\O)\backslash\O^{gen})\subset\hat\cf$|)}
\subsection{Some examples}\label{sub:example}
\subsubsection{The Steinberg representation}
A basic example is the Steinberg representation.
We obtain a well known expression for the Poincar\' e series of $W$.

This result was first (for equal labels, using Morse theory)
derived by Bott \cite{Bott}, and by elementary means by
Steinberg \cite{Stein}.
Macdonald \cite{Ma2} observed that the arbitrary parameter
case can be obtained by Steinberg's method.
Macdonald proved formula (\ref{poin}) below,
expressing the Poincar\'e polynomial in terms of the roots,
in an elementary way using case-by-case
verifications.
In \cite{Ma3} Macdonald reproved the formula in a uniform way.
Also note that the Steinberg representation is a representation
of $e_-\H e_-$ (cf. \ref{sub:range}).
Hence its formal degree can also be computed
by means of the (simpler) techniques of \cite{HOH}.

We assume that $Q\subset X\subset P$.
Let $\pi_{St}$ be the Steinberg representation, which is the
representation defined by $\pi_{St}(T_w)=(-1)^{l(w)}$. This is
a one dimensional discrete series representation
provided that $r_{St}\in T_{rs}^-$ (see Proposition \ref{prop:onediml}).
Recall that  $r_{St}\in T_{rs}$ is defined by $\forall\a\in F_0:
\a(r_{St})=q_{\a^\vee/2}^{-1/2}q_{\a^\vee}^{-1}$.
Generically this residual point is regular.
In this regular case,
the residual codimension $1$ cosets containing
$wr_{St}$ form a normal crossing divisor
$D_w$ locally at $wr_{St}$. By Proposition \ref{prop:antidual} we find
that $[\xi_{wr_{St}}]\in H_n(U\backslash D_w)$ (with $U$ a small ball
around $wr$) is zero if $w\not=e$, and for $w=e$ it is
straightforward to see that $[\xi_{r_{St}}]=(-1)^n.C$, where $C$ is the
positive generator of $H_n(U\backslash D_e)$. In view of
(\ref{eq:tochhandig}) we find that (still in the regular case)
\begin{equation}
\ka_{wr_{St}}=\frac{(-1)^n\d_{w,e}}{|X:Q|},
\end{equation}
and thus that ${\overline \ka}_{W_0r_{St}}=(-1)^n|W_0|^{-1}|X:Q|^{-1}$.

Hence if $\pi$ is a discrete series representation with central
character $r_{St}$ then, assuming that $r_{St}$ is regular,
equation (\ref{eq:rats}) implies that $\pi$ can only have a
nonzero weight space for the weight $r_{St}$. But this weight space is
one dimensional, so that $\pi=\pi_{St}$. Hence in the regular case,
the Steinberg representation is the only member of
$\Delta_{W_0r_{St}}$, and thus $d_{\pi_{St}}=1$ in this case.

Inserting the values of these constants, above identity thus
specializes to (using the Hilbert bases $(N_w)_{w\in W}$)
Macdonald's product formula for the Poincar\'e series of $W$:
\begin{equation}\label{poin}
\sum_{w\in W}q(w)^{-1}=\frac{(-1)^n|X:Q|}{m_{\{r_{St}\}}(r_{St})}
\end{equation}
By continuity, this formula holds in general provided that
$r_{St}\in T_{rs}^-$.
When we take $q(s)=\q$ for all $s\in S^{\mathrm{aff}}$, and
$X=Q$, then we obtain
\begin{equation}
\sum_{w\in W}\q^{-l(w)}=
\prod_{i=1}^n\frac{(\q^{m_i+1}-1)}{(\q^{m_i}-1)(\q-1)}
\end{equation}
where $(m_i)$ is the list of exponents of $W_0$.
\subsubsection{The subregular unipotent orbit of $Sp_{2n}$}
Let $F$ be a
nonarchimedean
local field, and
let $\q$ be the
cardinality of its residue field. Consider the group
$\mathcal{G}=\operatorname{SO}_{2n+1}(F)$ for $n\geq 3$.
The Langlands dual group of $\mathcal{G}$ is
$\hat{G}=\operatorname{Sp}_{2n}(\mathbb{C})$, whose
root datum (with basis) we write as
$\Ri=(C_n,\Z^n,B_n,\Z^n,F_0)$,
with $F_0=(e_1-e_2,\dots,e_{n-1}-e_n,2e_n)$.
We normalize the Haar measure of $\mathcal{G}$
by $\operatorname{Vol}(\mathcal{I})=1$, where
$\mathcal{I}\subset\mathcal{G}$ is an Iwahori subgroup.
Let us compute the formal dimensions of the
irreducible square integrable, Iwahori-spherical
representations of $\mathcal{G}$ whose
Kazhdan-Lusztig parameters $(r_u,u,\rho)$ (cf.
Appendix \ref{KL}) are such that $u$ is the subregular
unipotent orbit of $\hat{G}$.
Take $r_u\in T_{rs}$ dominant in its
$W_0$-orbit. It follows from the
discussion in Appendix \ref{KL} that the value $\a(r_u)$
with $\a\in F_0$ is given by $\q^{D_u(\a)/2}$, where
$D_u(\a)$ is the weight of $\a$ in the Bala-Carter diagram
of $u$ (cf. \cite{C}). In our case, the vector
of values $\a(r_u)$ with $\a\in F_0$ is $(\q,\dots,\q,1,\q)$,
which is a residual point for $(\Ri,q_1)$,
where $q_1$ denotes the length multiplicative function
$q_1(w)=\q^{l(w)}$.

The Springer correspondence for all classical types has been
computed explicitly in \cite{Sho}. We use the description
of \cite{Lu5} (see \cite{C}). In our particular case,
the partition $\l\vdash 2n$ of elementary divisors of $u$
is $\l=(2,2n-2)$. Thus the Springer representation corresponding
with $(u,1)$ ($1$ denoting the trivial representation of the component
group $A(r_u,u)$ (see Appendix \ref{KL})) is the representation
$\phi_{(n-1,1)}$ of $W_0$ labeled by the double partition
$(n-1,1)$ of $n$.
This is the reflection representation of $W_0$.

The component group is equal to $A(r_u,u)\simeq C_2$
(Chapter 13, loc. cit.). Both representations
$\pm 1$ of $A(r_u,u)$ are geometric, and one easily finds that
the Springer correspondent of $(u,-1)$ is the
representation $\phi_{(-,n)}$ of $W_0$.
This is the $1$-dimensional
representation in which $s_i$ acts by $1$ for $i=1,\dots,n-1$,
and in which $s_n$ acts by $-1$.

Let us denote by $\pi_{\pm 1}$ the irreducible
square integrable $\I$-spherical representations of $\mathcal{G}$
with the Kazhdan-Lusztig parameters $(r_u,u,\pm 1)$,
and put $\rho_{\pm 1,1}:=\pi_{\pm 1}^\I$.
The Kazhdan-Lusztig model \cite{KL}, and the explicit results
of \cite{Lu0} imply the following:
$\rho_{1,1}$ is an $(n+1)$-dimensional
discrete series representation of $\H(\Ri,q_1)$, with
central character $W_0r_u$, and with restriction to
$\H(W_0,q_1)$
whose limit for $\q\to 1$ is equal to
$\phi_{(-,1^n)}\otimes(\phi_{(n-1,1)}\oplus\phi_{(n,-)})$
(here $\phi_{(-,1^n)}$ is the sign representation, and
$\phi_{(n,-)}$ is the trivial representation of $W_0$).
The representation $\rho_{-1,1}$ is $1$-dimensional,
and has $\H(W_0,q_0)$-type corresponding to
$\phi_{(-,1^n)}\otimes\phi_{(-,n)}=\phi_{(1^n,-)}$
in the limit $\q\to 1$.

According to Corollary \ref{cor:fdim} (also see
Subsection \ref{sub:uni}) we have
\begin{equation}
\operatorname{fdim}(\pi_{\pm 1})=|W_0r_u|{\overline \ka}_{W_0r_u}
d_{\rho_{\pm 1}}m_{\{r_u\}}(r_u),
\end{equation}
with $m_{\{r_u\}}(r_u)$ equal to the rational function
(\ref{eq:m_L}), ${\overline \ka}_{W_0r_u}\in\mathbb{Q}^\times$ and
$d_{\rho_{\pm 1}}\in\mathbb{R}_+$,
subject to the condition
$(n+1)d_{\rho_{+1}}+d_{\rho_{-1}}=1$.

In general I do not know how to compute the constants
${\overline{\ka}}_{W_0r_u}$ and $d_{\rho_{\pm 1}}$
(there is a tedious ``algorithm'' for ${\overline \ka}_{W_0r_u}$
(analogous to \cite{HOH0}), and for $d_{\rho_{\pm 1}}$ not
even that). However, in the case of {\it regular} central characters
these constants are easy to determine. In the situation at hand we
we are able to determine the constants by slightly deforming $q$,
since the orbits of residual points that ``emerge'' from $r_u$
(there are two of them, corresponding to the two representations
$\rho_{\pm 1,1}$) under such a deformation are regular.
Moreover, one can show in the current example that the formal
dimensions are continuous under this deformation.

So let us consider {\it generic} root labels $q_f$
(cf. \cite{HOH0}, \cite{Slooten}) defined by
$q_f(s_i)=\q$ ($i=1,\dots,n-1$) and
$q_f(s_n)=\q^{f}$, where $0<f<2$, $f\not=1$.
There are two generic orbits $W_0r_{\pm 1,f}$
of residual points such that $W_0r_{\pm 1,1}=W_0r_u$.
By the generic parametrization of \cite{HOH0} of
orbits of residual points of the graded affine Hecke
algebra (which, by Theorem \ref{thm:lieres}, can also be
used for $(\Ri,q_f)$ residual points) of type $C_n$
such a generic orbit corresponds to a partition of $n$.
In this case the partitions are $\xi_1=(n-1,1)$
and $\xi_{-1}=1^n$.
The (standard basis) coordinates
of these residual points (suitably chosen within
their $W_0$-orbits) are (by \cite{HOH0})
$r_{1,f}=(\q^{2-n-f/2},\q^{3-n-f/2},\dots,\q^{1-f/2})$
and
$r_{-1,f}=(\q^{1-n+f/2},\q^{2-n+f/2},\dots,\q^{f/2})$.
In particular, these are regular orbits of residual points.

By Theorem \ref{thm:support}, for each of these central characters
there exists at least $1$ irreducible square integrable
representation of $\H$.
In addition, it is not difficult to see (cf. \cite{Slooten})
that the residual
Hilbert algebra of a {\it regular} orbit of residual points is
in fact simple. Thus for $f\not=1$, we find precisely two
irreducible square
integrable representations $\rho_{\pm 1,f}$, with central
characters $W_0r_{\pm 1,f}$.

One checks directly that $r_{-1,f}$ is the $\A$-weight space
of a $1$ dimensional (square integrable, by Casselman's criterion)
representation where $T_i$ ($i<n$) acts by $-1$, and $T_n$ by $\q^f$.
This is a continuous family of square integrable
representations in the parameter $f$ (if $f$ is in the range
$0<f<2$ and $n\geq 3$). We call this parameter family $\rho_{-1,f}$.

The other orbit $W_0r_{1,f}$ also carries a continuous parameter family
of square integrable representations $\rho_{1,f}$, the twist by
the automorphism $i$
(see Subsection \ref{sub:range})
of the affine reflection representation
of $\H$ (a representation of dimension $(n+1)$).
To see this, we give the following model for the representation
(there are several possible constructions one could invoke here,
but none of these is obvious (as far as I know)). Our approach
here is based on the simplifying circumstance that the
representation contains the sign representation of $\H(W_0)$
(is ``anti-spherical'').

We will use the spherical function $\phi(\mu,k)$ of the Yang system
(cf. \cite{HOH0}),
with $R=C_n$, $k_\a=\log(q(s_\a))$, and $\k=\log(\q)$.
Recall that this function depends {\it analytically} on $(\mu,k)$.
First we consider $-\k<0$ (the attractive case in \cite{HOH0}),
and we consider the residual point
$\mu=\log(w_0r_{1,f})=-\k(f/2+n-2,f/2+n-3,\dots,f/2,f/2-1)$.
The list of positive roots $\a$ such that $\a(\mu)=-k_\a$
is $\L=(e_1-e_2,\dots,e_{n-1}-e_n,2e_{n-1})$. In order to compute the
dimension (in the regular case $f\not=1$) of the graded Hecke
module generated by $\phi(\mu,k)$ we have to count the number of
exponentials $e^{w\mu}$ which have a nonzero coefficient in
$\phi(\mu,k)$. Assuming that $0<f<2, f\not=1$,
we see that $\mu$ satisfies the condition of Lemma 3.3 of
\cite{HOH0},
and by Remark 3.4 of \cite{HOH0}
this shows that $\mu$ is an exponent of
$\phi(\mu,k)$. Then $w^{-1}\mu$ is also an exponent iff
$w\d=w(n,n-1,\dots,1)$ satisfies $w\d(\a)>0$ for all $\a\in\L$.
One easily verifies that this is satisfied iff
$w\d=(n,n-1,\dots,2,\pm 1)$ or
$w\d=(n,n-1,\dots,\hat{j},\dots,1,-j)$ ($j=n,n-1,\dots,2$).
Hence the module generated by $\phi(\mu,-k)$ (with $\mu$ as above,
and $0<f<2$) is a spherical discrete series module of
the graded Hecke algebra, which is irreducible (this is always true,
see the discussion above Section 3, loc. cit.), and of dimension $n+1$ if
$f\not=1$. Now apply the involution $i$ (see Section 5, loc. cit.,
and also Subsection \ref{sub:range}) to replace $-\k$ by $\k$, and
then integrate the representation (as in \cite{Lu}, Section 9)
so obtained to get a representation of $\H$. We obtain a
parameter family
(depending on $f$ with $0<f<2$) of irreducible square
integrable representations generated by an anti-spherical vector,
of dimension $n+1$ if $f\not=1$. Now observe that for $f=1$ this
representation has to be irreducible of dimension $n+1$ as well,
by the classification of the square integrable representations
with central character $W_0r_u$ as described above. We call this
$n+1$-dimensional family $\rho_{1,f}$. It follows easily from the
above discussion that the characters of $\rho_{\pm 1,f}$ are
continuous in $0<f<2$, and uniformly square integrable.

Hence we can compute the formal dimension of both representations
by taking the limit for $f\to 1$ of the corresponding
generic formal dimensions. It is easy to see that
\begin{equation}
\lim_{f\to 1}m_{\{r_{\pm 1,f}\}}(r_{\pm 1,f})
=\pm\frac{1}{2}m_{\{r_u\}}(r_u).
\end{equation}
For $f\not=1$ one obviously has $d_{\rho_{\pm}}=1$
and $\ka_{wr_{\pm 1,f}}=\pm(-1)^n|X:Q|^{-1}=\pm(-1)^n/2$
for all $w$ such that
$wr_{\pm 1,f}$ is a weight in $\rho_{\pm 1,f}$, and $=0$ else.

Combining these facts, we find in the limit $f\to1$ that
$|W_0r_u|{\overline \ka}_{W_0r_u}=(-1)^n(n+2)/4$,
and $d_{\rho_{\pm 1}}=1/(n+2)$. Hence both constants
$\l_{\rho_{\pm 1}}$ are equal to $(-1)^n/4$, which is in accordance
with Reeder's conjectural formula (\cite{Re}, equation (0.5))
for the formal dimension (up to a sign).
A computation yields:
\begin{equation}\label{eq:fdimsubreg}
\operatorname{fdim}(\pi_{\pm 1})=\frac{1}{4}
\frac{\q(\q-1)^{n+2}(\q^{n-2}-1)\prod_{i=1}^{n-2}(\q^{2i+1}-1)}
{(\q^2-1)(\q^n-1)\prod_{i=1}^{n-1}(\q^{2i}-1)}
\end{equation}
\begin{rem} It would be interesting to work out the product
formula (\ref{eq:fdim}) for formal dimensions
(without the precise analysis of the constants $\l_\rho$)
for classical root systems
in general (for ``special parameters'', see \cite{Slooten}),
and to express the answer (in the case of real
central characters) in terms of the symbol of the Springer
correspondent according to the conjecture in \cite{Slooten}.
\end{rem}
\subsection{$K$-types}\label{sub:types}
We touch superficially upon the analogue of the problem of the
``$K$-type decomposition'' of admissible representations of a reductive
group for tempered representations of the affine Hecke algebra $\H$.
We refer to \cite{Re} for a deep connection between the ``$K$-types''
of an irreducible discrete series representation, and its formal
dimension.
We refer to \cite{Slooten} for precise conjectures on the
$K$-types of the irreducible tempered modules
with real central character for affine Hecke algebras of
classical type (and general root labels).

The role of $K$ can be played by any maximal finite type
Hecke subalgebra of the form $\H(W_J)\subset\H$, with
$J\subset F^{\mathrm{aff}}$  a maximal proper subset.
Such a subalgebra is a finite dimensional $*$-subalgebra.
The restriction of $\tau$
to $\H(W_J)$ is equal to the usual trace of the finite
type Hecke algebra $\H(W_J)$, normalized in such a way that
$\tau(T_e)=1$.

For $\sigma\in\hat{W_J}$ we denote by $d_{J,\sigma}(q)$ its
generic degree
with respect to $\H(W_J)$ with label $q|_{W_J}$.
Thus we have (\ref{eq:finite})
\begin{equation}
\tau|_{\H(W_J)}=(P_{W_J}(q))^{-1}
\sum_{\sigma\in\hat{W_J}}d_{J,\sigma}(q)\chi_\sigma,
\end{equation}
where
$P_{W_J}(q)$ denotes the Poincar\'e polynomial of
$W_J$ with respect
to the label function $q$ (restricted to $W_J$).

Now observe that the restriction to $\H(W_J)$ of $\pi(\om)$
is independent of $\om\in \O$. We denote the multiplicities
by $n_\O(\sigma)$, thus
\begin{equation}
\chi_{\pi(\om)}|_{\H(W_J)}=\sum_{\sigma\in\hat{W_J}}
n_\O(\sigma)\chi_\sigma.
\end{equation}
We introduce for $\g=(\Ri_P,W_Pr_P)\in\Gamma$
the following rational functions of $q$:
\begin{equation}
M_\g:=\int_{t\in T^P_u} m^P(r_Pt)d^Pt.
\end{equation}
Notice that for all orbits of the form
$\O=K_\d\backslash\Lambda_\g\times\d$,
\begin{equation}
\int_\O m^P(\om)d^\O\om=M_\g.
\end{equation}
From the Plancherel decomposition of $\H$ (Theorem
\ref{thm:mainp}) we thus obtain the following identities:
For all $J\subset F^{\mathrm{aff}}$ a maximal proper subset, and
each $\sigma\in\hat{W_J}$,
\begin{equation}\label{eq:gendgs}
d_{J,\sigma}(q)=P_{W_J}(q)
\sum_{\g\in\Gamma_a}M_\g
\sum_{\d\in\Delta_{\g,a}}
|W(\O)|^{-1}n_\O(\sigma)\mu_{\Ri_P,Pl}(K_{P(\g)}\d),
\end{equation}
where $\Gamma_a$ is a complete set of representatives
for the association classes (=$\W$-orbits) in $\Gamma$,
$\Delta_{\g,a}$ is a complete set of representatives in
$\Delta_\g$ for the action of $\W(\g)$, and where
$\O$ denotes the orbit $\O=K_\d\backslash\Lambda_\g\times\d$ of
cuspidal representations of $\H^{P(\g)}$ (for a given pair
$(\g,\d)\in\Gamma_a\times\Delta_{\g,a}$).
\begin{ex}
It is instructive to verify (\ref{eq:gendgs}) for
$R_0=B_2$ (equal label case), both for $X=Q$ and $X=P$,
using the discussion in Example \ref{ex:b2}.
The residual point $(\q,-1)$ for $X=Q$
(notation of Example \ref{ex:b2}) is the most complicated part.
This orbit of residual points carries $2$
one-dimensional discrete series representations which are
exchanged by the nontrivial affine diagram automorphism.
Their direct sum lifts to the two-dimensional irreducible
discrete series
representation which is carried by the (regular) orbit
$(\q^{1/2},-1)$
of residual points for the extended affine Hecke algebra with
$X=P$. Using Corollary \ref{cor:normeq} one concludes that
the formal dimension of this two-dimensional representation
(which is easily computed, since the underlying central character
is regular) is equal to the formal dimension of each of the
two one-dimensional discrete series in which it decomposes upon
restriction to the case $X=Q$.
\end{ex}
\subsection{A remark on the residual degrees $d_{\pi}$}
\label{sub:remarks}
We mention one further consequence of Corollary \ref{cor:normeq}
regarding the constants $d_{\pi}\in\mathbb{R}_+$.
\begin{cor}
Assume that the constants $f_s$ in Convention \ref{eq:scale} are
integers.
Let $r$ be a residual point, and let $\pi\in\Delta_{W_0r}$. Assume
that the character values of $\pi$ on $T_w$ are contained in
$k[\q^{1/2},\q^{-1/2}]$, where $k$ is a subfield of $\C$. Then
$d_{\pi}\in\mathbb{R}_+\cap k$.
\end{cor}
See also Conjecture \ref{rem:ell}; we expect that the
$d_{\pi}\in \Q$.
\begin{proof}
The main step is to show that Casselman's bound of Corollary
\ref{rem:cas} becomes uniform in $\q$ under the assumption. Let
$r=sc$ and choose $n\in\N$ such that $s^n=s$. Let us first fix
$\q>1$. Consider the isomorphism of localized Hecke algebras
\begin{equation}
j_{1/n}:\H^{an}_{q^n}(U)\to\H^{an}(U_{1/n}),
\end{equation}
where $U=W_0sc^n\exp{B}=W_0r^n\exp{B}$, with $B$ a suitably small
ball around the origin in $\mathfrak{t}_\mathbb{C}$. We have, by
the assumption that $s^n=s$:
\begin{equation}
j_{1/n}(\theta_{x,q^n})=\theta_{nx}.
\end{equation}
On the other hand, for all $s\in F_0$, the eigenvalues of the self
adjoint operator $\pi(j_{1/n}(N_{s,\q^n}))$ are of the form
$\pm\q^{\pm f_sn/2}$. Hence the operator norms of the operators
$\pi(j_{1/n}(N_{s,\q^n}^{\pm 1}))$ are bounded by $\q^{Mn}$, for a
suitable constant $M$.

Given $w\in W$ we can write $w=uxv$ with $x\in X^+$, $u\in W^x$
and $v\in W_0$, where $W^x$ denotes the set of shortest length
representatives of the left cosets of the stabilizer $W_x$ of $x$
in $W_0$. If we write $u=s_{i_1}\dots s_{i_k}$ and $v=s_{j_1}\dots
s_{j_l}$, we can thus choose signs $e_i$ and $d_j$ such that
\begin{equation}
N_w=N_{i_1}^{e_1}\dots N_{i_k}^{e_k}\theta_x N_{j_1}^{d_1}\dots
N_{j_l}^{d_l}.
\end{equation}
Let us simply denote this decomposition by $N_w=N_u^e\theta_x
N_v^d$.

Now by Theorem \ref{thm:normunif}(ii), and the remark that
$j_{1/n}$ intertwines the action of $\H^{an}_{\q^n}(U)$ on
$\pi|_{\q^n}$ with that of $\H^{an}(U_{1/n})$ on $\pi|_{\q}$, we
have ($M$ is a constant, not necessarily the same as above):
\begin{equation}
\begin{split}
|\chi_{\pi,\q^n}(N_{w,\q^n})|^2&=|\chi_{\pi,\q}(j_{1/n}(N_{w,\q^{n}}))|^2\\
&=|\chi_{\pi,\q}(j_{1/n}(N_{u,\q^n}^e
N_{v,\q^n}^d)\theta_{nx,\q})|^2\\
&\leq\operatorname{dim}(\pi)^2\Vert\pi(j_{1/n}(N_{u,\q^n}^e
N_{v,\q^n}^d))\Vert^2_o \Vert \theta_{nx,\q}\Vert_{ds}^2\\ &\leq
C\q^{2n(M-\e l(x))}\\
\end{split}
\end{equation}
where $C$ is independent of $w$ and $n$. In particular, this
implies that the highest power of $\q$ in $|\chi_\pi(N_w)|^2=
\chi_\pi(N_w)\chi_\pi(N_{w^{-1}})\in
k_\mathbb{R}[\q^{1/2},\q^{-1/2}]$ tends to $-\infty$ with $l(w)$
(with $k_\mathbb{R}:=k\cap\mathbb{R}$). Hence the left hand side
of the equality Corollary \ref{cor:normeq} is a Laurent series in
$\q^{-1/2}$ with coefficients in $k_\mathbb{R}$.

On the other hand, according to Proposition \ref{prop:par}(iv),
$m_{\{r\}}(r)^{-1}$ can be expanded as a Laurent series in
$\q^{-1/2}$ with coefficients in $\mathbb{Q}$. The desired result
follows.
\end{proof}
\section{Appendix: Residual Cosets}\label{sub:defn}
\subsection{Introduction and quick guide}\label{sub:quick}
Our approach to the spectral resolution is through residues of
certain rational $n$-forms on a complex torus $T$.
In order for our method to work well, we need to have a
certain a priori knowledge on the geometric and combinatorial
properties of the set of poles of these rational forms.
The present section serves to collect such facts about the
set of poles, and to classify the collection of ``residual
cosets'', the sets of maximal pole order, which will eventually
turn out to constitute the projection of support of the
Plancherel measure to ${W_0}\backslash T$.

Recall that we have chosen a rational, positive definite,
$W_0$-invariant symmetric form on $X$.
This defines an isomorphism
between $X\otimes_\Z \Q$ and $Y\otimes_\Z \Q$, and thus also a
rational, positive definite symmetric form on $Y$.
We extend this
form to a positive definite Hermitian form on
$\mathfrak{t}_\C:=\operatorname{Lie}(T)=Y\otimes_\Z\C$, where $T$
is the complex torus $T=\operatorname{Hom}(X,\C^\times)$. Via the
exponential covering map $\exp:\mathfrak{t}_\C\to T$ this
determines a distance function on $T$.

Let $q$ be a set of root labels. If $2\alpha\not\in
R_{\mathrm{nr}}$ we formally put $q_{\alpha^\vee/2}=1$, and
always $q_{\alpha^\vee/2}^{1/2}$ denotes the positive square root
of $q_{\alpha^\vee/2}$. Let $L$ be a coset of a subtorus
$T^L\subset T$ of $T$. Put $R_L:=\{\alpha\in R_0\mid
\alpha(T^L)=1\}$. This is a parabolic subsystem of $R_0$. The
corresponding parabolic subgroup of $W_0$ is denoted by $W_L$.
Define
\begin{equation}
R_L^p:=\{\alpha\in R_L\mid \alpha(L)=-q_{\alpha^\vee/2}^{1/2}
\ \mathrm{or} \  \alpha(L)=q_{\alpha^\vee/2}^{1/2}q_{\alpha^\vee}\}
\end{equation}
and
\begin{equation}
R_L^z:=\{\alpha\in R_L\mid \alpha(L)=\pm 1\}.
\end{equation}
We write $R_L^{p,ess}=R_L^p\backslash R_L^z$ and
$R_L^{z,ess}=R_L^z\backslash R_L^p$. We define an index $i_L$ by
\begin{equation}
i_L:=|R_L^p|-|R_L^z|.
\end{equation}
As a motivation for the somewhat more technical definition in the next
subsection,
we remark that this index $i_L$ computes the order of the pole
along $L$ of the rational $(n,0)$-form
\begin{equation}
\omega:=\frac{dt}{c(t,q)c(t^{-1},q)},
\end{equation}
\index{0y@$\omega=\frac{dt}{c(t,q)c(t^{-1},q)}$, $(n,0)$-form on $T$}
which plays a main role in this paper (cf. equation
(\ref{eq:basic})). We will find (cf. Corollary
\ref{cor:simpledefres}) that for each coset $L$ of a subtorus of
$T$,
\begin{equation}\label{eq:eq}
i_L\leq\operatorname{codim}(L).
\end{equation}
Suppose that
$L$ is a coset such that $i_L<\operatorname{codim}(L)$, and let
$T_L$ denote the subtorus orthogonal to $T^L$. Let $C_L$ be a
cycle of dimension
$\operatorname{dim}(C_L)=\operatorname{codim}(L)$ in a
sufficiently small neighborhood of $e$ in $T_L$, and let $C^L$ be
any compact cycle in $L\backslash\cup_{L\not\subset
L^\prime}L^\prime$ of dimension
$\operatorname{dim}(C^L)=\operatorname{dim}(L)$. Then for every
homolorphic function $f$ on $T$,
\begin{equation}\label{eq:f}
\int_{C^L\times C_L}f\omega=0.
\end{equation}
We call a coset $L$ {\it residual} if $i_L=\operatorname{codim}(L)$.
It will turn out that the support of the spectral measure of the
restriction of the trace $\tau$ to the center of the Hecke algebra is precisely
equal to the union of all the ``tempered forms'' of the residual cosets
(see Theorem \ref{thm:support}). The spectral measure arises as a sum of
integrals of the form (\ref{eq:f}).

For technical convenience, the Definition \ref{dfn:ressub} of the notion
``residual coset'' in the next subsection is slightly more complicated.
We will define the
residual cosets by induction on their codimension in $T$, in such
a way that the collection of residual cosets is easily amenable
to classification. In the next subsection we discuss their elementary
properties and show how the classification can be reduced to the
case of residual subspaces in the sense of \cite{HOH0}. These residual
subspaces were already
classified in the paper \cite{HOH0}. By this classification we
verify equation (\ref{eq:eq}) (cf. Corollary \ref{cor:simpledefres}).
Using Lemma \ref{lem:ats} this implies that the following are equivalent
for a coset $L\subset T$:
\begin{enumerate}
\item[(i)] $L$ is residual (in the sense of Definition \ref{dfn:ressub}).
\item[(ii)] $i_L\geq\operatorname{codim}(L)$.
\item[(iii)] $i_L=\operatorname{codim}(L)$.
\end{enumerate}
\subsection{Definition and Classification of Residual Cosets}
We give the following recursive definition of the notion {\it
residual coset}\label{sub:dfnser}.
\begin{dfn}\label{dfn:ressub}
A coset $L$ of a subtorus of $T$ is called residual if either
$L=T$, or else if there exists a residual coset $M\supset L$
such that $\mathrm{dim}(M)=\mathrm{dim}(L)+1$ and
\begin{equation}\label{eq:tech}
i_L\geq i_M+1.
\end{equation}
\end{dfn}
\begin{cor}\label{cor:fininv}
The collection of residual cosets is a nonempty, finite collection of
cosets of algebraic subtori of $T$, closed for the
action of the group of automorphisms of the root system preserving
$q$ (in particular the elements of $W_0$, but also for example
$-\operatorname{Id}$).
\end{cor}
\begin{proof}
By induction on the codimension.
In a residual coset $M$ of
codimension $k-1$ we find only finitely many cosets $L\subset M$
of codimension $1$ in $M$ with $i_L>i_M$. The invariance is
obvious from the invariance of the index function $i_L$.
\end{proof}
\begin{prop}\label{prop:red} If $L$ is residual, then
\begin{itemize}
\item[(i)] $R_L^{p,ess}$ spans a subspace
$V_L$ of dimension $\operatorname{codim}(L)$ in the $\Q$
vectorspace $V=X\otimes\Q$.
\item[(ii)] We have $R_L=V_L\cap R_0$, and the rank of $R_L$ equals
$\operatorname{codim}(L)$.
\item[(iii)] Put ${}_LX:=V_L\cap X$ and $X^L:=X/{}_LX$.
Then $T^L=\{t\in T\mid x(t)=1
\ \forall x\in{}_LX \}=
\operatorname{Hom}(X^L,\C^\times)=(T^{W_L})^0$.
\item[(iv)] Put $Y_L:=Y\cap\Q R^\vee_L$ and
${}^LX:=Y_L^\perp\cap X$. Let $X_L:=X/{}^LX$. We identify $R_L$
with its image in $X_L$. Let $F_L$ be the basis of $R_L$ such that
$F_L\subset R_{0,+}$. Then $\Ri_L:=(X_L,Y_L,R_L,R_L^\vee,F_L)$ is
a root datum.
\item[(v)] Put
$T_L:=\operatorname{Hom}(X_L,\C^\times)\subset T$ (we identify
$T_L$ with its canonical image in $T$). Then $T_L$ is the subtorus
in $T$ orthogonal to $L$. Define $K_L:=T^L\cap T_L=
Hom(X/({}_LX+{}^LX),\C^\times)\subset T$, a finite subgroup of
$T$. The intersection $L\cap T_L$ is a $K_L$-coset consisting of
{\it residual points} in $T_L$ with respect to the root datum
$\Ri_L$ and the root label $q_L$ obtained from $q$ by restriction
to $R_{L,\mathrm{nr}}^\vee\subset R_\mathrm{nr}^\vee$. When
$r_L\in T_L\cap L$, we have $L=r_LT^L$. Such $r_L$ is determined
up to multiplication by elements of $K_L$.
\end{itemize}
\end{prop}
\begin{proof}
By induction on $\operatorname{codim}(L)$ we may assume
that the assertions of (i) and (ii) hold true for $M$
in (\ref{eq:tech}).
From the definition we see that
$R_L^{p,ess}\backslash R_M^{p,ess}$ is not the empty set.
An element $\alpha$ of $R_L^{p,ess}\backslash R_M^{p,ess}$
can not be constant on $M$,
and hence $\alpha\not\in R_M=V_M\cap R_0$.
Thus
$$\operatorname{dim}(V_L)\geq\operatorname{dim}(V_M)+1
=\operatorname{codim}(M)+1=\operatorname{codim}(L).$$
Since also
$$V_L\subset \operatorname{Lie}(T^L)^\perp,$$
equality has to hold.
Hence $R_L\subset V_L$ and $R_L$ spans $V_L$. Since $R_L$ is
parabolic, we conclude that $R_L= V_L\cap R_0$. This proves
(i) and (ii). The subgroup $\{t\in T\mid x(t)=1
\ \forall x\in{}_LX \}\subset T$ is isomorphic to
$\operatorname{Hom}(X^L,\C^\times)$, which is a torus
because $X^L$ is free.
By (ii) then, its dimension equals
$\operatorname{dim}(T^L)$.  It contains $T^L$, hence is equal
to $T^L$. It follows that $T^L$ is the connected component
of the group of fixed points for $W_L$, proving (iii).
The statements (iv) and (v) are trivial.
\end{proof}

For later reference we introduce the following notation. A
residual coset $L$ determines a parabolic subsystem $R_L\subset
R_0$, and we associated with this a root datum $\Ri_L$. When
$\Sigma\subset R_0$ is any root subsystem, {\it not necessarily
parabolic}, we associate to $\Sigma$ two new root data, namely
$\Ri^\Sigma:=(X,Y,\Sigma,\Sigma^\vee,F_\Sigma)$ with $F_\Sigma$
determined by the requirement $F_\Sigma\subset R_{0,+}$, and
$\Ri_\Sigma:=(X_\Sigma,Y_\Sigma,\Sigma,\Sigma^\vee,F_\Sigma)$
where the lattice $X\to X_\Sigma$ is the quotient of $X$ by the
sublattice orthogonal to $\Sigma^\vee$, and $Y_\Sigma\subset Y$
is the sublattice of elements of $Y$ which are in the $\R$-linear
span of $\Sigma^\vee$.

There is an obvious converse to Proposition \ref{prop:red}:
\begin{prop}\label{prop:conv}
Let $R^\prime\subset R_0$ be a parabolic subsystem of roots, and
let $T^L\subset T$ be the subtorus such that $R^\prime=R_L$. Let
$T_L\subset T$ be the subtorus whose Lie algebra
$\operatorname{Lie}(T_L)$ is spanned by $R_L^\vee$. Let $r\in T_L$
be a residual point with respect to $(\Ri_L,q_L)$ as in
Proposition \ref{prop:red}(v). Then $L:=rT^L$ is a residual coset
for $(\Ri,q)$.
\end{prop}

The recursive nature of the definition of residual cosets makes it feasible
to give a complete classification of them. By Lemma \ref{prop:red}, this
classification problem reduces to the classification of the {\it residual
points}. In turn, Lusztig \cite{Lu} indicates how the classification of
residual points reduces to the classification of residual points in the sense
of \cite{HOH0} for certain graded affine Hecke algebras. This classification
is known by the results in \cite{HOH0}. Let us explain this in detail.
Following \cite{Lu} we call a root datum $\Ri=(X,Y,R_0,R_0^\vee,F_0)$ {\it
primitive} if one of the following conditions is satisfied:
\begin{enumerate}
\item[(1)] $\forall\a\in R_0:\ \a^\vee\not\in 2Y$.
\item[(2)] There is a unique $\a\in F_0$ with $\a^\vee\in 2Y$
and $\{w(\a)\mid w\in W_0\}$ generates $X$.
\end{enumerate}
A primitive root datum $\Ri$ satisfying (2) is of the type
$C_n^{\operatorname{aff}}$ ($n\geq 1$), by which we mean that
$$\Ri=(Q(B_n)=\Z^n,P(C_n)=\Z^n,
B_n,C_n,\{e_1-e_2,\dots,e_{n-1}-e_n,e_n\}).$$
By \cite{Lu} we know that every root datum is a direct sum
of primitive summands.
\begin{prop}\label{prop:lusz}
Let $r\in T$ be a residual point, and write $r=sc\in T_uT_{rs}$
for its polar decomposition (with
$T_u=\operatorname{Hom}(X,S^1)$ and
$T_{rs}=\operatorname{Hom}(X,\R_+)$).
The root system
\[
R_{s,1}:=\{\alpha\in R_1\mid \alpha(s)=1\}
\]
\index{R4@$R_{s,0}$ ($R_{s,1}$), roots of $R_0$ ($R_1$) vanishing
in $s\in T$}
has rank $\operatorname{dim}(T)$. The system
\[
R_{s,0}:=
\{\a\in R_0\mid k\a\in R_{s,1}\ \mathrm{for\ some}\ k\in\N\}
\]
contains both $R_{\{r\}}^{p,ess}$ and $R_{\{r\}}^{z,ess}$, and $r$
is residual with respect to the affine Hecke subalgebra
$\H^s\subset \H$ whose root datum is given by
$\Ri^{s}:=(X,Y,R_{s,0},R_{s,0}^\vee,F_{s,0})$ (with $F_{s,0}$ the
basis of $R_{s,0}$ contained in $R_{0,+}$).
\end{prop}
\begin{proof}
It is clear from the definitions that $R_{s,0}$ contains
$R_{\{r\}}^{p,ess}$ and $R_{\{r\}}^{z,ess}$, and hence has maximal
rank. Given a full flag of $\Ri$-residual subspaces
$\{c\}=L_0\subset L_1\subset \dots\subset L_n=T$, satisfying
(\ref{eq:tech}) at each level, we see that the sets $R_{L_i}^p$,
$R_{L_i}^z$ are contained in $R_{s,0}$. It follows by reverse
induction on $i$ (starting with $L_n=T$) that each element of the
flag is $\Ri^{s}$-residual.
\end{proof}
\begin{lem}\label{lem:order2}
Given a residual point $r=sc$, let $s_0\in
T_u=\operatorname{Hom}(X,S^1)$ be the element which coincides with
$s$ on each primitive summand of type $C_n^{\operatorname{aff}}$
and is trivial on the complement of these summands. Then $s_0$ has
at most order $2$.
\end{lem}
\begin{proof}
To see this we may assume that
$\Ri$ is of type $C_n^{\operatorname{aff}}$. Then $R_1$ is of
type $C_n$, $s=s_0$, and
$R_{s,1}=\{\alpha\in R_1\mid \alpha(s_0)=1\}$, being
of maximal rank in $R_1$, is of type $C_k+C_{n-k}$ for some $k$.
In particular, $\pm 2e_i\in R_{s,1}$ for all $i=1,\dots,n$. Moreover
the index of $\Z R_{s,1}$ in $\Z R_1$ is at most $2$. Thus $s_0$
takes values in $\{\pm 1\}$ on $R_1$,
and is trivial on elements of the form $\pm 2 e_i$. It follows
that $s_0$ is of order at most $2$ on $X=\Z^n$.
\end{proof}
Denote by $h\in\operatorname{Hom}(Q,S^1)$ the image of $s_0$ in
$\operatorname{Hom}(Q,S^1)$. Choose root labels $k_\a=k_{s,\a}\in
\R$ with $\a\in R_{s,0}$ by the requirement ($k_\a$ depends on the
image of $s$ in $\operatorname{Hom}(Q,S^1)$, but we suppress this
in the notation if there is no danger of confusion)
\begin{equation}\label{eq:log}
\begin{split}
\mathrm{e}^{k_\a}=&q_{\a^\vee}^{h(\a)/2}q_{\a^\vee+1}^{1/2}\\
=&\begin{cases}
q_{\a^\vee/2}^{1/2}q_{\a^\vee}&\mathrm{\ if\ }h(\a)=+1\\
q_{\a^\vee/2}^{1/2}&\mathrm{\ if\ }h(\a)=-1\\
\end{cases}
\end{split}
\end{equation}
\begin{thm}\label{thm:lieres}
Let $r=sc$ be a $(\Ri,q)$-residual point. Then
$\gamma:=\log(c)\in\mathfrak{t}:=\operatorname{Lie}(T_{rs})$ is a
residual point in the sense of \cite{HOH0} for the graded Hecke
algebra $H^{s}=\C[W(R_{s,0})]\otimes
\operatorname{Sym}(\mathfrak{t})$ with root system $R_{s,0}$ and
root labels $k_{s}:=(k_{s,\a})_{\a\in R_{s,0}}$. This means
explicitly that there exists a full flag of affine linear
subspaces $\{\gamma\}=\mathfrak{l}_n\subset\mathfrak{l}_{n-1}
\subset\dots \subset\mathfrak{l}_0=\mathfrak{t}$ such that the
sequence
\begin{equation}
i_{s,\mathfrak{l}_i}:=|R_{s,0,i}^p|-|R_{s,0,i}^z|
\end{equation}
is strictly increasing, where
\begin{equation}
R_{s,0,i}^p=\{\a\in R_{s,0}\mid \a(\mathfrak{l}_i)=k_{s,\a}\},
\end{equation}
and
\begin{equation}
R_{s,0,i}^z=\{\a\in R_{s,0}\mid \a(\mathfrak{l}_i)=0\}.
\end{equation}
Conversely, given a $s\in T_u$ such that $R_{s,1}\subset R_1$
has rank equal to $\operatorname{rank}(X)$, and a residual point
$\gamma\in\mathfrak{t}$ for the root system $R_{s,0}$ with
labels $(k_{s,\a})$ defined by (\ref{eq:log}), the point
$r:=s\exp{\gamma}$ is $(\Ri,q)$-residual. This sets up a
$1-1$ correspondence between $W_0$-orbits of $(\Ri,q)$-
residual points and the collection of pairs $(s,\gamma)$
where $s$ runs over the $W_0$-orbits of elements of
$T_u$ such that $R_{s,1}$ has rank equal
to $\operatorname{rank}(X)$, and
$\gamma\in\mathfrak{t}$ runs over the
$W(R_{s,0})$-orbits of
residual points (in the sense of
\cite{HOH0}) for $R_{s,0}$ with the labels $k_{s}$.
\end{thm}
\begin{proof}
Straightforward from the definitions.
\end{proof}
For convenience we include the following lemma:
\begin{lemma}\label{lem:conj}
If the rank of $R_0$ equals the rank of $X$ (a necessary condition
for existence of residual points!), the $W_0$-orbits of points
$s\in T_u$ such that $R_{s,1}\subset R_1$
has maximal rank correspond
$1-1$ to the $\operatorname{Hom}(P(R_1)/X,S^1)\simeq
Y/Q(R_1^\vee)$-orbits on the affine Dynkin diagram
$R^{(1)}_1$. In particular, $R_{s,1}$ only depends
on the corresponding $P(R_1^\vee)/Q(R_1^\vee)$-orbit of vertices of
$R^{(1)}_1$.
\end{lemma}
\begin{proof}
In the compact torus $\operatorname{Hom}(P(R_1),S^1)$,
the $W_0$-orbits of such points
correspond to the vertices of the fundamental alcove
for the action of the affine Weyl group
$W_0\ltimes 2\pi iQ(R_1^\vee)$ on $Y\otimes 2\pi i\R$.
Now we have to restrict to $X\subset P(R_1)$.
\end{proof}
With the results of this subsection at hand,
the classification of residual cosets is
now reduced to the classification of residual subspaces
as was given in \cite{HOH0}.
\begin{ex}\label{ex:b2}
Let $R_0=B_2=\{\pm e_1,\pm e_2,\pm e_1\pm e_2\}$ with basis
$\a_1=e_1-e_2$, $\a_2=e_2$, and let $X=Q=\Z^2$
(this is $C_2^{\mathrm{aff}}$).
Assume that $q(s_i)=\q$ for $i=0,1,2$. Then
$R_1=\{\pm 2e_1,\pm 2e_2,\pm e_1\pm e_2\}$ and thus $X=P(R_1)$.
We use $(\a_1,\a_2)$ as a basis of $X$ (so a point $t\in T$
is represented by $(\a_1(t),\a_2(t))$). The orbits
of points $s\in T_u$ such that $R_{s,1}$ has rank 2 are
represented by $(1,1),(1,-1)$ and $(-1,1)$. The latter
point corresponds to $R_{s,0}=\{\pm e_1,\pm e_2\}$, but since it
has value $-1$ on $\pm e_1$, the labels of $\pm e_1$
are equal to $1$ (by (\ref{eq:log})). Therefore there are no residual
points associated with $(-1,1)$. The other two points each give
one orbit of residual points, namely $(\q,\q)$ and $(\q,-1)$.

In addition we have $2$
orbits of one-dimensional residual cosets (with $K_{\{1\}}=C_2$
and $K_{\{2\}}=1$), and finally the principal two-dimensional
one, $T$.

Let us now consider $R_0=B_2$ with the lattice $X=P$
and again $q(s_i)=\q$ for $i=0,1,2$.
We take $(\a_1/2,\a_2)$ as a basis for $P$.
Now
$R_1=R_0$, and thus $X=P(R_1)$. So again we have $3$ orbits
of points $s$ for which the rank of $R_{s,0}$ is $2$, namely
$(1,1),(1,-1)$ and $(-1,1)$. Each corresponds to a
(regular) orbit of residual points:
$(\q^{1/2},\q),(\q^{1/2},-1)$, and $(-\q^{1/2},\q)$.

In addition there are $3$ one-dimensional residual cosets, $2$
associated with $P=\{1\}$ (with $K_{\{1\}}=1$) and $1$
with $P=\{2\}$ (with $K_{\{2\}}=C_2$). Finally we have the
principal residual coset $T$.
\end{ex}
\subsection{Properties of residual and tempered cosets}\label{sub:resiprop}
In the derivation of the Plancherel formula of the affine Hecke
algebra, the following properties of residual cosets will play a
crucial role.
\begin{theorem}\label{thm:equal}
For each residual coset $L\subset T$ we have
\begin{equation}
i_L=\operatorname{codim}(L).
\end{equation}
In other words, for every inclusion $L\subset M$ of residual
cosets with $\operatorname{dim}(L)=\operatorname{dim}(M)-1$, the
inequality (\ref{eq:tech}) is actually an equality.
\end{theorem}
\begin{proof}
Unfortunately, I have no classification free proof of this fact.
With the classification of residual subspaces at hand it can be
checked on a case-by-case basis. By the previous subsection
(Proposition \ref{prop:red} and Theorem \ref{thm:lieres}) the
verification reduces to the case of residual points for graded
affine Hecke algebras. In \cite{HOH0} (cf. Theorem 3.9) this
matter was verified.
\end{proof}
Theorem \ref{thm:equal} has important consequences, as we will see
later. At this point we show that the definition of residual
cosets can be simplified as a consequence of Theorem
\ref{thm:equal}. We begin with a simple lemma:
\begin{lem}\label{lem:ats}
Let $V$ be a complex vector space of dimension $n$, and suppose
that $\L$ is the intersection lattice of a set $\P$ of linear
hyperplanes in $V$. Assume that each hyperplane $H\in\P$ comes
with a multiplicity $m_H\in\Z$, and define the multiplicity $m_L$
for $L\in\L$ by $m_L:=\sum m_H$, where the sum is taken over the
hyperplanes $H\in \P$ such that $L\subset H$. Assume that
$\{0\}\in\L$ and that $m_{\{0\}}\geq n$. Then there exists a full
flag of subspaces $V=V_0\supset V_1\dots\supset V_n=\{0\}$ such
that $m_k:=m_{V_k}\geq k$.
\end{lem}
\begin{proof}
We construct the sequence inductively, starting with $V_0$.
Suppose we already constructed the flag up to $V_k$, with $k\leq
n-2$. Let $\P_k\subset\L$ denote the set of elements of $\L$ of
dimension $n-k-1$ contained in $V_k$, and let $N_k$ denote the
cardinality of $\P_k$. By assumption, $N_k\geq n-k\geq 2$. Since
every $H\in\P$ either contains $V_k$ or intersects $V_k$ in an
element of $\P_k$, we have
\begin{equation}
\sum_{L\in\P_k}(m_L-m_k)=m_n-m_k.
\end{equation}
Assume that $\forall L\in\P_k:\ m_L\leq k$. Then, because $m_k\geq
k$ and $N_k\geq 2$,
\begin{equation}
m_n\leq kN_k+(1-N_k)m_k\leq k\leq n-2,
\end{equation}
contradicting the assumption $m_n\geq n$. Hence there exists a
$L\in\P_k$ with $m_L\geq k+1$, which we can define to be
$V_{k+1}$.
\end{proof}
\begin{cor}\label{cor:simpledefres} For every coset $L\subset T$
one has $i_L\leq \operatorname{codim}(L)$, and $L$ is residual if
and only if $i_L=\operatorname{codim}(L)$.
\end{cor}
\begin{proof}
Define $\P$ to be the (multi-)set of codimension $1$ cosets of $T$
arising as connected components of the following codimension $1$
sets:
\begin{equation}
\begin{split}
&L^+_{\a,1}:=\{t\in T\mid\a(t)=q_{\a^\vee}q_{\a^\vee/2}^{1/2}\}\\
&L^+_{\a,2}:=\{t\in T\mid\a(t)=-q_{\a^\vee/2}^{1/2}\}\\
&L^-_{\a,1}:=\{t\in T\mid\a(t)=1\}\\ &L^-_{\a,2}:=\{t\in
T\mid\a(t)=-1\}\\
\end{split}
\end{equation}
Here $\a\in R_0$, and $q_{\a^\vee/2}=1$ when $2\a\not\in R_1$. We
give the components of $L^+_{\a,1}$, $L^+_{\a,1}$ the index $+1$,
and we give the components of $L^-_{\a,1}$, $L^-_{\a,1}$ index
$-1$.

Suppose that $L$ is any coset of a subtorus $T^L$ in $T$. Then
$i_L$ is equal to the sum of the indices the elements of $\P$
containing $L$.

Assume that $i_L\geq\operatorname{codim}(L)=k$. By Lemma
\ref{lem:ats} there exists a sequence $L\subset L_{k-e}\subset
L_{k-e-1}\dots\subset L_0=T$ of components of intersections of
elements of $\P$ such that $i_{L_{k-e}}=i_L\geq k$ and
$i_{L_j}\geq j=\operatorname{codim}(L_j)$ (we did not assume that
$L$ is a component of an intersection of elements in the multiset
$\P$, hence $e>0$ may occur). If $k(0)$ is the smallest index such
that $i_{L_{k(0)}}>k(0)$, then $L_{k(0)}$ is by definition
residual, and thus violates Theorem \ref{thm:equal}. Hence such
$k(0)$ does not exist and we conclude that $i_{L_k}=k$ for all
$k$. This proves that $e=0$ and that $L$ is residual.
\end{proof}
\begin{rem}
This solves the question raised in Remark 3.11 of \cite{HOH0}.
\end{rem}
\begin{theorem}\label{thm:ster}
\begin{enumerate}
\item[(i)] Let $R_0$ be indecomposable, and let $r=c$ be a real residual point in
$\overline{T_{rs,+}}$.
If $\omega:T_{rs}\to T_{rs}$ is a homomorphism
which acts on the root system $R_0$ by means of a diagram automorphism of $F_0$,
then $\omega(r)=r$.
\item[(ii)]
Define $*:T\to T$ by $x(t^*)=\overline{x(t)^{-1}}$
\index{*@$*$!$t\to t^*$, anti-holomorphic involution on $T$}.
If $r=cs\in T$
is a residual point, then $r^*\in W(R_{s,0})r$.
\item[(iii)] If $r=sc$ is a residual point, then the values $\a(c)$ of the roots
$\a\in R_0$ on $c$ are in the subgroup of $\mathbb{R}_+$ generated by the positive
square roots of the root labels $q_{\a^\vee}$, with $\a\in R_{nr}$.
\end{enumerate}
\end{theorem}
\begin{proof}
(i). If $R_0$ allows a nontrivial diagram automorphism then $R_0$ is
simply laced. So we are in the situation of the Bala-Carter classification
of distinguished weighted Dynkin diagrams.
A glance at the tables of section 5.9
of \cite{C} shows that this fact holds true.

(ii). This is a consequence of (i), since $*:sT_{rs}\to sT_{rs}$ acts on
$R_{s,0}$ by means of an automorphism (see also \cite{HOH0}, Theorem
3.10) which acts trivially on the set of indecomposable summands of $R_{s,0}$.

(iii). For this fact I have also no other proof to offer than a case-by-case
checking, using the results of this section and the
list of real residual points from \cite{HOH0}.  The amount
of work reduces a lot by the remark that it is well known in the simply laced
cases (see Corollary \ref{cor:cor} of the appendix Section \ref{KL}).

In the classical cases other than $C_n^{\text{aff}}$,
it follows from a well known theorem of Borel and de
Siebenthal \cite{BS} that the index of $Q(R_{s,0})\subset Q$ is at most $2$.
Hence the desired result
follows if we verify that for {\it real} residual points of the classical
root systems, the values $\a(c)$ are
in the subgroup of $\mathbb{R}_+$ generated by the root labels,
which is direct
from the classification lists in \cite{HOH0}.

For the real points of $C_n^{\text{aff}}$ it is also immediate from the
above and
(\ref{eq:log}). For nonreal points $r=sc$ we look at (the proof of)
Lemma \ref{lem:order2}. If $s$ has order $2$, then $R_{s,1}$ is of type
$C_{k}+C_{n-k}$.
We need to check the values of the roots $e_n$ and $e_k-e_{k+1}$ on $c$
in this case.
But the roots $2e_i$ are in $R_{s,1}$, and take rational values in the
labels $q_{\a^\vee}$ ($\a\in R_{nr}$) on $c$.

The real residual points of $F_4$ are all rational in the root labels (see
\cite{HOH0}). Again using the Theorem of Borel and de Siebenthal, we need
to check in addition the nonreal residual points $r=sc$ with
$R_{s,0}=A_2\times A_2\subset F_4$ (generating a lattice of index $3$ in
$Q(F_4)$)
and $R_{s,0}=A_3\times A_1\subset F_4$ (index $4$). These cases can be
checked without difficulty.

In the case of $G_2$, there are
generically $3$ real residual points, two of which have rational coordinates
and one has rational coordinates only in the square roots of the labels.
In addition there are two nonreal residual points $r=sc$ for $G_2$, which
are easily checked.
(We need to check only the case with $R_{s,0}=A_2$
(index $3$ in $Q(G_2)$)).
\end{proof}
\begin{rem}
In fact the result (ii) of the previous Theorem will also turn out to
be a consequence of Theorem \ref{thm:support}, in view of Proposition
\ref{prop:im}.
\end{rem}
\begin{dfn}
Let $L$ be a residual coset, and write
$L=r_LT^L$ with $r_L\in T_L\cap L$.
This is determined up to multiplication of $r_L$ by elements of
the finite group $K_L=T_L\cap T^L$. Write $r_L=s_Lc_L$ with
$s_L\in T_{L,u}$ and $c_L\in T_{L,rs}$.
We call $c_L$ the ``center'' of $L$, and we call
$L^{temp}:=r_LT^L_u$ the tempered compact form of
$L$ (both notions are independent of the choice of $r_L$, since
$K_L\subset T^L_u$). The cosets of the form $L^{temp}$ in $T$ will
be called ``tempered residual cosets''.
\end{dfn}
\index{c1@$c_L$, center of $L$}
\index{L@$L^{temp}$, tempered residual coset}
\begin{theorem}\label{thm:nonnest}
Suppose that $L\subset M$ are two residual cosets. Write
$L=r_LT^L=s_Lc_LT^L$ and
$M=r_MT^M=s_Mc_MT^M$ as before.
If $c_L=c_M$ then $L=M$.
\end{theorem}
\begin{proof}
According to Proposition \ref{prop:regequiv},
$c_L=c_M\Leftrightarrow e\in
c_L^{-1}T_{M,rs}:=M_L\Leftrightarrow L^{temp}\subset M^{temp}$.
Hence the proof reduces to Remark 3.14 of \cite{HOH0},
or can be proved directly in our setup in the same way, cf.
Remark \ref{rem:smoothnest}.
\end{proof}
Theorem \ref{thm:nonnest} shows that a tempered coset can not be
a subset of a strictly larger tempered coset.
In fact even more is true:
\begin{thm}\label{thm:nonint}(Slooten \cite{Klaas} (cf. \cite{Slooten}
for the classical cases))
Let $L_1$ and $L_2$ be residual subspaces.
If $L_1^{temp}\cap L_2^{temp}\not=\emptyset$ then $L_1=w(L_2)$
for some $w\in W_0$.
\end{thm}
We will not use this result in this paper, but it is
important for the combinatorial fine structure of the
spectrum of $\cf$. We note that the proof of this statement reduces
easily to the case of two residual subspaces (in the sense of
\cite{HOH0}) with the same center. This reduces the
statement of the theorem to the problem in Remark 3.12 of
\cite{HOH0}. This problem was solved by Slooten \cite{Klaas}.
\section{Appendix: Kazhdan-Lusztig parameters}\label{KL}
Let $F$ be a p-adic field.
Let $\mathcal G$ be the group of $F$-rational points of a
split semisimple algebraic group of adjoint type over
$F$, and let $\mathcal{I}$ be an Iwahori subgroup of
$\mathcal G$. The  centralizer algebra of the representation of
$\mathcal G$ induced from the trivial representation of
$\mathcal{I}$ is
isomorphic to an affine Hecke algebra $\H$ with ``equal labels'',
that is, the labels are given as in Convention \ref{eq:scale}
with $\q$ equal to
the cardinality of the residue field of $F$, and the exponents
$f_s$ all equal to $1$. Moreover, the lattice $X$ is equal to the
weight lattice of $R_0$ in this case. The Langlands dual group $G$
is the simply connected semisimple group with root system $R_0$, and
the torus $T$ can be viewed as a maximal torus in $G$.

In this situation Kazhdan and Lusztig \cite{KL} have given a
complete classification of the irreducible representations of
$\H$, and also of the tempered and square integrable irreducible
representations. Let us explain the connection with
residual cosets explicitly.

We assume that we are in the ``equal label case'' in this
subsection, unless stated otherwise. We put $k=\log(\q)/2$. Let $G$
be a connected semisimple group over $\C$, with fixed maximal
torus $T=\operatorname{Hom}(X,\C^\times)$. We make no assumption
on the isogeny class of $G$ yet.
\begin{prop}\label{prop:cor}
\begin{enumerate}
\item[(i)] If $r$ is a residual point with
polar decomposition $r=sc=s\exp(\gamma)\in T_uT_{rs}$ and $\gamma$
dominant, then the centralizer $C_{\mathfrak{g}}(s)$ of $s$ in
$\mathfrak{g}:=\operatorname{Lie}(G)$ is a semisimple subalgebra
of $\mathfrak{g}$ of rank equal to
$\operatorname{rk}(\mathfrak{g})$, and $\gamma/k$ is the
weighted Dynkin diagram (cf. \cite{C}) of a distinguished nilpotent
class of $C_\mathfrak{g}(s)$.
\item[(ii)] Conversely, let $s\in T_u$ be such
that the centralizer algebra $C_{\mathfrak{g}}(s)$ is semisimple
and let $e\in C_{\mathfrak{g}}(s)$ be a distinguished nilpotent
element. If $h$ denotes the weighted Dynkin diagram of $e$
then $r=sc$ with $c:=\exp(kh)$ is a residual point.
\item[(iii)]
The above maps define a $1-1$ correspondence between $W_0$-orbits
of residual points on the one hand, and conjugacy classes of pairs
$(s,e)$ with $s\in G$ semisimple such that  $C_{\mathfrak{g}}(s)$
is semisimple, and $e$ a distinguished nilpotent element in
$C_{\mathfrak{g}}(s)$.
\item[(iv)]
Likewise there is a $1-1$ correspondence between $W_0$-orbits of
residual points and conjugacy classes of pairs $(s,u)$ with
$C_G(s)$  semisimple and $u$ a distinguished unipotent element of
$C_G(s)^0$.
\end{enumerate}
\end{prop}
\begin{proof}
(i). We already saw in Appendix \ref{sub:defn} that the rank of
$C_\mathfrak{g}(s)$ is indeed maximal. So we are reduced to the
case $s=1$. Let $\langle\q\rangle$ denote the group of integer
powers of $\q$, and denote by $R_{\q}\subset R_0$ the root subsystem
of roots $\a\in R_0$ such that $\a(c)\in\langle\q\rangle$. Now
$R_{\q}$ is a root subsystem of rank equal to
$\operatorname{rk}(R_0)$, with the property that $\forall \a,\b\in
R_{\q}$ such that $\a+\b\in R_0$ we have $\a+\b\in R_{\q}$. Of course,
$c$ is a residual point of $R_{\q}$. By an elementary result of
Borel and De Siebenthal there exists a finite subgroup $Z\subset
T_u$ such that $C_\mathfrak{g}(Z)$ is semisimple with root system
$R_{\q}$.

We claim that for every simple root $\a$ of $R_{\q}$ we have
$\a(c)=1$ or $\a(c)={\q}$. To see this, observe that all the roots
$\a\in R_{\q}$ with $\a(c)={\q}$ are in the parabolic system obtained
from $R_{\q}$ by omitting the simple roots $\a$ such that
$\a(c)=\q^{l}$ with $l>1$. If this would be a proper parabolic
subsystem, $c$ would violate Theorem \ref{thm:equal} in this
parabolic. This proves the claim.

Define the element $h:=\gamma/k$. Note that $h$  belongs to
$2P(R_{\q}^\vee)$ by the previous remarks. Consider the grading of
$R_{\q}$ given by this element, and define a standard parabolic
subalgebra $\mathfrak{p}$ of $C_\mathfrak{g}(Z)$ by
\[
\mathfrak{p}:= \mathfrak{t}\oplus\sum_{\{\a\in R_{\q}:\a(h)\geq
0\}}\mathfrak{g}_\a= \sum_{i\geq 0}C_\mathfrak{g}(Z)(i).
\]
Its nilpotent radical $\mathfrak{n}$ is
\[
\mathfrak{n}:=\sum_{i\geq 2}C_\mathfrak{g}(Z)(i),
\]
and by the definition of residual points we see that $P\subset
C_\mathfrak{g}(Z)$ is a {\it distinguished parabolic subalgebra}
(see \cite{C}, Corollary 5.8.3.). According to (\cite{C},
Proposition 5.8.8.) we can choose $e\in\mathfrak{n}(2)$ in the
Richardson class associated with $\mathfrak{p}$, and $f\in
C_\mathfrak{g}(Z)(-2)$, such that $(f,h,e)$ form a
$\mathfrak{sl}_2$-triple in $C_\mathfrak{g}(Z)$. By
$\mathfrak{sl}_2$ representation theory it is now clear that $h\in
P(R^\vee_0)$. Consider the grading of $\mathfrak{g}$ and $R_0$ induced
by $h$. By definition of $Z$ we see that
$\mathfrak{g}(0)=C_\mathfrak{g}(Z)(0)$ and
$\mathfrak{g}(2)=C_\mathfrak{g}(Z)(2)$. Hence $e$ is distinguished
in $\mathfrak{g}$ by (\cite{C}, Proposition 5.7.5.), proving the
desired result. Note also that, by (\cite{C}, Proposition 5.7.6.),
in fact $\mathfrak{g}(1)=0$, and hence that $R_{\q}=R_0$.

(ii). Is immediate from the defining property
\[
\operatorname{dim}(C_\mathfrak{g}(s)(0))=
\operatorname{dim}(C_\mathfrak{g}(s)(2))
\]
of the grading with respect to the Dynkin diagram of a
distinguished class.

(iii). Is clear by the well known $1-1$ correspondence between
distinguished classes and their Dynkin diagrams.

(iv). The result follows from the well known $1-1$ correspondence
between unipotent classes and nilpotent classes for connected
semisimple groups over $\C$.
\end{proof}
\begin{cor}\label{cor:cor}
From the proof of Proposition \ref{prop:cor}(i) we see that if
$r=sc$ is a residual point, then $\a(c)\in\langle\q^{1/2}\rangle$
for all $\a\in R_0$. If $s=1$ we have $\a(c)\in\langle\q\rangle$ for
all $\a\in R_0$.
\end{cor}
Let $M\subset T$ be a residual coset. Write $M=rT^M\subset
T\subset G$ with $r\in T_M$ as in Proposition \ref{prop:red}. Let
$r=sc=s\exp(kh/2)$ be the polar decomposition of $r$ in $T_M$. Let
$L_M\subset G$ be the Levi subgroup $L_M:=C_G(T^M)$ and let
$L^\prime_M$ denote its semisimple part. By Proposition
\ref{prop:red} we see that the root system of $L^\prime_M$ is
$R_M$, $T_M$ is a maximal torus of $L^\prime_M$, and the connected
center of $L_M$ is $T^M$. Moreover, $r\in T_M$ is a residual point
with respect to $R_M$. Thus by Proposition \ref{prop:cor},
$C_{L^\prime_M}(s)$ is semisimple, and there exists a
distinguished unipotent element $u=\exp(e)$ in
$C_{L^\prime_M}(s)^0$ such that $[h,e]=2e$. This implies that the
set $N=N_u$ of all elements $t\in G$ such that
\begin{equation}\label{eq:KL}
tut^{-1}=u^{\q}.
\end{equation}
is of the form $N=rC_G(u)$. The centralizer $C_G(r,u)=C_G(s,c,u)$
is known to be maximal reductive in $C_G(s,u)$, and it contains
$T^M$. Its intersection with $L^\prime_M$ is also reductive but,
since $u$ is distinguished in $C_{L^\prime_M}(s)^0$, the rank of
this intersection is $0$. Hence $L^\prime_M\cap C_G(r,u)$ is
finite. We conclude that $T^M$ is a maximal torus in $C_G(s,u)$.
Let $u^\prime$ be another unipotent element in $G$ such that
$M\subset N^\prime=N_{u^\prime}$ and such that $T^M$ is a maximal
torus of $C_G(s,u^\prime)$. We see that $u^\prime\in
C_{L^\prime_M}(s)^0$ is distinguished and associated to the Dynkin
diagram $h$. Hence $u^\prime$ is conjugate to $u$ in
$C_{L^\prime_M}(s)^0$ by an element of $C_{L_M^\prime}(r)$. We
have shown:
\begin{prop}\label{prop:heen}
For each residual coset $M=rT^M=scT^M\subset T$ there exists a
unipotent element $u$ such that $tut^{-1}=u^{\q}$ for all $t\in M$,
and such that $T^M$ is a maximal torus of $C_G(s,u)$. This $u$ is
an element of $C_{L_M^\prime}(s)$ with $L_M:=C_G(T^M)$, and is
distinguished in this semisimple group. It is unique up to
conjugation by elements of the reductive group
$C_{L_M^\prime}(r)$.
\end{prop}
Let us consider the converse construction. From now in this
subsection we assume that $G$ is simply connected. We will be
interested in conjugacy classes of pairs $(t,u)$ with $t$
semisimple and $u$ unipotent, satisfying (\ref{eq:KL}). We choose an
element $(t,u)$ in the conjugacy class. By Jacobson-Morozov's
theorem there exists a homomorphism
\begin{equation}\label{JM}
\phi:SL_2(\C)\mapsto G
\end{equation}
such that
\[
u=\phi \left(
\begin{array}{cc}
1&1\\ 0&1\\
\end{array}
\right)
\]
We put
\[
c:= \phi \left(
\begin{array}{cc}
\q^{1/2}&0\\ 0&\q^{-1/2}\\
\end{array}
\right),\
h:=d\phi \left(
\begin{array}{cc}
1&0\\ 0&-1\\
\end{array}
\right),\ e:= d\phi \left(
\begin{array}{cc}
0&1\\ 0&0\\
\end{array}
\right).
\]
Denote by $C_G(\phi)$ the centralizer of the image of $\phi$. We
have $C_G(\phi)=C_G(d\phi)$, and by $\mathfrak{sl}_2$
representation theory we see that $C_G(d\phi)=C_G(h,e)$. Hence
$C_G(\phi)=C_G(c,u)$. By \cite{KL}, Section 2, this is a maximal
reductive subgroup of $C_G(u)$, and we can choose $\phi$ in such a
way that $t\in cC_G(\phi)$. In this case $t$ commutes with $c$,
and thus $t_1:=tc^{-1}\in C_G(\phi)$ commutes with $c,t$, and is
semisimple. It follows that $C_G(t_1,\phi)=C_G(\phi)\cap C_G(t_1)$
is reductive in $C_G(t_1)$, and contains $t_1$ in its center.
According to \cite{KL}, the choice of $\phi$ such that $t_1\in
C_G(\phi)$ is unique up to conjugation by elements in $C_G(t,u)$.

By conjugating $(t,u)$ and $\phi$ suitably we can arrange that
$\overline{T}:=(T\cap C_G(t_1,\phi))^0$ is a maximal torus of
$C_G(t_1,\phi)$. Put $L=C_G(\overline{T})$, a Levi group of $G$.
We claim that $L$ is minimal among the Levi groups of $G$
containing $\phi$ and $t_1$. Indeed, if $N$ would be a strictly
smaller Levi group of $G$ also containing $\phi$ and $t_1$, then
its connected center ${T^N}$ would be a torus contained in
$C_G(t_1,\phi)$ on the one hand, but strictly larger than
$\overline{T}$ on the other hand. This contradicts the choice of
$\overline{T}$, proving the claim. In particular, since the
connected center $T^L$ of $L$ satisfies $\overline{T}\subset
T^L\subset C_G(t_1,\phi)$, we have the equality
$\overline{T}=T^L$.

Note that maximal tori of $L$ are the maximal tori of $G$
containing $T^L$, and these are conjugate under the action of $L$.
The derived group $L^\prime$ is simply connected, because the
cocharacter lattice $Y_L$ of its torus $T_L$ equals
$Y_L=Q(R^\vee_0)\cap\Q R^\vee_L=Q(R^\vee_L)$. Hence, by a well known
result of Steinberg, $C:=C_L(t_1)\subset L$ is connected, and
reductive. This implies that there exist maximal tori of $C$
containing $t$. Thus there exist maximal tori of $L$ containing
both the commuting semisimple elements $t_1$ and $t$. Therefore we
may and will assume (after conjugation of $(t,u)$ and $\phi$ by a
suitable element of $L$) that $T^L$ and the elements $t_1, t$ are
inside $T$.

Both the image of $\phi$ and $t_1$ are contained in $C$. Let
$C^\prime\subset L^\prime$ denote its derived group. If the
semisimple rank of $C$ would be strictly smaller than that of $L$,
there would exist a Levi group $N$ such that $C\subset
N\subsetneqq L$, a contradiction. Hence $C^\prime$ has maximal
rank in $L^\prime$.

Choose $s_{L}$ in the intersection $t_1T^L\cap L^\prime$. By the
above, $s_L$ is in $T_{L,u}$, the compact form of the maximal
torus $T_L:=(L^\prime\cap T)^0$ of $L^\prime$. We put
$r_{L}=s_Lc\in L^\prime$, and we claim that this is a
$R_L$-residual point of $T_L$. By Proposition \ref{prop:cor} this
is equivalent to showing that $u$ is a distinguished unipotent
element of $C^\prime=C_{L^\prime}(s_L)$. This means that we have
to show that  $C_{L^\prime}(s_L,\phi)$ does not contain a
nontrivial torus. But $L=C_G(T^L)$ with $T^L$ a maximal torus in
$C_G(t_1,\phi)$. Hence $C_G(t_1,\phi)^0\cap L=T^L$, and thus
\begin{equation}\label{eq:heus}
C_G(s_L,\phi)^0\cap L^\prime= C_G(t_1,\phi)^0\cap
L^\prime={T_L\cap T^L},
\end{equation}
proving the claim.

This proves that $M:=tT^L=r_LT^L\subset T$ is a residual coset, by
application of Proposition \ref{prop:conv}.

Notice that (\ref{eq:heus}) shows that $T^L$ is also a maximal torus
of $C_G(s_L,\phi)$, and thus of $C_G(s_L,u)$.

Finally notice that the $W_0$-orbit of the pair $(t,M)$ is uniquely
determined by the conjugacy class of $(t,u)$ by the above
procedure.
%Notice that $C_G(t_L,u)\cap L^\prime=
%C_G(t_L,\phi)\cap L^\prime=
%C_G(t,\phi)\cap L^\prime\subset C_G(t,u)\cap L^\prime$.
We have shown:
\begin{prop}\label{prop:terug}
%For every conjugacy class $\mathfrak{C}$ of pair
%$(t^\prime,u^\prime)$ with
%$t^\prime$ semisimple and $u^\prime$ unipotent
%satisfying (\ref{eq:KL}), we can choose a residual
%subspace $M=s_LcT^L \subset T$ and a $t\in M$ such that
%the pair $(t,u)\in \mathfrak{C}$, and where $u$ is the unipotent
%element associated with $M$ as in Proposition \ref{prop:heen}.
For every pair $(t,u)$ with $t$ semisimple and $u$ unipotent
satisfying (\ref{eq:KL}), we can find a homomorphism $\phi$ as in
\ref{JM} such that $t$ commutes with $c$. Let $T^L$ be a maximal
torus of $C_G(t_1=tc^{-1},u)$ and put $M=tT^L$. By suitable
conjugation we can arrange that $t$, $c$ and $M$ are in $T$. Then
$M\subset T$ is a residual coset. If we write $t=rt^L$ with
$r=sc\in T_{L,u}T_{L,rs}$ and $t^L\in T^L$, then $T^L$ is also a
maximal torus of $C_G(s,c,u)$. The $W_0$-orbit of the pair $(t,M)$
is uniquely determined by $(t,u)$.
\end{prop}
\begin{cor}
There is a one-to-one correspondence between conjugacy classes of
pairs $(t,u)$ satisfying (\ref{eq:KL}) and $W_0$-orbits of pairs
$(t,M)$ with $M\subset T$ a residual coset, and $t\in M$.
\end{cor}
\begin{proof} The maps between these two sets
as defined in Proposition
\ref{prop:heen} and Proposition
\ref{prop:terug} are clearly
inverse to each other.
\end{proof}
\begin{rem}
Let $(c,u)$ (with $c\in T_{rs}$) be a pair satisfying (\ref{eq:KL}),
with $u$ a distinguished unipotent element of $G$. Then $u$ will
be distinguished in $C_G(s)$ for each $s$ in the finite group
$C_G(c,u)$. In particular, $C_G(s)$ is semisimple. Hence $s$ gives
rise to a residual point $cs^\prime$ in $T$ where $s^\prime\in T$
is conjugate with $s$ in $G$. This defines a one-to-one
correspondence between the orbits in $C_G(c,u)$ with respect to
the normalizer $N_G(C_G(c,u))$ and the residual points in $T$ with
split part $c$.
\end{rem}
%It is interesting to see how this relates the $\H(W_0)$ characters
%that occur in the residual characters associated with the residual
%points of the form $cs^\prime$ with the families of $W_0$ characters
%and the exotic Fourier transform for that family (by restriction
%of the Plancherel decomposition to $\H(W_0)$. I should think about
%the case $G_2$!!
The Kazhdan-Lusztig parameters for irreducible representations of
$\H$ consist of triples $(t,u,\rho)$ where $(t,u)$ is as above,
and $\rho$ is an irreducible representation of the
finite group
\[
A(t,u)=C_G(t,u)/(Z_GC_G(t,u)^0),
\]
where $Z_G$ is the center of $G$. However, not all the irreducible
representations of $A(t,u)$ arise, but only those representations
of $A(t,u)$ which appear in the natural action of $A(t,u)$
on the homology of the variety of Borel subgroups of $G$
containing $t$ and $u$.

Moreover, Kazhdan and Lusztig
show that the irreducible representation $\pi(t,u,\rho)$ is
tempered if and only if $t\in M^{temp}$, where $M$ is the residual
subspace associated to the pair $(t,u)$. In this way we obtain a
precise geometric description of the set of minimal central
idempotents $\{e_i\}_{i=1}^{l_t}$ of the residue  algebra $\H^t$
for $R_M$-generic $t\in M^{temp}$.
%$\H=\H({\mathcal R},q)$
%\newpage
%\include{notation}
% references
%% ref.tex
%%%%%%%%%%%%%%%%%%%%%%%%%%%%%%%%%%%%%%%%%%%%%%%%%%%%%
%                     References                    %
%%%%%%%%%%%%%%%%%%%%%%%%%%%%%%%%%%%%%%%%%%%%%%%%%%%%%
\index{0kl@${\overline \ka}_{W_0r}$(=${\overline \ka}_{\Ri,W_0r}$),
rational factor in $\nu(\{r\})$|see{${\overline \ka}_{W_LL}$}}
\index{m@$m_{\{r\}}$(=$m_{\Ri,\{r\}}$)|see{$m_L$}}

\printindex
\end{document}